\documentclass[10pt, reqno]{amsart}
\usepackage{amssymb, latexsym}
\usepackage{amsmath}
\usepackage{amsthm}
\usepackage{graphicx}

\usepackage{titletoc}
\usepackage{pdfsync}
\usepackage{color, xcolor}
\usepackage{amsthm, amsmath, amssymb}
\usepackage{amsthm, amsmath, amssymb}
\usepackage{mathrsfs}

\usepackage{hyperref}
\hypersetup{
    colorlinks=true,
    linkcolor=blue,
    citecolor=red,
    urlcolor=blue
}

\usepackage{mathrsfs}

\numberwithin{equation}{section}
\newtheorem{theorem}{Theorem}[section]
\newtheorem{proposition}[theorem]{Proposition}
\newtheorem{lemma}[theorem]{Lemma}

\usepackage{multirow}
\usepackage[font=bf, aboveskip=15pt]{caption}
\usepackage[toc, page]{appendix}

\theoremstyle{definition}

\newcommand{\R} {\mathbb R}

\newcommand{\cuad}{{\sqcap\kern-.68em\sqcup}}

\newcommand{\e}{\epsilon}

\newcommand{\be}{\begin{equation}}
\newcommand{\ee}{\end{equation}}

\newtheorem{remark}[theorem]{Remark}

\def\G{{\mathfrak G}}

\def\R{{\mathfrak R}}

\def\begeq{\begin{equation}}
\def\endeq{\end{equation}}
\def\p{\partial}

\def\lf{\left}
\def\ri{\right}

\def\R{\Bbb R}

\def\G{\Gamma}

\allowdisplaybreaks

\begin{document}

\title
[Schr\"odinger equation   involving mass subcritical and supercritical exponents]
{Solutions with one dimensional concentration for a two dimensional Schr\"odinger model with general potential and involving mass subcritical and supercritical exponents}

\author{Lipeng Duan}
\address{Lipeng Duan,
\newline\indent School of Mathematics and Information Science, Guangzhou University,
\newline\indent Guangzhou 510006, P. R. China.
}
\email{lpduan777@sina.com}

\author{Suting Wei}
\address{Suting Wei,
\newline\indent Department of Mathematics, South China
Agricultural University,
\newline\indent Guangzhou, 510642, P. R. China.
}
\email{stwei@scau.edu.cn}

\author{Jun Yang$^\S$}
\address{School of Mathematics and Information Science,
Guangzhou University,
\newline\indent Guangzhou 510006, P. R. China.
}
\email{jyang2019@gzhu.edu.cn}

\begin{abstract}
We  prove the existence of  solution  $(\tilde u, \lambda)$  which concentrates on   stationary and non-degenerate   curve, as $p_\epsilon \to 6$,  to the following mass-constraint nonlinear Schr\"odinger equation
\begin{align*}
\begin{cases}
 - \Delta \tilde u  +   W({\bf y}) \tilde u   -  |\tilde u|^{p_\epsilon-2} \tilde u = \lambda \tilde u,
\quad {\bf y} \in \mathbb{R}^2,
\\[2mm]
 \int_{\R^2} |\tilde u|^2 {\mathrm d} {\bf y} =a.
\end{cases}
\end{align*}
Here
$\lambda$ is a Lagrange multiplier, $ W({\bf y})$ is a real-valued potential, $a \in \mathbb{R}_+$ is a constant, $p_\epsilon = 6 \pm \epsilon, \epsilon \ll 1$ and $6$ is the $L^2$-critical exponent   in the sense of concentration on a codimension 1 submanifold.
When  constructing the solution concentrating  on a  curve, we  utilize  some  delicate transformations to overcome the challenges posed by the non-uniform boundedness of  length of  concentrating  curves.
 In addition, we   employ some refined  analytical techniques to overcome the strong resonance phenomena caused by the $L^2$ constraint.
Our result partially answers the Ambrosetti-Malchiodi-Ni conjecture in the setting of the  mass-constrained nonlinear elliptic equation.


\vspace{4mm}

\vspace{2mm}
{\textbf{Keywords:}  Schr\"odinger  equation, Mass  subcritical,  Mass  supercritical, General potential, Concentration on curves}
\vspace{2mm}

{\textbf{AMS Subject Classification:} 35A01, 35B09, 35B25.}
\end{abstract}


\date{\today}


\maketitle

\tableofcontents

\medskip\section{Introduction} \label{Section1}

\subsection{ Background and problem}\

Consider the time dependent Schr\"odinger equation
\begin{align} \label{000original}
i \frac{\partial \Psi}{\partial t} = - \Delta
\Psi + W({\bf y} )\Psi -   |\Psi|^{p_\epsilon-2}\Psi, \quad (t, {\bf y}) \in {\mathbb R}\times \Omega,  ~ \Omega \subset \mathbb{R}^N, \ N\geq 2
\end{align}
with a mass constraint
\begin{align} \label{constraint00}
\int_{\R^N} |\Psi|^2  {\mathrm d} {\bf y} =a,
\end{align}
where $2<p_\epsilon<\frac{2N} {N-2} $,  $i$ is the imaginary unit, $a>0$  is a constant and $W( {\bf y} )$ is the external trap potential.
This kind of equation \eqref{000original} with the mass constraint \eqref{constraint00} is a
well-known mathematical model to describe Bose-Einstein condensates
\cite{aftalionriviere, aftalionbook, Cornell, Davis, Gross, LS, Lieb2, Pethick-Smith, Pitaevskii}.
If we look for a standing wave solution of \eqref{000original}-\eqref{constraint00} with the form
$$
\Psi({\bf y}, t) = \tilde u({\bf y}) e^{-i \lambda t},
$$
where $\lambda \in \R$ is   a frequency parameter to be determined and $\tilde u$ is a real function,
then the pair $(\lambda, \tilde u)$ must satisfy the following time-independent equation
\begin{align}  \label{ori0}
 - \Delta \tilde u + W({\bf y}) \tilde u -  |\tilde u|^{p_\epsilon-2} \tilde u = \lambda \tilde u,
\quad {\bf y} \in   \Omega\subset \mathbb{R}^N, \quad N\ge 2,
\end{align}
with the mass constraint
\begin{align}  \label{mass 0}
 \int_{\R^N} |\tilde u|^2 {\mathrm d} {\bf y} =a.
\end{align}
We also know that a ground state solution of \eqref{eq1} is the minimizer of the following problem
\begin{align}\label{minimzation}
 E_a : = \min_{ u \in \mathscr M} J_a (u),
\end{align}
where
\[ J_a (u) = \int_{\R^N} \Big( \frac{1}{2}| \nabla u|^2 + \frac{1}{2} W({\bf y}) u^2 -\frac 1 {p_\e} |u|^{p_\e} \Big) {\mathrm d} {\bf y}, \]
and
\[\mathscr M_a = \Big\{ u \in H^1(\R^N)\, :\, \int_{\R^N} |u|^2 {\mathrm d} {\bf y} =a\quad \text{and} \quad \int_{\R^N} W({\bf y}) |u|^2 {\mathrm d} {\bf y} <+\infty\Big\}. \]
From   \eqref{minimzation},
 $\lambda$  can be viewed as the Lagrange multiplier.
By  scaling
\begin{align}   \label{tilde scaling U}
    \tilde w({\bf y}) = \lambda^{\frac{1}{p_\epsilon-1}} V(\lambda^{\frac{1}{2}} {\bf y}),
\end{align}
    we have
\begin{align}
    \int_{\R^N} |\tilde w|^2 {\mathrm d} {\bf y}  =   \lambda^{-\frac{2}{p_\epsilon-1}+\frac N2}\varrho_N  = a,  \quad \text{with}~ \varrho_N  =     \int_{\R^N}V^2,
    \end{align}
    where  $U$  solves
\[
-\Delta V + V = V^p, \quad V \in H^1(\mathbb{R}^N), \quad V > 0 \text{ in } \mathbb{R}^N.
\]
We know that $\tilde w({\bf y})$  given in \eqref{tilde scaling U} solves   problem   \eqref{ori0}-\eqref{mass 0}  in autonomous case (with $W\equiv 0$) and $\Omega =\R^N$.
If
\begin{align}
-\frac{2}{p_\epsilon-1}+\frac N2=0,\ \text{i.e.}\ p_\epsilon = 1+ \frac 4N:= \bar p,
\end{align}
we   refer to constraint problem  \eqref{ori0}-\eqref{mass 0} as mass-critical problem;  if $p_\epsilon >1+ \frac 4N$ (or $p_\epsilon <1+ \frac 4N$),  then  problem  \eqref{ori0}-\eqref{mass 0}  is said to be mass-supercritical   (or mass-subcritical) problem.
 The  criticality of  such exponent  crucially determines the existence of normalized  ground state solutions for the problem \eqref{ori0}-\eqref{mass 0}.  Analysis in this regard can be found in \cite{Pistoia}.
 On the other hand,
 this  exponent's criticality carries implications for other features of  \eqref{ori0}-\eqref{mass 0}  as well, especially concerning dynamical behaviors such as orbital stability and blow-up, while also linking to the exponents appearing in the Gagliardo-Nirenberg inequality; see \cite{Ca, weinstein1983}.

 \medskip
 Note that because of the role played by the mass constraint in   \eqref{ori0}-\eqref{mass 0}, the existence of normalized solutions becomes    a nontrivial issue, and many techniques developed for the classical Schr\"odinger equation
 cannot be directly applied in this setting.
    As for these  reasons, problem \eqref{ori0}-\eqref{mass 0} has significant mathematical interest and  garnered considerable attention from researchers, see \cite{bartsch2013, bartsch2018, bartsch2016, bartsch2017,  bartsch2019,   Guo-luo-peng, KG, Luo} and references therein.
 Variational methods focus more on minimizers or saddle point solutions of problem  \eqref{ori0}-\eqref{mass 0} while excited state solutions are rarely studied.
 On the other hand,   the minimization problem \eqref{minimzation}  is  not achieved in some cases, for example when the  exponent $p_\epsilon$  is supercritical.  Therefore, it is necessary to adopt  non-variational methods  to study  problem  \eqref{ori0}-\eqref{mass 0}.
    In \cite{Pistoia},  problem \eqref{ori0}-\eqref{mass 0}  is viewed  as a singular perturbation problem. For fixed $p_\epsilon$,  by applying the Lyapunov-Schmidt reduction method, the authors constructed  a solution  $(\lambda, \tilde u)$, with $\lambda \to - \infty$ and $\tilde u$     concentrating  on a non-degenerate critical point of the potential $W$, as $a \to a^*$, where
\begin{align*}
  a^*= \begin{cases}
 +\infty, \quad\text{if $p_\epsilon$ is subcritical},
 \\[2mm]
  0, \qquad\,\text{if $p_\epsilon$ is supercritical},
  \\[2mm]
  2  \varrho_N, \quad\text{if $p_\epsilon =  1+ \frac4N$}.
\end{cases}
\end{align*}
Recently,   Q. Guo, R. He, B. Li, and S. Yan \cite{Guoqidong1}  studied problem  \eqref{ori0}-\eqref{mass 0} in $\R^N$  in both the mass-supercritical and mass-subcritical cases
   and  constructed    multiple-peak blow-up solutions as $p_\epsilon \to \bar p$.  In  their study, they showed that  it is the number of the critical points of the potential $W$ that affects the existence and   number of solutions  to  problem \eqref{eq1}-\eqref{constraint0}.
    Moreover,    they  proved   local uniqueness for the solutions they constructed.
        On the subject of concentration phenomena of solutions to  problem \eqref{eq1}-\eqref{constraint0} via singular perturbation methods,  see \cite{DWY, Guoqing-chunhua,   Zhouyang, Guoyang,  Huangxiaomeng,   Luo} and the references cited therein.

     We also note  the mathematical similarity between   problem  \eqref{ori0}-\eqref{mass 0} and  the following singular perturbation problem:
\begin{align} \label{schrodinger}
 - \epsilon^2 \Delta u + W({\bf y}) u - | u|^{p-1} u = 0,\quad u \in H^1(\Omega),   \quad u>0~\text{in}~\Omega \subseteq {\mathbb R}^N,
\end{align}
where $1<p<\frac{N+2}{N-2}, $ $\epsilon >0$ is a small parameter, and $W$ (independent of $\epsilon$) is the potential.
Problem \eqref{schrodinger} appears in many applications, for  example,  in chemotaxis \cite{LinNiTak1988} and biological pattern formation \cite{Gierer-Meinhardt}.
Problem \eqref{schrodinger} with Dirichlet or Neumann boundary conditions  has been  studied
extensively in recent years, involving  the necessary and sufficient conditions for concentration,  boundary concentration, multi-peak phenomena, and high-dimensional concentration (i.e., on submanifolds). See  \cite{AMNconjecture1, AMNconjecture2, Ambrosetti, AmbMalSecc2001,
BP, BP1, BP2, CPY2019, DY10,
delPFelmer1996, delPFelmer1997, delPKowWei1, FloerWein1986, GuoYang, Malchodi1, WangWeiYang, Wangzhao, WeiYang, Wei-xu-yang, SWei-Yang}
for  background and studies on \eqref{schrodinger}.

        High-dimensional concentration for  problem \eqref{schrodinger}  refers to    the concentration of solutions on  \( k \)-dimensional manifolds for \( k = 1, \dots, N - 1 \).
        In   \cite{AMNconjecture1},   the authors  formulated   the following  Ambrosetti-Malchiodi-Ni   conjecture: if the limit set $\Gamma$ satisfies the stationary  condition
\begin{align}  \label{station in Nk}
 \theta_k \nabla^\perp W= W\mathbf{H} \quad \text{with} \quad \theta_k = \frac{p+1}{p-1} - \frac{1}{2}(N-k),
\end{align}
then there exists a solution to \eqref{schrodinger} which concentrates  near $\Gamma$, at least along a subsequence  $\epsilon_j \to 0$, as $j \to \infty$,
where $\mathbf{H}$ is the mean curvature of $\Gamma$
and $\nabla^\perp$ denotes the component of $W$ normal to $\Gamma$.
This  conjecture was substantiated in  \cite{delPKowWei1, machiodi,  machiodi2, WangWeiYang} for different dimensions   $N$ and various values of $k~(1\le k \le N-1$).
     Recently,  we  studied    high-dimensional concentration  of   the following  mass constraint problem  in   \cite{DWY},
\begin{align} \label{eq1-gp}
 - \Delta \tilde u + W({\bf y}) \tilde u -a |\tilde u|^2 \tilde u = \lambda \tilde u,
\quad  \int_{{\mathbb R}^N} |\tilde u|^2 {\mathrm d} {\bf y} =1,
\end{align}
      and demonstrated that condition  \eqref{station in Nk}  is necessary for high-dimensional concentrating solutions  in the special case $N=2$  and $k=1$.

       Motivated by the study of the Schr\"odinger equation in higher dimensions, we  are led to ask the following question:
       {\em can the    Ambrosetti-Malchiodi-Ni conjecture be generalized to the mass  constraint problem \eqref{ori0}-\eqref{mass 0}?}
 More precisely, we investigate whether there exist solutions concentrating on stationary sets as \( p_\e \to \bar{p}' \), where $\bar p'$  denotes  the  $L^2$ mass critical exponent  in the sense of high-dimensional concentration for solutions  to \eqref{ori0}-\eqref{mass 0}.   In this paper, we  provide a partial answer to this problem.

\medskip
       Unlike the case of pointwise concentration, the  $L^2$ mass critical exponent $\bar p'$  in the sense of high-dimensional concentration   is no longer   the usual
           $1+ \frac 4N$.
     In order to construct high-dimensional concentrated solutions with
 $L^2$
  mass constraint  problem \eqref{ori0}-\eqref{mass 0},  we  first pose  the following problem:
 \[
\textit{What is  the mass critical exponent $\bar p'$  for problem \eqref{ori0}-\eqref{mass 0} in the sense of high-dimensional concentration?}
\]
 For the above question, we partially answer it by taking the concentrating submanifold to be of codimension 1.
      For   simplicity, we consider the situation where the concentrating solution is radial and equation \eqref{ori0} is autonomous (i.e., $W\equiv 0$), thereby illustrating the above problem.

 \medskip
\noindent$\clubsuit${\em  Analysis of  higher concentration phenomena for radial case.}
   For  \eqref{ori0},
we consider the  autonomous   case   $W\equiv 0$  in  $ \R^N$
\begin{align}  \label{ori0W0}
 - \Delta \tilde u  -  |\tilde u|^{p_\epsilon-2} \tilde u = \lambda \tilde u,
\quad {\bf y} \in \mathbb{R}^N, \quad N\ge 2,
\end{align}
with the mass constraint
\begin{align}  \label{mass 0-1}
 \int_{\R^N} |\tilde u|^2 {\mathrm d} {\bf y} =a.
\end{align}
Following the approach in \cite{AMNconjecture1,  AMNconjecture2}, we present   radial solution   to problem \eqref{ori0W0}   which  concentrates  on the sphere of radius $r_0$ with the following  asymptotic profile
\begin{align}
\tilde u({\bf y}) \sim \tilde u_1({\bf y})  := (-\lambda)^{\frac 1{p_\epsilon-2}} U_\epsilon\Big( \sqrt{-\lambda\,} {\bf y}\Big),
\end{align}
where   $U_\epsilon(r)$  is  the unique  positive solution of
    $$
    - U_\epsilon'' + U_\epsilon = U_\epsilon^{p_\e-1}, \quad   \text{in}~ \R.
    $$
Moreover,  we know that
\begin{align}  \label{sphere scaling}
  \int_{\R^N} |\tilde u|^2 {\mathrm d} {\bf y}  \sim &  \int_{\R^N} |\tilde u_1|^2 {\mathrm d} {\bf y}   =    {(-\lambda)}^{\frac 2{p_\epsilon-2}} \int_{\R^N}\left|U_\epsilon\left(\sqrt{-\lambda\,} (| {\bf y}|-r_0) \right)\right|^2{\mathrm d} {\bf y}
     \nonumber\\[2mm]
  &  = {(-\lambda)}^{\frac 2{p_\epsilon-2}}   \int_{\mathbb{S}^{N-1}} {\mathrm d} \sigma \int_0^{+\infty} \left|U_\epsilon\left(\sqrt{-\lambda\,} (r-r_0) \right)\right|^2 r^{N-1}  {\mathrm d} r
    \nonumber\\[2mm]
  &  = {(-\lambda)}^{\frac 2{p_\epsilon-2}-\frac 12}   \int_{\mathbb{S}^{N-1}} {\mathrm d} \sigma \int_{-\sqrt{-\lambda\,} r_0} ^{+\infty} |U_\epsilon\left(r \right)|^2 \left( {(-\lambda)}^{-\frac12} r + r_0\right)^{N-1}  {\mathrm d} r
  \nonumber\\[2mm]
& =  |\mathbb{S}^{N-1}| {(-\lambda)}^{\frac 2{p_\epsilon-2}- \frac 12}
\int_{-\sqrt{-\lambda\,} r_0}^{+\infty} |U_\epsilon( r)|^2  \sum_{j=0}^{N-1} C_{N-1}^j \Big[{(-\lambda)}^{-\frac12} r\Big]^j (r_0)^{N-1-j}     {\mathrm d} r,
  \end{align}
    where $|\mathbb{S}^{N-1}|$   denotes Lebesgue measure of the unit sphere $\mathbb{S}^{N-1}$.
            For the high-dimension concentration solution,  it's known  that scaling invariance  does not hold.
            However, from  \eqref{sphere scaling},  we may reasonably define the critical exponent $\bar p'$ 	
  in this case by the following
     \[\frac 2{\bar p'-2}- \frac 12 =0, \quad \text{i.e}, \quad \bar p' =6.  \]
     On the other hand, from \eqref{sphere scaling}, we can also obtain   the asymptotic behavior of the solution $(\tilde u, \lambda)$  of   \eqref{ori0W0}-\eqref{mass 0-1} concentrating on sphere  as $p_\epsilon \to \bar p'$.
      Under the condition  that  $-\lambda$  is   large,  \eqref{sphere scaling} yields the following,
\begin{itemize}

\item[(i).]
{\em
  if $a > |\mathbb{S}^{N-1}| r_0  \int_\R |U_\epsilon( r)|^2    {\mathrm d} r $,
    then  $ -\lambda \to \infty$  as  $p_\epsilon \to 6^-$;
  }
   \\[0.5mm]
   \item[(ii).]
{\em
  if $a < |\mathbb{S}^{N-1}| r_0  \int_\R |U_\epsilon( r)|^2    {\mathrm d} r $,   then  $-\lambda \to \infty$ as  $p_\epsilon \to 6^+$.
  }

     \end{itemize}

\smallskip
\subsection{ Constructions of solutions    concentrating on a codimension-one submanifold as $p_\epsilon\to\bar p'$}\
 For simplicity, we consider only
 $N=2$ in this paper.   We set $p_\e = 6+ \e$ or $p_\e = 6-\e$,  with $\epsilon \to 0$.
Then we consider the problem
\begin{align} \label{eq1}
 - \Delta \tilde u + W({\bf y}) \tilde u -  |\tilde u|^{p_\epsilon-2} \tilde u = \lambda \tilde u,
\quad {\bf y} \in \mathbb{R}^2,
\end{align}
with the mass constraint
\begin{align} \label{constraint0}
 \int_{\R^N} |\tilde u|^2 {\mathrm d} {\bf y} =a.
\end{align}
 In this paper, we study the existence of non-radial  solution  $\tilde u_{p_\e} $ of    \eqref{eq1}-\eqref{constraint0},   concentrating on a closed curve, with $p_\e \to \bar p'$.
Throughout this paper,  we impose the following assumption:
\begin{align} \label{lambda a to bf a}
-\lambda \rightarrow + \infty, \quad \text{as}~\e \to 0.
\end{align}
The quantitative relationship between
$\lambda$
and $\e$ will be given later.
Letting
\begin{align}
\tilde u =\left( -\lambda \right)^{\frac1{p_\epsilon-2}}u, \quad \iota= \frac 1{\sqrt{-\lambda\,}},
\end{align}
then we get that $ u $ satisfies
\begin{align} \label{equ power 6}
 \iota^2 \Delta u - \Big(1 + \iota^2 W({\bf y}) \Big) u + | u|^{p_\epsilon-2} u = 0,
\quad {\bf y} \in \mathbb{R}^2,
\end{align}
and
\begin{align} \label{constraint1}
\int_{\R^N} | u|^2 {\mathrm d} {\bf y}  \,=\, a \left( {-\lambda}  \right)^{-\frac2{p_\epsilon-2} +1}=  a\iota^{-2+ \frac4{p_\epsilon-2}}.
\end{align}
  The substitution
\begin{equation}\label{definitionofy}
 u ({\bf y}) \,=\, v(y)\quad\mbox{with}\quad {\bf y} = \iota y,
\end{equation}
in problem \eqref{equ power 6}-\eqref{constraint1} will imply that
\begin{align} \label{eq3}
 -\Delta v + \big[1 + \iota^2  W ( \iota y ) \big] v - |v|^{p_\epsilon-2} v \,=\, 0
 \quad\mbox{and}\quad v>0, \quad \forall\, y\in \R^2,
\end{align}
and
\begin{align} \label{constraint2}
 \int_{\R^2} |v|^2 {\mathrm d} y = a\iota^{-2+ \frac4{p_\epsilon-2}}.
\end{align}

\medskip
Since the potential in  \eqref{eq3}  depends on  $\iota$ and from  \eqref{station in Nk}  together with the necessary conditions for concentrating submanifolds established in  \cite{DWY},
 it follows  that    the concentrating curve is   related to $\iota$.
 To this end, we denote these  concentrating   curves by  $\G_\iota$.
In this paper, we intend  to construct a solution for  problem  \eqref{eq3}-\eqref{constraint2} which concentrates on the curve $\Gamma_\iota$, as  $p_\epsilon \to \bar p'$.
For  this  purpose,  we adopt the following settings  to avoid tedious computations and make the  mathematical analysis rigorous.
 We first impose assumptions on the curves. These curves are the concentration sets of the solutions

\smallskip\noindent {\textbf{(A1).}}
{\it
Assume that $\Gamma_\iota$ are simple smooth closed curves in ${\mathbb R}^2$.
Consider the natural parameterization $\gamma_\iota(\theta)$ of $\Gamma_\iota$ with positive orientation, where $\theta$ denotes the arclength parameter measured from a fixed point of $\Gamma_\iota$.
We denote
\begin{align}
\ell = |\Gamma_\iota|,
\label{ell-epsilon}
\end{align}
and the curvature of $\G_\iota$ by $\kappa(\theta)$.
We also assume that there exists a positive constant $C_0$ independent of $\epsilon$ such that
\begin{align} \label{curvature bounded}
|\kappa(\theta)|, \, |\kappa'(\theta)|\leq C_0,
\quad\forall\, \theta\in (0, \ell).
\end{align}
}

\smallskip
The assumption $-C_0\leq\kappa(\theta)\leq C_0$ in \eqref{curvature bounded} implies that
there exists a small positive number $M_0$ (independent of $\epsilon$) in such a way that any point ${\bf y} $ near $\Gamma_\iota$ in $\R^2 $ can be represented in {\textbf{Fermi coordinates}}
\begin{align}\label{fermi}
{\bf y}= \gamma_\iota(\theta) +t \nu_\iota(\theta), \quad \forall\, \theta\in [0, \ell), \ \forall\,t\in(-M_0, M_0),
\end{align}
where $ \nu_\iota(\theta) $ denotes the unit outer normal to $\Gamma_\iota$.
By a slight abuse use of notation, for ${\bf y}$ in \eqref{fermi} and any function $W$, we denote $W(t, \theta)$ to actually mean $W({\bf y})$.
 We now state the assumptions for the setting of the potential function $W$.

\smallskip
\noindent {\textbf{(A2).}}
{\it
The smooth potential function $W$ is positive on ${\mathbb R}^2$.
Moreover, there exist two positive constants $C$ and $\varrho<1$ independent of $\iota$, such that
\begin{align}
\int_0^{\ell}\big[1+ \iota^2 W (0, \theta) \big]\,{\mathrm d}\theta \,\leq\, C\iota^{\varrho-1},
\label{length-constraint}
\end{align}
\begin{align}
\max\big\{\, \iota^2|W_t (t, \theta)|, \ \iota^2|W_{tt} (t, \theta)| \, \big\}
\leq
C \big[1+\iota^2 W(0, \theta)\big], \quad \forall\, \theta \in [0, \ell),\ \forall\,t\in(-M_0, M_0),
\label{H0}
\end{align}
and
\begin{align} \label{Wthetaestimate}
& \max\Big\{\, \iota^2|W_{\theta} (0, \theta)|, \iota^2|W_{t\theta} (0, \theta)|, \iota^2|W_{t\theta} (0, \theta)|,  \iota^2|W_{t\theta\theta} (0, \theta)|,  \iota^2|W_{tt\theta\theta} (0, \theta)|,  \nonumber
\\[2mm]
 & \quad \quad \quad \iota^2|W_{\theta\theta} (0, \theta)|, \ \iota^2|W_{\theta\theta\theta} (0, \theta)|,\, \iota^2|W_{\theta\theta\theta\theta} (0, \theta)|\Big\} \leq C \big[1+\iota^2 W(0, \theta)\big],
\quad \forall\, \theta \in [0, \ell).
\end{align}
}

\smallskip
The following definitions are adapted from \cite{delPKowWei1}.
Recall the Fermi coordinates defined as in \eqref{fermi}. For $\iota$ small, any curve $\G^g_\iota$ close to $\G_\iota$ can be parameterized as
\begin{align*}
\gamma_\iota^g (\theta) = \gamma_\iota(\theta) + g (\theta) \nu_\iota (\theta), \quad \theta \in (0, \ell).
\end{align*}
Denote
\begin{align}
\Lambda_\epsilon= \frac{p_\epsilon} {p_\epsilon-2} - \frac12 =\frac 12 \frac{p_\epsilon+2} {p_\epsilon-2} = \frac 12 \frac{8\pm \e} {4\pm \e}.
\end{align}
Following the analysis and computations as in \cite{delPKowWei1}, we  introduce the following  weight length functional and deduce that
\begin{align*}
J_\iota (g) & = \int_{\G_\iota^g} \big[1+ \iota^2 W (g, \theta) \big]^{\Lambda_\epsilon} {\mathrm d} \sigma
\\[2mm]
&= \int_0^{\ell} \Big[1+ \iota^2 W (g, \theta) \Big]^{\Lambda_\epsilon}\, \Big[ (1+\kappa g)^2 + |g'|^2\Big] {\mathrm d} \theta.
\end{align*}
In the present paper, for any fixed small $\iota$, we just make the {\emph {stationary assumption} for the curve $\Gamma_\iota$, which means that the first variation of the functional $J_\iota (g)$ at $g=0$ equals to $0$.
In other words, for any smooth, $\ell$-periodic function $g(\theta)$, there must be
\[
\displaystyle\int_0^{\ell} \Bigg\{ \big[1+ \iota^2 W (0, \theta) \big]^{\Lambda_\epsilon-1}  \iota^2 W_t (0, \theta)\, g+ \big[1+ \iota^2 W (0, \theta) \big]^{\Lambda_\epsilon} \kappa g\Bigg\} {\mathrm d} \theta =0,
\]
which is equivalent to the following condition
\begin{align} \label{stationary assumption}
\kappa(\theta) \big[1+ \iota^2 W (0, \theta) \big] + \Lambda_\epsilon \iota^2 W_t(0, \theta) = 0,
\quad\forall\, \theta\in (0, \ell).
\end{align}
We say that a stationary curve $\Gamma_\iota$ is {\emph {non-degenerate}} if
the following problem
\begin{align}
\begin{cases}
\Big(\big[1+ \iota^2 W (0, \theta) \big]^{\Lambda_\epsilon}g' \Big)'
- \Bigg( \Big( \big[1+ \iota^2 W \big]^{\Lambda_\epsilon} \Big)_{tt} (0, \theta)
- 2 \big[1+ \iota^2 W (0, \theta) \big]^{\Lambda_\epsilon} \kappa^2 \Bigg) g=0,
\\[5mm]
g(0) = g(\ell), \quad g'(0) = g'(\ell),
\label{nondeneracy}
\end{cases}
\end{align}
only has the trivial solution $g\equiv 0$.

By setting
\begin{align}
\mathbb{F}_\iota(h)
=
\Big(\big[1+ \iota^2 W(0, \theta) \big]^{\Lambda_\epsilon}h' \Big)'
\ -\
\Bigg\{ \Big( \big[1+ \iota^2 W(t, \theta) \big]^{\Lambda_\epsilon} \Big)_{tt}\Big|_{t=0}
\,-\, 2 \big[1+ \iota^2 W(0, \theta) \big]^{\Lambda_\epsilon} \kappa^2 \Bigg\} h,
\label{mathbbF}
\end{align}
(see \eqref{weighted-reduced tilde f} for the expression of $\mathbb{F}_\iota$ in other coordinates),
we say that a stationary curve $\Gamma_\iota$ is {\emph {non-degenerate}} if
the following problem
\begin{align}
\mathbb{F}_\iota(h)=0,
\quad
h(0) = h(\ell), \quad h'(0) = h'(\ell),
\end{align}
only has the trivial solution $h\equiv 0$.
Inspired by this, to obtain a precise {\em a priori} estimate for the corresponding in-homogeneous problem in the reduction procedure, we propose the following assumptions

\smallskip\noindent {\bf{(B1).}}
{\it{
For any small $\epsilon$, $\Gamma_\iota$ is a stationary curve.
Moreover, we assume that
\smallskip
\begin{align}\label{resonace-f}
{\mathrm d}_\iota \geq \hat{C}_2 \iota^{\alpha}
\quad \text{with}\ {\mathrm d}_\iota = \min_{j=1,2,\cdots} |\Xi_{\iota, j}|,
\end{align}
where $\alpha\in(0, 1)$ is a constant with constraint in \eqref{constraint-alpha-varrho}, and $\{\Xi_{\iota, j}\} _{j=1}^\infty$ are the real eigenvalues of the eigenvalue problem
\begin{align}
\mathbb{F}_\iota(h) = \Xi_{\iota, j}\big[1+ \iota^2 W(0, \theta) \big]^{\Lambda_\epsilon+2}h,
\qquad
h(0) = h(\ell), \qquad h'(0) = h'(\ell).
\label{mathbbF-boundary}
\end{align}
}}

\medskip
On the other hand, we will also encounter delicate resonance phenomena in the construction of solutions, which were analyzed in \cite{delPKowWei1}.
For the similar purpose, the following gap condition for $\epsilon$ small will be useful.

\smallskip\noindent {\bf{(B2).}}
{\it{
For given $c>0$ small, there holds
\begin{align} \label{gap condition1}
\Big| j^2 \frac{\iota^2}{{\tilde\ell}^2} - \lambda_* \Big| \geq \frac{c \iota}{{\tilde\ell}}, \quad \forall\,j \in \mathbb N,
\end{align}
where by recalling the unique positive eigenvalue $\lambda_0$ of the eigenvalue problem \eqref{eigenvalue}, we have used the conventions
\begin{align}
\lambda_* = \frac{\lambda_0} {4\pi^2}
\qquad \text{and}\qquad
{\tilde\ell} = \int_0^{\ell} \big[1+ \iota^2 W(0, \theta)\big]^{\frac{1}{2}} {\mathrm d} \theta.
\label{lambda*}
\end{align}
}}

\medskip
To guarantee that the reduction process is successfully implemented, the following conditions are needed to ensure that the reduced equation  is solvable, see Section \ref{section10}.

\smallskip \noindent {\bf{(B3).}}
{\it{
For $\varrho, \alpha\in (0, 1)$ given in \eqref{length-constraint} and \eqref{resonace-f} respectively, they will satisfy
\begin{align}
 \varrho> \frac 13+ \frac43 \alpha.
\label{constraint-alpha-varrho}
\end{align}
}
}

The  constructed solutions  are required to satisfy the constraint \eqref{constraint2}, hence we obtain the following  identity:
\begin{align} \label{mass identity-section1}
1&=     \iota^{1- \frac4{p_\epsilon-2}} \frac{\hat \ell   \varrho_0}a \left(1+ O\left(\iota\right) \right),
  \end{align}
  where $\hat \ell =\int_0^{\ell} \beta^{\frac{6-p_\epsilon} {p_\epsilon-2}} \,{\mathrm d}\theta,  \varrho_0 = \int_\R |U_\epsilon(r)|^2\,{\mathrm d}r. $
  For the detailed derivation of \eqref{mass identity-section1}, see Section \ref{Section9}.
    In order to find a solution $\iota (\epsilon)$   to   \eqref{mass identity-section1} satisfying \eqref{lambda a to bf a}, we need the following conditions.

 \smallskip \noindent {\bf{(B4).}}
{\it{
 We assume that  there exists $c_i(\e), i =1,2,3$  independent of $\iota$ such that
\begin{itemize}

\item[(i)]
if  $\frac{\hat \ell   \varrho_0}a < c_1(\e)<1, $   and we take $p_\epsilon = 6- \e$,   with $\frac \iota {1-c_1(\epsilon)} \in o\left(1 \right), $ as $\e \to 0$;
   \\[0.5mm]
   \item[(ii)]
  if  $1<c_2 (\epsilon) <\frac{\hat \ell   \varrho_0}a <c_3(\e), $  we take $p_\epsilon = 6+ \e$,  with $\frac \iota {c_2(\epsilon)-1} \in o\left(1 \right).$
 \end{itemize}
}
}

\begin{remark}\label{remark2-4}
{\it
Following the approach in   \cite{DWY}, we can show that the lengths of the concentrated  curve $\Gamma_\iota$ possess a uniform lower bound, i.e, .
\begin{align}\label{ellzero}
 \ell\,\geq\, \ell_0.
\end{align}
From \eqref{ellzero}, we  know that,   if we take $p_\epsilon = 6-\epsilon$,  there holds
\begin{align*}
 \hat\ell\geq \ell_0>0.
\end{align*}
Since the variable
 $\hat \ell$
  is related to
 $\iota$ and little quantitative information is available about the relation between them, it is difficult to give a direct solution of equation \eqref{mass identity-section1}. To ensure its solvability, we impose conditions in  {\bf{(B4)}}.
}
\end{remark}

\medskip
\begin{theorem}
\label{theorem 1.1} \label{theorem1-1}
Assume that $N=2$, the smooth simple closed curves $\Gamma_\iota$ and the potential $W$ satisfy the assumptions in
{\bf{(A1)}}-{\bf{(A2)}} and {\bf{(B1)}}-{\bf{(B4)}}.
  Then there  exists a solution $u_\epsilon$  to problem  \eqref{equ power 6}-\eqref{constraint1}.
  \\
\smallskip
\noindent $\clubsuit$
For ${\bf y}$ given by \eqref{fermi} in the neighbourhood of $\Gamma_\iota$,
$u_{\epsilon}$ has the local profile
\begin{equation}
\label{taketheform}
 u_{\epsilon}({\bf y})\,=\, \big[1 + \iota^2 W(0, \theta) \big]^{\frac 12}\,  U_\epsilon\Bigg(\big[1  + \iota^2  W(0, \theta) \big]^{\frac 12}\, \frac{t}  {\iota} \Bigg)\big(1+o(1)\big),
\end{equation}
where $U_\epsilon$ is the  radial solution of  \eqref{Schrodiner Profile}. \\
\smallskip
\noindent$\clubsuit$
The solution $u_{\epsilon}$ satisfies globally
$$u_{\epsilon}({\bf y})\,\leq\, C\exp \Big(-  \frac{c_0}\iota{\mathrm{dist}} ({\bf y}, \Gamma_\iota)\Big),$$
where $c_0$ and $C$ are two positive constants independent of $\iota$.

 Moreover,  we have  $\iota \in \left(e^{\frac{6} \e \ln  {\frac{\ell_0   \varrho_0}a}},   e^{\frac{2} \e \ln  c_1(\e)}   \right)$ ( for case i)  and  $\iota \in \left( e^{-\frac{2} \e \ln c_2(\e)},   e^{-\frac{6} \e \ln  c_3(\e)}  \right)$ ( for case ii), where $\ell_0$ is the   uniform lower bound of the length of
 concentrated  curve $\Gamma_\iota$, and there holds
$$\iota\to 0 \quad \text{as}~ \epsilon \to 0 .$$

\qed
\end{theorem}

\smallskip
To prove Theorem \ref{theorem1-1}, we will employ the infinite-dimensional Lyapunov-Schmidt reduction method developed in \cite{delPKowWei1}, together with some clever coordinate transformations.
The paper is organized as follows.
In Sections \ref{Section3.1}-\ref{Section3.2}, we will introduce the Fermi coordinates and write down equation \eqref{eq3} in local forms,
see \eqref{eq4} as well as \eqref{expression in xtau}.
After that, some preliminary facts on the profile function $U_\epsilon$, and some conventions will be given in Section \ref{section3-3}.
 Next, we will present the  construction and the reduction process in Sections   \ref{section-approximatesolution}-\ref{section10}.
 In order to fulfill the mass constraint \eqref{constraint2}, in Section \ref{Section9} we shall substitute the solution constructed in the above steps into it and then find the relation between the parameters $\iota$ and $\epsilon$, in such a way that $\iota$ can be settle down for fixed small $\epsilon$.

We conclude this section with additional remarks on Theorems \ref{theorem 1.1},  and a discussion of the difficulties in constructing higher-dimensional concentrating solutions.

\begin{remark}   {\em
The eigenfunction  corresponding to the positive eigenvalue  $\lambda_0$  of eigenvalue problem  \eqref{eigenvalue} represent unstable mode.
Similar to the situation in \cite{delPKowWei1, DWY},    the resonance phenomena arise from these following  eigenvalues near zero
\begin{align}
\Big(\frac{2j\pi\iota}{\ell}\Big)^2,
\qquad
\lambda_0\,-\,\Big(\frac{2j\pi\iota}{\ell}\Big)^2,
\qquad
j\in{\mathbb N}^+,
\label{eigenvalue-problem}
\end{align}  of the linearized problem
$$
\phi_{zz}+\phi_{ss}-\phi+p_\epsilon U_\epsilon^{p_\e-1}(s)\phi=0, \quad s\in {\mathbb R},\ z\in \big[0, \ell/\iota\big).
$$
In \cite{delPKowWei1}, in order to deal with the resonance, they introduced two free parameters $f$ and $e$ in the Lyapunov-Schmidt reduction procedure.
A novelty of this work  compared to \cite{delPKowWei1}  is that the  resonance effects are amplified.
First, it follows from the eigenvalues $\lambda_0\,-\,\Big(\frac{2j\pi\iota}{\ell}\Big)^2, j  \in{\mathbb N}^+$ that the  length of the concentrating curve significantly affects the resonance phenomena; furthermore, the  length is not uniformly bounded, thereby  making the problem more complicated.
On the other hand, the curve length $\ell$ may depend on $\iota$, and the term $\iota^2 W(0,\theta)$ appears in the operator $\mathbb F_\iota$ in \eqref{mathbbF}. Consequently, the eigenvalues $\Xi_{\iota, j}$ of \eqref{mathbbF-boundary} depend on $\iota$ and may approach zero, potentially causing resonance in the linear problem \eqref{problem-auxil}-\eqref{problem-auxil-boundary}.
To ensure the solvability of problems \eqref{liner hd}-\eqref{boundaryhatf-e}, obtain a good a priori estimate, and ultimately guarantee that the reduction procedure can be carried out smoothly, we impose conditions {\bf{(B1)}},  {\bf{(B2)}},  {\bf{(B3)}}.

}
\end{remark}

\begin{remark}   {\em
The curve length is not uniformly bounded  gives rise to difficulties in making a prior estimates to functions involving $\theta$ with $\theta \in [0, \ell). $
To resolve this issue, we perform the following transformation,
\begin{align*}
 \tau \,=\, \frac{2\pi}{\tilde\ell} \int_0^{\theta} \beta (\xi) \,{\mathrm d} \xi
\qquad\mbox{with}\quad
 {\tilde\ell} \,=\, \int_0^{\ell} \beta (\xi) \,{\mathrm d} \xi,
\end{align*}
leading to a new expression
 of the error function
 $S(w)$  with  $\tilde \ell$
  explicitly exhibited. This technique was employed by  authors  paper in \cite{DWY}.
  Furthermore, in contrast to \cite{DWY}, the nonlinearity  of equation  \eqref{equ power 6}  differs from that in \cite{DWY}. Consequently, we need to employ different transformations when   defining  the norms (e.g., $\|f\|_*$ and $\|e\|_{**}$ defined  in \eqref{constraint-f}, \eqref{constraint-e}, and the weighted Sobolev norm $\|\cdot\|_{H^2_*(\mathfrak S)} $ defined  in \eqref{norm h2 mathfrak S} ) and choosing the approximate solution (see \eqref{first-approximation}).
  We emphasize that the choice of the approximate solution and the norms, as specified above, is essential for the reduction procedure to be carried out successfully.

}
\end{remark}

\medskip\section{Some preliminaries and preparing works} \label{Section2}

\subsection{  The local forms of the equations in Fermi coordinates}\label{Section3.1}\

In the Fermi coordinates $(t, \theta)$ given in \eqref{fermi}, the metric of ${\mathbb R}^2$ can be locally written as
\begin{equation}\label{eqmetric}
\mathfrak{g} \,=\,{\mathrm d}t^2 \ +\ \big(1 + \kappa(\theta) t\big)^2 {\mathrm d}\theta^2.
\end{equation}
The Laplace-Beltrami operator has the form
\begin{equation}
\begin{split}
\Delta_\mathfrak{g}
& \,=\, \frac{1}{1+\kappa(\theta) t} \frac{\partial}{\partial t} \Bigg[ \big(1+\kappa(\theta) t\big) \frac{\partial}{\partial t} \Bigg]
\ +\
\frac{1}{1+\kappa(\theta) t}\frac{\partial}{\partial \theta} \Bigg[ \frac{1}{1+\kappa(\theta) t} \frac{\partial}{\partial \theta} \Bigg]
\\[2mm]
& \,=\, \frac{\partial^2}{\partial t^2}+ \frac{\kappa(\theta)}{1+\kappa(\theta) t} \frac{\partial }{\partial t}
\ +\
\frac{1}{\big(1+\kappa(\theta) t\big)^2 }\frac{\partial^2}{\partial \theta^2}
\ -\
\frac{ t \kappa'(\theta)}{\big(1+ \kappa(\theta) t\big)^3 }\frac{\partial}{\partial \theta}.
\end{split}
\end{equation}
It is trivial that equation \eqref{equ power 6} takes the local form
\begin{align}
&  \frac{\iota^2}{1+\kappa(\theta) t} \frac{\partial}{\partial t} \Bigg[ \big(1+\kappa(\theta) t\big) \frac{\partial}{\partial t} u\Bigg]
\ +\
\frac{\iota^2}{1+\kappa(\theta) t}\frac{\partial}{\partial \theta} \Bigg[ \frac{1}{1+\kappa(\theta) t} \frac{\partial}{\partial \theta} u\Bigg]
\nonumber\\[2mm]
&\ -\
\big[1+ \iota^2 W (t, \theta) \big] u
\ +\
 | u|^{p_\epsilon-2} u
 =0.
\label{eq3.1-0}
\end{align}

\medskip
 In the sequel, we will denote
\begin{equation}\label{eqDepsilon}
\tilde{\Gamma}_\epsilon \,=\, \iota^{-1}\Gamma_\iota,
\end{equation}
and set the natural stretched coordinates in the neighbourhood of the curve $\tilde{\Gamma}_\epsilon$,
\begin{equation}\label{eqszfinal}
(s, z) \,=\,
\iota^{-1} \big(t,\, \theta\big),
\end{equation}
which are defined for
\begin{equation}\label{eqFermi-stretched}
s\in \Big(- \frac{ M_0}{\iota}, \frac{ M_0}{\iota}\Big), \ \ \ z\in
\Big[0, \frac{\ell}{\iota} \Big).
\end{equation}
Hence, equation \eqref{eq3} has the local form
\begin{align} \label{eq3.1}
& \frac{1}{1+\iota \kappa s} \frac{\partial}{\partial s} \Big[ \big(1+\iota \kappa s\big) \frac{\partial}{\partial s} v \Big]
 +
 \frac{1}{1+\iota \kappa s}\frac{\partial}{\partial z} \Big[ \frac{1}{1+\iota \kappa s} \frac{\partial}{\partial z} v \Big]
 -
 \big[1+ \iota^2 W (\iota s, \iota z) \big] v +  |v|^{p_\epsilon-2} v
 =0,
\end{align}
which can also be written in the following way
\begin{align}
v_{ss} + v_{zz} + B_1 (v) - \big[1+ \iota^2 W (\iota s, \iota z) \big] v +   |v|^{p_\epsilon-2} v  \,=\, 0,
\label{equationlocal}
\end{align}
 where
\begin{align}
B_1(v)
& \,=\,
\frac{\iota \kappa (\iota z) }{1+\iota \kappa (\iota z) s} v_{s}
+\Big[ \frac{1}{(1+\iota \kappa(\iota z) s)^2 } -1 \Big]v_{zz}
+\frac{\iota^2 s \kappa' (\iota z) }{(1+\iota \kappa (\iota z) s)^3 }v_z.
\label{B1}
\end{align}

\medskip
In the neighborhood of $\tilde \Gamma_\iota$, by taking the Taylor expansion,
we have
\begin{align} \label{expansion of W}
1+ \iota^2 W (\iota s, \iota z)
 \,=\, & 1+ \iota^2 W (0, \iota z)+ \iota^3s W_t(0, \iota z) + \frac 12\iota^4s^2 W_{tt}(0, \iota z)
\nonumber\\[2mm]
&+ \frac 1 6 \iota^5s^3 W_{ttt}(0, \iota z)
+\frac{1}{24} \iota^6s^4 W_{tttt}\big(\hbar(\iota z)\iota s, \, \iota z\big),
\end{align}
where $\hbar(\theta)$ is a smooth function with $|\hbar(\theta)|\leq M_0, \, \forall \, \theta\in (0, \ell)$.
On the other hand, it is trivial that
\begin{align}
B_1(v) & \,=\, \big(\iota \kappa -\iota^2 \kappa^2 s \big) v_{s} + B_0(v),
\end{align}
where
\begin{align}
\label{B0v}
B_0 (v) \,=\, \iota^3 s^2 \kappa^3(\iota z) a_3(\iota s, \iota z)v_{s}
\,+\,
\iota s \kappa a_2(\iota s, \iota z) v_{zz}
\,+\,
\iota^2s \kappa'(\iota z) a_1(\iota s, \iota z) v_z,
\end{align}
with
\begin{align}
a_3(t, \theta)=\frac{1}{1+t \kappa (\theta) },
\label{a1a2a3}
\qquad
a_2(t, \theta)=\frac{-2-t \kappa (\theta)}{(1+t \kappa (\theta) )^2 },
\qquad
a_1(t, \theta)&=\frac{1}{(1+t \kappa (\theta) )^3 }.
\end{align}
Then we can rewrite $a_1, a_2, a_3$ as
\begin{align}
a_1(t, \theta) =  1- \iota \frac{3\kappa}{\beta}(x+f) + \iota^2  {\tilde a}_1,
\quad
&{\tilde a}_1=\frac{6\beta^2\kappa^2(x+f)^2 \,+\,8\iota\beta\kappa^3(x+f)^3
\,+\,3\iota^2\kappa^4(x+f)^4}{\beta^4\big[1+\iota\beta^{-1}\kappa(x+f)\big]^3}, \label{a1}
 \\[2mm]
a_2(t, \theta) =  -2 + \iota \frac{3\kappa}{\beta}(x+f) + \iota^2  {\tilde a}_2,
\quad
&{\tilde a}_2=\frac{-4\beta^{-2}\kappa^2(x+f)^2 \,-\,3\iota\beta^{-3}\kappa^3(x+f)^3}{\big[1+\iota\beta^{-1}\kappa(x+f)\big]^2}, \label{a2}
 \\[2mm]
a_3(t, \theta) =  1- \iota \frac{\kappa}{\beta}(x+f) + \iota^2 {\tilde a}_3,
\quad
&{\tilde a}_3=\frac{\beta^{-2}\kappa^2(x+f)^2}{1+\iota\beta^{-1}\kappa(x+f)}. \label{a3}
\end{align}
The assumptions in \eqref{curvature bounded} imply that
$a_1, a_2$ and $a_3$ are three smooth bounded functions on the following domain
$$
\Big\{\, (t, \theta)\, :\, |t|<M_0, \quad \theta\in (0, \ell)\, \Big\}
$$
provided that the constant $M_0$ in \eqref{fermi} is chosen sufficiently small.

By using the above calculations, we derive from \eqref{equationlocal} to get
\begin{align} \label{eq4}
{\mathcal D}(v):=\,&v_{ss} + v_{zz} - \big[1+ \iota^2 W (0, \iota z) \big] v +  |v|^{p_\epsilon-2} v
+\big(\iota \kappa -\iota^2 \kappa^2 s \big) v_{s}
+ B_0(v)
- \iota^3s W_t(0, \iota z) v
\nonumber\\[2mm]
 &  - \frac 12 \iota^4s^2 W_{tt}(0, \iota z) v- \frac 16 \iota^5s^3 W_{ttt}(0, \iota z)v
 -\frac{1}{24} \iota^6s^4 W_{tttt}\big(\hbar(\iota z)\iota s, \, \iota z\big)v \,=\, 0.
\end{align}

\medskip
\subsection{  Setting up near the curves}\label{Section3.2}\

\medskip
\noindent
{\textbf{Step 1.}}
We denote
\begin{align}  \label{beta}
\beta(\theta) = \big[1+ \iota^2 W (0, \iota z) \big]^{\frac 12},
\end{align}
and
  define a function ${\tilde v}({\tilde s}, z) $ by the relation
\begin{align} \label{coordinate-tildes}
v(s, z) \,=\, \beta^{\frac 2{p_\epsilon-2}}(\iota z) {\tilde v}({\tilde s}, z)
\qquad \text{with}\quad
{\tilde s} \,=\, \beta(\iota z)s,
\end{align}

We will express equation (\ref{eq4}) in terms of new coordinates $({\tilde s}, z)$ by the following preparations.
It's easy to get that
\begin{align*}
\p_z {\tilde s}
 \,=\,
 \iota\beta's \,=\, \iota\frac{\beta'}{\beta}{\tilde s},
\qquad
\p_{zz} {\tilde s} \,=\, \iota^2\beta''s
\,=\,
\iota^2\frac{\beta''}{\beta}{\tilde s}.
\end{align*}
We employ the notation
\begin{align}
 {\mathfrak{p}_\epsilon} =\frac 2{p_\epsilon-2}.
\end{align}
After that, there hold
\begin{align*}
v_{s}& \,=\, \beta^{{\mathfrak{p}_\epsilon}+1}{\tilde v}_{\tilde s},
\qquad\qquad
v_{ss} \,=\,  \beta^{{\mathfrak{p}_\epsilon}+2}{\tilde v}_{{\tilde s}{\tilde s}},
 \end{align*}
\begin{align*}
v_z  & \,=\, \iota \beta^{{\mathfrak{p}_\epsilon}-1}{\mathfrak{p}_\epsilon} \beta' {\tilde v}
+\iota  \beta^{{\mathfrak{p}_\epsilon}-1}  \beta' {\tilde s}{\tilde v}_{\tilde s}
+\beta^{{\mathfrak{p}_\epsilon}} {\tilde v}_z
\\
\nonumber
& =\beta^{{\mathfrak{p}_\epsilon}}\left[  \iota {\mathfrak{p}_\epsilon} \frac{\beta'}{\beta} {\tilde v}
+\iota\frac{\beta'}{\beta}   {\tilde s}{\tilde v}_{\tilde s} + {\tilde v}_z \right],
 \end{align*}
 and
\begin{align*}
v_{zz}&
\,=\,\iota \beta^{{\mathfrak{p}_\epsilon}-2}\beta'  ({\mathfrak{p}_\epsilon} -1)\left[  \iota{\mathfrak{p}_\epsilon}  \beta' {\tilde v}
+  \iota\beta' {\tilde s}{\tilde v}_{\tilde s} +\beta {\tilde v}_z \right]
\\
\nonumber
& \quad     +  \beta^{{\mathfrak{p}_\epsilon}-1} \left[ \iota^2
{\mathfrak{p}}_\epsilon\beta'' {\tilde v}
+  \iota^2( {\mathfrak{p}_\epsilon}+1)\,  \frac{|\beta'|^2}{\beta} {\tilde s} {\tilde v}_{\tilde s}
+ \iota (\mathfrak{p}_\epsilon+1) \beta' {\tilde v}_z
+ 2\iota \beta' {\tilde s} {\tilde v}_{z{\tilde s}}
+  \iota^2 \frac{|\beta'|^2}{\beta} {\tilde s}^2 {\tilde v}_{{\tilde s}{\tilde s}}
+  \iota^2 \beta'' {\tilde s} {\tilde v}_{\tilde s}
+ \beta {\tilde v}_{zz}\right].
\\
\nonumber
&
\,=\,  \beta^{{\mathfrak{p}_\epsilon}}\Bigg\{\iota^2  \Big[{\mathfrak{p}_\epsilon} ({\mathfrak{p}_\epsilon} -1)\Big(\frac{\beta'}{\beta}\Big)^2 + {\mathfrak{p}}\frac{\beta''}{\beta}\Big]{\tilde v}
+  \iota^2 \Big[2{\mathfrak{p}_\epsilon} \Big(\frac{\beta'}{\beta}\Big)^2 +\frac{\beta''}{\beta}\Big]  {\tilde s}{\tilde v}_{\tilde s}
\\
\nonumber
&
 \qquad   \qquad
  +  \iota^2 \Big(\frac{\beta'}{\beta}\Big)^2 {\tilde s}^2 {\tilde v}_{{\tilde s}{\tilde s}}
+ 2\iota \frac{\beta'}{\beta} {\tilde s} {\tilde v}_{z{\tilde s}}
+2\iota \mathfrak{p}_\epsilon \frac{\beta'}{\beta} {\tilde v}_z
+  {\tilde v}_{zz} \Bigg\} .
\end{align*}

Furthermore, for later use, we also compute
\begin{align}
\label{beta'}
\beta' (\theta)
\,=\,
\frac{ \iota^2 W_\theta(0, \theta)}{ 2 \beta(\theta)}
\,=\,
\frac{ \iota^2 W_\theta(0, \theta)}{ 2 \big[1+ \iota^2 W (0, \theta) \big]^{\frac 12}},
\end{align}
\begin{align}
\label{beta''}
\beta'' (\theta)
\,=\,
\frac{ \, \iota^2 W_{\theta \theta}(0, \theta)\, }{2\beta(\theta)}
\,-\,
\frac{\, \iota^4|W_{\theta}(0, \theta)|^2\, }{4\beta^3(\theta)}
\,=\,
\frac{ \, \iota^2 W_{\theta \theta}(0, \theta)\, }{2\big[1+ \iota^2 W (0, \theta) \big]^{\frac{1}{2}}}
\,-\,
\frac{\, \iota^4|W_{\theta}(0, \theta)|^2\, }{4\big[1+ \iota^2 W (0, \theta) \big]^{\frac 3 2}},
\end{align}
and
\begin{align}
\label{beta'''}
\beta'''(\theta)
\,=\,
\frac{ \, \iota^2 W_{\theta \theta \theta}(0, \theta)\, }{2\big[1+ \iota^2 W (0, \theta) \big]^{\frac{1}{2}}}
\,-\,\frac{ 3\iota^4 W_{\theta \theta}(0, \theta)\,W_{\theta }(0, \theta) }{4\big[1+ \iota^2 W (0, \theta) \big]^{\frac 3 2}}
\,+\,
\frac{\,3 \iota^6|W_{\theta}(0, \theta)|^3\, }{8\big[1+ \iota^2 W (0, \theta) \big]^{\frac 5 2}}.
\end{align}
The combination of the above formulas and the assumptions in \eqref{Wthetaestimate} will imply that
\begin{align}
|\beta'|\,\leq\, C\beta,
\qquad
|\beta''|\,\leq\, C\beta, \qquad |\beta'''|\,\leq\, C\beta.
\label{beta-bound}
\end{align}

\medskip
Finally,  dividing ${\mathcal D}(v)$ by $\beta^{ {\mathfrak{p}_\epsilon} +2}$, we conclude that problem \eqref{eq4} would be transformed to
\begin{align} \label{eq5}
{\tilde S}({\tilde v}) :\,=\, {\tilde v}_{{\tilde s}{\tilde s}} - {\tilde v}+ |\tilde v|^{p_\epsilon-2} \tilde v  +\frac{1}{\beta^2} {\tilde v}_{zz} + B_3({\tilde v}) \,=\, 0,
\end{align}
with
\begin{align} \label{B3w}
B_3({\tilde v})& \,=\,
  \Bigg[\iota \frac{\kappa}{\, \beta\, }
\,-\, \iota^2\frac{\kappa^2}{\beta^2} {\tilde s}\Bigg] {\tilde v}_{\tilde s}
+\iota^2  \Big[{\mathfrak{p}_\epsilon} ({\mathfrak{p}_\epsilon} -1)\frac{(\beta')^2}{\beta^4} + {\mathfrak{p}}\frac{\beta''}{\beta^3}\Big]{\tilde v}
+  \iota^2 \Big[2{\mathfrak{p}_\epsilon}\frac{(\beta')^2}{\beta^4} +\frac{\beta''}{\beta^3}\Big]  {\tilde s}{\tilde v}_{\tilde s}
\nonumber
\\
&\quad
  +  \iota^2 \frac{(\beta')^2}{\beta^4} {\tilde s}^2 {\tilde v}_{{\tilde s}{\tilde s}}
+ 2\iota \frac{\beta'}{\beta^3} {\tilde s} {\tilde v}_{z{\tilde s}}
+2\iota \mathfrak{p}_\epsilon \frac{\beta'}{\beta^3} {\tilde v}_z
 - \iota \frac{ \iota^2 W_t(0, \iota z)}{\beta^3} {\tilde s} {\tilde v}
-\frac{1}{2}\, \iota^2\frac{\iota^2 W_{tt}(0, \iota z)}{\beta^4} {\tilde s}^2 {\tilde v}
+ B_2 ({\tilde v}).
\end{align}
In the above formula, we have denoted
\begin{align} \label{B2w}
B_2 ({\tilde v}) \,=\,&
 - \frac{\iota^5}{6} \frac{\, \, W_{ttt}(0, \iota z)}{\beta^5} {\tilde s}^3{\tilde v}
 -\frac{\iota^6}{24} \frac{\, \, W_{tttt}\big(\hbar(\iota z)\iota s, \, \iota z\big)}{\beta^6} {\tilde s}^4{\tilde v}
 + \frac{1}{\,  \beta^{{\mathfrak{p}_\epsilon}+2}\, }B_0 (v)
\nonumber\\[2mm]
 \,=\, &
 - \frac{\iota^5}{6} \frac{\, \, W_{ttt}(0, \iota z)}{\beta^5} {\tilde s}^3{\tilde v}
 -\frac{\iota^6}{24} \frac{\, \, W_{tttt}\big(\hbar(\iota z)\iota s, \, \iota z\big)}{\beta^6} {\tilde s}^4{\tilde v}
+\iota^3  \kappa^3(\iota z) a_3(\iota s, \iota z) \frac{1}{\, \beta^3\, }{\tilde s}^2{\tilde v}_{\tilde s}
\nonumber\\[2mm]
&\,+\,
\iota \kappa(\iota z) a_2(\iota s, \iota z)
\Bigg\{\iota^2\Big[{\mathfrak{p}_\epsilon}(\mathfrak{p}_\epsilon-1) \frac{|\beta'|^2}{\beta^5}
+ {\mathfrak{p}_\epsilon}\frac{\beta''} {\beta^4} \Big]{\tilde s} {\tilde v}
+ \iota^2\Big[2{\mathfrak{p}_\epsilon}\,\frac{|\beta'|^2}{\beta^5} +\frac{\beta''} {\beta^4} \Big]  {\tilde s}^2 {\tilde v}_{\tilde s}
\nonumber\\[2mm]
&\qquad\qquad\qquad\qquad\quad
+  \iota^2 \frac{|\beta'|^2}{\beta^5} {\tilde s}^3{\tilde v}_{{\tilde s}{\tilde s}}
+ 2 \iota \frac{\beta'}{\beta^4} {\tilde s}^2 {\tilde v}_{z{\tilde s}}
+ 2\iota \mathfrak{p}_\epsilon\frac{\beta'} {\beta^4} {\tilde s}{\tilde v}_z
+\frac{1}{\, \beta^3\, } {\tilde s}{\tilde v}_{zz}
\Bigg\}
\nonumber\\[2mm]
&\,+\,
\iota^2 \kappa'(\iota z) a_1(\iota s, \iota z)
\Bigg[
\iota\mathfrak{p}_\epsilon\frac{\beta'}{\beta^4} {\tilde s}{\tilde v}
+ \iota\frac{\beta'}{\beta^4} {\tilde s}^2 {\tilde v}_{\tilde s}
+ \frac{1}{\, \beta^3\, } {\tilde s}{\tilde v}_z
\Bigg],
\end{align}
where the linear operator $B_0 (v)$ and the function $\hbar$ are defined in \eqref{B0v}, \eqref{expansion of W} respectively.
In the above, the terms $\beta'$ and $\beta''$ have been calculated in \eqref{beta'} and \eqref{beta''}.

\medskip
\noindent{\textbf{Step 2.}}
We take the transformation
\begin{align}
 \tau (z) \,=\, \frac{2\pi}{\tilde\ell} \int_0^{\iota z} \beta (\xi) \,{\mathrm d} \xi
\qquad\mbox{with}\quad
 {\tilde\ell} \,=\, \int_0^{\ell} \beta (\xi) \,{\mathrm d} \xi,
\label{ztau}
\end{align}
 and
\begin{align}\label{definition of hat v}
{\tilde v}({\tilde s}, z)
\,=\, {\hat v}({\tilde s}, \tau).
\end{align}
The above setting will imply that
\begin{align}
 \frac{{\mathrm d} \tau} {{\mathrm d}z}
 \,=\, \iota\frac{ 2\pi} {\tilde\ell} \beta (\iota z),
 \qquad
 \frac{{\mathrm d}^2 \tau} {{\mathrm d} z^2}
 \,=\, \iota^2\frac{ 2\pi} { {\tilde\ell}} \beta' (\iota z),
\end{align}
\begin{align}
{\tilde v}_{\tilde s}({\tilde s}, z) \,=\,{\hat v}_{\tilde s} ({\tilde s}, \tau),\qquad
{\tilde v}_{{\tilde s}{\tilde s}}({\tilde s}, z) \,=\,{\hat v}_{{\tilde s}{\tilde s}} ({\tilde s}, \tau),
\end{align}
\begin{align}
{\tilde v}_z ({\tilde s}, z)
\,=\,\iota\frac{2\pi}{{\tilde\ell}}\beta(\iota z){\hat v}_\tau ({\tilde s}, \tau),
\qquad
 {\tilde v}_{z{\tilde s}} ({\tilde s}, z)
\,=\, \iota\frac{2\pi}{{\tilde\ell}}\beta(\iota z){\hat v}_{\tau {\tilde s}} ({\tilde s}, \tau),
\end{align}
and
\begin{align}
{\tilde v}_{zz}({\tilde s}, z)
\,=\,
&\iota\partial_z\Big[\frac{2\pi}{{\tilde\ell}}\beta(\iota z){\hat v}_\tau ({\tilde s}, \tau)\Big]
\,=\,
\iota^2\frac{2\pi}{{\tilde\ell}}\beta'{\hat v}_\tau ({\tilde s}, \tau)
\,+\,\iota^2\frac{ 4\pi^2}{ {\tilde\ell^2}}\beta^2{\hat v}_{\tau\tau }({\tilde s}, \tau).
\end{align}

\medskip
We substitute the above results into \eqref{eq5} and then obtain
\begin{align}
0= {\tilde S}({\tilde v})
\,=\,&
{\tilde v}_{{\tilde s}{\tilde s}} - {\tilde v}+ |\tilde v|^{p_\epsilon-2} \tilde v
+\frac{1}{\beta^2} {\tilde v}_{zz} +  \Bigg[\iota \frac{\kappa}{\, \beta\, }
\,-\, \iota^2\frac{\kappa^2}{\beta^2} {\tilde s}\Bigg] {\tilde v}_{\tilde s}
\nonumber\\[2mm]
&
+\iota^2  \Big[{\mathfrak{p}_\epsilon} ({\mathfrak{p}_\epsilon} -1)\frac{(\beta')^2}{\beta^4} + {\mathfrak{p}}\frac{\beta''}{\beta^3}\Big]{\tilde v}
+  \iota^2 \Big[2{\mathfrak{p}_\epsilon}\frac{(\beta')^2}{\beta^4} +\frac{\beta''}{\beta^3}\Big]  {\tilde s}{\tilde v}_{\tilde s}
+  \iota^2 \frac{(\beta')^2}{\beta^4} {\tilde s}^2 {\tilde v}_{{\tilde s}{\tilde s}}
\nonumber \\[2mm]
&
+ 2\iota \frac{\beta'}{\beta^3} {\tilde s} {\tilde v}_{z{\tilde s}}
+2\iota \mathfrak{p}_\epsilon \frac{\beta'}{\beta^3} {\tilde v}_z
 - \iota \frac{ \iota^2 W_t(0, \iota z)}{\beta^3} {\tilde s} {\tilde v}
-\frac{1}{2}\, \iota^2\frac{\iota^2 W_{tt}(0, \iota z)}{\beta^4} {\tilde s}^2 {\tilde v}
+ B_2 ({\tilde v})
\nonumber\\[2mm]
\,=\,
&{\hat v}_{{\tilde s}{\tilde s}} - {\hat v} +   |\hat v|^{p_\epsilon-2} \hat v
+\iota^2\frac{1}{\beta^2}\Bigg[\frac{4\pi^2}{{\tilde\ell}^2} \beta^{2} {\hat v}_{\tau\tau}
+\frac{2\pi}{ {\tilde\ell}}\beta'{\hat v}_\tau
\Bigg]
+ \Big[\iota \frac{\kappa}{\, \beta\, }\,-\, \iota^2\frac{\kappa^2}{\beta^2} {\tilde s}\Big] {\hat v}_{\tilde s}
\nonumber\\[2mm]
& +\iota^2  \Big[{\mathfrak{p}_\epsilon} ({\mathfrak{p}_\epsilon} -1)\frac{(\beta')^2}{\beta^4} + {\mathfrak{p}}\frac{\beta''}{\beta^3}\Big]{\hat v}
+  \iota^2 \Big[2{\mathfrak{p}_\epsilon} \frac{(\beta')^2}{\beta^4} +\frac{\beta''}{\beta^3}\Big]  {\tilde s}{\hat v}_{\tilde s}
+  \iota^2 \frac{(\beta')^2}{\beta^4} {\tilde s}^2 {\hat  v}_{{\tilde s}{\tilde s}}
\nonumber \\[2mm]
&
+4\iota^2  \frac{\pi}{{\tilde\ell}} \frac{\beta'}{\beta^2} {\tilde s}{\hat v}_{\tau {\tilde s}}+4\iota^2\mathfrak{p}_\epsilon \frac{\pi}{{\tilde\ell}} \frac{\beta'}{\beta^2} {\hat v}_\tau
 - \iota \frac{ \iota^2 W_t(0, \iota z)}{\beta^3} {\tilde s} {\hat  v}
-\frac{1}{2}\, \iota^2\frac{\iota^2 W_{tt}(0, \iota z)}{\beta^4} {\tilde s}^2 {\hat  v}
+ B_2 ({\tilde v})
\nonumber\\[2mm]
\,=\,&{\hat v}_{{\tilde s}{\tilde s}} - {\hat v} +  |\hat v|^{p_\epsilon-2} \hat v
+\iota^2\frac{4\pi^2}{{\tilde\ell}^2}  {\hat v}_{\tau\tau}
\,+\,\iota \frac{\kappa}{\beta} \Bigg[ {\hat v}_{\tilde s} +\frac 1{\Lambda_\epsilon}{\tilde s}{\hat v}\Bigg]
+(4 \mathfrak{p}_\epsilon+2) \iota^2\frac{\pi}{{\tilde\ell}}\frac{\beta'} {\beta^2}{\hat v}_\tau
+4\iota^2  \frac{\pi}{{\tilde\ell}} \frac{\beta'}{\beta^2} {\tilde s}{\hat v}_{\tau {\tilde s}}
\nonumber\\[2mm]
&
+ \iota^2\left[ \mathfrak{p}_\epsilon(\mathfrak{p}_\epsilon-1)   \frac{|\beta'|^2} {\beta^4}  +{\mathfrak{p}_\epsilon} \frac{\beta''}{\beta^3} \right] {\hat v}
+\iota^2\Big[  \frac{\beta''}{\beta^3}\,+\,2{\mathfrak{p}_\epsilon} \frac{|\beta'|^2}{\beta^4} \,-\,\frac{\kappa^2}{\beta^2}\Big]{\tilde s} {\hat v}_{\tilde s}
\nonumber\\[2mm]
&
+ \iota^2 \frac{|\beta'|^2}{\beta^4}{\tilde s}^2 {\hat v}_{{\tilde s}{\tilde s}}
-\frac{1}{2}\, \iota^2\frac{\iota^2W_{tt}(0, \iota z)}{\beta^4} {\tilde s}^2 {\hat v}
+B_2 ({\tilde v}).
\label{eqn-hatv}
\end{align}
Here, we have used the stationary condition in \eqref{stationary assumption}.
In the coordinates $(\tilde{s}, \tau)$, the expression of  $B_2 ({\tilde v})$ in \eqref{eqn-hatv} has the form
\begin{align}
B_2 ({\tilde v})
\,=\, &
 - \frac{\iota^3}{6} \frac{\,\iota^2 W_{ttt}(0, \iota z)}{\beta^5} {\tilde s}^3{\hat v}
 -\frac{\iota^4}{24} \frac{\,\iota^2 W_{tttt}\big(\hbar(\iota z)\iota s, \, \iota z\big)}{\beta^6} {\tilde s}^4{\hat v}
 +\iota^3  \kappa^3(\iota z) a_3(\iota s, \iota z) \frac{1}{\, \beta^3\, }{\tilde s}^2{\hat v}_{\tilde s}
\nonumber\\[2mm]
&\,+\,
\iota \kappa(\iota z) a_2(\iota s, \iota z)
\Bigg\{\iota^2\Big[{\mathfrak{p}_\epsilon}(\mathfrak{p}_\epsilon-1) \frac{|\beta'|^2}{\beta^5}
+ {\mathfrak{p}_\epsilon}\frac{\beta''} {\beta^4} \Big]{\tilde s} {\hat v}
+ \iota^2\Big[2{\mathfrak{p}_\epsilon}\,\frac{|\beta'|^2}{\beta^5} +\frac{\beta''} {\beta^4} \Big]  {\tilde s}^2 {\tilde v}_{\hat s}
\nonumber\\[2mm]
&\qquad\qquad\qquad\qquad\quad
+  \iota^2 \frac{|\beta'|^2}{\beta^5} {\tilde s}^3{\hat v}_{{\tilde s}{\tilde s}}
+ 4 \iota^2 \frac{\pi}{{\tilde\ell}}\frac{\beta'}{\beta^3} {\tilde s}^2 {\hat v}_{\tau {\tilde s}}
+ 4\iota^2 \mathfrak{p}_\epsilon\frac{\pi}{{\tilde\ell}}\frac{\beta'} {\beta^3} {\tilde s}{\hat v}_\tau
\nonumber\\[2mm]
&\qquad\qquad\qquad\qquad\quad
+\iota^2\frac{2\pi}{{\tilde\ell}}\frac{\beta'}{\, \beta^3\, }{\tilde s} {\hat v}_\tau
\,+\,\iota^2\frac{ 4\pi^2}{ {\tilde\ell^2}}\frac{1}{\, \beta\, }{\tilde s}{\hat v}_{\tau\tau }
\Bigg\}
\nonumber\\[2mm]
&\,+\,\iota^3 \kappa'(\iota z)a_1(\iota s, \iota z)
\Bigg\{
\mathfrak{p}_\epsilon\frac{\beta'}{\beta^4} {\tilde s}{\hat v}
+ \frac{\beta'}{\beta^4} {\tilde s}^2 {\hat v}_{\tilde s}
+ \frac{1}{\, \beta^3\, } {\tilde s}\frac{ 2\pi}{ {\tilde\ell}}\beta{\hat v}_\tau
\Bigg\}
\nonumber\\[2mm]
\,=\, &
 - \frac{\iota^3}{6} \frac{\, \iota^2 W_{ttt}(0, \iota z)}{\beta^5} {\tilde s}^3{\hat v}
 -\frac{\iota^4}{24} \frac{\,\iota^2 W_{tttt}\big(\hbar(\iota z)\iota s, \, \iota z\big)}{\beta^6} {\tilde s}^4{\hat v}
 +\iota^3  \kappa^3(\iota z) a_3(\iota s, \iota z) \frac{1}{\, \beta^3\, }{\tilde s}^2{\tilde v}_{\tilde s}
\nonumber\\[2mm]
&\,+\,
\iota^3 \kappa(\iota z) a_2(\iota s, \iota z)
\Bigg\{
\left[ \mathfrak{p}_\epsilon\frac{\beta''} {\beta^4} +( \mathfrak{p}^2- \mathfrak{p}_\epsilon)  \frac{|\beta'|^2}{\beta^5} \right] {\tilde s} {\hat v}
+\, \left[ 2{\mathfrak{p}_\epsilon} \frac{|\beta'|^2}{\beta^5} +  \frac{\beta''}{\beta^4}  \right] {\tilde s}^2 {\hat v}_{\tilde s}
\nonumber\\[2mm]
&\qquad\qquad\qquad\qquad\quad
+\big(4\mathfrak{p}_\epsilon+2\big)\frac{ \pi}{ {\tilde\ell}}\frac{\beta'} {\beta^3}  {\tilde s}{\hat v}_\tau
+  4 \frac{\pi}{ {\tilde\ell}}\frac{\beta'} {\beta^3}  {\tilde s}^2{\hat v}_{\tau \tilde s }  + \frac{|\beta'|^2}{\beta^5} {\tilde s}^3{\hat v}_{{\tilde s}{\tilde s}}
+\frac{ 4\pi^2} {{\tilde\ell}^2}\frac{1}{\beta} {\tilde s}{\hat v}_{\tau\tau}\Bigg\}
\nonumber\\[2mm]
&\,+\,
\iota^3a_1(\iota s, \iota z) \kappa'(\iota z)
\Bigg\{
\mathfrak{p}_\epsilon\frac{\beta'}{\beta^4} {\tilde s}{\hat v}
+ \frac{\beta'}{\beta^4} {\tilde s}^2 {\hat v}_{\tilde s}
+\frac{ 2\pi}{ {\tilde\ell}}\frac{1}{ \beta^2} {\tilde s}{\hat v}_\tau
\Bigg\}.
\label{B2tildev}
\end{align}

\noindent{\textbf{Step 3.}}
We define a function $w (x, \tau) $ by the relation
\begin{align} \label{definition of x 2.14}
{\hat v}({\tilde s}, \tau) \,=\,  w(x, \tau)
\qquad \text{with}\quad
x \,=\, {\tilde s}-f(\tau)  = \beta s- f(\tau),
\end{align}
where $f :(0, 2\pi) \mapsto\R$ with the constraint \eqref{constraint-f} is a function  to be chosen later and satisfies the boundary conditions
$$
\quad f(0) \,=\, f(2\pi), \quad f'(0) \,=\, f'(2\pi).
$$

It's easy to get that
\begin{align*}
\p_\tau x & \,=\, -f',
\qquad
\p_{\tau\tau} x \,=\, -f''.
\end{align*}
After that, there hold
\begin{align*}
&{\hat v}_{{\tilde s}} \,=\, w_x,
\qquad\qquad
{\hat v}_{\tilde s \tilde s} \,=\, w_{xx},
\\[2mm]
&{\hat v}_\tau \,=\, - f' w_x  +w_\tau,
\qquad\qquad
{\hat v}_{{\tilde s}\tau} \,=\, - f' w_{xx}  +w_{x\tau},
\\[2mm]
&{\hat v}_{\tau\tau} \,=\,  |f'|^2 w_{xx}  -f''w_x +w_{\tau\tau} -2f' w_{x\tau}.
\end{align*}
Therefore, we write \eqref{eqn-hatv} in the form
\begin{align}
&\, w_{xx} - w + w^{{\mathfrak{p}_\epsilon} -1}
+ \iota^2 \frac{4\pi^2}{{\tilde\ell}^2} w_{\tau\tau}
+\iota^2\frac{4\pi^2}{{\tilde\ell}^2}\Big[|f'|^2 w_{xx}  -f''w_x  -2f' w_{x\tau}\Big]
\nonumber\\[2mm]
&
+\iota \frac{\kappa}{\, \beta\, } \Big[ w_x +{\frac 1{\Lambda_\epsilon}} (x+f)w\Big]
+\big(4{\mathfrak{p}_\epsilon} +2 \big) \iota^2\frac{\pi}{{\tilde\ell}}\frac{\beta'} {\beta^2}\big[- f' w_{x}  +w_{\tau}\big]
\nonumber\\[2mm]
&
+  4 \iota^2  \frac{\pi}{{\tilde\ell}} \frac{\beta'}{\beta^2} (x+f)\big[- f' w_{xx}  +w_{x\tau}\big]
+\iota^2 \Big[  \frac{\beta''}{\beta^3}\,+\, 2{\mathfrak{p}_\epsilon}  \frac{|\beta'|^2}{\beta^4} \,-\,\frac{\kappa^2}{\beta^2}\Big](x+f) w_{x}
\nonumber\\[2mm]
&+\iota^2 \left[{\mathfrak{p}_\epsilon} \frac{\beta''}{\beta^3} + \mathfrak{p}_\epsilon(\mathfrak{p}_\epsilon-1)   \frac{|\beta'|^2} {\beta^4}  \right]w
+ \iota^2  \frac{|\beta'|^2}{\beta^4}(x+f)^2 w_{xx}
\nonumber\\[2mm]
&-\frac{1}{2}\, \iota^2\frac{\iota^2W_{tt}(0, \iota z)}{\beta^4} (x+f)^2 w
+B_2 ({\tilde v})
=0,
\label{S(w)}
\end{align}
where $(x, \tau)\in{\mathbb R}\times(0, 2\pi)$.

\medskip
For the simplification of notation in \eqref{S(w)}, we denote
\begin{align}   \label{new opreator B4}
 B_4(w)
 =\iota \frac{\kappa}{\, \beta\, } \Big[ w_x +\frac 1{\Lambda_\epsilon} (x+f)w\Big],
\end{align}
and
\begin{align}   \label{new opreator B5}
 B_5(w)
\,=\,&\iota^2\frac{4\pi^2}{{\tilde\ell}^2}\Big[|f'|^2 w_{xx}  -f''w_x  -2f' w_{x\tau}\Big]
+\iota^2\big(4{\mathfrak{p}_\epsilon} +2 \big)\frac{\pi}{{\tilde\ell}}\frac{\beta'} {\beta^2} \big[- f' w_{x}  +w_{\tau}\big]
\nonumber\\[2mm]
&
+ \iota^2  4  \frac{\pi}{{\tilde\ell}} \frac{\beta'}{\beta^2} (x+f)\big[- f' w_{xx}  +w_{x\tau}\big]
+
 \iota^2 \Big[  \frac{\beta''}{\beta^3}\,+\, 2 {\mathfrak{p}_\epsilon}  \frac{|\beta'|^2}{\beta^4} \,-\,\frac{\kappa^2}{\beta^2}\Big](x+f) w_{x}
\nonumber\\[2mm]
&+\iota^2\left[{\mathfrak{p}_\epsilon} \frac{\beta''}{\beta^3} + \mathfrak{p}_\epsilon(\mathfrak{p}_\epsilon-1)   \frac{|\beta'|^2} {\beta^4}  \right]w
+ \iota^2 \frac{|\beta'|^2}{\beta^4}(x+f)^2 w_{xx}
-\frac{1}{2}\, \iota^2\frac{\iota^2W_{tt}(0, \iota z)}{\beta^4} (x+f)^2 w.
\end{align}
For the last term in \eqref{S(w)}, by \eqref{B2tildev} and the relation $\hat v (\tilde s, \tau)=w(x, \tau)$, we derive
\begin{align}
 B_2 ({\tilde v})
 \,=\,&  - \frac{\iota^3}{6} \frac{\, \, \iota^2 W_{ttt}(0, \iota z)}{\beta^5} (x+f)^3 w
 -\frac{\iota^4}{24} \frac{\, \, \iota^2 W_{tttt}\big(\hbar(\iota z)\iota s, \, \iota z\big)}{\beta^6} (x+f)^4w
\nonumber\\[2mm]
&
 +\iota^3  \kappa^3(\iota z) a_3(\iota s, \iota z) \frac{1}{\, \beta^3\, }(x+f)^2w_{x}
\nonumber\\[2mm]
&\,+\,
\iota^3 \kappa(\iota z) a_2(\iota s, \iota z)
\Bigg\{
 \left[ \mathfrak{p}_\epsilon\frac{\beta''} {\beta^4} +( \mathfrak{p}^2- \mathfrak{p}_\epsilon)  \frac{|\beta'|^2}{\beta^5} \right] (x+f) w
+ \left[ 2{\mathfrak{p}_\epsilon} \frac{|\beta'|^2}{\beta^5} +  \frac{\beta''}{\beta^4}  \right](x+f)^2 w_x
\nonumber\\[2mm]
&\qquad\qquad\qquad\qquad\qquad
+ \big[4\mathfrak{p}_\epsilon +2\big] \frac{\pi}{ {\tilde\ell}}\frac{\beta'} {\beta^3}  (x+f)\big[-f'w_x +w_\tau\big]
\nonumber\\[2mm]
&\qquad\qquad\qquad\qquad\qquad
+ 4\frac{\pi}{ {\tilde\ell}}\frac{\beta'} {\beta^3}  (x+f)^2[- f' w_{xx}  +w_{x\tau}]
+  \frac{|\beta'|^2}{\beta^5} (x+f)^3w_{xx}
\nonumber\\[2mm]
&\qquad\qquad\qquad\qquad\qquad
+\frac{ 4\pi^2} {{\tilde\ell}^2}\frac{1}{\beta} (x+f)\big[ |f'|^2 w_{xx}  -f''w_x +w_{\tau\tau} -2f' w_{x\tau}\big]
\Bigg\}
\nonumber\\[2mm]
&\,+\,
\iota^3 \kappa'(\iota z)a_1(\iota s, \iota z)
\Bigg\{
\mathfrak{p}_\epsilon\frac{\beta'}{\beta^4} (x+f)w
+  \frac{\beta'}{\beta^4} (x+f)^2 w_{x}
+\frac{ 2\pi}{ {\tilde\ell}}\frac{1}{ \beta^2} (x+f)\big[-f'w_x +w_\tau\big]
\Bigg\}
\nonumber\\[2mm]
: \,=\,&  B_6(w).
\label{B2w-2222new}
\end{align}
We finally use \eqref{S(w)} to obtain the equation
\begin{align}  \label{expression in xtau}
S(w):=&\, \frac{\iota^2 4\pi^2}{{\tilde\ell}^2} w_{\tau\tau} + w_{xx} - w +|w|^{p_\epsilon-2} w  + B_4(w)+ B_5(w)+ B_6(w)=0,
\end{align}
where $(x, \tau)\in{\mathbb R}\times(0, 2\pi)$.

\subsection{ The block solution $U_\epsilon$ and conventions}\label{section3-3}\

We provide the following known facts about the non-linear Schr\"odinger equation
\begin{align} \label{Schrodiner Profile}
{U_\epsilon}'' - U_\epsilon + U_\epsilon^{p_\e -1}= 0 \quad\mbox{and}\quad U_\epsilon>0\quad \text{for}\,x\in \R,
\qquad {U_\epsilon}'(0) =0, \qquad U_\epsilon(\pm \infty) =0.
\end{align}
It is well-known that $U_\epsilon$ is a unique solution, which is even in $x$, and also enjoys the asymptotic behaviours
\begin{align}
 \lim_{|x| \rightarrow +\infty} U_\epsilon(x) e^{|x|} = C<+\infty,
\qquad\text{and} \qquad
\lim_{|x| \rightarrow +\infty} \frac{U_\epsilon(x)}{{U_\epsilon}'(x)}=-1.
\label{U-asymp}
\end{align}
Moreover, $U_\epsilon$ satisfies the non-degeneracy condition, in the sense that the kernel of the linear operator $-\Delta + 1 -(p_\epsilon-1)U_{\epsilon}^{p_\epsilon-2} $ in $H^1(\mathbb{R}) $ is spanned by the function ${U_\epsilon}'(x)$.
For the associated linearized eigenvalue problem
\begin{align}
\label{eigenvalue}
h{''}-h+(p_\e-1){U_\epsilon}^{p_\epsilon-2}h\,=\, \lambda h\quad {\rm in}\,\Bbb R, \,\quad h(\pm \infty)\,=\, 0,
\end{align}
it is also well-known that it possesses a unique positive eigenvalue $\lambda_0$ in $H^1(\Bbb R)$,
with associated eigenfunction ${\mathcal Z}$ (even and positive) which can be normalized so that
$$
\int_{\R}{\mathcal Z}^2(x) \,{\mathrm d}x=1.
$$
In fact, a simple computation shows that
\begin{equation}
\label{lambda0}
\lambda_0\,=\, \frac14(p_\epsilon-2) (p_\epsilon +2),
\qquad
\mathcal Z(x)\,=\, \frac{{(U_\e)}^{\frac{p_\epsilon}2}(x)}{\sqrt{\displaystyle\int_{\mathbb R} {(U_\e)}^{p_\epsilon} \,{\mathrm d}x}\, }.
\end{equation}

Suppose $U$ satisfies
\begin{align}  \label{Ue}
{U}'' - U + {U}^{p-1}= 0,  \text{in}~ \R^2, \quad \text{with}~p>2.
\end{align}
Multiplying $U$  on  \eqref{Ue} and integrating, we have
\begin{align}   \label{eq1-}
\int_\R |U'|^2 + \int_\R|U|^2 = \int_\R|U|^{p}.
\end{align}
Similarly,  multiplying $r U'$   and  integrating, we have
\begin{align}  \label{eq2-}
-\frac12\int_\R |U'|^2 +\frac 12 \int_\R|U|^2 =  \frac 1{p} \int_\R|U|^{p}.
\end{align}
By \eqref{eq1-}-\eqref{eq2-}, we have
\begin{align}
\Big(1+  \frac {p}2\Big) \int_\R |U'|^2 + \Big(  1 -  \frac {p}2 \Big)\int_\R|U|^2   = 0,
\end{align}
i.e,
\begin{align}   \label{identity  for sigma}
 \int_\R |U'|^2 = \frac {p-2}{p+2} \int_\R |U|^2 \, =\,  - \frac {2 (p-2)}{p+2}    \int_\R x UU' = -\frac1{ \left(\frac p {p-2} - \frac12\right)}  \int_\R x UU'.
\end{align}
 In particular, if $p = 6\pm \e$, we have
\begin{align}
 \int_\R |U_\epsilon'|^2 = \frac {4\pm \e}{8\pm \e} \int_\R |U_\epsilon|^2  =  (-2)  \frac {4\pm \e}{8\pm \e}    \int_\R x U_\epsilon U_\epsilon' = - \frac 1{\Lambda_\e}  \int_\R x U_\epsilon U_\epsilon'.
 \end{align}

\medskip
We end this section by providing a collection of notation.
{Throughout the paper, ${\mathfrak S}$ represents the strip in $\mathbb{R}^2$ of the form
\begin{equation}
\label{huaxies}
{\mathfrak S}=\Bigg\{ (x, \tau): \quad x\in \R, \quad 0<\tau< 2\pi\Bigg\}.
\end{equation}
For the introduction of the parameters $f$ and $e$, the reader can refer to \eqref{definition of x 2.14}
and \eqref{first-approximation}.
In all what follows, we will make the following assumptions on the parameters $f$ and $e$
\begin{align} \label{constraint-f}
\|f\|_* :\,=\,
\left\|\beta^{\mathfrak p}   f\right\|_{L^2(0, 2\pi)}
+\left\|\frac{\beta^{\mathfrak p}}{\tilde\ell} f_\tau\right\|_{L^2(0, 2\pi)}
+\left\|\frac{\beta^{\mathfrak p}}{{\tilde\ell}^2}  f_{\tau\tau}\right\|_{L^2(0, 2\pi)}
\leq \epsilon^{1-\alpha},
\end{align}
\begin{align} \label{constraint-e}
\|e\|_{**} : \,=\,
\|e\|_{L^2(0,  2\pi)}
 + \frac{\iota}{{\tilde\ell}}  \|e'\|_{L^2 (0, 2\pi)}
 + \frac{\iota^2} {{\tilde\ell}^2}\|e''\|_{L^2(0, 2\pi)}
 \leq \iota^{\frac 32 \varrho + \frac 12- 2\alpha}.
\end{align}
By the norm defined in \eqref{constraint-f}, we have
\begin{align} \label{Linfty of f}
\left\|\beta^{\mathfrak p}  f\right\|_{L^\infty(0, 2\pi)}  \le {\tilde\ell}^{\frac 12}  \|f\|_*, \quad  \left\|\frac{\beta^{\mathfrak p}}{\tilde\ell} f_\tau\right\|_{L^\infty(0, 2\pi)} \le  {\tilde\ell}^{\frac 12}  \|f\|_*.
\end{align}
The proof of \eqref{Linfty of f} can be found in \cite{DWY}.  We omit the proof here for concise.

For simplicity, we also denote
{
{\it
\begin{align}
\label{mathcalF}
\hat{\mathcal F}\,=\, \Big\{\, (f, e)\, :\, \mbox{ the functions } f \mbox{ and } e \mbox{ satisfy } (\ref{constraint-f})
\mbox{ and } (\ref{constraint-e}) \mbox{ respectively}\, \Big\}.
\end{align}
}
By \eqref{constraint-alpha-varrho} and combining  the conditions  $\alpha, \varrho \in (0, 1)$,  we conclude that
\begin{align}
\|e\|_{**}  \ll   \epsilon^{1-\alpha}.
\end{align}

\medskip\section{The approximate solutions}\label{section-approximatesolution}

\subsection{  The first approximate solution}
\subsubsection{First approximate solution and its error}
In the following,  for simplicity,  we will drop the subscript  and supscript
  from
 $U_\epsilon, \mathfrak{p}_\epsilon, \Lambda_\epsilon$  and denote  these  by  $U, \mathfrak{p}, \Lambda$ respectively.
We take
\begin{align}
w_0(x, \tau)=U(x)+\iota \frac {e(\tau)}{\beta^{\mathfrak{p}} (\iota z)} \mathcal{Z}(x),
\qquad (x, \tau)\in{\mathbb R}\times(0, 2\pi)
\label{first-approximation}
\end{align}
as {\textbf{the first approximate solution}},  with $\mathfrak{p} =\frac2{p_\epsilon-2}$.
The parameter $e$ with the constraint \eqref{constraint-e} will be chosen by the reduction method in section \ref{section10}.

It is obvious that
\begin{align}
\Big(\frac{e(\tau)}{\beta^{\mathfrak{p}}(\iota z)} {\mathcal Z}(x)\Big)_x= \frac {e(\tau)}{\beta^{\mathfrak{p}}(\iota z)}{\mathcal Z}',
\qquad\qquad
\Big(\frac{e(\tau)}{\beta^{\mathfrak{p}}(\iota z)} {\mathcal Z}(x)\Big)_{xx}=\frac {e(\tau)}{\beta^{\mathfrak{p}}(\iota z)} {\mathcal Z}'',
\end{align}
\begin{align}
\left(\frac{e(\tau)}{\beta^{\mathfrak{p}}(\iota z)} {\mathcal Z}(x)\right)_{\tau}
= \left(e'-\mathfrak{p} \frac{\beta'} {\beta^{2}} \frac{\tilde \ell} {2\pi } e\right)\frac{{\mathcal Z}}{\beta^{\mathfrak{p}}},
\end{align}
\begin{align}
 \left(\frac{e(\tau)}{\beta^{\mathfrak{p}}(\iota z)} {\mathcal Z}(x)\right)_{x\tau}
=\left(e'-\mathfrak{p} \frac{\beta'} {\beta^{2}} \frac{\tilde \ell} {2\pi } e\right)\frac{{\mathcal Z}'}{\beta^{\mathfrak{p}}},
\end{align}
\begin{align}
 \left(\frac{e(\tau)}{\beta^{\mathfrak{p}}(\iota z)} {\mathcal Z}(x)\right)_{\tau\tau}
 = \left[e''- 2\mathfrak{p}\frac{\beta'} {\beta^{2} } \frac{\tilde \ell} {2\pi } e' +\left(- \mathfrak{p}\frac{\beta''} {\beta^{3} } + \mathfrak{p}( {\mathfrak{p}+2})\frac{|\beta'|^2} {\beta^{4} } \right) \frac{\tilde \ell^2} {4\pi^2 }e\right] \frac{{\mathcal Z}}{\beta^{\mathfrak{p}}}.
\end{align}
Therefore, there holds
\begin{align}
S(w_0)=&\, w_{0,xx} - w_{0} + w_{0}^{p_\e-1}
+ \iota^2\frac{4\pi^2}{{\tilde\ell}^2} w_{0,\tau\tau}
+ B_4(w_0)+B_5(w_0)+ B_6\big(w_0\big)
\nonumber\\[2mm]
=&\,
\iota^3\frac{4\pi^2}{{\tilde\ell}^2}  \big( \frac {e}{\beta^{\mathfrak{p}}} \big)''\mathcal Z
 +\frac{\iota\lambda_0 e{\mathcal Z}}{\beta^{\mathfrak{p}}}  +  \iota^2 \frac{(p_\epsilon-1) (p_\epsilon-2)}{2\beta^{2\mathfrak{p}}}e^2U^{p_\epsilon- 3} {\mathcal Z}^2
\nonumber\\[2mm]
& + B_4\big(w_0\big)
+ B_5\big(w_0\big)
+ B_6\big(w_0\big) + N_{0\e} \left(\iota\frac{e} {\beta^{\mathfrak{p}}} \mathcal Z\right)
\label{S(w_0)}
\nonumber\\[2mm]
=&\, \iota^3\frac{4\pi^2}{{\tilde\ell}^2}  \left[e''- 2\mathfrak{p}\frac{\beta'} {\beta^2} \frac{\tilde \ell} {2\pi} e' +\left(-\mathfrak{p} \frac{\beta''} {\beta^3} + \mathfrak{p}(\mathfrak{p}+2) \frac{|\beta'|^2} {\beta^4} \right) \frac{\tilde \ell^2} {4\pi^2}e\right]\frac{\mathcal Z}{\beta^{\mathfrak{p}}}
\nonumber\\[2mm]
&  +\frac{\iota\lambda_0 e{\mathcal Z}}{\beta^{\mathfrak{p}}} +\iota^2 \frac{(p_\epsilon-1) (p_\epsilon-2)}{2\beta^{2\mathfrak{p}}}e^2U^{p_\epsilon- 3} {\mathcal Z}^2
 +N_{0\e} \left(\iota\frac{e} {\beta^{\mathfrak{p}}} \mathcal Z\right)
\nonumber\\[2mm]
& + B_4\big(w_0\big)
+ B_5\big(U\big) + B_5 \left(\iota\frac{e} {\beta^{\mathfrak{p}}} \mathcal Z\right)
+ B_6\big(U\big) + B_6 \left(\iota\frac{e} {\beta^{\mathfrak{p}}} \mathcal Z\right),
\end{align}
where
\begin{align}   \label{N0epsilonZ}
N_{0\e} \left(\iota\frac{e} {\beta^{\mathfrak{p}}} \mathcal Z\right)
=  \left( U + \iota\frac{e} {\beta^{\mathfrak{p}}} \mathcal Z\right)^{p_\epsilon-1}
-  U^{p_\epsilon-1}  - \iota (p_\epsilon-1) U^{p_\epsilon-2} \frac{e} {\beta^{\mathfrak{p}}} \mathcal Z
-  \iota^2\frac{(p_\epsilon-1) (p_\epsilon-2)}{2\beta^{2\mathfrak{p}}}e^2U^{p_\epsilon- 3} {\mathcal Z}^2.
\end{align}
The last five terms in \eqref{S(w_0)} will be calculated on details.

\smallskip
\noindent$\clubsuit$
From the expressions of $B_4$ and $B_5$ given in \eqref{new opreator B4}-\eqref{new opreator B5}, we have
\begin{align}
B_4(w_0)\,=\,
\iota\frac{\kappa}{\, \beta\, } \Big[ U'+\frac 1{\Lambda} (x+f)  U\Big]
+\iota^2 \frac{\kappa}{\, \beta^{1+ \mathfrak{p}}\, } \left[  e {\mathcal Z}' +\frac 1{\Lambda} (x+f) e{\mathcal Z}\right],
\end{align}
\begin{align}
B_5(U)
= &
-\iota^2\frac{\iota^2W_{tt}(0, \iota z)}{\beta^4}f xU
+\iota^2\Big[\frac{|\beta'|^2}{\beta^4} 2f
-\frac{4\pi}{{\tilde\ell}}\frac{\beta'}{\beta^2}f'
\Big]xU''
\nonumber\\[2mm]
&
+\iota^2\Big[
-\frac{4\pi^2}{{\tilde\ell}^2}f''
-(4\mathfrak{p} + 2) \frac{\pi}{\tilde\ell}  \frac{\beta'}{\beta^2} f'
+2\mathfrak{p}\frac{|\beta'|^2}{\beta^4}f
+\frac{\beta''} {\beta^3}f
\,-\,\frac{\kappa^2}{\beta^2}f
\Big]U'
\nonumber\\[2mm]
&
+ \iota^2\Big[
-\frac{1}{2}\, \frac{\iota^2 W_{tt}(0, \iota z)}{\beta^4}f^2  + \mathfrak{p}  \frac{\beta''} {\beta^3} + \mathfrak{p} (\mathfrak{p}-1)  \frac{|\beta'|^2} {\beta_\iota^4}
\Big]U
\nonumber\\[2mm]
&+\iota^2\Big[\frac{4\pi^2}{{\tilde\ell}^2}|f'|^2
-\frac{4\pi}{{\tilde\ell}}\frac{\beta'}{\beta^2}ff'
+\frac{|\beta'|^2}{\beta^4} f^2
\Big]U''
-\iota^2\frac{1}{2}\,  \frac{\iota^2W_{tt}(0, \iota z)}{\beta^4} x^2U
\nonumber\\[2mm]
&
+\iota^2\Big[2\mathfrak{p}\frac{|\beta'|^2}{\beta^4}
+\frac{\beta''} {\beta^3}
\,-\,\frac{\kappa^2}{\beta^2}
+\iota^2 \frac{|\beta'|^2}{\beta^4} x^2 U'',
\end{align}
and
\begin{align}
 B_5 \left(\iota\frac{e} {\beta^{\mathfrak{p}}} \mathcal Z\right)
\,=\,&
\iota^3\frac{4\pi^2}{{\tilde\ell}^2}
\left[
|f'|^2  \frac{e}{\beta^{\mathfrak{p}}} {\mathcal Z}''
-f'' \frac{e}{\beta^{\mathfrak{p}}} {\mathcal Z}'
-2f'   \left(e'-\mathfrak{p} \frac{\beta'} {\beta^{2}} \frac{\tilde \ell} {2\pi } e\right)\frac{{\mathcal Z}'}{\beta^{\mathfrak{p}}}
\right]
\nonumber\\[2mm]
&
+\iota^3\big[ 4\mathfrak{p} +2 \big] \frac{\pi}{ {\tilde\ell}}\frac{\beta'} {\beta^2}
\left[- f'  \frac{e}{\beta^{\mathfrak{p}}} {\mathcal Z}'
+ \left(e'-\mathfrak{p} \frac{\beta'} {\beta^{2}} \frac{\tilde \ell} {2\pi } e\right)\frac{{\mathcal Z}}{\beta^{\mathfrak{p}}}  \right]
\nonumber\\[2mm]
&
+\iota^3\frac{4\pi}{{\tilde\ell}} \frac{\beta'}{\beta^2} (x+f)
\left[- f'  \frac{e}{\beta^{\mathfrak{p}}} {\mathcal Z}''
 \left(e'-\mathfrak{p} \frac{\beta'} {\beta^{2}} \frac{\tilde \ell} {2\pi } e\right)\frac{{\mathcal Z}'}{\beta^{\mathfrak{p}}}  \right]
\nonumber\\[2mm]
&+
\iota^3\Big( \frac{\beta''}{\beta^3} \,+\, 2\mathfrak{p} \frac{|\beta'|^2}{\beta^4} \,-\,\frac{\kappa^2}{\beta^2} \Big) (x+f)  \frac{e}{\beta^{\mathfrak{p}}} {\mathcal Z}'
+\iota^3\mathfrak{p}\left[ \frac{\beta''}{\beta^3} + (\mathfrak{p}-1) \frac{|\beta'|^2} {\beta^4}  \right] \frac{e}{\beta^{\mathfrak{p}}}{\mathcal Z}
\nonumber\\[2mm]
&+\iota^3 \frac{|\beta'|^2}{\beta^4}(x^2+2xf+f^2)  \frac{e}{\beta^{\mathfrak{p}}} {\mathcal Z}''
  -\iota^3\frac{1}{2}\,\frac{\iota^2W_{tt}(0, \iota z)}{\beta^4} (x^2+2xf+f^2)  \frac{e}{\beta^{\mathfrak{p}}}{\mathcal Z}.
\end{align}

\smallskip
\noindent$\clubsuit$
On the other hand, the expression of $B_6$ given in \eqref{B2w-2222new} will give
\begin{align}
B_6(U)
\,=\,&  - \frac{\iota^3}{6} \frac{\, \iota^2 W_{ttt}(0, \iota z)}{\beta^5} (x+f)^3 U
 -\frac{\iota^4}{24} \frac{\,\iota^2 W_{tttt}\big(\hbar(\iota z)\iota s, \, \iota z\big)}{\beta^6} (x+f)^4U
\nonumber\\[2mm]
&
\,+\,\iota^3  \kappa^3(\iota z) a_3(\iota s, \iota z)
\frac{1}{\, \beta^3\, }(x^2U'+2fxU'+f^2U')
\nonumber\\[2mm]
&\,+\,
\iota^3 \kappa(\iota z) a_2(\iota s, \iota z)
\Bigg\{
\left[ \mathfrak{p}\frac{\beta''} {\beta^4}  + (\mathfrak{p}^2- \mathfrak{p})  \frac{|\beta'|^2} {\beta^5} \right] (xU+fU)
+ \left[ 2\mathfrak{p}\frac{|\beta'|^2}{\beta^5} + \frac{\beta''}{\beta^4}  \right] (x^2U'+2fxU'+f^2U')
\nonumber\\[2mm]
&\qquad\qquad\qquad\qquad\qquad
+(4\mathfrak{p}+2) \frac{\pi}{{\tilde\ell}}\frac{\beta'} {\beta^3} \big[- f' xU'  - ff' U'  \big]
- \frac{4\pi}{{\tilde\ell}} \frac{\beta'} {\beta^3} f' x^2U''
-\frac{ 8\pi}{ {\tilde\ell}}\frac{\beta'} {\beta^3} f' fxU''
\nonumber\\[2mm]
&\qquad\qquad\qquad\qquad\qquad
-\frac{4\pi}{{\tilde\ell}} \frac{\beta'} {\beta^3}  f' f^2U''
+ \frac{|\beta'|^2}{\beta^5} \big[x^3U''+3fx^2U''+3f^2xU''+f^3U''\big]
\nonumber\\[2mm]
&\qquad\qquad\qquad\qquad\qquad
+\frac{4\pi^2}{{\tilde\ell}^2}\frac{1}{\beta} \big[ |f'|^2 xU''  -f''xU'\big]
+\frac{4\pi^2}{{\tilde\ell}^2}\frac{1}{\beta} \big[ f|f'|^2 U''  -ff''U' \big]
\Bigg\}
\nonumber\\[2mm]
&\,+\,
\iota^3 \kappa'(\iota z)a_1(\iota s, \iota z)
\Bigg\{
\mathfrak{p}\frac{\beta'}{\beta^4} (xU+fU)
+\frac{\beta'}{\beta^4} (x^2U'+2fxU'+f^2U')
\nonumber\\[2mm]
&\qquad\qquad\qquad\qquad\qquad
-\frac{2\pi}{{\tilde\ell}}\frac{1}{ \beta^2}f' xU'
-\frac{2\pi}{{\tilde\ell}}\frac{1}{ \beta^2} ff' U'
\Bigg\}.
\end{align}
  The expressions of $a_1, a_2, a_3$ can be found in \eqref{a1}, \eqref{a2}, \eqref{a3}.
Then we can  rewritten $B_6(U)$  as
\begin{align}
B_6(U)= &   - \iota^3 \frac{\,\iota^2 W_{ttt}(0, \iota z)}{6\beta^5}\, x^3 U
\ -\
\iota^3\frac{\, \iota^2W_{ttt}(0, \iota z)}{2\beta^5}f\,x^2 U
\nonumber\\[2mm]
&
\ +\ \iota^3\Big[- \frac{\, \, \iota^2W_{ttt}(0, \iota z)}{2\beta^5}f^2 - 2\kappa \big[\mathfrak{p}\frac{\beta''}{\beta^4}   +  (\mathfrak{p}^2-\mathfrak{p})\frac{|\beta'|^2}{\beta^5}\big]   +\mathfrak{p}\frac{\kappa'\beta'} {\beta^4} \Big] x U
\nonumber\\[2mm]
& + \iota^3\Big[-\frac{\, \, \iota^2 W_{ttt}(0, \iota z)}{6\beta^5}f^3 -2\kappa\big[\mathfrak{p}\frac{\beta''}{\beta^4} +  (\mathfrak{p}^2-\mathfrak{p})\frac{|\beta'|^2}{\beta^5}  \big]   f  + \mathfrak{p} \frac{\kappa'\beta'} {\beta^4} f\Big] U
\nonumber\\[2mm]
 &
\ +\ \iota^3\Big[ \frac{\kappa^3}{\beta^3}-4 \mathfrak{p} \frac{\kappa|\beta'|^2}{\beta^5}  -\frac{2\kappa\beta''}{\beta^4} + \frac{\kappa'\beta'} {\beta^4} \Big]x^2U'
\nonumber\\[2mm]
 & + \iota^3\Big[
 \frac{2\kappa^3}{\beta^3}f
 -8\mathfrak{p}\frac{\kappa|\beta'|^2}{\beta^5}  f
 -\frac{4\kappa\beta''}{\beta^4}f
 -
 +\frac{8\pi^2}{{\tilde\ell}^2}\frac{\kappa}{\beta}   f''
+(8\mathfrak{p} +4) \frac{\pi }{ {\tilde\ell}}\frac{\kappa\beta'}{\beta^3} f'
+\frac{2\kappa'\beta'}{\beta^4}f
- \frac{2\pi}{{\tilde\ell}} \frac{\kappa'}{\beta^2}  f'\Big]xU'
\nonumber\\[2mm]
 &  + \iota^3 \Big[
 \frac{\kappa^3}{\beta^3}f^2
 -4\mathfrak{p}\frac{\kappa|\beta'|^2}{\beta^5} f^2
 - \frac{2\kappa\beta''}{\beta^4}f^2
 +\frac{8\pi^2}{{\tilde\ell}^2}\frac{\kappa}{\beta}  f''f
 + (8\mathfrak{p} +4) \frac{\pi }{ {\tilde\ell}}\frac{\kappa\beta'}{\beta^3} f' f
 +  \frac{\kappa'\beta'} {\beta^4}f^2
  - \frac{2\pi}{{\tilde\ell}}  \frac{\kappa'}{\beta^2}  f'f\Big]U'
\nonumber\\[2mm]
 & -\iota^3 \frac{2\kappa|\beta'|^2}{\beta^5}  x^3 U''
 \ +\
  \iota^3\Big[ -\frac{6\kappa|\beta'|^2}{\beta^5}f
 +\frac{8\pi}{ {\tilde\ell}}\frac{\kappa\beta'}{\beta^3}   f'\Big] x^2 U''
\nonumber\\[2mm]
 &  +\iota^3\Big[ -\frac{6\kappa|\beta'|^2}{\beta^5}f^2
 + \frac{16\pi}{ {\tilde\ell}} \frac{\kappa\beta'}{\beta^3} f' f
-\frac{8\pi^2} {{\tilde\ell}^2}\frac{\kappa}{\beta} |f'|^2\Big] x U''
\nonumber\\[2mm]
 & +\iota^3\Big[ -\frac{2\kappa|\beta'|^2}{\beta^5}f^3
 +\frac{8\pi}{ {\tilde\ell}} \frac{\kappa}{\beta^3} \big[  \beta' f' f^2  -  f |f'|^2 \big]
  \Big] U''
\nonumber\\[2mm]
 & +  \frac{\iota^4}{\beta^{1+\mathfrak{p}}} \mathcal F_1\left[\tau, x, \beta^{\mathfrak p} f, \frac{\beta^{\mathfrak p} f'}{\tilde\ell}\right]
  +   \frac{\iota^4}{\beta^{1+\mathfrak{p}}} \hat{\mathcal F}_1\left[\tau, x, \beta^{\mathfrak p} f, \frac{\beta^{\mathfrak p} f''}{\tilde\ell^2} \right],
\end{align}
 where
\begin{align}  \label{mathcal F1}
 & \frac{\iota^4}{\beta^{1+\mathfrak{p}}} \mathcal F_1\left[\theta, x, \beta^{\mathfrak p} f, \frac{\beta^{\mathfrak p} f'}{\tilde\ell} \right]
\nonumber\\[2mm]
 =&  -\frac{\iota^4}{24} \frac{\,\iota^2 W_{tttt}\big(\hbar(\iota z)\iota s, \, \iota z\big)}{\beta^6} (x+f)^4U
  +\iota^4  \kappa^3(\iota z)  \Big(  -\frac{\kappa}{\beta}(x+f) + \iota {\tilde a}_3 \Big)
\frac{1}{\, \beta^3\, }(x^2U'+2fxU'+f^2U')
\nonumber\\[2mm]
&\,+\,
\iota^4 \kappa(\iota z)  \Big( \frac{3\kappa}{\beta}(x+f) + \iota  {\tilde a}_2 \Big)
\nonumber\\[2mm]
&\qquad \quad\times
\Bigg\{
\left[ \mathfrak{p}\frac{\beta''} {\beta^4}  + (\mathfrak{p}^2- \mathfrak{p})  \frac{|\beta'|^2} {\beta^5} \right] (xU+fU)
+ \left[ 2\mathfrak{p}\frac{|\beta'|^2}{\beta^5} + \frac{\beta''}{\beta^4}  \right] (x^2U'+2fxU'+f^2U')
\nonumber\\[2mm]
&\qquad\qquad\quad
+(4\mathfrak{p}+2) \frac{\pi}{{\tilde\ell}}\frac{\beta'} {\beta^3} \big[- f' xU'  - ff' U'  \big]
- \frac{4\pi}{{\tilde\ell}} \frac{\beta'} {\beta^3} f' x^2U''
-\frac{ 8\pi}{ {\tilde\ell}}\frac{\beta'} {\beta^3} f' fxU'' -\frac{4\pi}{{\tilde\ell}} \frac{\beta'} {\beta^3}  f' f^2U''
\nonumber\\[2mm]
&\qquad\qquad\quad
+ \frac{|\beta'|^2}{\beta^5} \big[x^3U''+3fx^2U''+3f^2xU''+f^3U''\big]
+\frac{4\pi^2}{{\tilde\ell}^2}\frac{1}{\beta} \big[ |f'|^2 xU''  + f|f'|^2 U''  \big]
\Bigg\}
\nonumber\\[2mm]
&\quad \,+\,
\iota^4 \kappa'(\iota z)\Big(-   \frac{3\kappa}{\beta}(x+f) + \iota  {\tilde a}_1\Big)
\nonumber\\[2mm]
&\qquad \qquad\times
\Bigg\{
\frac{\beta'}{\beta^4} (xU+fU)
+\frac{\beta'}{\beta^4} (x^2U'+2fxU'+f^2U')
-\frac{2\pi}{{\tilde\ell}}\frac{1}{ \beta^2}f' xU'
-\frac{2\pi}{{\tilde\ell}}\frac{1}{ \beta^2} ff' U'
\Bigg\},
\end{align}
and
\begin{align}  \label{hat mathcal F1}
 \frac{\iota^4}{\beta^{1+\mathfrak{p}}} \hat{\mathcal F}_1\left[\theta, x, \beta^{\mathfrak p} f, \frac{\beta^{\mathfrak p} f''}{\tilde\ell^2} \right]
  =  \iota^4 \kappa(\iota z)  \Big( \frac{3\kappa}{\beta}(x+f) + \iota  {\tilde a}_2 \Big)  \left\{\frac{4\pi^2}{{\tilde\ell}^2}\frac{1}{\beta} \big[  ff''U'   +f''xU'\big]
\right\}.
\end{align}

\smallskip
\noindent$\clubsuit$
From the expression of $B_6$ given in \eqref{B2w-2222new},
we can select from $B_6 \left(\iota\frac{e} {\beta^{\mathfrak{p}}} \mathcal Z\right)  $ the terms involving $e''$
\begin{align}
 &\iota^3 \kappa a_2 \frac{ 4\pi^2} {{\tilde\ell}^2} \frac{1}{\, \beta\, }
 (x+f)  \left(\iota\frac{e} {\beta^{\mathfrak{p}}} \mathcal Z\right) _{\tau\tau}
 \nonumber\\[2mm]
& =\iota^4\frac{ 4\pi^2} {{\tilde\ell}^2}  \frac{\kappa a_2}{\, \beta\, } (x+f)
 \left[e''- 2\mathfrak{p}\frac{\beta'} {\beta^2} \frac{\tilde \ell} {2\pi} e' +\left(-\mathfrak{p} \frac{\beta''} {\beta^3} + \mathfrak{p}(\mathfrak{p}+2) \frac{|\beta'|^2} {\beta^4} \right) \frac{\tilde \ell^2} {4\pi^2}e\right]
 \frac{\mathcal Z}{\beta^{\mathfrak{p}} } ,
\end{align}
and the term from $B_6 \left(\iota\frac{e} {\beta^{\mathfrak{p}}} \mathcal Z\right)  $  involving $f''$ is   $-\iota^4 \kappa a_2 \frac{ 4\pi^2} {{\tilde\ell}^2} \frac{1}{\, \beta^{1+ \mathfrak{p}} \, }  (x+f) f'' e   \mathcal Z' $.
 Then we can rewrite $B_6 \left(\iota\frac{e} {\beta^{\mathfrak{p}}} \mathcal Z\right) $ as
\begin{align}\label{expression of B6mathcal Z}
B_6 \left(\iota\frac{e} {\beta^{\mathfrak{p}}} \mathcal Z\right) &
 =
\iota^4 \frac{4\pi^2} {{\tilde\ell}^2} a_2\frac{\kappa}{\, \beta^{1+ \mathfrak{p}} \, }  e'' (x+f)\mathcal Z -\iota^4 \kappa a_2 \frac{ 4\pi^2} {{\tilde\ell}^2} \frac{1}{\, \beta^{1+\mathfrak{p}} \, }  (x+f) f'' e   \mathcal Z'
\nonumber\\[2mm]
&
\quad +\frac{\iota^4} {\beta^3}  \mathcal G_1\left[\theta, x, \beta^{\mathfrak p} f, \frac{\beta^{\mathfrak p} f'}{\tilde\ell},   e, \frac{e'}{\tilde\ell}\right],
\end{align}
where
\begin{align}  \label{mathcal G1}
&  \frac{\iota^4}{\beta^{1+\mathfrak{p}}}  \mathcal G_1\left[\theta, x, \beta f, \frac{\beta^{\mathfrak p} f'}{\tilde\ell},   e, \frac{e'}{\tilde\ell}\right]
\nonumber\\[2mm]
\,=\,&  - \frac{\iota^3}{6} \frac{\, \, \iota^2 W_{ttt}(0, \iota z)}{\beta^5} (x+f)^3  \left(\iota\frac{e} {\beta^{\mathfrak{p}}} \mathcal Z\right)
 -\frac{\iota^4}{24} \frac{\, \, \iota^2 W_{tttt}\big(\hbar(\iota z)\iota s, \, \iota z\big)}{\beta^6} (x+f)^4 \left(\iota\frac{e} {\beta^{\mathfrak{p}}} \mathcal Z\right)
\nonumber\\[2mm]
&+\iota^3  \kappa^3(\iota z) a_3(\iota s, \iota z) \frac{1}{\, \beta^3\, }(x+f)^2 \left(\iota\frac{e} {\beta^{\mathfrak{p}}} \mathcal Z\right)_{x}
\nonumber\\[2mm]
&+ \iota^3 \kappa(\iota z) a_2(\iota s, \iota z)
 \nonumber\\[2mm]
&\qquad\times
\Bigg\{
\left[\mathfrak{p} \frac{\beta''} {\beta^4} +  \mathfrak{p}(\mathfrak{p}-1)  \frac{|\beta'|^2} {\beta^5} \right]  (x+f)  \left(\iota\frac{e} {\beta^{\mathfrak{p}}} \mathcal Z\right)
+ \Big( 2\mathfrak{p} \frac{|\beta'|^2}{\beta^5} +  \frac{\beta''}{\beta^4} \Big)  (x+f)^2  \left(\iota\frac{e} {\beta^{\mathfrak{p}}} \mathcal Z\right)_x
\nonumber\\[2mm]
&\qquad\qquad+ \big(4 \mathfrak{p} +2 \big) \frac{\pi}{ {\tilde\ell}}\frac{\beta'} {\beta^3}  (x+f)\left[-f' \left(\iota\frac{e} {\beta^{\mathfrak{p}}} \mathcal Z\right)_x + \left(\iota\frac{e} {\beta^{\mathfrak{p}}} \mathcal Z\right)_\tau\right]
\nonumber\\[2mm]
&\qquad\qquad
+\frac{ 4\pi}{ {\tilde\ell}}\frac{\beta'} {\beta^3}  (x+f)^2\Bigg[- f'  \left(\iota\frac{e} {\beta^{\mathfrak{p}}} \mathcal Z\right)_{xx}  + \left(\iota\frac{e} {\beta^{\mathfrak{p}}} \mathcal Z\right)_{x\tau}\Bigg]
+ \frac{|\beta'|^2}{\beta^5} (x+f)^3 \left(\iota\frac{e} {\beta^{\mathfrak{p}}} \mathcal Z\right)_{xx}
\nonumber\\[2mm]
&\qquad\qquad
+\frac{ 4\pi^2} {{\tilde\ell}^2}\frac{1}{\beta} (x+f)\left[ |f'|^2  \left(\iota\frac{e} {\beta^{\mathfrak{p}}} \mathcal Z\right)_{xx}    -2f'  \left(\iota\frac{e} {\beta^{\mathfrak{p}}} \mathcal Z\right)_{x\tau}\right]
\Bigg\}
\nonumber\\[2mm]
&+
\iota^3 \kappa'(\iota z)a_1(\iota s, \iota z)
\nonumber\\[2mm]
&\qquad\times
\Bigg\{
\mathfrak{p}\frac{\beta'}{\beta^4} (x+f) \left(\iota\frac{e} {\beta^{\mathfrak{p}}} \mathcal Z\right)
+\frac{\beta'}{\beta^4} (x+f)^2  \left(\iota\frac{e} {\beta^{\mathfrak{p}}} \mathcal Z\right)_{x}
\nonumber\\[2mm]
&\qquad\qquad+\frac{ 2\pi}{ {\tilde\ell}}\frac{1}{ \beta^2} (x+f)\left[-f' \left(\iota\frac{e} {\beta^{\mathfrak{p}}} \mathcal Z\right)_x + \left(\iota\frac{e} {\beta^{\mathfrak{p}}} \mathcal Z\right)_\tau\right]
\Bigg\}
\nonumber\\[2mm]
&+
\iota^4\frac{ 4\pi^2} {{\tilde\ell}^2}  \frac{\kappa a_2}{\, \beta\, } (x+f)
\left[ - 2\mathfrak{p}\frac{\beta'} {\beta^2} \frac{\tilde \ell} {2\pi} e' +\left(-\mathfrak{p} \frac{\beta''} {\beta^3} + \mathfrak{p}(\mathfrak{p}+2) \frac{|\beta'|^2} {\beta^4} \right) \frac{\tilde \ell^2} {4\pi^2}e
\right]\frac{\mathcal Z}{\beta^{\mathfrak{p}}} .
\end{align}

\medskip
 \subsubsection{More precise expression of $S(w_0)$}
For later use, we will rearrange all terms in $S(w_0)$ based on oddness and evenness of the terms in the variable $x$.
We denote
\begin{align}
\iota S_1 = \iota \frac{ \kappa}{\, \beta\, } \Big[ U'+\frac 1{\Lambda} x  U\Big],
\qquad
\iota S_2 = \iota\frac 1{\Lambda}   \frac{ \kappa}{\, \beta\, } f  U,
\end{align}
\begin{align}\label{expression of s3}
 \iota^2 S_3&  = \iota^2
 \Bigg\{\frac{\kappa}{\, \beta^{1+ \mathfrak{p}}\, }\left[ e {\mathcal Z}'+\frac 1{\Lambda} e x\mathcal Z \right]
 -\,\frac{\iota^2W_{tt}(0, \iota z)}{\beta^4}f xU
+\Big[\frac{|\beta'|^2}{\beta^4} 2f
-\frac{4\pi}{{\tilde\ell}}\frac{\beta'}{\beta^2}f'
\Big]xU''
\nonumber\\[2mm]
&\qquad\quad
+ \Big[
-\frac{4\pi^2}{{\tilde\ell}^2}f''
-\big( 4\mathfrak{p}+2 \big) \frac{\pi}{\tilde\ell}  \frac{\beta'}{\beta^2} f'
+2  \mathfrak{p}\frac{|\beta'|^2}{\beta^4}f
+\frac{\beta''} {\beta^3}f
\,-\,\frac{\kappa^2}{\beta^2}f
\Big]U'
\Bigg\},
\end{align}
\begin{align}
 \iota^2 S_4&  =\iota^2
 \Bigg\{ \Big[
-\frac{1}{2}\,  \frac{\iota^2W_{tt}(0, \iota z)}{\beta^4}f^2  +  \mathfrak{p}  \frac{\beta''} {\beta^3}  +\mathfrak{p}( \mathfrak{p}-1)  \frac{|\beta'|^2} {\beta^4} \Big]U
+\Big[\frac{4\pi^2}{{\tilde\ell}^2}|f'|^2
-\frac{4\pi}{{\tilde\ell}}\frac{\beta'}{\beta^2}ff'
+\frac{|\beta'|^2}{\beta^4} f^2
\Big]U''
\nonumber\\[2mm]
&\qquad\quad -\frac{1}{2}\, \frac{\iota^2W_{tt}(0, \iota z)}{\beta^4} x^2U
+\Big[2  \mathfrak{p} \frac{|\beta'|^2}{\beta^4}
+\frac{\beta''} {\beta^3}
\,-\,\frac{\kappa^2}{\beta^2}
\Big]xU'
+\frac{|\beta'|^2}{\beta^4} x^2 U''
\nonumber\\[2mm]
&\qquad \quad
+\frac 1{\Lambda} \frac{\kappa}{\, \beta^{1+ \mathfrak{p}} \, } fe{\mathcal Z}  +\frac{(p_\epsilon-1) (p_\epsilon-2)}{2\beta^{2\mathfrak{p}}}e^2U^{p_\epsilon- 3} {\mathcal Z}^2
\Bigg\},
\end{align}
and
\begin{align}\label{S5}
 \iota^3 S_5
 =  &
 \iota^3 \Bigg\{-\frac{1}{6} \frac{\, \iota^2 W_{ttt}(0, \iota z)}{\beta^5} x^3 U
 + \Big[- \frac{1}{2} \frac{\iota^2 W_{ttt}(0, \iota z)}{\beta^5}f^2-2 \kappa \mathfrak{p} \big[ \frac{\beta''}{\beta^4} + (\mathfrak{p}-1)  \frac{|\beta'|^2} {\beta^5} \big]      +\frac{\kappa'\beta'} {\beta^4} \Big] x U
\nonumber\\[2mm]
 & \quad
 + \Big[ \frac{\kappa^3}{\beta^3}- 4\mathfrak{p}\frac{\kappa|\beta'|^2}{\beta^5}-\frac{2\kappa\beta''}{\beta^4} + \frac{\kappa'\beta'} {\beta^4} \Big]x^2U'
  -\frac{2\kappa|\beta'|^2}{\beta^5}  x^3 U''
 \nonumber\\[2mm]
 & \quad  + \Big[ -\frac{6\kappa|\beta'|^2}{\beta^5}f^2
 + \frac{16\pi}{ {\tilde\ell}} \frac{\kappa\beta'}{\beta^3} f' f
-\frac{8\pi^2} {{\tilde\ell}^2}\frac{\kappa}{\beta} |f'|^2\Big]x U''
\nonumber\\[2mm]
 & \quad  + \Big[
 \frac{\kappa^3}{\beta^3}f^2
 -4 \mathfrak{p}\frac{\kappa|\beta'|^2}{\beta^5} f^2
 - \frac{2\kappa\beta''}{\beta^4}f^2
 +\frac{8\pi^2}{{\tilde\ell}^2}\frac{\kappa}{\beta}  f''f
 + \big[8\mathfrak{p} +4 \big] \frac{\pi }{ {\tilde\ell}}\frac{\kappa\beta'}{\beta^3} f' f
 +  \frac{\kappa'\beta'} {\beta^4}f^2
  - \frac{2\pi}{{\tilde\ell}}  \frac{\kappa'}{\beta^2}  f'f\Big]U'
\Bigg\} ,
\end{align}
\begin{align}\label{S6}
 \iota^3 S_6
  =  & \iota^3\Bigg\{
 - \frac{1}{2} \frac{\,\iota^2 W_{ttt}(0, \iota z)}{\beta^5}f x^2 U
\ +\
\Big[-\frac{\, \, \iota^2 W_{ttt}(0, \iota z)}{6\beta^5}f^3 -2\kappa\mathfrak{p}\big[\frac{\beta''}{\beta^4} + (\mathfrak{p}-1)  \frac{|\beta'|^2}{\beta^5}  \big]   f  + \mathfrak{p} \frac{\kappa'\beta'} {\beta^4} f\Big] U
\nonumber\\[2mm]
 & \qquad +  \Big[
 \frac{2\kappa^3}{\beta^3}f
 -4\mathfrak{p} \frac{\kappa|\beta'|^2}{\beta^5}  f
 -\frac{2\kappa\beta''}{\beta^4}f
 +\frac{8\pi^2}{{\tilde\ell}^2}\frac{\kappa}{\beta}   f''
+\big[8\mathfrak{p} +4 \big]\frac{\pi }{ {\tilde\ell}}\frac{\kappa\beta'}{\beta^3} f'
+\frac{2\kappa'\beta'}{\beta^4}f
- \frac{2\pi}{{\tilde\ell}} \frac{\kappa'}{\beta^2}  f'\Big]xU'
\nonumber\\[2mm]
 & \qquad     +  \Big[ -\frac{6\kappa|\beta'|^2}{\beta^5}f
 +\frac{8\pi}{ {\tilde\ell}}\frac{\kappa\beta'}{\beta^3}   f'\Big]  x^2 U''
  + \Big[ -\frac{2\kappa|\beta'|^2}{\beta^5}f^3
 +\frac{8\pi}{ {\tilde\ell}} \frac{\kappa\beta'}{\beta^3} f' f^2
 -\frac{8\pi^2} {{\tilde\ell}^2}\frac{\kappa}{\beta} |f'|^2 f \Big]U''
 \Bigg\},
\end{align}
\begin{align}\label{S7}
\iota^3S_7
 = &\frac{\iota^3}{\beta^{\mathfrak{p}} } \Bigg\{ \Big[
-\frac{4\pi^2}{{\tilde\ell}^2}f''
-  2\frac{\pi}{\tilde\ell}  \frac{\beta'}{\beta^2} f'
+\frac{\beta''} {\beta^3}f
\,-\,\frac{\kappa^2}{\beta^2}f
\Big]e {\mathcal Z}'
+\Big[\frac{2|\beta'|^2}{\beta^4} f
-\frac{4\pi}{{\tilde\ell}}\frac{\beta'}{\beta^2}f'
\Big]e x {\mathcal Z}''
\nonumber\\[2mm]
&\qquad\quad -\,\frac{\iota^2W_{tt}(0, \iota z)}{\beta^4}f e x{\mathcal Z}
+
\Big[ -\frac{8\pi^2}{{\tilde\ell}^2}f'  +\frac{4\pi}{{\tilde\ell}}\frac{\beta'}{\beta^2}f\Big]e' {\mathcal Z}'
\Bigg\},
\end{align}
\begin{align}\label{S8}
  \iota^3S_8
\,=\,&  \frac{\iota^3}{\beta^{\mathfrak{p}}}\Bigg\{
\frac{2\pi}{\tilde\ell} \frac{\beta'}{\beta^2}     e'{\mathcal Z}
+ \left[ \frac{4\pi}{{\tilde\ell}}\frac{\beta'}{\beta^2} e'   + \frac{\beta''}{\beta^3}e - \frac{\kappa^2}{\beta^2}e  \right]
 x\mathcal Z'
 \nonumber\\[2mm]
&\qquad
  +
 \left[ -\frac{1}{2}\,  \frac{\iota^2W_{tt}(0, \iota z)}{\beta^4}f^2  -  \mathfrak{p} (\mathfrak{p}+2) \frac{|\beta'|^2}{\beta^4}
+ \mathfrak{p} \frac{\beta''} {\beta^3}  \right]  e\mathcal Z
\nonumber\\[2mm]
&\qquad
+\Big[\frac{4\pi^2}{{\tilde\ell}^2}|f'|^2
-\frac{4\pi}{{\tilde\ell}}\frac{\beta'}{\beta^2}ff'
+\frac{|\beta'|^2}{\beta^4} f^2
\Big]e {\mathcal Z}''
-\frac{1}{2}\, \frac{\iota^2W_{tt}(0, \iota z)}{\beta^4}  e x^2{\mathcal Z}
+
\frac{|\beta'|^2}{\beta^4} e x^2{\mathcal Z}''
\Bigg\}.
\end{align}
Finally, we can rewrite $S(w_0)$ in \eqref{S(w_0)} as
\begin{align}
 S(w_0)\,=\, &
 \iota S_1 + \iota S_2 + \iota^2 S_3 + \iota^2 S_4 + \iota^3 S_5
 + \iota^3 S_6 + \iota^3 S_7+\iota^3 S_8
 \nonumber\\[2mm]
 &
 +\iota^3\frac{4\pi^2}{{\tilde\ell}^2}  \frac{e''}{\beta^{\mathfrak{p}}}\mathcal Z
  +\frac1{\beta^{\mathfrak{p}}} \iota\lambda_0 e{\mathcal Z}+N_{0\e} \left(\iota\frac{e} {\beta^{\mathfrak{p}}} \mathcal Z\right)
\nonumber\\[2mm]
 &
  \underbrace{ +\iota^4 \frac{4\pi^2} {{\tilde\ell}^2} a_2\frac{\kappa}{\, \beta^{1+\mathfrak{p}}\, }  e'' (x+f)\mathcal Z -\iota^4 \kappa a_2 \frac{ 4\pi^2} {{\tilde\ell}^2} \frac{1}{\, \beta^{1+\mathfrak{p}}\, }  (x+f) f'' e   \mathcal Z'    + \frac{\iota^4}{\beta^{2} } \hat{\mathcal F}_1\left[\tau, x, \beta^{\mathfrak p} f, \frac{\beta^{\mathfrak p} f''}{\tilde\ell^2} \right] }_{\text{terms involving}~f'', e''}
\nonumber\\[2mm]
 &
 +  \frac{\iota^4}{\beta^{1+\mathfrak{p}}} \mathcal F_1\left[\tau, x, \beta^{\mathfrak p} f, \frac{\beta^{\mathfrak p} f'}{\tilde\ell}\right]
 +\frac{\iota^4}{\beta^{1+\mathfrak{p}}}  \mathcal G_1\left[\tau, x, \beta^{\mathfrak p} f, \frac{\beta^{\mathfrak p} f'}{\tilde\ell},   e, \frac{e'}{\tilde\ell}\right].
\label{error1}
\end{align}
We note that the expressions of the following terms
$$
\frac{\iota^4}{\beta^{1+\mathfrak{p}}} \mathcal F_1\left[\tau, x, \beta^{\mathfrak p} f, \frac{\beta^{\mathfrak p} f'}{\tilde\ell}\right]
\quad
\frac{\iota^4}{\beta^{1+\mathfrak{p}}}  \mathcal G_1\left[\tau, x, \beta^{\mathfrak p} f, \frac{\beta^{\mathfrak p} f'}{\tilde\ell},   e, \frac{e'}{\tilde\ell}\right],
\quad
\frac{\iota^4}{\beta^{1+\mathfrak{p}}} \hat{\mathcal F}_1\left[\tau, x, \beta^{\mathfrak p} f, \frac{\beta^{\mathfrak p} f''}{\tilde\ell^2} \right]
$$
are given in \eqref{mathcal F1},  \eqref{hat mathcal F1},  \eqref{mathcal G1}, respectively.

\medskip
\subsection{  The second approximate solution via first correction}\
\label{section7.3}
\subsubsection{Second approximate solution and its error}
In this part, we want to construct a correction function to eliminate the terms $\iota S_1$ and $\iota S_2$
in the error $S(w_0)$.
By the same method in \cite{delPKowWei1}, we first choose correction terms $\phi_{11}, {\phi_{12}}$ to cancel the terms
$ S_1,   S_2. $
That is to say, for fixed $\tau\in (0, 2\pi)$,
we need solutions of
\begin{align} \label{correction1}
\phi_{11, xx} - \phi_{11} + (p_\epsilon-1) U^{p_\epsilon-2} \phi_{11} \,=\, -  S_1, \quad \phi_{11}(\pm \infty, \tau) \,=\, 0,
\end{align}
\begin{align} \label{correction2}
\phi_{12, xx} - {\phi_{12}} + (p_\epsilon-1) U^{p_\epsilon-2} {\phi_{12}} \,=\, -  S_2, \quad{\phi_{12}}(\pm \infty, \tau) \,=\, 0.
\end{align}
\medskip
As it is well known, those problems are solvable provided that
\begin{equation*}
 \int_{\R}  S_1\, U'  \,{\mathrm d}x \,=\, 0,
\qquad
 \int_{\R} S_2\, U'  \,{\mathrm d}x \,=\, 0.
\end{equation*}
Using the identity  \eqref{identity  for sigma},
we can get
\begin{align*}
  \int_{\R} S_1\, U'  \,{\mathrm d}x
\,=\, &    \frac{\kappa}{\beta}\int_{\R}\Big[ U' +\frac 1{\Lambda}xU \Big] \, U' \,{\mathrm d}x=0.
\end{align*}
Thus, for $\phi_{11}$, we have
\begin{align} \label{definition of phi11}
\phi_{11} (x, \tau) \,=\,      \frac{ \kappa}{\beta}\varpi_1(x),
\end{align}
and $\varpi_1$ is the unique solution of
\begin{align}\label{varphi1equation new}
\varpi_{1, xx} - \varpi_1 + (p_\epsilon-1) U^{p_\epsilon-2} \varpi_1 \,=\, -\Big( U' + \frac 1{\Lambda} x U\Big), \quad \int_{\R} \varpi_1 U'  \,{\mathrm d}x \,=\, 0.
\end{align}
In fact, the function $\varpi_1$ is an odd function and enjoys the exponential decay
\[ |\varpi_1 |\leq C |x| e^{-c|x|}, \quad \text{when $|x|$ is large}. \]
On the other hand, we get
\begin{align}
{\phi_{12}} (x, \tau)
 \,=\, \frac 1{\Lambda}    \frac{ \kappa}{\, \beta\, } f \varpi_2(x),
\end{align}
with
$\varpi_2 (x) \,=\, -\frac1{p_\epsilon-2} \big(U+ x U'\big), $
and $\varpi_2$ satisfies
\begin{align} \label{varphi2equation new}
\varpi_{2, xx} - \varpi_2 + (p_\epsilon-1) U^{p_\epsilon-2} \varpi_2 \,=\, -U.
\end{align}

Therefore, we define the correction as
\begin{align}\label{definition of phi1}
\iota\phi_1(x, \tau)
 \,=\, \iota{\phi_{11}}(x, \tau) + \iota{\phi_{12}}(x, \tau)
 =\iota{\mathbf b}_1(\iota z)\varpi_1(x)
+
 \iota\frac 1{\Lambda} {\mathbf b}_1(\iota z)f(\tau) \varpi_2(x),
\end{align}
where the relation between $z$ and $\tau$ is given in \eqref{ztau}
and
$$
{\mathbf b}_1(\theta)=\frac{\kappa(\theta)}{\beta(\theta)}.
$$
Then we take
\begin{align}\label{secondapproximatesolution}
w_1(x, \tau) \,=\,  w_0(x, \tau)\,+\, \iota\phi_1(x, \tau)
\end{align}
as the {\textbf{second approximate solution}}.
Then we use \eqref{error1} to derive
\begin{align}\label{sw1new}
S (w_0+ \iota\phi_1)
\,=\, &S(w_0)  + \iota \Big[\phi_{1, xx} -  \phi_1  +(p_\epsilon-1) U^{p_\epsilon-2}\phi_1 \Big]
+ \frac{\iota^2 4\pi^2}{{\tilde\ell}^2} \iota\phi_{1, \tau\tau}
\nonumber\\[2mm]
& +\left[ \left( w_0+\iota\phi_1\right)^{p_\epsilon-1}  -w_0^{p_\epsilon-1}  -\iota (p_\epsilon-1) U^{p_\epsilon-2} \phi_1 \right]
+ B_4(\iota\phi_1)+ B_5(\iota\phi_1)+ B_6(\iota\phi_1)
\nonumber\\[2mm]
\,=\,&
\iota^2 S_3+\iota^2 S_4   + \iota^3 S_5 +\iota^3 S_6
+\iota^3 S_7+\iota^3 S_8
+\frac1{\beta^{\mathfrak{p}}} \iota^3\frac{4\pi^2}{{\tilde\ell}^2}  e''\mathcal Z
+ \frac1{\beta^{\mathfrak{p}} } \iota\lambda_0 e{\mathcal Z} +N_{0\e} \left(\iota\frac{e} {\beta^{\mathfrak{p}}} \mathcal Z\right)
\nonumber\\[2mm]
&
  \underbrace{ +\iota^4 \frac{4\pi^2} {{\tilde\ell}^2} a_2\frac{\kappa}{\, \beta^{1+ \mathfrak{p}} \, }  e'' (x+f)\mathcal Z -\iota^4 \kappa a_2 \frac{ 4\pi^2} {{\tilde\ell}^2} \frac{1}{\, \beta^{1+ \mathfrak{p}} \, }  (x+f) f'' e   \mathcal Z'    + \frac{\iota^4}{\beta^{1+\mathfrak{p}}} \hat{\mathcal F}_1\left[\tau, x, \beta^{\mathfrak p} f, \frac{\beta^{\mathfrak p} f''}{\tilde\ell^2} \right] }_{\text{terms involving}~f'', e''}
\nonumber\\[2mm]
 &
 +  \frac{\iota^4}{\beta^{1+\mathfrak{p}}} \mathcal F_1\left[\tau, x, \beta^{\mathfrak p} f, \frac{\beta^{\mathfrak p} f'}{\tilde\ell}\right]
 +\frac{\iota^4}{\beta^{1+\mathfrak{p}}}  \mathcal G_1\left[\tau, x, \beta^{\mathfrak p} f, \frac{\beta^{\mathfrak p} f'}{\tilde\ell},   e, \frac{e'}{\tilde\ell}\right]
\nonumber\\[2mm]
&+ \frac{\iota^3 4\pi^2}{{\tilde\ell}^2}\phi_{1, \tau\tau}
+  \mathbb{N}_0(\iota\phi_1)     +  B_4(\iota\phi_1)+ B_5(\iota\phi_1)+ B_6(\iota\phi_1),
\end{align}
where  the  linear operators $ B_4, B_5, B_6 $ are given in \eqref{new opreator B4}, \eqref{new opreator B5}, \eqref{B2w-2222new},
and the nonlinear term is
\begin{align}
\mathbb{N}_0(\iota\phi_1)
=
\left( w_0+\iota\phi_1\right)^{p_\epsilon-1}  -w_0^{p_\epsilon-1}  -\iota (p_\epsilon-1) U^{p_\epsilon-2} \phi_1.
\end{align}

By some calculations, we provide the following facts
\begin{align}
\phi_{1, x}
={\mathbf b}_1\varpi_{1,x} +\frac 1{\Lambda}{\mathbf b}_1 f \varpi_{2,x},
\qquad
\phi_{1, xx}
={\mathbf b}_1\varpi_{1,xx} +\frac 1{\Lambda}{\mathbf b}_1 f \varpi_{2,xx},
\end{align}
\begin{align}
\phi_{1, \tau}
= & {\mathbf b}_1'\frac{{\tilde\ell}}{2\pi \beta} \varpi_{1}
+
\frac 1{\Lambda}{\mathbf b}_1'\frac{{\tilde\ell}}{2\pi\beta} f\varpi_{2}
+
 \frac 1{\Lambda}{\mathbf b}_1 f'\varpi_{2}
\nonumber\\[2mm]
= & \frac{{\tilde\ell}}{2\pi} {\mathbf b}_2 \varpi_{1}
+
\frac 1{\Lambda}\frac{{\tilde\ell}}{2\pi}{\mathbf b}_2f \varpi_{2}
+
\frac 1{\Lambda}{\mathbf b}_1 f'\varpi_{2},
\end{align}
where
\begin{align}
{\mathbf b}_2(\theta)=\frac{{\mathbf b}_1'(\theta)}{\beta(\theta)},
\label{mathbfb2}
\end{align}
and
\begin{align}
\phi_{1, x\tau}
= &
\frac{{\tilde\ell}}{2\pi} {\mathbf b}_2 \varpi_{1,x}
+
\frac 1{\Lambda}\frac{{\tilde\ell}}{2\pi}{\mathbf b}_2 f\varpi_{2,x}
+
\frac{1}{\Lambda} {\mathbf b}_1 f'\varpi_{2,x},
\end{align}

\begin{align}
\phi_{1, \tau\tau}
= &
\frac{{\tilde\ell}}{2\pi} {\mathbf b}_2'\frac{{\tilde\ell}}{2\pi\beta} \varpi_{1}
+
\frac{1}{\Lambda} \frac{{\tilde\ell}}{2\pi}{\mathbf b}_2'\frac{{\tilde\ell}}{2\pi\beta}f \varpi_{2}
+
\frac{1}{\Lambda} \frac{{\tilde\ell}}{2\pi}{\mathbf b}_2f' \varpi_{2}
+
\frac{1}{\Lambda} {\mathbf b}_1'\frac{{\tilde\ell}}{2\pi\beta} f'\varpi_{2}
 +
\frac{1}{\Lambda} {\mathbf b}_1 f''\varpi_{2}
\nonumber\\
= &
\frac{{\tilde\ell}^2}{4\pi^2} {\mathbf b}_3 \varpi_{1}
+
\frac{1}{\Lambda} \frac{{\tilde\ell}^2}{4\pi^2}{\mathbf b}_3 f \varpi_{2}
+
\frac{2}{\Lambda} \frac{{\tilde\ell}}{2\pi}{\mathbf b}_2f' \varpi_{2}
+
\frac{1}{\Lambda} {\mathbf b}_1 f''\varpi_{2},
\label{expression of phi1tautau}
\end{align}
where
\begin{align}
{\mathbf b}_3(\theta)=\frac{{\mathbf b}_2'(\theta)}{\beta(\theta)}.
\label{mathbfb3}
\end{align}

\medskip
\subsubsection{More precise expression of $S(w_1)$}

We will do the following analysis.

\smallskip\noindent $\clubsuit$
We can further decompose the term $\iota^2S_4$ as follows
\begin{align}  \label{decompose S3+S4}
\iota^2 S_4  =\iota^2 S_{41} +\iota^2 S_{42},
\end{align}
where
\begin{align}\label{barS41}
\iota^2 S_{41} =& \iota^2\Bigg\{
 \Big[ \mathfrak{p}  \frac{\beta''} {\beta^3}  +  \mathfrak{p}( \mathfrak{p}-1)  \frac{|\beta'|^2} {\beta^4}  \Big] U
-
\frac{1}{2}\, \frac{\iota^2W_{tt}(0, \iota z)}{\beta^4} x^2U
\nonumber
\\
&\qquad
+\Big[ 2\mathfrak{p}\frac{|\beta'|^2}{\beta^4}
+\frac{\beta''} {\beta^3}
\,-\,\frac{\kappa^2}{\beta^2}
\Big]xU'
+
\frac{|\beta'|^2}{\beta^4} x^2 U''
\Bigg\},
\end{align}

\begin{align}\label{barS42}
\iota^2 S_{42}&  =\iota^2\Bigg\{
-\frac{1}{2}\, \frac{\iota^2W_{tt}(0, \iota z)}{\beta^4}|f|^2 U
+
\Big[\frac{4\pi^2}{{\tilde\ell}^2}|f'|^2
-\frac{4\pi}{{\tilde\ell}}\frac{\beta'}{\beta^2}ff'
+\frac{|\beta'|^2}{\beta^4} |f|^2
\Big]U''
\nonumber\\[2mm]
 & \qquad\qquad
+
\frac1{\Lambda}  \frac{\kappa}{\, \beta^{1+ \mathfrak{p}} \, } fe{\mathcal Z}
+
 \frac{(p_\epsilon-1) (p_\epsilon-2)}{2\beta^{2\mathfrak{p}}}e^2U^{p_\epsilon- 3} {\mathcal Z}^2
\Bigg\}.
\end{align}

\smallskip\noindent $\clubsuit$
Moreover,  by the definition of $B_4$ in  \eqref{new opreator B4}, we  can decompose $  B_4(\iota\phi_1)$ as
\begin{align}\label{B4phi1}
 B_4(\iota\phi_1)
 =&\iota^2 \frac{\kappa}{\, \beta\, } \Big[ \phi_{1,x} + \frac{1}{\Lambda} x\phi_1+\frac{1}{\Lambda_\e} f\phi_1\Big]
\nonumber\\[2mm]
=&\iota^2 {\mathbf b}_1^2\Big[
\varpi_{1,x}
+ \frac{1}{\Lambda} f\varpi_{2,x}
+\frac{1}{\Lambda} x\varpi_1
+ \frac{1}{\Lambda^2}f x\varpi_2
+\frac{1}{\Lambda} f\varpi_1
+\frac{1}{\Lambda^2} f^2\varpi_2
  \Big]
\nonumber\\[2mm]
\,=\,&\iota^2 {\bf B}_{41} +\iota^2{\bf B}_{42},
\end{align}
with
\begin{align}   \label{B41}
{\bf B}_{41} =
 {\mathbf b}_1^2\Big[
\varpi_{1,x}
+\frac{1}{\Lambda_\e} x\varpi_1
\Big],
\end{align}
\begin{align}   \label{B42}
{\bf B}_{42} =
 {\mathbf b}_1^2\Big[
\frac{1}{\Lambda} f\varpi_{2,x}
+ \frac{1}{\Lambda^2} f x\varpi_2
+\frac{1}{\Lambda} f\varpi_1
+ \frac{1}{\Lambda^2} f^2\varpi_2
\Big].
\end{align}

\smallskip\noindent $\clubsuit$
  From the expression of $B_5$ given in \eqref{new opreator B5}, we have
\begin{align} \label{expression of B5phi1}
B_5(\iota\phi_1)
\, =\,&\iota^3
\Bigg\{\big[4\mathfrak{p}+2 \big]\frac{\pi}{{\tilde\ell}}\frac{\beta'} {\beta^2}\Big[- f' \phi_{1,x} +\phi_{1,\tau}\Big]
 +\frac{4\pi^2}{{\tilde\ell}^2}\Big[|f'|^2 \phi_{1,xx}  -f''\phi_{1,x}  -2f' \phi_{1,x\tau}\Big]
\nonumber\\[2mm]
&\qquad
+ \frac{4\pi}{{\tilde\ell}} \frac{\beta'}{\beta^2} (x+f)\big[- f' \phi_{1,xx}  +\phi_{1,x\tau}\big]
+\Big[\frac{\beta''}{\beta^3} + 2\mathfrak{p} \frac{|\beta'|^2}{\beta^4}  \,-\,\frac{\kappa^2}{\beta^2} \Big]
(x+f) \phi_{1,x}
\nonumber\\[2mm]
&\qquad
+\Big[\mathfrak{p} \frac{\beta''}{\beta^3}  +  \mathfrak{p}  (\mathfrak{p} -1) \frac{|\beta'|^2}{\beta^4} \Big] w+ \frac{|\beta'|^2}{\beta^4}  (x+f)^2 \phi_{1,xx}
-\frac{1}{2}\, \frac{\iota^2W_{tt}(0, \iota z)}{\beta^4} (x+f)^2\phi_{1}
\Bigg\}.
\end{align}

\smallskip\noindent $\clubsuit$
By the definition of $B_6$  given in \eqref{B2w-2222new},  we
 select  the terms from  $  B_6(\iota\phi_1)$    involving $f''$
\begin{align}
&\iota^4 \kappa(\iota z) a_2(\iota s, \iota z)  \frac{ 4\pi^2} {{\tilde\ell}^2}\frac{1}{\beta} (x+f)\big[  -f''\phi_{1, x} +\phi_{1, \tau\tau}\big]
\nonumber\\[2mm]
& =\iota^4 \kappa(\iota z) a_2(\iota s, \iota z)  \frac{ 4\pi^2} {{\tilde\ell}^2}\frac{1}{\beta} (x+f)
\Bigg[  -f'' \big({\mathbf b}_1\varpi_{1,x} +\frac{1}{\Lambda} {\mathbf b}_1 f \varpi_{2,x} \big)
\nonumber\\[2mm]
&  \qquad\qquad
+ \Big(\frac{{\tilde\ell}^2}{4\pi^2} {\mathbf b}_3 \varpi_{1}
+
\frac{1}{\Lambda} \frac{{\tilde\ell}^2}{4\pi^2}{\mathbf b}_3 f \varpi_{2}
+
\frac{2}{\Lambda} \frac{{\tilde\ell}}{2\pi}{\mathbf b}_2f' \varpi_{2}
+
\frac{1}{\Lambda} {\mathbf b}_1 f''\varpi_{2} \Big)\Bigg]
\nonumber\\[2mm]
& = \underbrace{\iota^4 \kappa(\iota z) a_2(\iota s, \iota z)  \frac{ 4\pi^2} {{\tilde\ell}^2}\frac{1}{\beta} (x+f)\Bigg[  -f'' \big({\mathbf b}_1\varpi_{1,x} +\frac{1}{\Lambda} {\mathbf b}_1 f \varpi_{2,x} \big)
 +
\frac{1}{\Lambda} {\mathbf b}_1 f''\varpi_{2} \Bigg]}_{:= \frac{\iota^4}{\beta^{1+\mathfrak{p}}} \hat{\mathcal F}_2\left[\tau, x, \beta^{\mathfrak p} f, \frac{\beta^{\mathfrak p} f''}{\tilde\ell^2} \right] }
\nonumber\\[2mm]
& \quad + \iota^4 \kappa(\iota z) a_2(\iota s, \iota z)  \frac{ 4\pi^2} {{\tilde\ell}^2}\frac{1}{\beta} (x+f) \Big(\frac{{\tilde\ell}^2}{4\pi^2} {\mathbf b}_3 \varpi_{1}
+
\frac{1}{\Lambda} \frac{{\tilde\ell}^2}{4\pi^2}{\mathbf b}_3 f \varpi_{2}
+
\frac{2}{\Lambda} \frac{{\tilde\ell}}{2\pi}{\mathbf b}_2f' \varpi_{2} \Big)
\nonumber\\[2mm]
& =  \frac{\iota^4}{\beta^{1+\mathfrak{p}}} \hat{\mathcal F}_2\left[\tau, x, \beta^{\mathfrak p} f, \frac{\beta^{\mathfrak p} f''}{\tilde\ell^2} \right]
\nonumber\\[2mm]
& \quad+\iota^4 \kappa(\iota z) a_2(\iota s, \iota z)  \frac{ 4\pi^2} {{\tilde\ell}^2}\frac{1}{\beta} (x+f) \Big(\frac{{\tilde\ell}^2}{4\pi^2} {\mathbf b}_3 \varpi_{1}
+
\frac{1}{\Lambda} \frac{{\tilde\ell}^2}{4\pi^2}{\mathbf b}_3 f \varpi_{2}
+
\frac{2}{\Lambda} \frac{{\tilde\ell}}{2\pi}{\mathbf b}_2f' \varpi_{2} \Big). \label{hat mathcal F2}
\end{align}
And  we can rewrite $  B_6(\iota\phi_1) $ as
\begin{align}\label{expression of B6phi1}
  B_6(\iota\phi_1) =   \frac{\iota^4}{\beta^{1+\mathfrak{p}}}  \mathcal F_2\left[\theta, x, \beta^{\mathfrak p} f, \frac{\beta^{\mathfrak p} f'}{\tilde\ell}\right]
  + \frac{\iota^4}{\beta^{1+\mathfrak{p}}} \hat{\mathcal F}_2\left[\tau, x, \beta^{\mathfrak p} f, \frac{\beta^{\mathfrak p} f''}{\tilde\ell^2} \right],
\end{align}
where
\begin{align}  \label{mathcal F2}
& \frac{\iota^4}{\beta^{1+\mathfrak{p}}}  \mathcal F_2\left[\theta, x, \beta^{\mathfrak p} F, \frac{\beta^{\mathfrak p} f'}{\tilde\ell}\right],
\nonumber\\[2mm]
\,=\,&
 - \frac{\iota^4}{6} \frac{\, \,\iota^2 W_{ttt}(0, \iota z)}{\beta^5} (x+f)^3 \phi_1
 -\frac{\iota^5}{24} \frac{\, \,\iota^2 W_{tttt}\big(\hbar(\iota z)\iota s, \, \iota z\big)}{\beta^6} (x+f)^4\phi_1
\nonumber\\[2mm]
&
 +\iota^4  \kappa^3(\iota z) a_3(\iota s, \iota z) \frac{1}{\, \beta^3\, }(x+f)^2\phi_{1, x}
\nonumber\\[2mm]
&\,+\,
\iota^4 \kappa(\iota z) a_2(\iota s, \iota z)
\Bigg\{
\Big[\mathfrak{p} \frac{\beta''}{\beta^3}  +  \mathfrak{p}  (\mathfrak{p} -1) \frac{|\beta'|^2}{\beta^4} \Big] (x+f) \phi_1
+\Big[ 2 \mathfrak{p}\frac{|\beta'|^2}{\beta^5}+  \frac{\beta''}{\beta^4}  \Big]  (x+f)^2 \phi_{1, x}
\nonumber\\[2mm]
&\qquad\qquad\qquad\qquad\qquad
+ \big[4\mathfrak{p} +2 \big] \frac{\pi}{ {\tilde\ell}}\frac{\beta'} {\beta^3}  (x+f)\big[-f'\phi_{1, x} +\phi_{1, \tau}\big]
\nonumber\\[2mm]
&\qquad\qquad\qquad\qquad\qquad
+\frac{ 4\pi}{ {\tilde\ell}}\frac{\beta'} {\beta^3}  (x+f)^2[- f' w_{xx}  +w_{x\tau}]
+ \frac{|\beta'|^2}{\beta^5} (x+f)^3\phi_{1, xx}
\nonumber\\[2mm]
&\qquad\qquad\qquad\qquad\qquad
+\frac{ 4\pi^2} {{\tilde\ell}^2}\frac{1}{\beta} (x+f)\big[ |f'|^2 \phi_{1, xx}    -2f' \phi_{1, x\tau}\big]
\Bigg\}
\nonumber\\[2mm]
&\,+\,
\iota^4 \kappa'(\iota z)a_1(\iota s, \iota z)
\Bigg\{
\mathfrak{p}\frac{\beta'}{\beta^4} (x+f)\phi_1
+\frac{\beta'}{\beta^4} (x+f)^2 \phi_{1, x}
+\frac{ 2\pi}{ {\tilde\ell}}\frac{1}{ \beta^2} (x+f)\big[-f'\phi_{1, x} +\phi_{1, \tau}\big]
\Bigg\}
\nonumber\\[2mm]
& + \iota^4 \kappa(\iota z) a_2(\iota s, \iota z)  \frac{ 4\pi^2} {{\tilde\ell}^2}\frac{1}{\beta} (x+f) \Bigg\{\frac{{\tilde\ell}^2}{4\pi^2} {\mathbf b}_3 \varpi_{1}
+
\frac{1}{\Lambda}\frac{{\tilde\ell}^2}{4\pi^2}{\mathbf b}_3 f \varpi_{2}
+
\frac{2}{\Lambda} \frac{{\tilde\ell}}{2\pi}{\mathbf b}_2f' \varpi_{2} \Bigg\}.
\end{align}

\smallskip\noindent $\clubsuit$
Similarly, for $\mathbb{N}_0(\iota\phi_1)$, we also have the decomposition
\begin{align}\label{decompose N0}
 \mathbb{N}_0(\iota\phi_1)
\,=\,&\left( w_0+\iota\phi_1\right)^{p_\epsilon-1}  -w_0^{p_\epsilon-1}  -\iota (p_\epsilon-1)  w_0^{p_\epsilon-2} \phi_1 + \iota (p_\epsilon-1)  \left( w_0^{p_\epsilon-2}- U^{p_\epsilon-2}\right)  \phi_1
 \nonumber\\[2mm]
\,=\,& \iota^2 \frac{(p_\epsilon-1) (p_\epsilon-2)} {2}  w_0^{p_\epsilon- 3} \phi_1^2 +  \iota^3 \frac{(p_\epsilon-1) (p_\epsilon-2)(p_\epsilon- 3)} {6}  w_0^{p_\epsilon- 4} \phi_1^3
 \nonumber\\[2mm]
\quad &+  \iota^2(p_\epsilon-1)(p_\epsilon-2) U^{p_\epsilon- 3}   \frac{e} {\beta^{\mathfrak{p}}} \mathcal Z\phi_1   + \iota^3 \frac{(p_\epsilon-1)(p_\epsilon-2) (p_\epsilon- 3)}2 U^{p_\epsilon- 4}  \frac{e^2} {\beta^{2\mathfrak{p}}} \mathcal Z^2\phi_1     + \iota^4 \mathcal N(x, f, e)
\nonumber\\[2mm]
\,=\,&    \iota^2 \frac{(p_\epsilon-1) (p_\epsilon-2)} {2}  U^{p_\epsilon- 3} \phi_1^2  +  \iota^3 \frac{(p_\epsilon-1) (p_\epsilon-2) (p_\epsilon- 3)} {2}  U^{p_\epsilon- 4} \phi_1^2  \frac{e} {\beta^{\mathfrak{p}}} \mathcal Z
  \nonumber\\[2mm]
\quad &
  + \iota^3 \frac{(p_\epsilon-1) (p_\epsilon-2)(p_\epsilon- 3)} {6}   U^{p_\epsilon- 4} \phi_1^3
+  \iota^2(p_\epsilon-1)(p_\epsilon-2) U^{p_\epsilon- 3}   \frac{e} {\beta^{\mathfrak{p}}} \mathcal Z\phi_1
 \nonumber\\[2mm]
\quad &
+ \iota^3 \frac{(p_\epsilon-1)(p_\epsilon-2) (p_\epsilon- 3)}2 U^{p_\epsilon- 4}  \frac{e^2} {\beta^{2\mathfrak{p}}} \mathcal Z^2\phi_1     + \iota^4{\bf N}_{04}
\nonumber\\[2mm]
\,=\,&   \iota^2 \frac{(p_\epsilon-1) (p_\epsilon-2)} {2}  U^{p_\epsilon- 3}  \left[{\mathbf b}_1^2\varpi_1^2+ \frac 1{\Lambda^2} {\mathbf b}_1^2f^2\varpi_2^2+ 2 {\mathbf b}_1^2 \frac 1{\Lambda} f \varpi_1 \varpi_2 \right]
 \nonumber\\[2mm]
\quad &
+\iota^3 \frac{(p_\epsilon-1) (p_\epsilon-2) (p_\epsilon- 3)} {2}  U^{p_\epsilon- 4} \left[{\mathbf b}_1^2\varpi_1^2+ \frac 1{\Lambda^2} {\mathbf b}_1^2f^2\varpi_2^2+ 2 {\mathbf b}_1^2 \frac 1{\Lambda} f \varpi_1 \varpi_2 \right]     \frac{e} {\beta^{\mathfrak{p}}} \mathcal Z
 \nonumber\\[2mm]
\quad &   + \iota^3 \frac{(p_\epsilon-1) (p_\epsilon-2)(p_\epsilon- 3)} {6}   U^{p_\epsilon- 4}  \left[{\mathbf b}_1\varpi_1+ \frac 1{\Lambda} {\mathbf b}_1f\varpi_2 \right]^3
 \nonumber\\[2mm]
\quad &+  \iota^2(p_\epsilon-1)(p_\epsilon-2) U^{p_\epsilon- 3}   \frac{e} {\beta^{\mathfrak{p}}} \mathcal Z\Big( {\mathbf b}_1\varpi_1+ \frac 1{\Lambda} {\mathbf b}_1f\varpi_2\Big)
 \nonumber\\[2mm]
\quad &
  + \iota^3 \frac{(p_\epsilon-1)(p_\epsilon-2) (p_\epsilon- 3)}2 U^{p_\epsilon- 4}  \frac{e^2} {\beta^{2\mathfrak{p}}} \mathcal Z^2\Big({\mathbf b}_1\varpi_1+ \frac 1{\Lambda} {\mathbf b}_1f\varpi_2\Big)      +  \iota^4{\bf N}_{04}
\nonumber\\[2mm]
=&\iota^2 {\bf N}_{01}+\iota^2 {\bf N}_{02}+\iota^3{\bf N}_{03} +\iota^4{\bf N}_{04},
\end{align}
 where
\begin{align}\label{N01}
{\bf N}_{01} =  \frac{(p_\epsilon-1) (p_\epsilon-2)} {2}{\mathbf b}_1^2 U^{p_\epsilon- 3}\varpi_1^2,
\end{align}
\begin{align}\label{N02}
{\bf N}_{02} =  (p_\epsilon-1)(p_\epsilon-2) U^{p_\epsilon- 3}  \Bigg[  \frac 1{2\Lambda^2} {\mathbf b}_1^2f^2\varpi_2^2+  {\mathbf b}_1^2 \frac 1{\Lambda} f \varpi_1 \varpi_2
+  \frac{e} {\beta^{\mathfrak{p}}} \mathcal Z\Big( {\mathbf b}_1\varpi_1+ \frac 1{\Lambda} {\mathbf b}_1f\varpi_2\Big)    \Bigg] ,
\end{align}
\begin{align}\label{N03}
{\bf N}_{03}\,=\,&
  \frac{(p_\epsilon-1) (p_\epsilon-2) (p_\epsilon- 3)} {2}  U^{p_\epsilon- 4} \left[{\mathbf b}_1^2\varpi_1^2+ \frac 1{\Lambda^2} {\mathbf b}_1^2f^2\varpi_2^2+ 2 {\mathbf b}_1^2 \frac 1{\Lambda} f \varpi_1 \varpi_2 \right]     \frac{e} {\beta^{\mathfrak{p}}} \mathcal Z
 \nonumber\\[2mm]
\quad &   +   \frac{(p_\epsilon-1) (p_\epsilon-2)(p_\epsilon- 3)} {6}   U^{p_\epsilon- 4}  \left[{\mathbf b}_1\varpi_1+ \frac 1{\Lambda} {\mathbf b}_1f\varpi_2 \right]^3
 \nonumber\\[2mm]
\quad &
+ \frac{(p_\epsilon-1)(p_\epsilon-2) (p_\epsilon- 3)}2 U^{p_\epsilon- 4}  \frac{e^2} {\beta^{2\mathfrak{p}}} \mathcal Z^2\big({\mathbf b}_1\varpi_1+ \frac 1{\Lambda} {\mathbf b}_1f\varpi_2\big) ,
\end{align}
\begin{align}
  \iota^4{\bf N}_{04}
\,=\,& \iota^2 \frac{(p_\epsilon-1) (p_\epsilon-2)} {2} \left[   w_0^{p_\epsilon- 3}
   -    U^{p_\epsilon- 3} \right]\phi_1^2
  -  \iota^3 \frac{(p_\epsilon-1) (p_\epsilon-2) (p_\epsilon- 3)} {2}  U^{p_\epsilon- 4} \phi_1^2  \frac{e} {\beta^{\mathfrak{p}}} \mathcal Z
    \nonumber\\[2mm]
& \quad        +  \iota^3 \frac{(p_\epsilon-1) (p_\epsilon-2)(p_\epsilon- 3)} {6} \left[   w_0^{p_\epsilon- 4}
   -    U^{p_\epsilon- 4} \right]  \phi_1^3     +  \iota^4 \mathcal   {\mathcal N}(x, f, e),
 \end{align}
 and
\begin{align}
\iota^4 \mathcal N(x, f, e)
\,=\,  &   \left( w_0+\iota\phi_1\right)^{p_\epsilon-1}  -w_0^{p_\epsilon-1}  -\iota (p_\epsilon-1)  w_0^{p_\epsilon-2} \phi_1
   \nonumber\\[2mm]
&  -(p_\epsilon-1) (p_\epsilon-2) \left[\iota^2 \frac{1} {2}  w_0^{p_\epsilon- 3} \phi_1^2 +  \iota^3 \frac{(p_\epsilon- 3)} {6}  w_0^{p_\epsilon- 4} \phi_1^3  \right]
   \nonumber\\[2mm]
& +   (p_\epsilon-1)\phi_1 \Bigg\{\iota    \left( w_0^{p_\epsilon-2}- U^{p_\epsilon-2}\right)   - \iota^2 (p_\epsilon-2) U^{p_\epsilon- 3}   \frac{e} {\beta^{\mathfrak{p}}} \mathcal Z
   \nonumber\\[2mm]
& \qquad\qquad\qquad\qquad\qquad- \iota^3 \frac{(p_\epsilon-2) (p_\epsilon- 3)}2 U^{p_\epsilon- 4}  \frac{e^2} {\beta^{2\mathfrak{p}}} \mathcal Z^2 \Bigg\}.
 \end{align}
 In fact, it's easy to verify that
\begin{align}
   |\iota^4{\bf N}_{04}| \le \frac{\iota^4} {\beta^{1+\mathfrak{p}}} e^{-c|x|}  \left( |f| + |e| \right).
 \end{align}

\medskip
   Combining \eqref{sw1new}, \eqref{decompose S3+S4}, \eqref{B4phi1}, \eqref{decompose N0},  we can get
\begin{align}
  S(w_1)=&\,\iota^2 S_{3}
  + \iota^2 S_{41} + \iota^3 S_{42}+ \iota^3 S_{5}
   +\iota^3 S_{6} +  \iota^3 S_{7} +\iota^3S_{8}
   \nonumber\\[2mm]
 &
   + \iota^3\frac1{\beta^{\mathfrak{p}}} \frac{4\pi^2}{{\tilde\ell}^2}  e''\mathcal Z
 + \frac1{\beta^{\mathfrak{p}}} \iota\lambda_0 e{\mathcal Z}   + N_{0\e} \left(\iota\frac{e} {\beta^{\mathfrak{p}}} \mathcal Z\right)
\nonumber\\[2mm]
 &
  \underbrace{ +\iota^4 \frac{4\pi^2} {{\tilde\ell}^2} a_2\frac{\kappa}{\, \beta^{1+\mathfrak{p}} \, }  e'' (x+f)\mathcal Z -\iota^4 \kappa a_2 \frac{ 4\pi^2} {{\tilde\ell}^2} \frac{1}{\, \beta^{1+\mathfrak{p}}\, }  (x+f) f'' e   \mathcal Z'    + \frac{\iota^4}{\beta^{1+\mathfrak{p}} } \hat{\mathcal F}_1\left[\tau, x, \beta^{\mathfrak p} f, \frac{\beta^{\mathfrak p} f''}{\tilde\ell^2} \right] }_{\text{terms involving}~f'', e''}
 \nonumber\\[2mm]
 &
+\underbrace{\iota^2 {\bf B}_{41} +\iota^2{\bf B}_{42}}_{B_4(\iota\phi_1)}
+\underbrace{\iota^2 {\bf N}_{01}+\iota^2 {\bf N}_{02} +\iota^3{\bf N}_{03} + \iota^4 {\bf N}_{04}}_{\mathbb{N}_0(\iota\phi_1)}
\nonumber\\[2mm]
 &
  +  \frac{\iota^4}{\beta^{1+\mathfrak{p}}} \mathcal F_1\left[\tau, x, \beta^{\mathfrak p} f, \frac{\beta^{\mathfrak p} f'}{\tilde\ell}\right]
 +\frac{\iota^4}{\beta^{1+\mathfrak{p}}}  \mathcal G_1\left[\tau, x, \beta^{\mathfrak p} f, \frac{\beta^{\mathfrak p} f'}{\tilde\ell},   e, \frac{e'}{\tilde\ell}\right]
\nonumber\\[2mm]
 & + \iota^3\frac{ 4\pi^2}{{\tilde\ell}^2} {\phi}_{1, \tau\tau}
 +B_{5}(\iota\phi_1)+B_{6}(\iota\phi_1).
\end{align}
 Note that, for  $B_{5}(\iota\phi_1), B_{6}(\iota\phi_1)$,    \eqref{expression of B5phi1}, \eqref{expression of B6phi1} hold true.

\medskip
\subsection{  The third approximate solution via second correction}\

In this part, we want to eliminate the terms
$\iota^2 S_{41}, \iota^2 {\bf N}_{01}, \iota^2 {\bf B}_{41} $
in $S(w_1)$ by finding a solution to the problem
\begin{align}  \label{new second croection}
 & \phi_{2, xx} - \phi_{2} + (p_\epsilon-1) U^{p_\epsilon-2} \phi_{2}
  \,=\, - \left[  S_{41} + {\bf N}_{01}+  {\bf B}_{41} \right],
\quad \phi_{2}(\pm \infty, \tau) \,=\, 0.
\end{align}
In order to solve \eqref{new second croection}, we need to consider the following linear problem
\begin{align}\label{linzezrized 3}
\varpi_{3, xx} - \varpi_{3} + (p_\epsilon-1) U^{p_\epsilon-2} \varpi_{3} \,=\, - xU',
\qquad \varpi_3(\pm \infty) \,=\, 0,
\end{align}
\begin{align}\label{linzezrized 4}
\varpi_{4, xx} - \varpi_{4} + (p_\epsilon-1) U^{p_\epsilon-2} \varpi_{4} \,=\, - x^2U,
\qquad \varpi_4(\pm \infty) \,=\, 0,
\end{align}
\begin{align}\label{linzezrized 6}
\varpi_{5, xx} - \varpi_{5} + (p_\epsilon-1) U^{p_\epsilon-2}\varpi_{5}
\,=\, -x^2U'',
\qquad \varpi_5(\pm \infty) \,=\, 0,
\end{align}
\begin{align}\label{linzezrized 7}
\varpi_{6, xx} - \varpi_{6} + (p_\epsilon-1) U^{p_\epsilon-2}\varpi_{6}
\,=\,-\varpi_{1,x} ,
\qquad \varpi_6(\pm \infty) \,=\, 0,
\end{align}
\begin{align}\label{linzezrized 8}
\varpi_{7, xx} - \varpi_{7} + (p_\epsilon-1) U^{p_\epsilon-2}\varpi_{7}
\,=\, -x\varpi_1,
\qquad \varpi_7(\pm \infty) \,=\, 0,
\end{align}
\begin{align}\label{linzezrized 10}
\varpi_{8, xx} - \varpi_{8} + (p_\epsilon-1) U^{p_\epsilon-2}\varpi_{8}
\,=\, -U^{p_\epsilon- 3}\varpi_1^2,
\qquad \varpi_{8}(\pm \infty) \,=\, 0,
\end{align}
Here we note that $\varpi_1, \varpi_2$ are solutions of problem \eqref{varphi1equation new} and \eqref{varphi2equation new}  respectively.
Recall the expressions of $\iota^2 S_{41}$, $ \iota^2 {\bf N}_{01}, \iota^2 {\bf B}_{41} $ given in \eqref{barS41}, \eqref{N01} and \eqref{B41}  respectively.
Then we can give the solution of \eqref{new second croection} as
\begin{align}\label{definition of e2phi2}
\iota^2\phi_2 = \iota^2\sum_{j=1}^{7}\phi_{2j},
\end{align}
where
\begin{align}\label{phi21-phi210}
& \phi_{21} = \Big[\mathfrak{p}\frac{\beta''} {\beta^3} + (\mathfrak{p}-1) \mathfrak{p} \frac{|\beta'|^2} {\beta^4}  \Big] \varpi_2,
\qquad
\phi_{22} =-  \frac{1}{2}\,  \frac{\iota^2W_{tt}(0, \iota z)}{\beta^4}\varpi_4,
\qquad
\phi_{23} =  \Big[2\mathfrak{p}\frac{|\beta'|^2}{\beta^4}
+\frac{\beta''} {\beta^3}\,-\,\frac{\kappa^2}{\beta^2}\Big]\varpi_3,
\nonumber\\[2mm]
&
\phi_{24} =  \frac{|\beta'|^2}{\beta^4} \varpi_5,
\qquad \qquad
\phi_{25} =   {\mathbf b}_1^2\varpi_6,
\qquad \qquad
\phi_{26} =\frac1{\Lambda}   {\mathbf b}_1^2\varpi_7,
\qquad
\phi_{27} = \frac{(p_\epsilon-1) (p_\epsilon-2)} {2}   {\mathbf b}_1^2\varpi_{8}.
\end{align}
We observe that  $\phi_2$  is an even function of $x$.

\medskip\medskip
As a conclusion, we take
\begin{align}\label{newthirdapproximatesolution}
w_2(x, \tau) \,=\,  w_1(x, \tau)\,+\,\iota^2 \phi_2(x, \tau)
\end{align}
as the {\textbf{third approximate solution}}.
Then we have the new error
\begin{align}\label{sw2new}
S (w_2)
\,=\,
&S(w_1)
+
\iota^2\Big[ \phi_{2, xx} -  \phi_2 + (p_\epsilon-1)U^{p_\epsilon-2} \phi_2 \Big]
+
\left[ \left( w_1+ \iota^2\phi_2\right)^{p_\epsilon-1}  -w_1^{p_\epsilon-1} -  \iota^2 (p_\epsilon-1)U^{p_\epsilon-2} \phi_2  \right]
\nonumber\\[2mm]
& + \iota^4\frac{4\pi^2}{{\tilde\ell}^2}\phi_{2, \tau\tau}
+ B_4(\iota^2\phi_2)+ B_5(\iota^2\phi_2)+ B_6(\iota^2\phi_2)
\nonumber\\[2mm]
=&
\iota^2 S_{3}  + \iota^2 S_{42}  + \iota^3 S_{5} +\iota^3 S_{6}
 +\iota^3 S_{7} +\iota^3 S_{8}
 +\frac 1{\beta^{\mathfrak{p}}} \iota^3\frac{4\pi^2}{{\tilde\ell}^2}  e''\mathcal Z
 +\frac 1{\beta^{\mathfrak{p}} } \iota\lambda_0 e{\mathcal Z}
 \nonumber\\[2mm]
&
+  N_{0\e} \left(\iota\frac{e} {\beta^{\mathfrak{p}}} \mathcal Z\right)
+\iota^2 {\bf N}_{02} +\iota^3{\bf N}_{03} + \iota^4 {\bf N}_{04} +  \iota^2 {\bf B}_{42}
\nonumber\\[2mm]
&
  \underbrace{ +\iota^4 \frac{4\pi^2} {{\tilde\ell}^2} a_2\frac{\kappa}{\, \beta^{1+\mathfrak{p}}\, }  e'' (x+f)\mathcal Z -\iota^4 \kappa a_2 \frac{ 4\pi^2} {{\tilde\ell}^2} \frac{1}{\, \beta^{1+\mathfrak{p}} \, }  (x+f) f'' e   \mathcal Z'    + \frac{\iota^4}{\beta^{1+\mathfrak{p}}} \hat{\mathcal F}_1\left[\tau, x, \beta^{\mathfrak p} f, \frac{\beta^{\mathfrak p} f''}{\tilde\ell^2} \right] }_{\text{terms involving}~f'', e''}
 \nonumber\\[2mm]
&
 + \iota^3 \frac{4\pi^2}{{\tilde\ell}^2} {\phi}_{1, \tau\tau}
+ \iota^4 \frac{4\pi^2}{{\tilde\ell}^2}\phi_{2, \tau\tau}
+B_{5}(\iota\phi_1)+B_{6}(\iota\phi_1)
\nonumber\\[2mm]
&
  +  \frac{\iota^4}{\beta^{1+\mathfrak{p}}} \mathcal F_1\left[\tau, x, \beta^{\mathfrak p} f, \frac{\beta^{\mathfrak p} f'}{\tilde\ell}\right]
 +\frac{\iota^4}{\beta^{1+\mathfrak{p}}}  \mathcal G_1\left[\tau, x, \beta^{\mathfrak p} f, \frac{\beta^{\mathfrak p} f'}{\tilde\ell},   e, \frac{e'}{\tilde\ell}\right]
 \nonumber\\[2mm]
&+ \mathbb{N}_1 (\iota^2 \phi_2)
+ B_4(\iota^2\phi_2)+  B_5 (\iota^2{\phi_2})+  B_6 (\iota^2{\phi_2}),
\end{align}
where  the  linear operators $ B_4, B_5, B_6 $ are given in \eqref{new opreator B4}, \eqref{new opreator B5}, \eqref{B2w-2222new}.

\smallskip\noindent $\clubsuit$In the above, the term $\mathbb{N}_1 (\iota^2 \phi_2)$ is given by
\begin{align}
\mathbb{N}_1 (\iota^2\phi_2) =& \left( w_1+ \iota^2\phi_2\right)^{p_\epsilon-1}  -w_1^{p_\epsilon-1}   - \iota^2 (p_\epsilon-1) w_1^{p_\epsilon-2} \phi_2 +   \iota^2 (p_\epsilon-1) (w_1^{p_\epsilon-2} -U^{p_\epsilon-2}) \phi_2
\nonumber\\[2mm]
=&   \left( w_1+ \iota^2\phi_2\right)^{p_\epsilon-1}  -w_1^{p_\epsilon-1}   - \iota^2 (p_\epsilon-1) w_1^{p_\epsilon-2} \phi_2
\nonumber\\[2mm]  &
+   \iota^3 (p_\epsilon-1) (p_\epsilon-2) U^{p_\epsilon- 3}\Big(\phi_1+ \frac{e} {\beta^{\mathfrak{p}}} \mathcal Z \Big)  \phi_2
\nonumber\\[2mm]  &
+     \iota^2 (p_\epsilon-1)\phi_2 \Big[  (w_1^{p_\epsilon-2} -U^{p_\epsilon-2})  -  (p_\epsilon-2) \iota U^{p_\epsilon- 3} \Big(\phi_1+ \frac{e} {\beta^{\mathfrak{p}}} \mathcal Z\Big) \Big]
\nonumber\\[2mm]
=&   \iota^3 {\bf N}_{11}+\iota^4 {\bf N}_{12},
\label{definition of mathbbN1}
\end{align}
where $\phi_1$ and $\phi_2$ are given in \eqref{definition of phi1} and \eqref{definition of e2phi2}, and
\begin{align}
 \iota^3 {\bf N}_{11}   =  \iota^3 (p_\epsilon-1) (p_\epsilon-2) U^{p_\epsilon- 3}\Big(\phi_1+ \frac{e} {\beta^{\mathfrak{p}}} \mathcal Z \Big)  \phi_2,
\end{align}
\begin{align}
 \iota^4 {\bf N}_{12} = &  \left( w_1+ \iota^2\phi_2\right)^{p_\epsilon-1}  -w_1^{p_\epsilon-1}   - \iota^2 (p_\epsilon-1) w_1^{p_\epsilon-2} \phi_2
 \nonumber\\[2mm]  &
+     \iota^2 (p_\epsilon-1)\phi_2 \Big[  (w_1^{p_\epsilon-2} -U^{p_\epsilon-2})  -  (p_\epsilon-2) \iota U^{p_\epsilon- 3} \Big(\phi_1+ \frac{e} {\beta^{\mathfrak{p}}} \mathcal Z \Big) \Big].
 \end{align}

\smallskip\noindent $\clubsuit$
We note that the terms in $B_4(\iota^2\phi_2)$  are independent of $f'', e''$.
By the definition of $B_5$ given in \eqref{new opreator B5}, we   have
\begin{align}   \label{decomposition of  B5(phi2)}
B_5 (\iota^2{\phi_2})
\,=\,&\iota^4\frac{4\pi^2}{{\tilde\ell}^2}\Big[|f'|^2 \phi_{2, xx}  -f''\phi_{2, x}  -2f' \phi_{2, x\tau}\Big]
+
\iota^4\big(4\mathfrak{p} +2\big) \frac{\pi}{{\tilde\ell}}\frac{\beta'} {\beta^2}\Big[- f' \phi_{2, x} +\phi_{2, \tau}\Big]
\nonumber\\[2mm]
&
+  \iota^4  \frac{4\pi}{{\tilde\ell}} \frac{\beta'}{\beta^2} (x+f)\big[- f' \phi_{2, xx}  +\phi_{2, x\tau}\big]
+
\iota^4\Big[
\frac{\beta''}{\beta^3}\,+\, 2\mathfrak{p} \frac{|\beta'|^2}{\beta^4} \,-\,\frac{\kappa^2}{\beta^2} \Big] (x+f) \phi_{2, x}
\nonumber\\[2mm]
&
+ \iota^4\Big[ \mathfrak{p}(\mathfrak{p}-1)   \frac{|\beta'|^2}{\beta^4} + \mathfrak{p} \frac{\beta''}{\beta^3} \Big] \phi_2
 +  \iota^4  \frac{|\beta'|^2}{\beta^4} (x+f)^2 \phi_{2, xx}
-\frac{1}{2}\, \iota^4\frac{\iota^2W_{tt}(0, \iota z)}{\beta^4} (x+f)^2 \phi_2.
\end{align}
Here we note that the less regular term in $B_5 (\iota^2{\phi_2})$  is  $$-\iota^4\frac{4\pi^2}{{\tilde\ell}^2} f''\phi_{2, x}. $$

\smallskip\noindent $\clubsuit$
We select  the  term  in  $B_6 (\iota^2{\phi_2})$ that    involves   $f''$,    which is
\begin{align}
-\iota^5 \kappa(\iota z) a_2(\iota s, \iota z) \frac{ 4\pi^2} {{\tilde\ell}^2}\frac{1}{\beta} (x+f)f'' \phi_{2, x}.
\end{align}
Furthermore, we  denote
\begin{align}  \label{decomposition of B6phi2}
\frac{\iota^5}{\beta^{{\mathfrak p}+1}} \mathcal F_3\left[\tau, x, \beta^{\mathfrak p} f, \frac{\beta^{\mathfrak p} f'}{\tilde\ell}\right]\,=\,   B_6(\iota^2\phi_2)   +\iota^5 \kappa(\iota z) a_2(\iota s, \iota z) \frac{ 4\pi^2} {{\tilde\ell}^2}\frac{1}{\beta} (x+f)f'' \phi_{2, x}.
\end{align}
Here  the terms  in $\frac{\iota^5}{\beta^{{\mathfrak p}+1}} \mathcal F_3\left[\tau, x, \beta^{\mathfrak p} f, \frac{\beta^{\mathfrak p} f'}{\tilde\ell}\right],$  are all independent of $f'', e''$ and satisfies
\begin{align} \label{estimate of mathcal F3-x}
\frac{\iota^5}{\beta^{{\mathfrak p}+1}} \mathcal F_3\left[\tau, x, \beta^{\mathfrak p} f, \frac{\beta^{\mathfrak p} f'}{\tilde\ell}\right]\, \le\, C e^{-c|x|} \left(1+ |\beta^{\mathfrak p} f| +  \Big|\frac{\beta^{\mathfrak p} f'}{\tilde\ell}\Big| \right).
\end{align}

\medskip
\section{The gluing procedure to derive the projected form of equation \eqref{eq3}}
\label{section8}

\subsection{  The gluing procedure} \label{section 6-1}
The gluing method from \cite{delPKowWei1} will be used to transform problem \eqref{eq3} in $\R^2$
into a projected version on the infinite strip ${\mathfrak S}$, in which we can use the local coordinates $(x, \tau)$.
We define a smooth cut-off function
$
\eta_{3\delta}^{\iota}(s) \,=\, \eta_{3\delta}(\epsilon |s|),
$
where $\eta_{\delta}(t)$ is given by
$$
\eta_{\delta}(t) \,=\, 1\ \mbox{ for }\ 0\leq t\leq \delta
\qquad \mbox{and}\qquad
\eta_{\delta}(t) \,=\, 0
\ \mbox{ for }\ t>2\delta,
$$
for any fixed $\delta\leq \frac{\delta_0}{100}$ with small $\delta_0=M_0$ given in (\ref{fermi}).
For problem \eqref{eq3}, the global approximation $\mathbf{U}(y)$ on the whole space ${\mathbb R}^2$ can be defined simply in the form
\begin{equation}\label{definition of bf U}
\mathbf{U}(y) \,=\,
\eta_{3\delta}^\epsilon(s) \beta^{\mathfrak{p}} (\iota z)w_2(x, \tau), \quad \text{for any}\,y \in \R^2,
\end{equation}
where the local approximate solution $w_2 (x, \tau)$ is defined in \eqref{newthirdapproximatesolution} and the relations between $y$, $(s, z)$ and $ (x,\tau)$
are given in \eqref{definitionofy}, \eqref{fermi},   \eqref{eqszfinal}, \eqref{coordinate-tildes}, \eqref{ztau}, and \eqref{definition of x 2.14}.

\medskip
For a perturbation term $\Phi(y) \,=\, \eta_{3\delta}^{\iota}(\epsilon |s|){\tilde\phi}(y) \,+\, {\tilde\psi}(y)$, the function $v(y) \,=\, \mathbf{U}(y)+\Phi(y)$
satisfies (\ref{eq3}) if the pair $({{\tilde\phi}}, {{ \tilde \psi}})$ satisfies the following coupled system:
\begin{align}
\eta_{3\delta}^{\iota}\mathbf{L}(\tilde \phi)
 \,=\,
\eta_{\delta}^{\epsilon}\Big[-{\mathcal E}
+ \mathbf{N}\big(\eta_{3\delta}^{\iota}{\tilde\phi}+{\tilde\psi}\big)
-(p_\epsilon-1)\, \mathbf{U}^{p_\epsilon-2}  \, {\tilde\psi}\Big],
\label{equivalent-system-1}
\end{align}
\begin{align} \label{equivalent-system-2}
 \Delta \tilde \psi- [1+ \iota^2 W] \tilde \psi + (1-\eta_{\delta}^{\epsilon}) (p_\epsilon-1)\, \mathbf{U}^{p_\epsilon-2}  \, {\tilde\psi} & \,=\, -\iota^2 (\Delta \eta_{3\delta}^{\iota}){ \tilde \phi}-2\epsilon \nabla \eta_{3\delta}^{\iota} \cdot
\nabla{\tilde\phi}
-(1-\eta_{\delta}^{\epsilon}){\mathcal E}
\nonumber
\\[2mm]
&\quad +(1-\eta_{\delta}^{\epsilon})\mathbf{N}\big(\eta_{3\delta}^{\iota}{\tilde\phi}+{\tilde\psi}\big),
\end{align}
 with
\begin{align}
 \mathbf{L}({\tilde\phi}) \,=\, \Delta{\tilde\phi} - [1+ \iota^2 W]{\tilde\phi} +(p_\epsilon-1)\, \mathbf{U}^{p_\epsilon-2} {\tilde\phi},
\end{align}
\begin{align}
 {\mathcal E} \,=\, \Delta {\bf U} -[1+\iota^2 W] {\bf U} + {\bf U}^{p_\epsilon-1},
\end{align}
 and
\begin{align} \label{Noliner bf N1}
 \mathbf{N}\big(\eta_{3\delta}^{\iota}{\tilde\phi}+{\tilde\psi}\big)
 & \,=\,
 -\Big[{\bf U} + \big(\eta_{3\delta}^{\iota}{\tilde\phi}+{\tilde\psi}\big)\Big]^{p_\epsilon-1}  + {\bf U}^{p_\epsilon-1} + (p_\epsilon-1)\, \mathbf{U}^{p_\epsilon-2} \big(\eta_{3\delta}^{\iota}{\tilde\phi}+{\tilde\psi}\big)
 \nonumber
 \\[2mm]
 & =  \left[- (p_\epsilon-1)\, \left(\mathbf{U} + t  \big(\eta_{3\delta}^{\iota}{\tilde\phi}+{\tilde\psi}\big)\right) ^{p_\epsilon-2}  + (p_\epsilon-1)\, \mathbf{U}^{p_\epsilon-2}  \right]  \big(\eta_{3\delta}^{\iota}{\tilde\phi}+{\tilde\psi}\big),
\end{align}
for some $t\in (0,1)$.
The term ${ {\mathcal E}}$ can be locally recast in $(x, \tau)$ coordinate system by the relation
\begin{equation}\label{errorrelationinterior}
{{\mathcal E}} \,=\,   \beta^{\mathfrak{p}+2}\, S(w_2),  \quad \text{for $|\iota s|\leq 10\delta$},
\end{equation}
where for the explicit expression of $S(w_2)$, the reader can refer to \eqref{sw2new}.

 First, given a small $ \tilde {\phi}$ satisfying the following decay property
\begin{equation}\label{decay-assumption}
|\nabla {\tilde\phi}|\,+\, |{\tilde\phi}|\,\leq\, e^{- { \upsilon_1\, \delta}/{\iota}}, \, \, \, \, {\rm for }\, \, \, |s|>\delta/\iota,
\end{equation}
where $\upsilon_1$ is a very small positive constant, we will solve (\ref{equivalent-system-2}).
Note that the nonlinear operator $\mathbf{N}$ has a power-like behaviour with power greater than one, and a direct application of contraction mapping principle yields
that (\ref{equivalent-system-2}) has a unique (small) solution ${\tilde\psi} \,=\, {\tilde\psi}({\tilde\phi})$ with
\begin{equation}\label{contraction-varphi}
\|{\tilde\psi}({ \tilde \phi})\|_{L^{\infty}}\,\leq \, C\, \iota \lf[\|{\tilde\phi}\|_{L^{\infty}(|s|>\frac{\delta}{\iota})}\,+\, \|\nabla {\tilde\phi}\|_{L^{\infty}(|s|>\frac{\delta}{\iota})}\ri]\,+\,
e^{-\frac{ \upsilon_2 \delta}{\iota}},
\end{equation}
where $|s|>\delta/\iota$ denotes the complement of $\delta/\iota$-neighbourhood of
$\Gamma_\iota$ and $\upsilon_2$ denotes some small constant.
Moreover, the nonlinear operator $\tilde {\psi}$ satisfies a Lipschitz condition of the form
\begin{equation}\label{varphi-lip}
\|{\tilde\psi}({\tilde\phi}_1)-{\tilde\psi}({ \tilde \phi}_2)\|_{L^{\infty}}
\leq
C\iota \lf[ \|{\tilde\phi}_1-{\tilde\phi}_2\|_{L^{\infty}(|s|>\delta/\iota)}+\|\nabla {\tilde\phi}_1-\nabla {\tilde\phi}_2\|_{L^{\infty}(|s|>
\delta/\iota)}\ri].
\end{equation}

Therefore, after solving (\ref{equivalent-system-2}),
we can concern (\ref{equivalent-system-1}) as a local nonlinear problem involving ${\tilde\psi} \,=\, { \tilde \psi}({\tilde\phi})$,
which can be expressed in local coordinates $(x, \tau)$.
Using Fermi coordinates  \eqref{fermi},  \eqref{definitionofy} and  the calculations in Section \ref{Section2}, especially   \eqref{eqszfinal}, \eqref{coordinate-tildes}, \eqref{ztau} and \eqref{definition of x 2.14},   we give the following lemma.
\begin{lemma}
By the relation
\begin{align}
{\tilde\phi}(y) =   \beta^{\mathfrak{p}}(\iota z) \phi^*(x, \tau),  \quad \text{for $|\iota s|\leq 10\delta$},
\label{relation betwen tildephi and phi}
\end{align}
we have
\begin{align}
\mathcal L \left(\phi^*\right)
\,=\,
\frac{1}{\beta^{\mathfrak{p}+2}}\mathbf{L}({\tilde\phi}),  \quad \text{for $|\iota s|\leq 10\delta$},
\end{align}
where
\begin{align}\label{expression of mathcal L}
\mathcal L (\phi^*) \,=\,   L_2(\phi^*) \ +\ \chi(\e|s|) \big[\, B_4(\phi^*)+ B_5(\phi^*)+ B_6(\phi^*)\,\big],
\end{align}
with
\begin{align}\label{definition of L2operator}
 L_2(\phi^*)  = \frac{\iota^2 4\pi^2}{{\tilde\ell}^2} \phi^*_{\tau\tau} + \phi^*_{xx} - \phi^*+  ({\mathfrak p}-1)w_2^{{\mathfrak p}-2} \phi^*.
\end{align}
In the above,  the smooth cut-off function $\chi (r)$  is given by
$$
\chi(t) \,=\, 1, \quad \text{if}~t<10\delta, \qquad \chi(t) \,=\, 0, \quad\text{if}\ t>20\delta.
$$
For the local form of the nonlinear part, we denote by the notation
\begin{align}\label{Nmathcal0}
{\mathcal N}(\phi^*)
\,=\,
\frac{1}{\beta^{\mathfrak{p}+2}}\left[\, \mathbf{N}\big(\eta_{3\delta}^{\iota}{\tilde\phi}+{\tilde\psi}\big) -(p_\epsilon-1)\, \mathbf{U}^{p_\epsilon-2}  \, {\tilde\psi}(\tilde\phi) \,\right],
 \quad \text{for $|\iota s|\leq 10\delta$}.
\end{align}
\qed
\end{lemma}

\medskip
Note that the approximate solution $\mathbf{U}$ has unknown parameter pair $(f, e)$,
see \eqref{constraint-f}, \eqref{constraint-e} and \eqref{mathcalF}.
We will deal with the following projected problem:
for given $(f, e)$ in $\hat{\mathcal F}$, finding
functions $\phi^*\in H^2({\mathfrak S})$, \, and $c(\tau), d(\tau) \in L^2(0, 2\pi)$,
such that
\begin{align} \label{project 6.13}
\mathcal L (\phi^*) \,=\, \eta_\delta^\iota(s) \Big[-\frac{1}{\beta^{\mathfrak{p}}}{\mathcal E} + {\mathcal N}(\phi^*)\Big]
+ c(\tau) \chi(\iota|s|) U'(x)+ d(\tau) \chi(\epsilon|s|){\mathcal Z}(x), \quad \text{for}\,(x, \tau)\in \mathfrak S,
\end{align}

\begin{equation}
\phi^*(x, 0) \,=\, \phi^*(x, 2\pi), \qquad \phi^*_\tau(x, 0) \,=\, \phi^*_\tau(x, 2\pi),
\quad
-\infty<x<+\infty,
\label{periodic condition}
\end{equation}

\begin{equation}
\label{projectedproblem3}
 \int_{-\infty}^{+\infty} \phi^*(x, \tau){\mathcal Z}(x)  \,{\mathrm d}x
 \,=\, \int_{-\infty}^{+\infty} \phi^*(x, \tau) U'(x)\,{\mathrm d}x
 \,=\, 0, \quad 0<\tau<2\pi,
\end{equation}
where the definition of $\mathfrak S$ has been given in \eqref{huaxies}.

\medskip
\subsection{  Solvability of the linear problem}

We first consider the linear problem
\begin{align}  \label{linearphi1}
 \frac{\iota^2 4\pi^2}{{\tilde\ell}^2} \phi_{\tau\tau} + \phi_{xx} - \phi+(p_\epsilon-1)U^{p_\epsilon-2}  \phi = h,
\quad \text{for}\,(x, \tau)\in \mathfrak S,
\end{align}
\begin{align}
 \phi(x, 0) = \phi(x, 2\pi), \quad \phi_\tau(x, 0) = \phi_\tau(x, 2\pi), \quad \text{for}~x\in \R,
\end{align}
\begin{align}  \label{linearphi3}
 \int_\R\phi(x, \tau)\,\mathcal Z(x)\,{\mathrm d}x =   \int_\R\phi(x, \tau)\,U'(x)\,{\mathrm d}x \,=\, 0,  \quad \text{for}~\tau\in (0, 2\pi),
\end{align}
and then have the following lemma.
\begin{lemma}\label{linear lemma6-1}
  Suppose $\phi$ solves \eqref{linearphi1}-\eqref{linearphi3}. Then it holds true
\begin{align} \label{norm h2 mathfrak S estimate}
  \|\phi\|_{H^2_*(\mathfrak S)}   \leq  C   \|h\|_{L^2(\mathfrak S)},
\end{align}
where the norm is given by
\begin{align}  \label{norm h2 mathfrak S}
\|\phi\|_{H^2_*(\mathfrak S)}
= \left({\displaystyle\int}_{\mathfrak S} \Bigg[\, |\phi|^2+ |\phi_{x}|^2+ |\phi_{xx}|^2 +\Big|\frac{\iota}{\tilde\ell}\phi_{x\tau}\Big|^2
+ \Big|\frac{\iota}{\tilde\ell}\phi_{\tau}\Big|^2 +  \Big|\frac{\iota^2}{\tilde\ell^2}\phi_{\tau\tau}\Big|^2 \, \Bigg] \,{\mathrm d}x  {\mathrm d} \tau \right)^{\frac12}.
\end{align}
Moreover,
we have
\begin{align} \label{betanorm h2 mathfrak S estimate}
  \|\beta^{\mathfrak{p}} \phi\|_{H^2_*(\mathfrak S)}   \leq  C   \|\beta^{\mathfrak{p}}  h\|_{L^2(\mathfrak S)}.
\end{align}

\end{lemma}

\begin{proof}
 We express $h$ and $\phi$ as
\begin{align}  \label{fourier-h}
 h(x, \tau) = \sum_{k=0}^\infty \big[h_{1k} (x)\cos k\tau +  h_{2k} (x)\sin k\tau\big],
 \\[2mm]
 \label{fourier-phi}\phi(x, \tau) = \sum_{k=0}^\infty \big[\phi_{1k} (x)\cos k\tau +  \phi_{2k} (x)\sin k\tau\big].
\end{align}
 Then we get
\begin{align}  \label{linear philk1}
 - \frac{\iota^2k^2 4\pi^2}{{\tilde\ell}^2}\phi_{\ell k}+ \phi_{\ell k, xx} -  \phi_{\ell k}+(p_\epsilon-1)U^{p_\epsilon-2} \phi_{\ell k} = h_{\ell k}\quad \mbox{in } {\mathbb R},
\end{align}
\begin{align} \label{linear philk2}
 \int_\R\phi_{\ell k}{\mathcal Z}  \,{\mathrm d}x=   \int_\R\phi_{\ell k} U'  \,{\mathrm d}x \,=\, 0.
\end{align}
 Consider the bilinear form
\begin{align*}
 B(\psi, \psi) = \int_\R\Big[|\psi_x|^2+ (1-(p_\epsilon-1)U^{p_\epsilon-2} )|\psi|^2 \Big]  \,{\mathrm d}x.
\end{align*}
Since \eqref{linear philk2},   we can get
\begin{align}  \label{bilinear 1}
 \|\phi_{\ell k}\|^2_{L^2(\R)} +  \|\phi_{\ell k, x}\|^2_{L^2(\R)} \leq C B(\phi_{\ell k}, \phi_{\ell k}).
\end{align}
Then  from \eqref{linear philk1}, we claim that
\begin{align}  \label{lin7-12}
\Bigg[1+\frac{\iota^4k^4}{{\tilde\ell}^4}\Bigg]   \|\phi_{\ell k}\|^2_{L^2(\R)}
+  \Bigg[1+\frac{\iota^2k^2}{{\tilde\ell}^2}\Bigg]  \|\phi_{\ell k, x}\|^2_{L^2(\R)} \leq C\|h_{\ell k}\|^2_{L^2(\R)}.
\end{align}
In fact, from \eqref{linear philk1}, we have
\begin{align*}
  \Bigg[1+\frac{\iota^2k^2}{{\tilde\ell}^2}\Bigg]   \|\phi_{\ell k}\|^2_{L^2(\R)} +  \|\phi_{\ell k, x}\|^2_{L^2(\R)} \leq C\|h_{\ell k}\|^2_{L^2(\R)}.
\end{align*}
  If $ \frac{\iota^2k^2}{{\tilde\ell}^2} \leq 1$, then we have   $  \frac{\iota^4k^4}{{\tilde\ell}^4} \leq \frac{\iota^2k^2}{{\tilde\ell}^2}$. We can directly obtain \eqref{lin7-12}.
  If   $ \frac{\iota^2k^2}{{\tilde\ell}^2} \geq 1$, then from \eqref{linear philk1},  we have
\begin{align*}
  \frac{\iota^2k^2}{{\tilde\ell}^2} B(\phi_{\ell k}, \phi_{\ell k}) + \frac{\iota^4k^4}{{\tilde\ell}^4} \|\phi_{\ell k}\|^2_{L^2(\R)} \leq  C\|h_{\ell k}\|^2_{L^2(\R)} + (1-c)\frac{\iota^4k^4}{{\tilde\ell}^4} \|\phi_{\ell k}\|^2_{L^2(\R)},
\end{align*}
  i.e,
\begin{align*}
  \frac{\iota^2k^2}{{\tilde\ell}^2} B(\phi_{\ell k}, \phi_{\ell k}) + c\frac{\iota^4k^4}{{\tilde\ell}^4} \|\phi_{\ell k}\|^2_{L^2(\R)} \leq  C\|h_{\ell k}\|^2_{L^2(\R)},
\end{align*}
with $c\in (0, 1)$.
 Then by \eqref{bilinear 1}, we get  \eqref{lin7-12}.
By \eqref{lin7-12} and \eqref{linear philk1}, we can directly get
\begin{align}\label{second order estiamte}
\|\phi_{\ell k, xx}\|^2_{L^2(\R)}  \leq  C\|h_{\ell k}\|^2_{L^2(\R)}.
\end{align}

 From the expression of $\phi$ as in \eqref{fourier-phi}, we have
\begin{align*}
\phi_x(x, \tau)  &= \sum_{k=0}^\infty \big[\phi_{1k, x} (x)\cos k\tau +  \phi_{2k, x} (x)\sin k\tau\big],
\\[2mm]
\phi_{xx}(x, \tau)  &= \sum_{k=0}^\infty \big[\phi_{1k, xx} (x)\cos k\tau +  \phi_{2k, xx} (x)\sin k\tau\big],
\\[2mm]
\phi_\tau(x, \tau) & = k\sum_{k=0}^\infty \big[-\phi_{1k} (x)\sin k\tau +  \phi_{2k} (x)\cos k\tau\big],
\\[2mm]
\phi_{\tau\tau}(x, \tau) & = -k^2\sum_{k=0}^\infty \big[\phi_{1k} (x)\cos k\tau +  \phi_{2k} (x)\sin k\tau\big],
\\[2mm]
\phi_{x\tau}(x, \tau) & = k\sum_{k=0}^\infty \big[-\phi_{1k, x} (x)\sin k\tau +  \phi_{2k,x} (x)\cos k\tau\big].
\end{align*}
Therefore, we can derive that
\begin{align*}
  \|\phi\|_{L^2(\mathfrak S)}^2  & =\pi  \sum_{k=0}^\infty \left[\|\phi_{1k}\|^2_{L^2(\R)}+  \|\phi_{2k}\|^2_{L^2(\R)}\right],
  \\[2mm]
   \|\phi_x\|_{L^2(\mathfrak S)}^2  & = \pi\sum_{k=0}^\infty \left[\|\phi_{1k, x}\|^2_{L^2(\R)}+  \|\phi_{2k, x}\|^2_{L^2(\R)}\right],
   \\[2mm]
 \|\phi_{xx}\|_{L^2(\mathfrak S)}^2   &= \pi \sum_{k=0}^\infty \left[\|\phi_{1k, xx}\|^2_{L^2(\R)}+  \|\phi_{2k, xx}\|^2_{L^2(\R)}\right],
 \\[2mm]
  \|\phi_\tau\|_{L^2(\mathfrak S)}^2   &= \pi \sum_{k=0}^\infty \left[k^2\|\phi_{1k}\|^2_{L^2(\R)}+  k^2\|\phi_{2k}\|^2_{L^2(\R)}\right],
  \\[2mm]
   \|\phi_{\tau\tau}\|_{L^2(\mathfrak S)}^2   &=\pi  \sum_{k=0}^\infty \left[k^4\|\phi_{1k}\|^2_{L^2(\R)}+  k^4\|\phi_{2k}\|^2_{L^2(\R)}\right],
   \\[2mm]
     \|\phi_{x\tau}\|_{L^2(\mathfrak S)}^2   &= \pi  \sum_{k=0}^\infty \left[k^2\|\phi_{1k, x}\|^2_{L^2(\R)}+  k^2\|\phi_{2k,x}\|^2_{L^2(\R)}\right].
\end{align*}
  Then from \eqref{lin7-12}, \eqref{second order estiamte},  we have
\begin{align*}
     &\|\phi\|_{L^2(\mathfrak S)}^2 +   \|\phi_x\|_{L^2(\mathfrak S)}^2  +  \|\phi_{xx}\|_{L^2(\mathfrak S)}^2+  \frac{\iota^2}{\tilde\ell^2} \|\phi_\tau\|_{L^2(\mathfrak S)}^2  +     \frac{\iota^4}{\tilde\ell^4} \|\phi_{\tau\tau}\|_{L^2(\mathfrak S)}^2  +   \frac{\iota^2}{\tilde\ell^2} \|\phi_{x\tau}\|_{L^2(\mathfrak S)}^2
\nonumber\\[2mm]
 & \leq C\sum_{k=0}^\infty  \|h_{\ell k}\|^2_{L^2(\R)} \leq \|h\|_{L^2(\mathfrak S)}^2.
\end{align*}

 Next, we consider the proof of \eqref{betanorm h2 mathfrak S estimate}.
By denoting \[ \hat\phi = \beta^{\mathfrak{p}} \phi,  \]  we have,
for $\,(x, \tau)\in \mathfrak S$,
\begin{align*}
 \frac{\iota^2 4\pi^2}{{\tilde\ell}^2} \Big(\frac{\hat\phi} {\beta^{\mathfrak{p}}} \Big)_{\tau\tau} +   \frac1{\beta^{\mathfrak{p}}} \Big[\hat\phi_{xx} - \hat\phi+(p_\epsilon-1)U^{p_\epsilon-2}  \hat\phi  \Big]= h,
\end{align*}
i.e,
\begin{align*}
 \frac1{\beta^{\mathfrak{p}}}  \Big[  \frac{\iota^2 4\pi^2}{{\tilde\ell}^2} \hat\phi_{\tau\tau}  + \hat\phi_{xx} - \hat\phi+(p_\epsilon-1)U^{p_\epsilon-2}  \hat\phi  \Big] +  \frac{\iota^2 4\pi^2}{{\tilde\ell}^2}\Big[2\hat\phi_\tau \Big(\frac1{\beta^{\mathfrak{p}}} \Big)_\tau + \hat\phi \Big(\frac1{\beta^{\mathfrak{p}}} \Big)_{\tau\tau} \Big] = h,
\end{align*}
i.e,
\begin{align*}
  \frac1{\beta^{\mathfrak{p}}}\Big[  \frac{\iota^2 4\pi^2}{{\tilde\ell}^2} \hat\phi_{\tau\tau}  + \hat\phi_{xx} - \hat\phi+(p_\epsilon-1)U^{p_\epsilon-2}  \hat\phi  \Big] +  \frac{\iota^2 4\pi^2}{{\tilde\ell}^2}\Big[2\hat\phi_\tau \Big(-\mathfrak{p}\frac{\beta'}{\beta^{1+\mathfrak{p}}}  \frac{\tilde \ell} {2\pi \beta}\Big) + \hat\phi   \Big(-\frac{\beta''} {\beta^{2+\mathfrak{p}} } +(2+ \mathfrak{p}) \frac{|\beta'|^2} {\beta^{3+  \mathfrak{p}}} \Big)\frac{\tilde \ell^2} {4\pi^2}  \Big] = h,
\end{align*}
with
\begin{align*}
\Big(\frac1{\beta^{\mathfrak{p}}} \Big)_\tau =  -\mathfrak{p}\frac{\beta'}{\beta^{1+\mathfrak{p}}}  \frac{\tilde \ell} {2\pi \beta},
\qquad
\Big(\frac1{\beta^{\mathfrak{p}}} \Big)_{\tau\tau}
=  \Big(-\frac{\beta''} {\beta^{2+\mathfrak{p}} } +(2+ \mathfrak{p}) \frac{|\beta'|^2} {\beta^{3+  \mathfrak{p}}} \Big)\frac{\tilde \ell^2} {4\pi^2}.
\end{align*}
Then we get
\begin{align*}
\frac{\iota^2 4\pi^2}{{\tilde\ell}^2} \hat\phi_{\tau\tau}  + \hat\phi_{xx} - \hat\phi+(p_\epsilon-1)U^{p_\epsilon-2}  \hat\phi    = \tilde h,
\end{align*}
with
\begin{align*}
 \tilde h & = \beta h -\frac{\iota^2 4\pi^2}{{\tilde\ell}^2}\Big[2\hat\phi_\tau \Big(- \mathfrak{p}\frac{\beta'}{\beta^{2}}  \frac{\tilde \ell} {2\pi}\Big) + \hat\phi  \Big(-\frac{\beta''} {\beta^{2} } +(2+ \mathfrak{p}) \frac{|\beta'|^2} {\beta^{3}} \Big)\frac{\tilde \ell^2} {4\pi^2}  \Big]
  \nonumber
\\[2mm]
& = \beta^{\mathfrak{p}}  h - \iota^2\Big[2\hat\phi_\tau  \Big(- \mathfrak{p}\frac{\beta'}{\beta^{2}}  \frac {2\pi} {\tilde \ell}\Big)  + \hat\phi\Big(-\frac{\beta''} {\beta^{2} } +(2+ \mathfrak{p}) \frac{|\beta'|^2} {\beta^{3}} \Big)\Big].
\end{align*}
Then from \eqref{norm h2 mathfrak S estimate}, we get   \eqref{betanorm h2 mathfrak S estimate}  directly.
\end{proof}

 We next consider the following  projected  problem
\begin{align}  \label{pro-linearphi1}
 \frac{\iota^2 4\pi^2}{{\tilde\ell}^2} \phi_{\tau\tau} + \phi_{xx} - \phi+(p_\epsilon-1)U^{p_\epsilon-2}  \phi = h+c(\tau) U'(x)+ d(\tau){\mathcal Z}(x),
 \quad \text{for}\,(x, \tau)\in \mathfrak S,
\end{align}
\begin{align}
 \phi(x, 0) = \phi(x, 2\pi), \quad \phi_\tau(x, 0) = \phi_\tau(x, 2\pi), \quad \text{for}~x\in \R,
\end{align}
\begin{align}  \label{pro-linearphi3}
 \int_\R\phi(x,\tau){\mathcal Z}(x) \,{\mathrm d}x =   \int_\R\phi(x,\tau)  U'(x)\,{\mathrm d}x \,=\, 0,  \quad \text{for}~\tau\in (0, 2\pi).
\end{align}

\begin{lemma} \label{pro-linear lemma6-1}
Problem \eqref{pro-linearphi1}-\eqref{pro-linearphi3}  has a unique solution.  And, it holds true
\begin{align} \label{pro estimate H2-l2}
 \|\phi\|_{H^2_*(\mathfrak S)}   \leq  C   \|h\|_{L^2(\mathfrak S)}.
\end{align}
Moreover, we have
\begin{align} \label{pro estimate H2-l2-beta}
  \|\beta^{\mathfrak{p}} \phi\|_{H^2_*(\mathfrak S)}   \leq  C   \|\beta^{\mathfrak{p}}  h\|_{L^2(\mathfrak S)}.
\end{align}
\end{lemma}
\begin{proof}
Similar to the analysis in Lemma \ref{linear lemma6-1},  by the expansions \eqref{fourier-h}-\eqref{fourier-phi}, we  transfer  \eqref{pro-linearphi1}-\eqref{pro-linearphi3}  to the following
\begin{align}  \label{pro-linear philk1}
 - \frac{\iota^2k^2 4\pi^2}{{\tilde\ell}^2}\phi_{\ell k}+ \phi_{\ell k, xx} -  \phi_{\ell k}+(p_\epsilon-1)U^{p_\epsilon-2} \phi_{\ell k} = h_{\ell k} +c_{\ell k}U'+ d_{\ell k}\mathcal{Z}
 \quad \mbox{in } \R,
\end{align}
\begin{align} \label{pro-linear philk2}
 \int_\R\phi_{\ell k}{\mathcal Z}  \,{\mathrm d}x=   \int_\R\phi_{\ell k} U'  \,{\mathrm d}x \,=\, 0.
\end{align}
By Fredholm's theory, we know that problem \eqref{pro-linear philk1}-\eqref{pro-linear philk2}  is solvable if
\begin{align*}
c_{\ell k} = - \frac{\displaystyle\int_\R  h_{\ell k}  U'  \,{\mathrm d}x} {\displaystyle\int_\R |U'|^2  \,{\mathrm d}x},
\qquad
d_{\ell k} = - \frac{\displaystyle\int_\R  h_{\ell k} {\mathcal Z}  \,{\mathrm d}x} {\displaystyle\int_\R{\mathcal Z}^2  \,{\mathrm d}x}.
\end{align*}
Then define
\begin{align}
\phi(x, \tau) = \sum_{k=0}^\infty \big[\phi_{1k} (x)\cos k\tau +  \phi_{2k} (x)\sin k\tau\big],
\end{align}
and correspondingly
\begin{align}
c(\tau) & = \sum_{k=0}^\infty \big[c_{1k}\cos k\tau +  c_{2k}\sin k\tau\big],
\\
d(\tau) & = \sum_{k=0}^\infty \big[d_{1k}\cos k\tau +  d_{2k}\sin k\tau\big].
\end{align}
 By Lemma \ref{linear lemma6-1}, we have
\begin{align} \label{es1}
 \|\phi\|_{H^2_*(\mathfrak S)}   \leq  C \Big(   \|h\|_{L^2(\mathfrak S)} + \|c(\tau) U'\|_{L^2(\mathfrak S)} +\| d(\tau){\mathcal Z}\|_{L^2(\mathfrak S)} \Big).
\end{align}
There holds
\begin{align}  \label{es2}
 \|c(\tau) U'(x)\|_{L^2(\mathfrak S)}^2
 &  \leq C \sum_{k=0}^\infty(c_{1k}^2+ c_{2k}^2)
 \leq C \sum_{k=0}^\infty \int_\R|h_{\ell k}|^2  \,{\mathrm d}x \leq C\|h\|^2_{L^2(\mathfrak S)}.
\end{align}
 Similarly, we have
\begin{align} \label{es3}
  \|d(\tau){\mathcal Z}(x)\|_{L^2(\mathfrak S)}^2   \leq C\|h\|^2_{L^2(\mathfrak S)}.
\end{align}
By \eqref{es1}, \eqref{es2}, \eqref{es3}, we get \eqref{pro estimate H2-l2}.

The proof of \eqref{pro estimate H2-l2-beta} is similar to \eqref{betanorm h2 mathfrak S estimate}. Here we omit for concise.
\end{proof}

\medskip
 Consider the projected problem
\begin{align}  \label{L2-pro-linearphi1}
  \mathcal L (\phi^*) = h+c(\tau) U'(x)+ d(\tau){\mathcal Z}(x),
 \quad \text{for}\,(x, \tau)\in \mathfrak S,
\end{align}
\begin{align}
 \phi^*(x, 0) = \phi^*(x, 2\pi), \qquad \phi^*_\tau(x, 0) = \phi^*_\tau(x, 2\pi), \quad \text{for}~x\in \R,
\end{align}
\begin{align}  \label{L2-pro-linearphi3}
 \int_\R\phi^*(x,\tau){\mathcal Z}(x) \,{\mathrm d}x =   \int_\R\phi^*(x,\tau)  U'(x)\,{\mathrm d}x \,=\, 0,  \quad \text{for}~\tau\in (0, 2\pi).
\end{align}
The operator $\mathcal L$ is given in \eqref{expression of mathcal L}.
The invertibility of $\mathcal L$   in the projection space can be easily obtained from Lemma \ref{pro-linear lemma6-1}
and the Contraction Mapping Theorem.

\begin{lemma} \label{pro-linear lemma6-1-mathcal L}
Problem \eqref{L2-pro-linearphi1}-\eqref{L2-pro-linearphi3}  has a unique solution with the estimate
\begin{align} \label{pro estimate H2-l2-mathcal}
 \|\phi^*\|_{H^2_*(\mathfrak S)}   \leq  C   \|h\|_{L^2(\mathfrak S)},
\end{align}
for some universal constant $C$.
Moreover,
we have
\begin{align} \label{pro estimate H2-l2-mathcal-beta}
  \|\beta^{\mathfrak{p}} \phi^*\|_{H^2_*(\mathfrak S)}   \leq  C   \|\beta^{\mathfrak{p}} h\|_{L^2(\mathfrak S)}.
\end{align}
\end{lemma}
\begin{proof}
By Lemma \ref{pro-linear lemma6-1}, combined with the Contraction Mapping Theorem, the proof can be completed.
\end{proof}

\subsection{  The $L^2$ estimate of $E_{11}$ in $\mathfrak S$}
Note that, if  $|s|\leq {10\delta}/\e$, there holds
\begin{align}
 \frac{1}{\beta^{\mathfrak{p}+2}} {{\mathcal E}} \,=\, S(w_2).
\end{align}
Denote
\begin{align*}
\mathfrak S_\delta = \mathfrak S \cap \Big\{ |s| < \frac {10\delta}\epsilon \Big\}.
\end{align*}
We decompose $S(w_2)$ as
\begin{align}
S(w_2) = E_{11} + E_{12},
\end{align}
where
\begin{align}  \label{E11+E12}
 E_{12} =   \iota^3\frac{4\pi^2}{{\tilde\ell}^2}  e''\mathcal Z+ \e\lambda_0 e{\mathcal Z}.
\end{align}
We give the following  lemma.
\begin{lemma}
By the expression of $S(w_2)$ in \eqref{sw2new}  and decomposition \eqref{E11+E12}, we   easily get
\begin{align} \label{estimate of L2mathfrak S-E11}
\|\beta^{\mathfrak p} E_{11}\|_{L^2(\mathfrak S_\delta)}   \leq C\left( \iota^3+ \iota^2\|f\|_* + \iota^2\|e\|_{**}\right),
\end{align}
and
\begin{align} \label{Lipschitz of E11-E12 in mathfrak S}
\|\beta^{\mathfrak p}\left(E_{11} (f_1, e_1) - E_{11} (f_2, e_2)\right)\|_{L^2(\mathfrak S_\delta)}  \leq C  \iota^2 \left[\|f_1- f_2\|_*+ \|e_1- e_2\|_{**} \right],
\end{align}
for any $(f_1, e_1)$ and $(f_2, e_2)$ satisfy \eqref{constraint-f}, \eqref{constraint-e}.
\end{lemma}

\begin{proof}
 The proof follows directly from \eqref{sw2new} and \eqref{E11+E12}, so it is omitted here.
\end{proof}

\medskip
\subsection{  Solving the nonlinear intermediate problem}\

In this part, by using the facts given in the above, we will solve problem \eqref{project 6.13}-\eqref{projectedproblem3}.
The term $E_{12}$ has precisely the form $d \chi{\mathcal Z}$ and can be absorbed in the term.
Therefore,  problem  \eqref{project 6.13}-\eqref{projectedproblem3}  can be replaced by the following
\begin{align}  \label{projectedproblem1-refined-e11}
\mathcal L (\phi^*) \,=\, \eta_\delta^\iota(s) \Big[- E_{11}+ {\mathcal N}(\phi^*)\Big]
+ c(\tau) \chi(\epsilon|x|) U'(x)
+ d(\tau) \chi(\epsilon|x|){\mathcal Z}(x), \quad \text{for}\,(x, \tau)\in \mathfrak S,
\end{align}

\begin{equation}
{\phi^* }(x, 0) \,=\, {\phi^* }(x, 2\pi), \quad \phi^*_\tau(x, 0) \,=\, \phi^* _\tau(x, 2\pi),
\quad
-\infty<x<+\infty,
\end{equation}

\begin{equation}
\label{projectedproblem3-refined-aa}
 \int_{-\infty}^{+\infty} \phi^*(x, \tau){\mathcal Z}(x)  \,{\mathrm d}x \,=\, \int_{-\infty}^{+\infty} \phi^*(x, \tau) U'(x)\,{\mathrm d}x \,=\, 0, \quad 0<\tau<2\pi.
\end{equation}

\medskip
Let ${\mathcal T}$ be the operator defined by  Lemma \ref{pro-linear lemma6-1-mathcal L}.
 Then problem \eqref{projectedproblem1-refined-e11}-\eqref{projectedproblem3-refined-aa} is equivalent to the fixed-point problem
\begin{align}  \label{mathbb-A}
 \phi^* ={\mathcal T}\big(-E_{11} +{\mathcal N}(\phi^*)\big) : =\mathbb A (\phi^*).
\end{align}
 By the effect of cut-off function $\chi(\e|x|)$, using the maximum principle and Sobolev embedding,  we can get
\begin{align}
 |\left(\beta^{\mathfrak{p}} \phi^*\right)(x, \tau)| \,+\, |\nabla\left(\beta^{\mathfrak{p}} \phi^*\right)(x, \tau)| \,\leq\,    \|\beta^{\mathfrak{p}} \phi^*\|_{H^2_*(\mathfrak S)}  e^{-\frac{2\delta} \iota}, \quad \text{for}~|s|\geq \frac{40\delta}{\iota},
\end{align}
  where  the norm $\|\cdot\|_{H^2_*(\mathfrak S)}  $ is defined in \eqref{norm h2 mathfrak S}.
  In fact, from  the expression of $\mathcal L (\phi^*)$ given in \eqref{expression of mathcal L}, we know that, if $|\iota s|>20 \delta$,  $\phi^*$ satisfies an equation of the form
\begin{align} \label{asymptotic of phi equation}
\frac{\iota^2 4\pi^2}{{\tilde\ell}^2} \phi^*_{\tau\tau} + \phi^*_{xx} - \big(1+ o(1)\big)\phi^* = 0.
\end{align}
Then
we get
\begin{align}
\frac1{\beta^{\mathfrak{p}}}\left[\frac{\iota^2 4\pi^2}{{\tilde\ell}^2} \left(\beta^{\mathfrak{p}} \phi^*\right)_{\tau\tau}
+2\iota^2\Big(- \frac{\beta'}{\beta^2}  \frac{2\pi}{\tilde \ell} \Big)\left(\beta^{\mathfrak{p}} \phi^*\right)_\tau+   \left(\beta^{\mathfrak{p}} \phi^*\right)_{xx} - \big(1+ o(1)\big) \left(\beta^{\mathfrak{p}} \phi^*\right) \right] = 0.
\end{align}
We note
\begin{align}
\left(\beta^{\mathfrak{p}} \phi^*\right)_{\tau\tau}    &= \beta^{\mathfrak{p}} \phi^*_{\tau\tau}  + 2 (\beta^{\mathfrak{p}})_{\tau}   \phi^*_{\tau}  +  (\beta^{\mathfrak{p}})_{\tau\tau}\phi^*
\nonumber
\\[2mm]  &=  \beta^{\mathfrak{p}} \phi^*_{\tau\tau}  + 2 \left( \mathfrak{p}\beta^{\mathfrak{p}-1}  \beta' \frac{\tilde \ell} {2\pi \beta} \right)    \phi^*_{\tau}  +  \left( \mathfrak{p}\beta^{\mathfrak{p}-2}  \beta' \frac{\tilde \ell} {2\pi} \right)_{\tau}\phi^*
\nonumber
\\[2mm]  &=  \beta^{\mathfrak{p}} \phi^*_{\tau\tau}  + 2 \left( \mathfrak{p}\beta^{\mathfrak{p}-1}  \beta' \frac{\tilde \ell} {2\pi \beta} \right)    \phi^*_{\tau}  +  \left[ \Big(\mathfrak{p} (\mathfrak{p}-2) \beta^{\mathfrak{p}-3}  |\beta'|^2 +  \mathfrak{p}\beta^{\mathfrak{p}-2}  \beta''  \Big)\frac{\tilde \ell^2} {4\pi^2 \beta}   \right]\phi^*.
\end{align}
For $|s|\geq \frac{20 \delta} {\iota}$, we use  a barrier function of the form
\begin{align}
\varphi(x, \tau) = C \|\beta^{\mathfrak{p}} \phi^*\|_\infty  e^{-\frac12\big(x- \frac{20\beta\delta}{\iota}\big)}.
\end{align}
It's easy to verify that $\varphi(x, \tau) $ is a sup-solution of \eqref{asymptotic of phi equation}
and satisfies
\begin{align}
\varphi(x, \tau)  >  \left|\left(\beta^{\mathfrak{p}} \phi^*\right)(x, \tau)\right|,
\quad
\mbox{at}\ s= \frac{20\delta}{\iota}
\end{align}
where the relation between $s$ and $x$ is given by \eqref{coordinate-tildes} and \eqref{definition of x 2.14}.
Then we have
\begin{align}
 |\left(\beta^{\mathfrak{p}} \phi^*\right)(x, \tau)| \leq  C \| \beta^{\mathfrak{p}} \phi^*\|_\infty  e^{-\frac12\big(x- \frac{20\beta\delta}{\iota}\big)}
  \leq C \|\beta^{\mathfrak{p}}  \phi^*\|_{H^2_*(\mathfrak S)} e^{-\frac{5 \beta \delta} \iota},
 \quad \text{for}~|s|\geq \frac{40\delta}{\iota}.
\end{align}
The estimate of $\nabla \phi^*$ follows by local elliptic estimate and we omit the proof for brevity.

\medskip
By defining the following closed and bounded set
\begin{align}
 \mathfrak B  =  &\Bigg\{\phi(x, \tau)  \, :\, \phi~\text{satisfies}~\eqref{periodic condition},
\qquad
\| \beta^{\mathfrak{p}} \phi\|_{H^2_*(\mathfrak S)}\leq D \left( \iota^3+ \iota^2\|f\|_* + \iota^2\|e\|_{**}\right),
\nonumber \\[1mm]
 & \qquad \qquad \qquad
   \Big\|\, |\left( \beta^{\mathfrak{p}}  \phi\right)(x, \tau) |+ | \nabla \left( \beta^{\mathfrak{p}}  \phi\right)(x, \tau)|\, \Big\|_{L^\infty(|s|\geq \frac{40\delta}{\iota})}
\leq
\|\left( \beta^{\mathfrak{p}}  \phi\right)\|_{H^2_*(\mathfrak S)}  e^{-\frac{2\delta} \iota}\Bigg\},
\end{align}
we claim that:
{\em the map $\mathbb A $ defined in \eqref{mathbb-A} is a contraction map from $ \mathfrak B$ to itself.
}

In fact,   recall
\begin{align}
{\mathcal N}(\phi^*) & \,=\,
\frac{1}{\beta^{\mathfrak{p}+2}}
\left[\mathbf{N}( \eta_{3\delta}^{\iota}{\tilde\phi} +{\tilde\psi}) -(p_\epsilon-1)\, \mathbf{U}^{p_\epsilon-2}  \, {\tilde\psi}(\tilde \phi) \right]
\end{align}
with the relation
\[ \tilde\phi(y) =    \beta^{\mathfrak{p}}(\iota z)\phi^*(x, \tau),  \quad \text{for $|\iota s|\leq 10\delta$}, \]
 and
\begin{equation*}
\mathbf{U}(y) \,=\,
\eta_{3\delta}^\epsilon(s)  \beta^{\mathfrak{p}}(\iota z)w_2(x, \tau), \quad \text{for any}\,y \in \R^2.
\end{equation*}
By \eqref{Noliner bf N1}, we have
\begin{align}
\mathbf{N}(\phi_1) -  \mathbf{N}(\phi_2)   &=  \left[- (p_\epsilon-1)\, \left(\mathbf{U} + t_1\phi_1\right) ^{p_\epsilon-2}  + (p_\epsilon-1)\, \mathbf{U}^{p_\epsilon-2}  \right] \phi_1
  \nonumber
 \\[2mm]
 & \quad +  \left[ (p_\epsilon-1)\, \left(\mathbf{U} + t_2\phi_2\right) ^{p_\epsilon-2}  - (p_\epsilon-1)\, \mathbf{U}^{p_\epsilon-2}  \right] \phi_2
 \nonumber\\[2mm]
 & \le   \left|- (p_\epsilon-1)\, \left(\mathbf{U} + t_1\phi_1\right) ^{p_\epsilon-2}  + (p_\epsilon-1)\, \mathbf{U}^{p_\epsilon-2}  \right|  |\phi_1- \phi_2|
  \nonumber\\[2mm]
 & \quad +      \left|- (p_\epsilon-1)\, \left(\mathbf{U} + t_1\phi_1\right) ^{p_\epsilon-2}  + (p_\epsilon-1)\, \left(\mathbf{U} + t_2\phi_2\right) ^{p_\epsilon-2}   \right| \left(|\phi_1| + |\phi_2| \right)
  \nonumber\\[2mm]
 & \le   C \left(\mathbf{U}^{p_\epsilon- 3}\big(|\phi_1|+ |\phi_2|\big) + |\phi_1|^{p_\epsilon-2}\right)   |\phi_1- \phi_2|
  \nonumber\\[2mm]
 & \quad +      \left|- (p_\epsilon-1)\, (p_\epsilon- 2) \left(\mathbf{U} + \hat{t}_1\phi_1\right) ^{p_\epsilon- 3}  t_1\phi_1
 + (p_\epsilon-1)(p_\epsilon- 2) \, \left(\mathbf{U} + \hat{t}_2\phi_2\right) ^{p_\epsilon- 3} t_2 \phi_2   \right| \left(|\phi_1| + |\phi_2| \right)
 \nonumber\\[2mm]
 & \le   C \left(\mathbf{U}^{p_\epsilon- 3} \big(|\phi_1|+ |\phi_2|\big) + |\phi_1|^{p_\epsilon-2} + |\phi_2|^{p_\epsilon-2} \right)   |\phi_1- \phi_2|.
\end{align}
 For any $\phi_1^*, \phi_2^*\in \mathfrak B$,  we have
\begin{align}\label{absolute vale N1-N2}
  &\Big|{\mathcal N}(\phi_1^*) - {\mathcal N}(\phi_2^*) \Big|
\nonumber\\[2mm]
 &  = \frac{1}{\beta^{\mathfrak{p}+2}}   \left[\mathbf{N}( \eta_{3\delta}^{\iota}{\tilde\phi_1} +{\tilde\psi_1}) - \mathbf{N}( \eta_{3\delta}^{\iota}{\tilde\phi_2} +{\tilde\psi_2}) -(p_\epsilon-1)\, \mathbf{U}^{p_\epsilon-2}  \, {\tilde\psi_1} + (p_\epsilon-1)\, \mathbf{U}^{p_\epsilon-2}  \, {\tilde\psi_2}  \right]
   \nonumber\\[2mm]
 & \le   C    \frac{1}{\beta^{\mathfrak{p}+2}} \left\{ \left[ \mathbf{U}^{p_\epsilon- 3}  \left( \left|\eta_{3\delta}^{\iota}{\tilde\phi_1} +{\tilde\psi_1}\right| + \left| \eta_{3\delta}^{\iota}{\tilde\phi_2} +{\tilde\psi_2} \right|  \right)   + \left|\eta_{3\delta}^{\iota}{\tilde\phi_1} +{\tilde\psi_1}\right|^{p_\epsilon-2} + \left| \eta_{3\delta}^{\iota}{\tilde\phi_2} +{\tilde\psi_2} \right|^{p_\epsilon-2}      \right]     +   \mathbf{U}^{p_\epsilon-2}\right\}  \nonumber\\[2mm]
 &\quad \times   \left\{
 \left|   \eta_{3\delta}^{\iota}{\tilde\phi_1} -  \eta_{3\delta}^{\iota}{\tilde\phi_2}   \right|   +   \left| {\tilde\psi_1}  -{\tilde\psi_2}  \right| \right\}.
\end{align}
 Here we have used the notation
\[  \tilde{\psi}_j = {\tilde\psi}({ \tilde \phi}_j),\quad j =1,2.\]
 By Lemma \ref{pro-linear lemma6-1-mathcal L}, we can get
\begin{align}  \label{mathbb A1-mathbb A2}
\|\mathbb A (\phi^*_1) - \mathbb A (\phi^*_2)\|_{H^2_*(\mathfrak S)}&= \|\mathcal T\left(  {\mathcal N}(\phi^*_1)-  {\mathcal N}(\phi^*_2) \right) \|_{H^2_*(\mathfrak S)}
  \leq \|\mathcal T\|\, \|{\mathcal N}(\phi_1^*) - {\mathcal N}(\phi_2^*)\|_{L^2(\mathfrak S)},
\end{align} where $\mathcal T$ is the linear operator given in Lemma \ref{pro-linear lemma6-1-mathcal L}.
Then  from \eqref{absolute vale N1-N2}, we can get
\begin{align} \label{mathcal N1-mathcal N2}
 \|{\mathcal N}(\phi_1^*) - {\mathcal N}(\phi_2^*)\|_{L^2(\mathfrak S)}
 &  \leq \Bigg(\|\phi_1^*\|_{L^4(\mathfrak S)}+  \|\phi_2^*\|_{L^4(\mathfrak S)}
 +\Big\|\frac{{\tilde\psi}_1}{\beta^{\mathfrak{p}} }\Big\|_{L^4(\mathfrak S_\delta)}
 +  \Big\|\frac{{\tilde\psi}_2}{\beta^{\mathfrak{p}}}\Big\|_{L^4(\mathfrak S_\delta)} +1
   \nonumber
 \\[2mm]
 & \qquad
  + \|\phi_1^*\|^{p_\epsilon -2}_{L^{4p_\epsilon -8} (\mathfrak S)}+  \|\phi_2^*\|^{p_\epsilon -2}_{L^{4p_\epsilon -8}(\mathfrak S)}
  +\Big\|\frac{{\tilde\psi}_1}{\beta^{\mathfrak{p}} }\Big\|^{p_\epsilon -2}_{L^{4p_\epsilon -8}(\mathfrak S_\delta)}
  + \Big\|\frac{{\tilde\psi}_2}{\beta^{\mathfrak{p}}}\Big\|^{p_\epsilon -2}_{L^{4p_\epsilon -8}(\mathfrak S_\delta)}\Bigg)
     \nonumber
 \\[2mm]
 & \quad    \times
    \left( \|\phi_1^*- \phi_2^*\|_{L^4(\mathfrak S)}   +  \Big\|\frac{{\tilde\psi}_1- {\tilde\psi}_2}{\beta}\Big\|_{L^4(\mathfrak S_\delta)}  \right)
\end{align}
Here we used the  fact
\begin{align*}
 \tilde \phi_\ell (y)  =  \beta^{\mathfrak{p}}  (\iota z) \phi^*_\ell  (x, \tau),
  \quad \text{for}~y\in  \mathfrak S_\delta = \mathfrak S \cap \Big\{ |s| < \frac {10\delta}\iota \Big\},
\end{align*}
  and the area of $\mathfrak S_\delta$ is of order $O(\frac{\beta\delta}\iota)$.
 We recall the estimate and Lipschitz property of $\tilde\psi$
\begin{equation}  \label{estimate of  psi-star}
\|{\tilde\psi}({\tilde\phi})\|_{L^{\infty}}\,\leq \, C\, \iota \lf[\|\tilde\phi\|_{L^{\infty}(|s|>\frac{20\delta}{\epsilon})}\,+\, \|\nabla  \tilde\phi\|_{L^{\infty}(|s|>\frac{20\delta}{\epsilon})}\ri]\,+\,
e^{-\frac{ \upsilon_2 \delta}{\iota}},
\end{equation}
and
\begin{equation}\label{estimeta of psistar12}
\|{\tilde\psi}(\tilde\phi_1)-{\tilde\psi}(\tilde\phi_2)\|_{L^{\infty}}
\leq
C\iota \lf[ \|\tilde\phi_1-\tilde\phi_2\|_{L^{\infty}(|s|>{20\delta}/\iota)}+\|\nabla\tilde\phi_1-\nabla\tilde\phi_2\|_{L^{\infty}(|s|>
{20\delta}/\iota)}\ri].
\end{equation}
Then  by \eqref{estimate of psi-star} and the relation \eqref{relation betwen tildephi and phi}}, we get
\begin{align}\label{estimate tilde psi1 tilde psi 2 L4}
\Big\|\frac{{\tilde\psi}_1}{\beta}\Big\|_{L^4(\mathfrak S_\delta)}
+  \Big\|\frac{{\tilde\psi}_2}{\beta}\Big\|_{L^4(\mathfrak S_\delta)}
 & \leq
\left\{ \,\lf[ \sum_{j =1}^2\left(\| \tilde\phi_j\|_{L^{\infty}(|s|>\frac{20\delta}{\iota})}\,+\, \|\nabla \tilde\phi_j \|_{L^{\infty}(|s|>\frac{20\delta}{\iota})} \right)\ri]\,+\,
e^{-\frac{ \upsilon_2 \delta}{\iota}} \,\right\}
\left( \int_{\mathfrak S_\delta}  \,{\mathrm d}x  {\mathrm d}\tau\right)^{\frac14}
  \nonumber
 \\[2mm]
 & \leq C\left\{ \,\lf[ \sum_{j =1}^2\left(\| \tilde\phi_j\|_{L^{\infty}(|s|>\frac{20\delta}{\iota})}\,+\, \|\nabla \tilde\phi_j \|_{L^{\infty}(|s|>\frac{20\delta}{\iota})} \right)\ri]\,+\,
e^{-\frac{ \upsilon_2 \delta}{\iota}} \,\right\}\left( \frac{\int_0^{2\pi} \beta   {\mathrm d}\tau} {\iota} \right)^{\frac14}
    \nonumber
 \\[2mm]
 & \leq C\left( \frac{\int_0^{2\pi} \beta   {\mathrm d}\tau} {\iota} \right)^{\frac14}
   \,\sum_{j=1}^2\|\left(\beta\phi_j^*\right)\|_{H^2_*(\mathfrak S)}  e^{-\frac{2\delta} \iota} +  O(e^{-\frac{\delta}{\iota}})
 \leq C e^{-\frac{\delta}{\iota}}.
\end{align}
Here we have used the fact
\begin{align} \label{02pi integrlestimate}
\left( \frac{\int_0^{2\pi} \beta   {\mathrm d}\tau} {\iota} \right)^{\frac14} e^{-\frac{2\delta} \iota} \leq e^{-\frac{\delta} \iota},
\end{align}
  for $\iota$ small enough.
  In fact,  by \eqref{length-constraint}, we have
\begin{align}   \label{02pi integrl1}
  \int_0^{2\pi} \beta   {\mathrm d}\tau =    \frac{2\pi} {\tilde\ell}\int_0^{\ell} \beta^2   {\mathrm d}\theta \leq C \iota^{\varrho-1}.
\end{align}
And  \eqref{02pi integrlestimate} can be directly obtained  by \eqref{02pi integrl1}.

Next, by \eqref{estimeta of psistar12}, we have
\begin{align}\label{tilde psi 1-tilde psi 2}
\|\tilde\psi_1- \tilde\psi_2\|_{L^4(\mathfrak S_\delta)}
 &  \leq \|\tilde\psi_1- \tilde\psi_2\|_{L^{\infty} (\mathfrak S_\delta)} \left( \int_{\mathfrak S_\delta}  \,{\mathrm d}x  {\mathrm d}\tau\right)^{\frac14}
   \nonumber
 \\[2mm]
 &  \leq C \left( \frac{\int_0^{2\pi} \beta   {\mathrm d}\tau} {\iota} \right)^{\frac14}
 \big[\, \|\beta^{\mathfrak{p}} \left({\phi}^*_1-{ \phi}^*_2\right)\|_{L^{\infty}(|s|>{20\delta}/\iota)}+\|\nabla \left(\beta^{\mathfrak{p}}{\phi}^*_1-\beta{\phi}^*_2 \right) \|_{L^{\infty}(|s|>{20\delta}/\iota)}\big]
    \nonumber
 \\[2mm]
  & \leq o(1) \|\beta^{\mathfrak{p}} \left({\phi}^*_1-{ \phi}^*_2\right)\|_{H^2_*(\mathfrak S)}.
\end{align}
 The last inequality holds due to
\begin{align*}
 \big[\, \|\beta^{\mathfrak{p}} \left({\phi}^*_1-{ \phi}^*_2\right)\|_{L^{\infty}(|s|>{20\delta}/\iota)}+\|\nabla \left(\beta^{\mathfrak{p}}{\phi}^*_1-\beta^{\mathfrak{p}}{\phi}^*_2 \right) \|_{L^{\infty}(|s|>{20\delta}/\iota)}\big] \leq  \|\beta^{\mathfrak{p}} \left({\phi}^*_1-{ \phi}^*_2\right)\|_{H^2_*(\mathfrak S)}  e^{-\frac{2\delta} \iota},
\end{align*}
 for any ${\phi}^*_1, {\phi}^*_2 \in \mathfrak B.$
 By \eqref{mathbb A1-mathbb A2}, \eqref{mathcal N1-mathcal N2},  \eqref{estimate tilde psi1 tilde psi 2 L4}, \eqref{tilde psi 1-tilde psi 2} and Sobolev embedding, we  prove that $\mathbb{A}$ is a contraction map.

\medskip
We will show that $\mathbb{A}$ maps $\mathfrak B$ into itself.
For any $\phi^*\in \mathfrak B$,  by Lemma \ref{pro-linear lemma6-1-mathcal L}, then
\begin{align} \label{estimate of mathbbAphi1}
\|\beta^{\mathfrak{p}} \mathbb A (\phi^*)\|_{H^2_*(\mathfrak S)}
 =\,& \|\beta^{\mathfrak{p}}{\mathcal T}\left(  - E_{11}+ {\mathcal N}(\phi^*)\right) \|_{H^2_*(\mathfrak S)}
    \nonumber
 \\[2mm]
& \leq\, \|{\mathcal T}\| \Bigg\{\|\beta^{\mathfrak{p}} E_{11}\|_{L^2(\mathfrak S)}
\ +\
\Bigg\|\frac{\beta^{\mathfrak{p}}}{\beta^{2+\mathfrak{p}} }    \Big[\mathbf{N}( \eta_{3\delta}^{\iota}{\tilde\phi_1} +{\tilde\psi_1}) - \mathbf{N}( \eta_{3\delta}^{\iota}{\tilde\phi_2} +{\tilde\psi_2})
    \nonumber
 \\[2mm]   & \qquad\qquad
-(p_\epsilon-1)\, \mathbf{U}^{p_\epsilon-2}  \, {\tilde\psi_1} + (p_\epsilon-1)\, \mathbf{U}^{p_\epsilon-2}  \, {\tilde\psi_2}  \Big]  \Bigg\|_{L^2(\mathfrak S)} \Bigg\}
    \nonumber
 \\[2mm]
& \leq\,  \|{\mathcal T}\|\Bigg\{\|\beta^{\mathfrak{p}}  E_{11}\|_{L^2(\mathfrak S)}
 +C\|\beta^{\mathfrak{p}} \phi^*\|^2_{L^4(\mathfrak S)}
  +C \|\tilde\psi\|^2_{L^4(\mathfrak S_\delta)}
     \nonumber
 \\[2mm] & \qquad \qquad
   +C  \|\beta^{\mathfrak{p}} \phi^*\|^3_{L^6(\mathfrak S)}
    + C  \|\tilde\psi\|^3_{L^6(\mathfrak S_\delta)}
    +C \|\tilde\psi\|_{L^2(\mathfrak S_\delta)} \Bigg\}.
\end{align}
Using \eqref{estimate of  psi-star} and similar to  \eqref{estimate tilde psi1 tilde psi 2 L4},  we have
\begin{align}  \label{estimate of tildepsi}
   \|\tilde\psi\|_{L^p(\mathfrak S_\delta)}
&  \leq C\left( \frac{\int_0^{2\pi} \beta   {\mathrm d}\tau} {\iota} \right)^{\frac1p}
\, \|\beta^{\mathfrak{p}} \phi^*\|_{H^2_*(\mathfrak S)}  e^{-\frac{2\delta} \iota} +  O(e^{-\frac{\delta}{\iota}})
     \nonumber
 \\[2mm]
 &
 \leq  Ce^{-\frac{\delta}{\iota}} \Big(\|\beta^{\mathfrak{p}}\phi^*\|_{H^2_*(\mathfrak S)}+1\Big),
\end{align}
 provided with
\begin{align}
\left( \frac{\int_0^{2\pi} \beta   {\mathrm d}\tau} {\iota} \right)^{\frac1p} e^{-\frac{2\delta} \iota} \leq e^{-\frac{\delta} \iota}.
\end{align}  By \eqref{estimate of mathbbAphi1}, \eqref{estimate of tildepsi} and Sobolev embedding, we get
\begin{align}\label{estimate of mathbbAphi2}
\|\mathbb A (\phi^*)\|_{H^2_*(\mathfrak S)}
& \leq C\Big\{\|\beta^{\mathfrak{p}}  E_{11}\|_{L^2(\mathfrak S)}
+ C\|\beta^{\mathfrak{p}} \phi^*\|^2_{H^2_*(\mathfrak S)}
+ C e^{-\frac{\delta}{\iota}} \Big(\|\beta^{\mathfrak{p}}\phi^*\|_{H^2_*(\mathfrak S)}+1\Big)  \Big\}
 \nonumber
 \\[2mm]
 &  \leq C \left( \iota^3+ \iota^2\|f\|_* + \iota^2\|e\|_{**}\right).
\end{align}
Next,  for any $\phi^*$, since  $\varphi^* :=\mathbb A (\phi^*) $ satisfies  the equation
\begin{align}
\frac{\iota^2 4\pi^2}{{\tilde\ell}^2} \varphi^*_{\tau\tau} + \varphi^*_{xx} - \big(1+ o(1)\big)\varphi^* =  h
\end{align}
with $h$ compactly supported. Hence $ \varphi^* $ satisfies   the estimates required for functions in the  $\mathcal B$.
Then we know that $\mathbb{A}$ map $\mathfrak B$ to itself if $D$ large enough.

\medskip
By summarising the above results, we present the following   resolution theory.
\begin{proposition}\label{proposition6.5}
There is a number $D>0$ such that for all $\epsilon$ small enough and all parameters $(f, e)$ in $\hat{\mathcal F}$,
problem \eqref{project 6.13}-\eqref{projectedproblem3} has a unique solution $\phi^* \,=\, \phi^*(f, e)$, which satisfies the estimates
\begin{align*}
 \|\beta \phi^*\|_{H_*^2(\mathfrak S)} + \| \phi^*\|_{H_*^2(\mathfrak S)} & \,\leq\, D \left( \iota^3+ \iota^2\|f\|_* + \iota^2\|e\|_{**}\right),
  \end{align*}
  and
\begin{align*}
\big\| |\beta^{\mathfrak{p}} \phi^*|+|\nabla  \left(\beta^{\mathfrak{p}} \phi^*\right)| \big\|_{L^{\infty}(|x|>40\delta/\iota)}\,\leq\, \|\beta^{\mathfrak{p}} \phi^*\|_{H_*^2(\mathfrak S)}e^{-\delta/\iota}.
\end{align*}
 Moreover, $\phi^*$ depends Lipschitz-continuously on the parameters $f$ and $e$ in the sense of the estimate
\begin{equation}
\label{characteriztion}
\|\phi^*(f_1, e_1)\,-\, \phi^*(f_2, e_2)\|_{H_*^2(\mathfrak S)}
\,\leq\, C\iota^2\big [\, \|f_1-f_2\|_{*}\,+\, \|e_1-e_2\|_{**}\, \big].
\end{equation}
\end{proposition}

\medskip\section{Deriving the reduced system of differential equations}\label{deriv-greduce-equa}

After solving problem \eqref{project 6.13}-\eqref{projectedproblem3} in Proposition \ref{proposition6.5},
we shall choose the unknown parameters $f$ and $e$ to make the Lagrange multipliers $c$ and $d$ vanish.
In fact, $c\equiv 0$ and $d\equiv 0$ are equivalent to
\begin{align} \label{project 6.13-1}
\int_\R\Big\{\,\mathcal L (\phi^*) +\eta_\delta^\iota(s)\frac{1}{\beta^{\mathfrak{p}+2}}{\mathcal E} - \eta_\delta^\iota(s){\mathcal N}(\phi^*)\Big\}
U'(x) \,{\mathrm d}x\,=\, 0,
\end{align}
and
\begin{align} \label{project 6.13-2}
\int_\R\Big\{\,\mathcal L (\phi^*) +\eta_\delta^\iota(s)\frac{1}{\beta^{\mathfrak{p}+2}}{\mathcal E} - \eta_\delta^\iota(s){\mathcal N}(\phi^*)\Big\}
{\mathcal Z}(x) \,{\mathrm d}x\,=\, 0.
\end{align}
The calculations of the above terms will be given in Sections \ref{section9.1}-\ref{section9.2}.
These will consequently derive the reduced system with differential equations involving the
parameters $f$ and $e$, see Section \ref{section9.3}.

\subsection{  Estimates for projections of the error}\label{section9.1}
We   denote by $\tilde{\mathcal D}_i(\theta),  i \in \mathbb{N}^+, $ the smooth functions that are uniformly bounded and independent of $f, f', f'', e, e', e''$.

\subsubsection{Computation for $ \beta^{-\mathfrak{p}-2}\int_{\R}\, \eta_\delta^\iota(s)\, {\mathcal E}\, U'  \,{\mathrm d}x$}\

We first concern
\begin{align}
 \frac{1}{\beta^{\mathfrak{p}+2}}\int_\R\eta_\delta^\iota(s){\mathcal E}U'(x) \,{\mathrm d}x
=&\frac{1}{\beta^{\mathfrak{p}+2}}\int_\R {\mathcal E}U'(x) \,{\mathrm d}x
\ +\
\frac{1}{\beta^{\mathfrak{p}+2}}\int_\R \big[\eta_\delta^\iota(s)-1\big]{\mathcal E}U'(x) \,{\mathrm d}x
\nonumber\\[2mm]
=& \int_{\R}\, S\big(w_2\big)\, U'(x) \,{\mathrm d}x
\ +\ \int_\R \big[\eta_\delta^\iota(s)-1\big]S\big(w_2\big)U'(x) \,{\mathrm d}x
\nonumber\\[2mm]
=&
\int_\R \Bigg[\iota^2 S_{3}  + \iota^2 S_{42}  + \iota^3 S_{5} +\iota^3 S_{6}
 +\iota^3 S_{7} +\iota^3 S_{8}
  +\frac 1{\beta^{\mathfrak{p}}} \iota^3\frac{4\pi^2}{{\tilde\ell}^2}  e''\mathcal Z
  \nonumber\\[2mm]
&\qquad    +\frac 1{\beta^{\mathfrak{p}}}  \iota\lambda_0 e{\mathcal Z}   + N_{0\e} \left(\iota\frac{e} {\beta^{\mathfrak{p}}} \mathcal Z\right)
   + \iota^2 {\bf N}_{02} +\iota^3{\bf N}_{03} + \iota^4 {\bf N}_{04}
     \nonumber\\[2mm]
&\qquad+  \iota^2 {\bf B}_{42} +  \iota^3 \frac{4\pi^2}{{\tilde\ell}^2} {\phi}_{1, \tau\tau}
+ \iota^4 \frac{4\pi^2}{{\tilde\ell}^2}\phi_{2, \tau\tau}
+B_{5}(\iota\phi_1)
  \nonumber\\[2mm]
&\qquad
  + \mathbb{N}_1 (\iota^2 \phi_2)
+ B_4(\iota^2\phi_2)+  B_5 (\iota^2{\phi_2})   \Bigg]U'(x) \,{\mathrm d}x  + \mathfrak G(\theta)
\nonumber\\[2mm]
& +\int_\R \eta_\delta^\iota(s) \Bigg\{ \iota^4 \frac{4\pi^2} {{\tilde\ell}^2} a_2\frac{\kappa}{\, \beta^{1+ \mathfrak{p}} \, }  e'' (x+f)\mathcal Z -\iota^4 \kappa a_2 \frac{ 4\pi^2} {{\tilde\ell}^2} \frac{1}{\, \beta^{1+\mathfrak{p}} \, }  (x+f) f'' e   \mathcal Z'
  \nonumber\\[2mm]
&\qquad\qquad \quad+ \frac{\iota^4}{\beta^{1+\mathfrak{p}}} \hat{\mathcal F}_1\left[\tau, x, \beta^{\mathfrak p} f, \frac{\beta^{\mathfrak p} f''}{\tilde\ell^2} \right]
  +  \frac{\iota^4}{\beta^{1+\mathfrak{p}}} \mathcal F_1\left[\tau, x, \beta^{\mathfrak p} f, \frac{\beta^{\mathfrak p} f'}{\tilde\ell}\right]
 \nonumber\\[2mm]
&\qquad\qquad \quad
+B_{6}(\iota\phi_1)
+\frac{\iota^4}{\beta^{1+\mathfrak{p}}}  \mathcal G_1\left[\tau, x, \beta^{\mathfrak p} f, \frac{\beta^{\mathfrak p} f'}{\tilde\ell},   e, \frac{e'}{\tilde\ell}\right]
 +  B_6 (\iota^2{\phi_2}) \Bigg\} U'(x) \,{\mathrm d}x,
\end{align}
due to the expression of $\mathcal{E}$ given in \eqref{errorrelationinterior},
where
\begin{align*}
 \mathfrak G(\theta)  & =   \int_\R\big[\eta_\delta^\iota(s)-1\big]  S(w_2) U'(x) \,{\mathrm d}x.
\end{align*}
Note that the terms  from
\begin{align}
  & \iota^4 \frac{4\pi^2} {{\tilde\ell}^2} a_2\frac{\kappa}{\, \beta^{1+\mathfrak{p}} \, }  e'' (x+f)\mathcal Z -\iota^4 \kappa a_2 \frac{ 4\pi^2} {{\tilde\ell}^2} \frac{1}{\, \beta^{1+\mathfrak{p}}\, }  (x+f) f'' e   \mathcal Z'    + \frac{\iota^4}{\beta^{1+\mathfrak{p}}} \hat{\mathcal F}_1\left[\tau, x, \beta^{\mathfrak p} f, \frac{\beta^{\mathfrak p} f''}{\tilde\ell^2} \right]
 \nonumber\\[2mm]
&
+B_{6}(\iota\phi_1)
  +  \frac{\iota^4}{\beta^{1+\mathfrak{p}}} \mathcal F_1\left[\tau, x, \beta^{\mathfrak p} f, \frac{\beta^{\mathfrak p} f'}{\tilde\ell}\right]
 +\frac{\iota^4}{\beta^{1+\mathfrak{p}}}  \mathcal G_1\left[\tau, x, \beta^{\mathfrak p} f, \frac{\beta^{\mathfrak p} f'}{\tilde\ell},   e, \frac{e'}{\tilde\ell}\right]
 +  B_6 (\iota^2{\phi_2})
 \end{align}
 are involving $a_\ell(\iota s, \iota z), \ell =1,2,3$.
By the exponential decay of the terms in $S\big(w_2\big)$ with respect to $x$,
 we  neglect  the calculations of the terms in $ \mathfrak G(\theta). $
Since the following terms in $S\big(w_2\big)$
$$
\iota^2 S_{42},\quad \iota^3 S_6,\quad \iota^3 S_8,\quad
\frac 1{\beta^{\mathfrak{p}}} \iota^3\frac{4\pi^2}{{\tilde\ell}^2}  e''\mathcal Z +\frac 1{\beta^{\mathfrak{p}}}  \iota\lambda_0 e{\mathcal Z},
\quad \iota^4 \frac{4\pi^2}{{\tilde\ell}^2}\phi_{2, \tau\tau}, \quad N_{0\e} \left(\iota\frac{e} {\beta^{\mathfrak{p}}} \mathcal Z\right)
$$
are even functions of $x$, then integration against $U'(x)$ just vanishes.
 And there  holds
\begin{align}
  \iota^4 \frac{4\pi^2}{{\tilde\ell}^2}\int_\R\phi_{1, \tau\tau} U'(x) \,{\mathrm d}x=0,
\end{align}
 due to the fact
$$
\int_\R\phi_{1}(x,\tau) U'(x)\,{\mathrm d}x=0.
$$
 Then we have
\begin{align}
& \frac{1}{\beta^{\mathfrak{p}+2}}\int_\R\eta_\delta^\iota(s){\mathcal E}U'(x) \,{\mathrm d}x
 \nonumber\\[2mm]
& =\int_{\R}\,   \iota^2 S_{3} U'(x)\,{\mathrm d}x
+ \int_{\R}\,  \iota^3 S_{5}U'(x)\,{\mathrm d}x
+\int_{\R}\, \iota^3 S_{7}U'(x)\,{\mathrm d}x
+ \int_{\R}\, \iota^2 {\bf N}_{02} U'(x)\,{\mathrm d}x
\nonumber\\[2mm]
& \quad
+ \int_{\R}\, \iota^3{\bf N}_{03}U'(x)\,{\mathrm d}x
+  \int_{\R}\, \iota^4 {\bf N}_{04} U'(x)\,{\mathrm d}x
 +  \int_{\R}\,\iota^2 {\bf B}_{42} U'(x)\,{\mathrm d}x
     +\int_{\R}\, B_{5}(\iota\phi_1) U'(x)\,{\mathrm d}x
     \nonumber\\[2mm]
& \quad +   \int_{\R}\,  \mathbb{N}_1 (\iota^2 \phi_2)U'(x)\,{\mathrm d}x
+ \int_{\R}\, B_4(\iota^2\phi_2) U'(x)\,{\mathrm d}x
+  \int_{\R}\, B_5 (\iota^2{\phi_2})  U'(x)\,{\mathrm d}x
 \nonumber\\[2mm]
& \quad
  + \int_\R \eta_\delta^\iota(s)  \Bigg[ \iota^4 \frac{4\pi^2} {{\tilde\ell}^2} a_2\frac{\kappa}{\, \beta^{1+ \mathfrak{p}} \, }  e'' (x+f)\mathcal Z
  -\iota^4 \kappa a_2 \frac{ 4\pi^2} {{\tilde\ell}^2} \frac{1}{\, \beta^{1+ \mathfrak{p}} \, } (x+f) f'' e   \mathcal Z'     \Bigg]  U'(x)\,{\mathrm d}x
   \nonumber\\[2mm]
& \quad
  +  \int_\R \eta_\delta^\iota(s)   \frac{\iota^4}{\beta^{1+\mathfrak{p}}} \hat{\mathcal F}_1\left[\tau, x, \beta^{\mathfrak p} f, \frac{\beta^{\mathfrak p} f''}{\tilde\ell^2} \right]   U'(x)\,{\mathrm d}x
   +  \int_\R \eta_\delta^\iota(s)   \frac{\iota^4}{\beta^{1+\mathfrak{p}}} \mathcal F_1\left[\tau, x, \beta^{\mathfrak p} f, \frac{\beta^{\mathfrak p} f'}{\tilde\ell}\right] U'(x)\,{\mathrm d}x
\nonumber\\[2mm]
& \quad
 + \int_\R \eta_\delta^\iota(s)  \frac{\iota^4}{\beta^{1+\mathfrak{p}}}  \mathcal G_1\left[\tau, x, \beta^{\mathfrak p} f, \frac{\beta^{\mathfrak p} f'}{\tilde\ell},   e, \frac{e'}{\tilde\ell}\right]  U'(x)\,{\mathrm d}x
 + \int_\R \eta_\delta^\iota(s)  B_{6}(\iota\phi_1)  U'(x)\,{\mathrm d}x
 \nonumber\\[2mm]
& \quad +  \int_\R \eta_\delta^\iota(s)   B_6 (\iota^2{\phi_2}) U'(x)\,{\mathrm d}x
 \nonumber\\[2mm]
& = \sum_{i=1}^{17} I_j.
\end{align}

\smallskip\noindent$\clubsuit$
 By the expressions of $ \iota^2 S_3$  given in \eqref{expression of s3}, we have
\begin{align}  \label{bar S3intergralUx}
 I_1   &=  \int_{\R}\,   \iota^2 S_{3} U'(x)\,{\mathrm d}x
 \nonumber\\[2mm]
& = -\, \iota^2 \frac{\iota^2W_{tt}(0, \iota z)}{\beta^4}f \int_{\R} xU U'(x)\,{\mathrm d}x
+\iota^2\Big(\frac{|\beta'|^2}{\beta^4} 2f
-\frac{4\pi}{ {\tilde\ell}}\frac{\beta'}{\beta^2}f'
\Big)\int_{\R} xU'' U'(x)\,{\mathrm d}x
\nonumber\\[2mm]
&\quad
+\iota^2 \Big[
-\frac{4\pi^2}{{\tilde\ell}^2}f''
- \big(4\mathfrak{p}+2 \big) \frac{\pi}{\tilde\ell}  \frac{\beta'}{\beta^2} f'
+2 \mathfrak{p}\frac{|\beta'|^2}{\beta^4}f
+\frac{\beta''} {\beta^3}f
\,-\,\frac{\kappa^2}{\beta^2}f
\Big]\int_{\R} |U'|^2 \,{\mathrm d}x
\nonumber\\[2mm]
& \quad  +\iota^2\frac{\kappa}{\, \beta^{1+\mathfrak{p}} \, }\int_{\R} \left[ e {\mathcal Z}'+\frac 1{\Lambda} e x\mathcal Z \right]U'(x) \,{\mathrm d}x
\nonumber\\[2mm]
& =  \iota^2 \varrho_1\Bigg\{ -  \frac{4\pi^2}{{\tilde\ell}^2}f''  - 4\mathfrak{p}\frac{\pi}{\tilde\ell}  \frac{\beta'}{\beta^2} f'
+ \Big[(2\mathfrak{p}-1) \frac{|\beta'|^2}{\beta^4}
+\frac{\beta''} {\beta^3}
\,-\,\frac{\kappa^2}{\beta^2} \Big]   f  +  \,  \Lambda\frac{\iota^2 W_{tt}(0, \iota z)}{\beta^4} f\Bigg\}
+ \frac {\iota^2 \kappa d_1e}{\beta^{\mathfrak{p}+1}},
\end{align}
where
\[\varrho_1= \int_\R |U'|^2{\mathrm d}x, \quad   d_1 = \int_\R\left[ {\mathcal Z}'+\frac{1}{\Lambda}  x\mathcal Z \right]U'(x) \,{\mathrm d}x. \]

On the other hand, \eqref{B42} and \eqref{N02} will imply the following estimates
\begin{align}
&\int_\R  \iota^2{\bf N}_{02}U'(x) \,{\mathrm d}x
 \nonumber\\[2mm]
& =  \iota^2  (p_\epsilon-1)(p_\epsilon- 2)  \int_\R \Bigg[  \frac 1{2\Lambda^2} {\mathbf b}_1^2f^2\varpi_2^2+  {\mathbf b}_1^2 \frac 1{\Lambda} f \varpi_1 \varpi_2
+  \frac{e} {\beta^{\mathfrak{p}}} \mathcal Z\Big( {\mathbf b}_1\varpi_1+ \frac 1{\Lambda} {\mathbf b}_1f\varpi_2\Big)    \Bigg] U^{p_\epsilon- 3}U'(x) \,{\mathrm d}x
\nonumber\\[2mm]
& =  \iota^2 (p_\epsilon-1)(p_\epsilon- 2)\int_\R\left(\frac{ {\mathbf b}_1e}{\beta^{\mathfrak{p}}} {\mathcal Z}\varpi_1
+ \frac 1{\Lambda}{\mathbf b}_1^2f \varpi_1\varpi_2 \right) U^{p_\epsilon- 3}U'(x) \,{\mathrm d}x
\nonumber\\[2mm]
& = \iota^2\frac{ (p_\epsilon-1)(p_\epsilon- 2)  {\mathbf b}_1e}{\beta^{\mathfrak{p}}} \int_\R  {\mathcal Z}\varpi_1   U^{p_\epsilon- 3}U'(x) \,{\mathrm d}x
 + \iota^2\frac{(p_\epsilon-1)(p_\epsilon- 2)} {\Lambda}{\mathbf b}_1^2f \int_\R   \varpi_1\varpi_2    U^{p_\epsilon- 3}U'(x) \,{\mathrm d}x
\nonumber\\[2mm]
& = \iota^2 \frac{   {\mathbf b}_1e}{\beta^{\mathfrak{p}}}  d_2
+  \iota^2 \frac{(p_\epsilon-1)(p_\epsilon- 2)} {\Lambda}{\mathbf b}_1^2f \int_\R   \varpi_1\varpi_2    U^{p_\epsilon- 3}U'(x) \,{\mathrm d}x,
\end{align}
where
\begin{align*}
d_2 =  (p_\epsilon-1)(p_\epsilon- 2) \int_\R  {\mathcal Z}\varpi_1   U^{p_\epsilon- 3}U'(x) \,{\mathrm d}x.
\end{align*}
 From \eqref{B42}, we have
\begin{align}
 \int_\R   \iota^2 {\bf B}_{42}  U'(x) \,{\mathrm d}x
 &  = \iota^2 {\mathbf b}_1^2\int_\R \Big[
 \frac{1}{\Lambda}f\varpi_{2,x}
+  \frac{1}{\Lambda^2} f x\varpi_2
+ \frac{1}{\Lambda}f\varpi_1
\Big]U'(x) \,{\mathrm d}x
\nonumber\\[2mm]
& =   \iota^2 {\mathbf b}_1^2   \frac{1}{\Lambda}f  \int_\R \Big[ \varpi_{2,x} +  \varpi_1\Big] U'(x) \,{\mathrm d}x   +  \iota^2 {\mathbf b}_1^2   \frac{1}{\Lambda^2}f   \int_\R x\varpi_{2} U'(x) \,{\mathrm d}x.
\end{align}
By the identities
\begin{equation}
 \int_{\R}  x \varpi_{2} U'(x) \,{\mathrm d}x  = - \int_{\R}   \varpi_{2} U \,{\mathrm d}x  - \int_{\R}  x \varpi_{2, x} U \,{\mathrm d}x
\end{equation}
\begin{equation}
\int_{\R}  \varpi_{2, x} U'(x) \,{\mathrm d}x = - \left[  \frac1{p_\epsilon-2} + \frac 14\right]  \int_\R |U'|^2 \,{\mathrm d}x,
\qquad
\int_{\R}  \varpi_{2} U \,{\mathrm d}x = -  \left[ \frac{2\Lambda}{p_\epsilon-2}   -   \frac 12\Lambda \right] \int_\R |U'|^2 \,{\mathrm d}x,
\end{equation}
\begin{align}
 (p_\epsilon-1)(p_\epsilon- 2) \int_\R   \varpi_1\varpi_2    U^{p_\epsilon- 3}U'(x) \,{\mathrm d}x     =   \int_\R\left(U'(x) + \frac 1{\Lambda} x U \right) \varpi_{2, x}   \,{\mathrm d}x
 -  \int_\R  \varpi_1 U'(x)\,{\mathrm d}x,
\end{align}
we can get
\begin{align} \label{B42+ N01 project on U'}
  &\int_\R\left[\iota^2 {\bf B}_{42}  + \iota^2{\bf N}_{02} \right]U'(x) \,{\mathrm d}x
\nonumber\\[2mm]
& = \iota^2   \frac{ {\mathbf b}_1e}{\beta^{\mathfrak{p}}}  d_2
+ \iota^2\frac{1} {\Lambda} {\mathbf b}_1^2f  \Bigg[(p_\epsilon-1)(p_\epsilon- 2)\int_\R   \varpi_1\varpi_2    U^{p_\epsilon- 3}U'(x) \,{\mathrm d}x  +  \int_\R \Big[ \varpi_{2,x} +  \varpi_1\Big] U'(x) \,{\mathrm d}x \Bigg]
\nonumber\\[2mm]
& \quad + \iota^2 {\mathbf b}_1^2   \frac{1}{\Lambda^2}f \Bigg[- \int_{\R}   \varpi_{2} U \,{\mathrm d}x  - \int_{\R}  x \varpi_{2, x} U \,{\mathrm d}x \Bigg]
 \nonumber\\[2mm]
& = \iota^2  \frac{  {\mathbf b}_1e}{\beta^{\mathfrak{p}}}  d_2
+\iota^2 \frac{1} {\Lambda} {\mathbf b}_1^2f  \Bigg[ \int_\R\left(U'(x) + \frac 1{\Lambda} x U \right) \varpi_{2, x}   \,{\mathrm d}x
 -  \int_\R  \varpi_1 U'(x)\,{\mathrm d}x  +  \int_\R \Big[ \varpi_{2,x} +  \varpi_1\Big] U'(x) \,{\mathrm d}x \Bigg]
\nonumber\\[2mm]
& \quad + \iota^2 {\mathbf b}_1^2   \frac{1}{\Lambda^2}f \Bigg[- \int_{\R}   \varpi_{2} U \,{\mathrm d}x  - \int_{\R}  x \varpi_{2, x} U \,{\mathrm d}x \Bigg]
  \nonumber\\[2mm]
& = \iota^2  \frac{  {\mathbf b}_1e}{\beta^{\mathfrak{p}}}  d_2
+\iota^2\frac{2} {\Lambda} {\mathbf b}_1^2f   \int_\R U'(x)   \varpi_{2, x}   \,{\mathrm d}x
+ \iota^2 {\mathbf b}_1^2   \frac{1}{\Lambda}f   \left[ \frac{2}{p_\epsilon-2}   -   \frac 12 \right] \int_\R |U'|^2 \,{\mathrm d}x
  \nonumber\\[2mm]
& = \iota^2 \frac{  {\mathbf b}_1e}{\beta^{\mathfrak{p}}}  d_2
  + \iota^2 {\mathbf b}_1^2   \frac{1}{\Lambda}f    \left[  -\frac2{p_\epsilon-2} - \frac 12 +  \frac{2}{p_\epsilon-2}   -   \frac 12 \right] \int_\R |U'|^2 \,{\mathrm d}x.
 \end{align}

Therefore,  combining \eqref{bar S3intergralUx} and \eqref{B42+ N01 project on U'}, we have
\begin{align}  \label{S32+N03+B43intUx}
&   \int_\R\left[\iota^2 S_3 +\iota^2 {\bf B}_{42}  + \iota^2{\bf N}_{02} \right]U'(x) \,{\mathrm d}x
\nonumber
 \\[2mm]
& = \iota^2 \varrho_1\Bigg\{ -  \frac{4\pi^2}{{\tilde\ell}^2}f''  - 4\mathfrak{p}\frac{\pi}{\tilde\ell}  \frac{\beta'}{\beta^2} f'
+ \Big[(2\mathfrak{p}-1) \frac{|\beta'|^2}{\beta^4}
+\frac{\beta''} {\beta^3}
\,-\, \big[\frac 1{\Lambda}+1\big]\frac{\kappa^2}{\beta^2}  +\Lambda\frac{\iota^2 W_{tt}(0, \iota z)}{\beta^4} \Big]   f  \Bigg\}
\nonumber
 \\[2mm]
&\quad
+\iota^2\frac { {\mathbf b}_1(d_1+d_2)e}{\beta^{\mathfrak{p}}}.
\end{align}

\smallskip\noindent$\clubsuit$
Moreover, we have
\begin{align}
\int_\R  \iota^3{\bf N}_{03}U'(x) \,{\mathrm d}x
\,=\, & \iota^3 \frac{(p_\epsilon-1) (p_\epsilon- 2) (p_\epsilon- 3)} {2}  \int_\R U^{p_\epsilon- 4} \left[{\mathbf b}_1^2\varpi_1^2+ \frac 1{\Lambda_\e^2} {\mathbf b}_1^2f^2\varpi_2^2+ 2 {\mathbf b}_1^2 \frac 1{\Lambda} f \varpi_1 \varpi_2 \right]     \frac{e} {\beta^{\mathfrak{p}}} \mathcal Z U'(x) \,{\mathrm d}x
 \nonumber\\[2mm]
&   +   \iota^3  \frac{(p_\epsilon-1) (p_\epsilon- 2)(p_\epsilon- 3)} {6} \int_\R  U^{p_\epsilon- 4}  \left[{\mathbf b}_1\varpi_1+ \frac 1{\Lambda} {\mathbf b}_1f\varpi_2 \right]^3 U'(x) \,{\mathrm d}x
 \nonumber\\[2mm]
\quad &  +  \iota^3   \frac{(p_\epsilon-1)(p_\epsilon- 2) (p_\epsilon- 3)}2 \int_\R U^{p_\epsilon- 4}  \frac{e^2} {\beta^{2\mathfrak{p}}} \mathcal Z^2\big({\mathbf b}_1\varpi_1+ \frac 1{\Lambda} {\mathbf b}_1f\varpi_2\big)  U'(x) \,{\mathrm d}x
\nonumber\\[2mm]
\,=\,& \iota^3 {(p_\epsilon-1) (p_\epsilon- 2) (p_\epsilon- 3)}   {\mathbf b}_1^2 \frac 1{\Lambda} f  \frac{e} {\beta^{\mathfrak{p}}}  \int_\R U^{p_\epsilon- 4}   \varpi_1 \varpi_2     \mathcal Z U'(x) \,{\mathrm d}x \nonumber\\[2mm]
&   +    \iota^3\frac{(p_\epsilon-1) (p_\epsilon- 2)(p_\epsilon- 3)} {6}  {\mathbf b}_1^3\int_\R  U^{p_\epsilon- 4}  \varpi_1^3  U'(x) \,{\mathrm d}x
 \nonumber\\[2mm]
&   +   \iota^3 \frac{(p_\epsilon-1) (p_\epsilon- 2)(p_\epsilon- 3)} {2}  \frac 1{\Lambda^2} {\mathbf b}_1^3f^2 \int_\R  U^{p_\epsilon- 4}  \varpi_2^2\varpi_1U'(x) \,{\mathrm d}x
 \nonumber\\[2mm]
&  +  \iota^3  \frac{(p_\epsilon-1)(p_\epsilon- 2) (p_\epsilon- 3)}2 \frac{e^2} {\beta^{2\mathfrak{p}}}    {\mathbf b}_1\ \int_\R U^{p_\epsilon- 4}   \mathcal Z^2 \varpi_1  U'(x) \,{\mathrm d}x
    \nonumber\\[2mm]
\,=\,&  \iota^3    {\mathbf b}_1^2 \frac 1{\Lambda} f  \frac{e} {\beta^{\mathfrak{p}}} d_3+   \iota^3\frac{(p_\epsilon-1) (p_\epsilon- 2)(p_\epsilon- 3)} {6}  {\mathbf b}_1^3 d_4
 +   \iota^3    \frac 1{\Lambda^2} {\mathbf b}_1^3f^2  (d_5 + d_6),
\end{align}
where
\begin{align}
& d_3={(p_\epsilon-1) (p_\epsilon- 2) (p_\epsilon- 3)}\int_\R U^{p_\epsilon- 4}   \varpi_1 \varpi_2     \mathcal Z U'(x) \,{\mathrm d}x, \\
&  d_4 = \frac{(p_\epsilon-1) (p_\epsilon- 2) (p_\epsilon- 3)} 6 \int_\R  U^{p_\epsilon- 4}  \varpi_1^3  U'(x) \,{\mathrm d}x,
\end{align}
\begin{align}
 & d_5=\frac{(p_\epsilon-1) (p_\epsilon- 2) (p_\epsilon- 3)} 2\int_\R  U^{p_\epsilon- 4}  \varpi_1^3  U'(x) \,{\mathrm d}x, \\
 & d_6 = \frac{(p_\epsilon-1) (p_\epsilon- 2) (p_\epsilon- 3)} 2  \int_\R U^{p_\epsilon- 4}   \mathcal Z^2 \varpi_1  U'(x) \,{\mathrm d}x.
\end{align}
Moreover, there holds
\begin{align}
 \int_\R  \iota^4{\bf N}_{04}U'(x) \,{\mathrm d}x  & =   \frac{\iota^4} {\beta^{1+\mathfrak{p}} } \tilde {\mathcal F}_1(\tau, \beta^{\mathfrak p} f, e).
 \end{align}

\smallskip\noindent$\clubsuit$
 On the other hand,  recall the expression of $\iota^3 S_{5}$ given in  \eqref{S5}, we can easily  get
\begin{align}  \label{bar S31+bar S51intergralUx}
\int_\R  \iota^3 S_{5} U'(x) \,{\mathrm d}x     =
\iota^3\Bigg\{ \frac1{|\beta|^3}\tilde{\mathcal D}_1(\theta)
+\frac1{\beta^{1+\mathfrak{p}} }\tilde{\mathcal F}_2\left[\tau, \beta^{\mathfrak p} f, \frac{\beta^{\mathfrak p} {f}'}{\tilde\ell}\right]
- \frac{8\pi^2}{{\tilde\ell}^2}\frac{\kappa}{\beta}  f''f \varrho_1 \Bigg\},
\end{align}
 where
\begin{align*}
\frac1{|\beta|^3}\tilde{\mathcal D}_1(\theta)
& =   - \frac{1}{6} \frac{ \iota^2W_{ttt}(0, \iota z)}{\beta^5}\int_\R x^3 U U'(x) \,{\mathrm d}x
+ \Big[ -2 \kappa \mathfrak{p} \Big[ \frac{\beta''}{\beta^4} +  (\mathfrak{p}-1)  \frac{|\beta'|^2} {\beta^5}  \Big]       +\frac{\beta' \kappa'} {\beta^4} \Big](-\Lambda)\varrho_1
\nonumber\\[2mm]
 & \quad
 +  \Big[ \frac{\kappa^3}{\beta^3}- 4\mathfrak{p}\frac{\kappa|\beta'|^2}{\beta^5}-\frac{2\kappa\beta''}{\beta^4} + \frac{\kappa'\beta'} {\beta^4} \Big]\int_\R x^2|U'|^2 \,{\mathrm d}x-2 \frac{|\beta'|^2\kappa}{\beta^5}\int_\R x^3 U'' U'(x) \,{\mathrm d}x,
\end{align*}
and
\begin{align}
  &\iota^3 \frac1{\beta^{1+\mathfrak{p}} }\tilde{\mathcal F}_2\left[\tau, \beta^{\mathfrak p} f, \frac{\beta^{\mathfrak p} {f}'}{\tilde\ell}\right]
\nonumber\\[2mm]
 & = - \frac{1}{2} \frac{\iota^2 W_{ttt}(0, \iota z)}{\beta^5}f^2 \int_\R x UU_x \,{\mathrm d}x
\nonumber\\[2mm]
 & \quad +
  \Big[
 \Big(\frac{\kappa^3}{\beta^3}
  - \frac{2\kappa\beta''}{\beta^4}
 -4 \mathfrak{p}\frac{\kappa|\beta'|^2}{\beta^5}
  +  \frac{\kappa'\beta'} {\beta^4} \Big)f^2
 + \big(8\mathfrak{p} +4 \big) \frac{\pi }{ {\tilde\ell}}\frac{\kappa\beta'}{\beta^3} f' f
  - \frac{2\pi}{{\tilde\ell}}  \frac{\kappa'}{\beta^2}  f'f\Big]\int_\R|U'|^2 {\mathrm d}x     \nonumber\\[2mm]
 & \quad +
  \Big[ -\frac{6\kappa|\beta'|^2}{\beta^5}f^2
 + \frac{16\pi}{ {\tilde\ell}} \frac{\kappa\beta'}{\beta^3} f' f
-\frac{8\pi^2} {{\tilde\ell}^2}\frac{\kappa}{\beta} |f'|^2\Big]\int_\R x U'' U'{\mathrm d}x.
\end{align}

\smallskip\noindent $\clubsuit$
 From the expression of $\iota^3 S_7$ as in \eqref{S7}, we have
\begin{align}  \label{intB4-new}
 \iota^3 \int_{\R} S_7 U'(x)\,{\mathrm d}x
= &    \frac{\iota^3}{\beta^{\mathfrak{p}} } \Bigg\{
\Big[
-\frac{4\pi^2}{{\tilde\ell}^2}f''
-\frac{2\pi}{\tilde\ell}  \frac{\beta'}{\beta^2} f'
+\frac{\beta''} {\beta^3}f
\,-\,\frac{\kappa^2}{\beta^2}f
\Big]e\int_{\R}{\mathcal Z}' U'(x)\,{\mathrm d}x
\nonumber\\[2mm]
&\qquad\quad+ \Big[2\frac{|\beta'|^2}{\beta^4} f
-\frac{4\pi}{ {\tilde\ell}}\frac{\beta'}{\beta^2}f'
\Big]e\int_{\R} x{\mathcal Z}''U'(x)\,{\mathrm d}x
\nonumber\\[2mm]
&\qquad\quad -\, \frac{\iota^2W_{tt}(0, \iota z)}{\beta^4}f e \int_{\R}x{\mathcal Z}U'(x)\,{\mathrm d}x
+
\Big[ -\frac{8\pi^2}{{\tilde\ell}^2}f'  +\frac{4\pi}{ {\tilde\ell}}\frac{\beta'}{\beta^2}  f\Big]e' \int_{\R}{\mathcal Z}' U'(x)\,{\mathrm d}x
\Bigg\}
\nonumber\\[2mm]
 = & \frac{\iota^3}{\beta^{2+\mathfrak{p}} }  \tilde{\mathcal F}_3\left[\tau, \beta^{\mathfrak p} f, \frac{\beta^{\mathfrak p} f'}{\tilde\ell} \right] e
 + \frac{\iota^3}{\beta^{\mathfrak{p}} }  \Big[ -\frac{8\pi^2}{{\tilde\ell}^2}f'  +\frac{4\pi}{ {\tilde\ell}}\frac{\beta'}{\beta^2}  f\Big]e'd_{7} -  \frac{\iota^3}{\beta^{\mathfrak{p}} }\frac{4\pi^2}{{\tilde\ell}^2}f'' e    d_{7},
\end{align}
where
\begin{align}
  \frac{\iota^3}{\beta^{2+\mathfrak{p}} } \tilde{\mathcal F}_3\left[\tau, \beta^{\mathfrak p} f, \frac{\beta^{\mathfrak p} f'}{\tilde\ell} \right]  e  &=  \frac{\iota^3}{\beta^{\mathfrak{p}} }  \Bigg\{\Big[
-\frac{2\pi}{\tilde\ell}  \frac{\beta'}{\beta^2} f'
+\frac{\beta''} {\beta^3}f
\,-\,\frac{\kappa^2}{\beta^2}f
\Big] \int_{\R}{\mathcal Z}' U'(x)\,{\mathrm d}x
\nonumber\\[2mm]
& \qquad \qquad
+\Big[2\frac{|\beta'|^2}{\beta^4} f
-\frac{4\pi}{ {\tilde\ell}}\frac{\beta'}{\beta^2}f'
\Big]\int_{\R} x{\mathcal Z}''U'(x)\,{\mathrm d}x
\nonumber\\[2mm]
& \qquad \qquad - \frac{\iota^2W_{tt}(0, \iota z)}{\beta^4}f  \int_{\R}x{\mathcal Z}U'(x)\,{\mathrm d}x \Bigg\}e,
\end{align}
 and  $d_{7} =  \int_{\R}{\mathcal Z}' U'(x)\,{\mathrm d}x. $

\smallskip\noindent $\clubsuit$
There holds
\begin{align}\label{e-second order higher terms}
I_{12}
\,=\,  &-\iota^4 \frac{8\pi^2} {{\tilde\ell}^2} \frac{{\mathbf b}_1}{\, \beta^{\mathfrak{p}} \, }  e''\int_{\R} x\mathcal Z  U'(x) \,{\mathrm d}x
     \nonumber\\[2mm]
  &+   \underbrace{ \iota^5 \int_{\R} \eta_\delta^\iota(s)\frac{4\pi^2} {{\tilde\ell}^2}  \left[ \frac{3\kappa}{\beta}(x+f) + \iota {\tilde a}_2\right] \frac{\kappa}{\, \beta^{1+ \mathfrak{p}}\, }   (x+f)\mathcal Z  U'(x) \,{\mathrm d}x\,  e''   + O(e^{\frac{- \delta} \iota} )  \frac{\iota^2 e''}{\tilde\ell^2}}_{:=  \frac{\iota^3}{\beta^{2+ \mathfrak{p}} }  \tilde{\mathcal F}_4\left[\tau, \beta^{\mathfrak{p}} f \right]  \frac{\iota^2 e''}{\tilde\ell^2}}
     \nonumber\\[2mm]
  &
-\iota^4 \frac{8\pi^2} {{\tilde\ell}^2} \frac{{\mathbf b}_1}{\, \beta^{ \mathfrak{p}} \, } f f'' e\int_{\R} \mathcal Z'  U'(x) \,{\mathrm d}x
     \nonumber\\[2mm]
  &
 + \underbrace{\iota^5 \int_{\R} \eta_\delta^\iota(s)\frac{4\pi^2} {{\tilde\ell}^2}  \left[ \frac{3\kappa}{\beta}(x+f) + \iota {\tilde a}_2\right] \frac{\kappa}{\, \beta^{1+ \mathfrak{p}}\, }   (x+f)  e\mathcal Z'  U'(x) \,{\mathrm d}x\, f''
 + O(e^{\frac{- \delta} \iota} )  \frac{f''}{\tilde\ell^2}}_{:=  \frac{\iota^5}{\beta^{2+ \mathfrak{p}} }  \tilde{\mathcal F}_5\left[\theta, \beta f, e\right]  \frac{\beta^{\mathfrak p} f''}{\tilde\ell^2}}
 \nonumber\\[2mm]
\,=\,&\iota^4  \frac{d_8} {\beta^{\mathfrak{p}}}{\mathbf b}_1 \frac{e''}{\tilde\ell^2}
 -  \iota^4 8 \pi^2 \frac{d_7} {\beta^{ \mathfrak{p}} }{\mathbf b}_1  f e\frac{f''}{\tilde\ell^2}
 + \frac{\iota^3}{\beta^{2+  \mathfrak{p}}} \tilde{\mathcal F}_4\left[\tau, \beta^{\mathfrak{p}} f \right]  \frac{\iota^2 e''}{\tilde\ell^2}
  +\frac{\iota^5}{\beta^{2+ \mathfrak{p}} }  \tilde{\mathcal F}_5\left[\tau, \beta^{\mathfrak{p}} f, e \right]  \frac{\beta^{\mathfrak p} f''}{\tilde\ell^2},
\end{align}
where
\begin{align}
d_8= -8 \pi^2\int_{\R} x\mathcal Z  U'(x) \,{\mathrm d}x.
\end{align}

\smallskip\noindent$\clubsuit$
By the expression of $B_5(\iota \phi_1)$  given in \eqref{expression of B5phi1}, we have
\begin{align}
 \int_{\R} B_{5}(\iota \phi_1)\, U'(x) \,{\mathrm d}x
 \, =\,&\iota^3
\int_{\R} \Bigg\{\big(4\mathfrak{p}+2 \big)\frac{\pi}{{\tilde\ell}}\frac{\beta'} {\beta^2}\Big[- f' \phi_{1,x} +\phi_{1,\tau}\Big]
 +\frac{4\pi^2}{{\tilde\ell}^2}\Big[|f'|^2 \phi_{1,xx}  -f''\phi_{1,x}  -2f' \phi_{1,x\tau}\Big]
\nonumber\\[2mm]
&\qquad\quad
+ \frac{4\pi}{{\tilde\ell}} \frac{\beta'}{\beta^2} (x+f)\big(- f' \phi_{1,xx}  +\phi_{1,x\tau}\big)
+\Big[\frac{\beta''}{\beta^3} + 2\mathfrak{p} \frac{|\beta'|^2}{\beta^4}  \,-\,\frac{\kappa^2}{\beta^2} \Big]
(x+f) \phi_{1,x}
\nonumber\\[2mm]
&\qquad\quad
+\Big[\mathfrak{p} \frac{\beta''}{\beta^3}  +  \mathfrak{p}  (\mathfrak{p} -1) \frac{|\beta'|^2}{\beta^4} \Big] \phi_{1}+ \frac{|\beta'|^2}{\beta^4}  (x+f)^2 \phi_{1,xx}
\nonumber\\[2mm]
&\qquad\quad-\frac{1}{2}\, \frac{\iota^2W_{tt}(0, \iota z)}{\beta^4} (x+f)^2\phi_{1}
\Bigg\} \, U'(x) \,{\mathrm d}x
 \nonumber\\[2mm]
 \, =\,& -  \iota^3 \frac{4\pi^2}{{\tilde\ell}^2} f'' \int_{\R} \phi_{1,x}\, U'(x) \,{\mathrm d}x
 +\frac{\iota^3}{\beta^{1+\mathfrak p}}\tilde{\mathcal F}_6\left[\tau, \beta^{\mathfrak p} f, \frac{\beta^{\mathfrak p} {f}'}{\tilde\ell} \right]
 \nonumber\\[2mm]
 \, =\,& -  \iota^3 \frac{4\pi^2}{{\tilde\ell}^2} f''  \frac{{\mathbf b}_1} {\Lambda}f  d_9
  +\frac{\iota^3}{\beta^{1+\mathfrak p}}\tilde{\mathcal F}_6\left[\tau, \beta^{\mathfrak p} f, \frac{\beta^{\mathfrak p} {f}'}{\tilde\ell} \right],
\end{align}
 where
\begin{align}\label{B5phi1 intergral Ux}
\frac{\iota^3}{\beta^{1+\mathfrak p}}\tilde{\mathcal F}_6\left[\tau, \beta^{\mathfrak p} f, \frac{\beta^{\mathfrak p} {f}'}{\tilde\ell} \right]
\,=\,&
 \iota^3
\int_{\R} \Bigg\{\big(4\mathfrak{p}+2 \big)\frac{\pi}{{\tilde\ell}}\frac{\beta'} {\beta^2}\Big[- f' \phi_{1,x} +\phi_{1,\tau}\Big]
 +\frac{4\pi^2}{{\tilde\ell}^2}\Big[|f'|^2 \phi_{1,xx}    -2f' \phi_{1,x\tau}\Big]
\nonumber\\[2mm]
& \qquad\quad
+ \frac{4\pi}{{\tilde\ell}} \frac{\beta'}{\beta^2} (x+f)\big[- f' \phi_{1,xx}  +\phi_{1,x\tau}\big]
+\Big[\frac{\beta''}{\beta^3} + 2\mathfrak{p} \frac{|\beta'|^2}{\beta^4}  \,-\,\frac{\kappa^2}{\beta^2} \Big]
(x+f) \phi_{1,x}
\nonumber\\[2mm]
& \qquad\quad
+\Big[\mathfrak{p} \frac{\beta''}{\beta^3}  +  \mathfrak{p}  (\mathfrak{p} -1) \frac{|\beta'|^2}{\beta^4} \Big] \phi_{1}+ \frac{|\beta'|^2}{\beta^4}  (x+f)^2 \phi_{1,xx}
\nonumber\\[2mm]
& \qquad\quad-\frac{1}{2}\, \frac{\iota^2W_{tt}(0, \iota z)}{\beta^4} (x+f)^2\phi_{1}
\Bigg\} \, U'(x) \,{\mathrm d}x
\end{align}
 and \[d_9 = \int_{\R} \varpi_{2,x}\, U'(x) \,{\mathrm d}x.  \]

\smallskip\noindent$\clubsuit$
Since $\phi_2$ is even of $x$, from the expression of $ B_4$ as in \eqref{new opreator B4},  then we have
\begin{align}\label{B4phi2 intergral Ux}
 \int_{\R} B_4(\iota^2\phi_2)\, U'(x) \,{\mathrm d}x
 \,=\,& \iota^3 \frac{\kappa}{\, \beta\, }\int_{\R} \Big[ \phi_{2,x} +\frac 1{\Lambda} (x+f)\phi_2\Big] U'(x) \,{\mathrm d}x
\nonumber\\[2mm]
 \,=\,&\iota^3 \frac{\kappa}{\, \beta\, }\int_{\R} \Big[ \phi_{2,x} +\frac 1{\Lambda} x\phi_2\Big] U'(x) \,{\mathrm d}x=  \frac{\iota^3 }{\beta^3} \tilde{\mathcal D}_2(\theta).
\end{align}

\smallskip\noindent$\clubsuit$
By the expression of $  B_{5} (\iota^2 \phi_2)$ given in    \eqref{decomposition of  B5(phi2)}, we get
\begin{align}\label{B5phi2 intergral Ux}
 \int_{\R}   B_5 (\iota^2{\phi_2}) \, U'(x) \,{\mathrm d}x
 \,=\,  &   \frac{\iota^4}{\beta^{1+\mathfrak p}}\tilde{\mathcal F}_7\left[\tau, \beta^{\mathfrak p} f, \frac{\beta^{\mathfrak p} {f}'}{\tilde\ell} \right]-\iota^4\frac{4\pi^2}{{\tilde\ell}^2} f'' \int_\R\phi_{2, x}\, U'(x) \,{\mathrm d}x
 \nonumber\\[2mm]
 \,=\,& \frac{\iota^4}{\beta^{1+\mathfrak p}}\tilde{\mathcal F}_7\left[\tau, \beta^{\mathfrak p} f, \frac{\beta^{\mathfrak p} {f}'}{\tilde\ell} \right]-\iota^4\frac{4\pi^2}{{\tilde\ell}^2} f''  \tilde{\mathcal D}_3(\theta),
\end{align}
where
\begin{align*}
\tilde{\mathcal D}_3 (\theta)= \int_\R\phi_{2, x}\, U'(x) \,{\mathrm d}x.
\end{align*}

\smallskip\noindent$\clubsuit$
By the definition of $ \mathbb{N}_1(\iota^2 \phi_2)$ given in \eqref{definition of mathbbN1}, we have
\begin{align}
\int_{\R}\,  \mathbb{N}_1(\iota^2 \phi_2)\, U'(x) \,{\mathrm d}x  =  \frac{\iota^3}{\beta}  \tilde{\mathcal D}_4 (\theta) e+\frac{\iota^4}\beta  \tilde{\mathcal F}_8\left[\tau,  e \right].
\end{align}

 \smallskip\noindent$\clubsuit$
By recalling the decomposition of $B_6 (\iota \phi_1)$ given in \eqref{expression of B6phi1}, we have
\begin{align}\label{B6phi1 intergral Ux}
\int_{\R} \eta_\delta^\iota(s)   B_6 (\iota \phi_1)\, U'(x) \,{\mathrm d}x
\,=\,    & \int_{\R} \eta_\delta^\iota(s)   \frac{\iota^4}{\beta^{1+\mathfrak{p}}}  \mathcal F_2\left[\tau, x, \beta^{\mathfrak p} f, \frac{\beta^{\mathfrak p} f'}{\tilde\ell}\right] \, U'(x) \,{\mathrm d}x
  \nonumber\\[2mm]
   &+ \int_{\R} \eta_\delta^\iota(s) \frac{\iota^4}{\beta^{1+\mathfrak{p}}} \hat{\mathcal F}_2\left[\tau, x, \beta^{\mathfrak p} f, \frac{\beta^{\mathfrak p} f''}{\tilde\ell^2} \right]  \, U'(x) \,{\mathrm d}x
  \nonumber\\[2mm]
\,=\,&
 \frac{\iota^4}{\beta^{1+\mathfrak{p}}}\tilde{\mathcal F}_9\left[\tau, \beta^{\mathfrak p} f, \frac{\beta^{\mathfrak p} {f}'}{\tilde\ell}\right]
 +  \frac{\iota^4}{\beta^{1+\mathfrak{p}}}\tilde{\mathcal F}_{10} \left[\tau, \beta^{\mathfrak p}  f \right] \frac{\beta^{\mathfrak p}  f''} {\tilde\ell^2}.
\end{align}

\smallskip\noindent$\clubsuit$
We recall the expressions of
$$
\frac{\iota^4}{\beta^{1+\mathfrak{p}}}  \mathcal F_1\left[\tau, x, \beta^{\mathfrak p} f, \frac{\beta^{\mathfrak p} f'}{\tilde\ell} \right],
\quad
\frac{\iota^4}{\beta^{1+\mathfrak{p}}}  \mathcal G_1\left[\tau, x, \beta^{\mathfrak p} f, \frac{\beta^{\mathfrak p} f'}{\tilde\ell},   e, \frac{e'}{\tilde\ell}\right],
\quad
\frac{\iota^4}{\beta^{1+\mathfrak{p}}} \hat{\mathcal F}_1\left[\tau, x, \beta^{\mathfrak p} f, \frac{\beta^{\mathfrak p} f''}{\tilde\ell^2} \right]
$$
given in \eqref{mathcal F1}, \eqref{mathcal G1}, \eqref{hat mathcal F1}.
 Then we denote
\begin{align}  \label{intergral Ux higher order}
  & \int_{\R}  \eta_\delta^\iota(s)\Bigg\{\, \frac{\iota^4}{\beta^{1+\mathfrak{p}}}  \mathcal F_1\left[\theta, x, \beta^{\mathfrak p}  f, \frac{\beta^{\mathfrak p} f'}{\tilde\ell}  \right]
+\frac{\iota^4}{\beta^{1+\mathfrak{p}}}  \mathcal G_1\left[\tau, x, \beta^{\mathfrak p} f, \frac{\beta^{\mathfrak p} f'}{\tilde\ell},   e, \frac{e'}{\tilde\ell}\right]
\nonumber\\[2mm]
   &\qquad\qquad\qquad
   + \frac{\iota^4}{\beta^{1+\mathfrak{p}}} \hat{\mathcal F}_1\left[\tau, x, \beta^{\mathfrak p} f, \frac{\beta^{\mathfrak p} f''}{\tilde\ell^2} \right]  \Bigg\} \, U'(x) \,{\mathrm d}x
\nonumber\\[2mm]
   &=
 \frac{\iota^3}{\beta^{1+\mathfrak{p}}}\tilde{\mathcal F}_{11}\left[\tau, \beta^{\mathfrak p} f, \frac{\beta^{\mathfrak p} {f}'}{\tilde\ell}, e, \frac{\iota e'}{\tilde\ell}\right]+  \frac{\iota^4}{\beta^{1+\mathfrak{p}}}\tilde{\mathcal F}_{12}\left[\tau, \beta^{\mathfrak p}  f \right] \frac{\beta^{\mathfrak p}  f''} {\tilde\ell^2},
 \end{align}
 where
\begin{align*}
   \int_{\R}  \eta_\delta^\iota(s)  \frac{\iota^4}{\beta^{1+\mathfrak{p}}} \hat{\mathcal F}_1\left[\tau, x, \beta^{\mathfrak p} f, \frac{\beta^{\mathfrak p} f''}{\tilde\ell^2} \right]\, U'(x) \,{\mathrm d}x
   =
  \frac{\iota^4}{\beta^{1+\mathfrak{p}}}\tilde{\mathcal F}_{12}\left[\tau, \beta^{\mathfrak p}  f \right] \frac{\beta^{\mathfrak p}  f''} {\tilde\ell^2}.
 \end{align*}

\medskip
At last, combining \eqref{S32+N03+B43intUx}, \eqref{bar S31+bar S51intergralUx}, \eqref{intB4-new}, \eqref{e-second order higher terms},
  \eqref{B5phi1 intergral Ux}, \eqref{B6phi1 intergral Ux}, \eqref{B4phi2 intergral Ux},  \eqref{B5phi2 intergral Ux},  and  \eqref{intergral Ux higher order},
we conclude that
\begin{align}
 &\int_{\R}\, S\big(w_2\big)\, U'(x) \,{\mathrm d}x
\nonumber\\[2mm]
\,=\,&
 \iota^2\varrho_1\Bigg[ -  \frac{4\pi^2}{{\tilde\ell}^2}f''  - 4\mathfrak{p}\frac{\pi}{\tilde\ell}  \frac{\beta'}{\beta^2} f'
+ \Big[(2\mathfrak{p}-1) \frac{|\beta'|^2}{\beta^4}
+\frac{\beta''} {\beta^3}
\,-\, \big[\frac 1{\Lambda}+1\big]\frac{\kappa^2}{\beta^2}  +\Lambda\frac{\iota^2 W_{tt}(0, \iota z)}{\beta^4} \Big]   f\Bigg]
\nonumber\\[2mm]
&\qquad
+\iota^2  \frac {{\mathbf b}_1(d_1+d_2)e}{\beta^{\mathfrak{p}}}
+ \iota^3 {\mathbf b}_1^2 \frac 1{\Lambda} f  \frac{e} {\beta^{\mathfrak{p}}} d_3
+ \iota^3 \frac{(p_\epsilon-1) (p_\epsilon- 2)(p_\epsilon- 3)} {6}  {\mathbf b}_1^3 d_4
+ \iota^3  \frac 1{\Lambda^2} {\mathbf b}_1^3f^2  (d_5 + d_6)
\nonumber\\[2mm]
&\qquad
+ \frac{\iota^4} {\beta^{1+\mathfrak{p}} } \tilde {\mathcal F}_1(\tau, \beta^{\mathfrak{p}} f, e)
+ \frac{ \iota^3}{|\beta|^3}\tilde{\mathcal D}_1(\theta)
+\frac{ \iota^3}{\beta^{1+\mathfrak{p}}}\tilde{\mathcal F}_2\left[\tau, \beta^{\mathfrak p} f, \frac{\beta^{\mathfrak p} {f}'}{\tilde\ell}\right]
- \iota^3\frac{8\pi^2}{{\tilde\ell}^2}\frac{\kappa}{\beta}  f''f \varrho_1
\nonumber\\[2mm]
&\qquad
+ \frac{\iota^3}{\beta^{2+\mathfrak{p}} }  \tilde{\mathcal F}_3\left[\tau, \beta^{\mathfrak{p}} f, \frac{\beta^{\mathfrak p} f'}{\tilde\ell} \right] e
+ \frac{\iota^3}{\beta^{\mathfrak{p}} }  \Big[ -\frac{8\pi^2}{{\tilde\ell}^2}f'  +\frac{4\pi}{ {\tilde\ell}}\frac{\beta'}{\beta^2}  f\Big]e'd_7
-  \frac{\iota^3}{\beta^{\mathfrak{p}} }\frac{4\pi^2}{{\tilde\ell}^2}f'' e    d_7
 \nonumber\\[2mm]
&\qquad
+ \iota^4  \frac{d_8} {\beta}{\mathbf b}_1 \frac{e''}{\tilde\ell^2}+ \frac{\iota^3}{\beta^{2+  \mathfrak{p}}} \tilde{\mathcal F}_4\left[\tau, \beta f \right]  \frac{\iota^2 e''}{\tilde\ell^2}
-8 \iota^4 \pi^2  \frac{d_7} {\beta^{ \mathfrak{p}} }{\mathbf b}_1  f e\frac{f''}{\tilde\ell^2}
+\frac{\iota^5}{\beta^{2+ \mathfrak{p}} }  \tilde{\mathcal F}_5\left[\tau, \beta f, e \right]  \frac{\beta^{\mathfrak p} f''}{\tilde\ell^2}
 \nonumber\\[2mm]
&\qquad
- \iota^3 \frac{4\pi^2}{{\tilde\ell}^2} f''  \frac{{\mathbf b}_1} {\Lambda}f  d_9
+\frac{\iota^3}{\beta^{1+\mathfrak{p}} }\tilde{\mathcal F}_6\left[\tau, \beta^{\mathfrak p} f, \frac{\beta^{\mathfrak p} {f}'}{\tilde\ell} \right]
+ \frac{\iota^3 }{\beta^3} \tilde{\mathcal D}_2(\theta)
\nonumber\\[2mm]
&\qquad
+ \frac{\iota^4}{\beta^{1+\mathfrak{p}} }\tilde{\mathcal F}_7\left[\tau, \beta^{\mathfrak p} f, \frac{\beta^{\mathfrak p} {f}'}{\tilde\ell} \right]
-\iota^4\frac{4\pi^2}{{\tilde\ell}^2} f''  \tilde{\mathcal D}_3(\theta)
+  \frac{\iota^3}{\beta}  \tilde{\mathcal D}_4 (\theta) e
+\frac{\iota^4}\beta  \tilde{\mathcal F}_8\left[\tau,  e \right]
+ \frac{\iota^4}{\beta^{1+\mathfrak{p}}}\tilde{\mathcal F}_9\left[\tau, \beta^{\mathfrak p} f, \frac{\beta^{\mathfrak p} {f}'}{\tilde\ell}\right]
 \nonumber\\[2mm]
&\qquad
+  \frac{\iota^4}{\beta^{1+\mathfrak{p}}}\tilde{\mathcal F}_{10} \left[\tau, \beta f \right] \frac{\beta^{\mathfrak{p}}  f''} {\tilde\ell^2}
+ \frac{\iota^3}{\beta^{1+\mathfrak{p}}}\tilde{\mathcal F}_{11}\left[\tau, \beta^{\mathfrak p} f, \frac{\beta^{\mathfrak p} {f}'}{\tilde\ell}, e, \frac{\iota e'}{\tilde\ell}\right]
+  \frac{\iota^4}{\beta^{1+\mathfrak{p}}}\tilde{\mathcal F}_{12}\left[\tau, \beta f \right] \frac{\beta^{\mathfrak{p}}  f''} {\tilde\ell^2}
 \nonumber\\[2mm]
& =
 \iota^2\varrho_1\Bigg[ -  \frac{4\pi^2}{{\tilde\ell}^2}f''  - 4\mathfrak{p}\frac{\pi}{\tilde\ell}  \frac{\beta'}{\beta^2} f'
+ \Big[(2\mathfrak{p}-1) \frac{|\beta'|^2}{\beta^4}
+\frac{\beta''} {\beta^3}
\,-\, \Big(\frac 1{\Lambda}+1\Big)\frac{\kappa^2}{\beta^2}  +\Lambda\frac{\iota^2 W_{tt}(0, \iota z)}{\beta^4} \Big]   f\Bigg]
  \nonumber\\[2mm]
&  \quad
+ \iota^2\frac { {\mathbf b}_1(d_1+d_2)e}{\beta^{\mathfrak{p}}}
+ \iota^4  \frac{d_8} {\beta^{\mathfrak{p}}}{\mathbf b}_1 \frac{e''}{\tilde\ell^2}
+   \frac{\iota^3}{\beta^{\mathfrak{p}} }  \Big[ -\frac{8\pi^2}{{\tilde\ell}^2}f'  +\frac{4\pi}{ {\tilde\ell}}\frac{\beta'}{\beta^2}  f\Big]e'd_7
   \nonumber\\[2mm]
&  \quad
+\frac{\iota^3}{\beta^{1+\mathfrak{p}}}\tilde{\mathcal F}^{*} \left[\tau, \beta^{\mathfrak p} f, \frac{\beta^{\mathfrak p} {f}'}{\tilde\ell},  e, \frac{\iota e'}{\tilde\ell}\right]
 + \frac{ \iota^3}{\beta^3}\tilde{\mathcal D}(\theta)
 + \frac{\iota^3}{\beta^{1+\mathfrak{p}}}\tilde{\mathcal F}^{**} \left[\tau, \beta^{\mathfrak p} f,  e, \frac{\iota^2 e''}{\tilde\ell^2},  \frac{\beta^{\mathfrak p} f''}{\tilde\ell^2}\right] ,
\end{align}
where
\begin{align}
\frac{ \iota^3}{\beta^3}\tilde{\mathcal D}(\theta)  =  \iota^3\frac{(p_\epsilon-1) (p_\epsilon- 2)(p_\epsilon- 3)} {6}  {\mathbf b}_1^3 d_4 +  \frac{ \iota^3}{|\beta|^3}\tilde{\mathcal D}_1(\theta) +   \frac{ \iota^3}{\beta^3}\tilde{\mathcal D}_2(\theta),
\end{align}
\begin{align}
\frac{\iota^3}{\beta^{1+\mathfrak{p}}}\tilde{\mathcal F}^*  \left[\tau, \beta^{\mathfrak p} f, \frac{\beta^{\mathfrak p} {f}'}{\tilde\ell},  e, \frac{\iota e'}{\tilde\ell}\right]
&=\iota^3    {\mathbf b}_1^2 \frac 1{\Lambda} f  \frac{e} {\beta^{\mathfrak{p}}} d_3
+ \iota^3\frac{(p_\epsilon-1) (p_\epsilon- 2)(p_\epsilon- 3)} {6}  {\mathbf b}_1^3 d_4
+  \iota^3    \frac 1{\Lambda^2} {\mathbf b}_1^3f^2  (d_5 + d_6)
\nonumber\\[2mm]
&\quad
+ \frac{\iota^4} {\beta^{1+\mathfrak{p}} } \tilde {\mathcal F}_1(\tau, \beta^{\mathfrak p} f, e)
+ \frac{ \iota^3}{|\beta|^3}\tilde{\mathcal D}_1(\theta)
+\frac{ \iota^3}{\beta^{1+\mathfrak{p}}}\tilde{\mathcal F}_2\left[\tau, \beta^{\mathfrak p} f, \frac{\beta^{\mathfrak p} {f}'}{\tilde\ell}\right]
\nonumber\\[2mm]
&\quad
+ \frac{\iota^3}{\beta^{2+\mathfrak{p}} }  \tilde{\mathcal F}_3\left[\tau, \beta^{\mathfrak p} f, \frac{\beta^{\mathfrak p} f'}{\tilde\ell} \right] e
+\frac{\iota^3}{\beta^{1+\mathfrak{p}}}\tilde{\mathcal F}_6\left[\tau, \beta^{\mathfrak p} f, \frac{\beta^{\mathfrak p} {f}'}{\tilde\ell} \right]
\nonumber\\[2mm]
&\quad
+    \frac{\iota^4}{\beta^{1+\mathfrak{p}}}\tilde{\mathcal F}_7\left[\tau, \beta^{\mathfrak p} f, \frac{\beta^{\mathfrak p} {f}'}{\tilde\ell}\right]
+  \frac{\iota^3}{\beta}  \tilde{\mathcal D}_4 (\theta) e+\frac{\iota^4}\beta  \tilde{\mathcal F}_8\left[\tau,  e \right]
\nonumber\\[2mm]
&\quad
+ \frac{\iota^4}{\beta^{1+\mathfrak{p}}}\tilde{\mathcal F}_9\left[\tau, \beta^{\mathfrak p} f, \frac{\beta^{\mathfrak p} {f}'}{\tilde\ell}\right]
+  \frac{\iota^3}{\beta^{1+\mathfrak{p}}}\tilde{\mathcal F}_{11}\left[\tau, \beta^{\mathfrak p} f, \frac{\beta^{\mathfrak p} {f}'}{\tilde\ell}, e, \frac{\iota e'}{\tilde\ell}\right],
\end{align}
and
\begin{align}
\frac{\iota^3}{\beta^{1+\mathfrak{p}}}\tilde{\mathcal F}^{**} \left[\tau, \beta^{\mathfrak p} f,  e, \frac{\iota^2 e''}{\tilde\ell^2},  \frac{\beta^{\mathfrak p} f''}{\tilde\ell^2}\right]
\,=\,& -  \iota^3\frac{8\pi^2}{{\tilde\ell}^2}\frac{\kappa}{\beta}  f''f \varrho_1   -  \frac{\iota^3}{\beta^{\mathfrak{p}} }\frac{4\pi^2}{{\tilde\ell}^2}f'' e    d_7
  + \frac{\iota^3}{\beta^{2+  \mathfrak{p}}} \tilde{\mathcal F}_4\left[\tau, \beta^{\mathfrak p} f \right]  \frac{\iota^2 e''}{\tilde\ell^2}
 \nonumber\\[2mm]
&\qquad  - \iota^4 8 \pi^2  \frac{d_7} {\beta^{ \mathfrak{p}} }{\mathbf b}_1  f e\frac{f''}{\tilde\ell^2}+\frac{\iota^5}{\beta^{2+ \mathfrak{p}} }  \tilde{\mathcal F}_5\left[\tau, \beta^{\mathfrak p} f, e \right]  \frac{\beta^{\mathfrak p} f''}{\tilde\ell^2}
  -  \iota^3 \frac{4\pi^2}{{\tilde\ell}^2} f''  \frac{{\mathbf b}_1} {\Lambda}f  d_9
  \nonumber\\[2mm]
&\qquad
+  \frac{\iota^4}{\beta^{1+\mathfrak{p}}}\tilde{\mathcal F}_{10} \left[\tau, \beta^{\mathfrak{p}} f \right] \frac{\beta^{\mathfrak p} f''} {\tilde\ell^2}
+  \frac{\iota^4}{\beta^{1+\mathfrak{p}}}\tilde{\mathcal F}_{12}\left[\tau, \beta^{\mathfrak{p}} f \right] \frac{\beta^{\mathfrak p}f''} {\tilde\ell^2}.
\end{align}
Moreover,  by the expressions of the following terms
$$
\frac{\iota^3}{\beta^{1+\mathfrak{p}}}\tilde{\mathcal F}^{**} \left[\tau, \beta^{\mathfrak p} f,  e, \frac{\iota^2 e''}{\tilde\ell^2},  \frac{\beta^{\mathfrak p} f''}{\tilde\ell^2}\right],
\qquad
\frac{\iota^3}{\beta^{1+\mathfrak{p}}}\tilde{\mathcal F}^{*} \left[\theta, \beta^{\mathfrak p} f, \frac{\beta^{\mathfrak p} {f}'}{\tilde\ell},  e, \frac{\iota e'}{\tilde\ell}\right],
$$
we can  verify that
\begin{align}
\left\|  \beta  \frac{ \iota^3}{\beta^3}\tilde{\mathcal D}(\theta)  \right\|_{L^2(0, 2\pi)}  \le  C\iota^3,
\end{align}
\begin{align}
\left\|  \beta \frac{\iota^3}{\beta^{1+\mathfrak{p}}}\tilde{\mathcal F}^{**} \left[\tau, \beta^{\mathfrak p}  f,  e, \frac{\iota^2 e''}{\tilde\ell^2},  \frac{\beta^{\mathfrak p} f''}{\tilde\ell^2}\right] \right\|_{L^2(0, 2\pi)}  \le  C (\iota^3+ \iota^3 \|f\|_* + \iota^3 \|e\|_{**} ),
\end{align}
\begin{align}
\left\|\beta  \frac{\iota^3}{\beta^{1+\mathfrak{p}}}\tilde{\mathcal F}^{*} \left[\theta, \beta^{\mathfrak p} f, \frac{\beta^{\mathfrak p} {f}'}{\tilde\ell},  e, \frac{\iota e'}{\tilde\ell}\right]   \right\|_{L^2(0, 2\pi)}
\le C (\iota^3+ \iota^3 \|f\|_* +\iota^3 \|e\|_{**} ),
\end{align}
 and
\begin{align}
& \left\| \beta \frac{\iota^3}{\beta}\left(\tilde{\mathcal F}^{*} \left[\theta, \beta^{\mathfrak p} f_1, \frac{\beta^{\mathfrak p} {f_1}'}{\tilde\ell},  e_1, \frac{\e e_1'}{\tilde\ell}\right]    - \tilde{\mathcal F}^{*} \left[\theta, \beta^{\mathfrak p} f_2, \frac{\beta^{\mathfrak p} {f_2}'}{\tilde\ell},  e_2, \frac{\iota e_2'}{\tilde\ell}\right]  \right) \right\|_{L^2(0, 2\pi)}
\nonumber\\[2mm]
&
\le  \iota^3 \left(\|f_1- f_2\|_*+ \|e_1- e_2\|_{**}\right).
\end{align}

\medskip
\medskip
\subsubsection{Computation for $ \beta^{-\mathfrak{p}-2}\int_{\R}\, \eta_\delta^\iota(s)\, {\mathcal E}\,{\mathcal Z}  \,{\mathrm d}x$}\

By \eqref{errorrelationinterior}, we have
\begin{align}
& \frac{1}{\beta^{\mathfrak{p}+2}}\int_\R\eta_\delta^\iota(s){\mathcal E} \,{\mathcal Z} \,{\mathrm d}x
\nonumber\\[2mm]
\,=\,&\frac{1}{\beta^{\mathfrak{p}+2}}\int_\R {\mathcal E} \,{\mathcal Z} \,{\mathrm d}x
\ +\
\frac{1}{\beta^{\mathfrak{p}+2}}\int_\R \big[\eta_\delta^\iota(s)-1\big]{\mathcal E} \,{\mathcal Z} \,{\mathrm d}x
\nonumber\\[2mm]
\,=\,& \int_{\R}\, S\big(w_2\big)\,  \,{\mathcal Z} \,{\mathrm d}x
\ +\ \int_\R \big[\eta_\delta^\iota(s)-1\big]S\big(w_2\big) \,{\mathcal Z} \,{\mathrm d}x
\nonumber\\[2mm]
\,=\,&
\int_\R \Bigg[\iota^2 S_{3}  + \iota^2 S_{42}  + \iota^3 S_{5} +\iota^3 S_{6}
 +\iota^3 S_{7} +\iota^3 S_{8}
  +\frac 1{\beta^{\mathfrak{p}}} \iota^3\frac{4\pi^2}{{\tilde\ell}^2}  e''\mathcal Z
  +\frac 1{\beta^{\mathfrak{p}}}  \iota\lambda_0 e{\mathcal Z}
  \nonumber\\[2mm]
&\qquad
+ N_{0\e} \left(\iota\frac{e} {\beta^{\mathfrak{p}}} \mathcal Z\right)
+ \iota^2 {\bf N}_{02} +\iota^3{\bf N}_{03} + \iota^4 {\bf N}_{04} +  \iota^2 {\bf B}_{42}
 +  \iota^3 \frac{4\pi^2}{{\tilde\ell}^2} {\phi}_{1, \tau\tau}
  \nonumber\\[2mm]
&\qquad
+ \iota^4 \frac{4\pi^2}{{\tilde\ell}^2}\phi_{2, \tau\tau}
+B_{5}(\iota\phi_1)    + \mathbb{N}_1 (\iota^2 \phi_2)
+ B_4(\iota^2\phi_2)+  B_5 (\iota^2{\phi_2})   \Bigg] \,{\mathcal Z} \,{\mathrm d}x  +  \mathfrak P(\theta)
\nonumber\\[2mm]
&
  +\int_\R \eta_\delta^\iota(s) \Bigg\{ \iota^4 \frac{4\pi^2} {{\tilde\ell}^2} a_2\frac{\kappa}{\, \beta^{1+ \mathfrak{p}} \, }  e'' (x+f)\mathcal Z -\iota^4 \kappa a_2 \frac{ 4\pi^2} {{\tilde\ell}^2} \frac{1}{\, \beta^{1+\mathfrak{p}} \, }  (x+f) f'' e   \mathcal Z'
   \nonumber\\[2mm]
&\qquad \qquad\qquad
+ \frac{\iota^4}{\beta^{1+\mathfrak{p}}} \hat{\mathcal F}_1\left[\tau, x, \beta^{\mathfrak p} f, \frac{\beta^{\mathfrak p} f''}{\tilde\ell^2} \right]
  +  \frac{\iota^4}{\beta^{1+\mathfrak{p}}} \mathcal F_1\left[\tau, x, \beta^{\mathfrak p} f, \frac{\beta^{\mathfrak p} f'}{\tilde\ell}\right]
 \nonumber\\[2mm]
&\qquad \qquad\qquad
+B_{6}(\iota\phi_1)
 +\frac{\iota^4}{\beta^{1+\mathfrak{p}}}  \mathcal G_1\left[\tau, x, \beta^{\mathfrak p} f, \frac{\beta^{\mathfrak p} f'}{\tilde\ell},   e, \frac{e'}{\tilde\ell}\right]
 +  B_6 (\iota^2{\phi_2}) \Bigg\} \,{\mathcal Z} \,{\mathrm d}x,
\end{align}
 where
\begin{align*}
 \mathfrak P(\theta)  & =    \int_\R \big[\eta_\delta^\iota(s)-1\big]S\big(w_2\big) \,{\mathcal Z} \,{\mathrm d}x.
\end{align*}

By the exponential decay of the terms in  $S\big(w_2\big)$ with respect to $x$,
 we  neglect  the calculations of the terms in $ \mathfrak P(\theta). $
Since $S_{3},  S_5, S_7$ are odd functions of $x$, then integration against $\mathcal Z$ just vanishes.
This gives that
\begin{align}\label{projected eqs for enew}
 \frac{1}{\beta^{\mathfrak{p}+2}}\int_\R\eta_\delta^\iota(s){\mathcal E}\,{\mathcal Z} \,{\mathrm d}x
=&
\int_\R \Bigg\{ \iota^2 S_{42}   +\iota^3 S_{6}
+\iota^3 S_{8}
+\iota^3 \frac 1{\beta^{\mathfrak{p}}} \frac{4\pi^2}{{\tilde\ell}^2}  e''\mathcal Z
+\iota \frac 1{\beta^{\mathfrak{p}}}  \lambda_0 e{\mathcal Z}   + N_{0\e} \left(\iota\frac{e} {\beta^{\mathfrak{p}}} \mathcal Z\right)
  \nonumber\\[2mm]
&\qquad
     + \iota^2 {\bf N}_{02} +\iota^3{\bf N}_{03}  + \iota^4 {\bf N}_{04} +  \iota^2 {\bf B}_{42}
 +  \iota^3 \frac{4\pi^2}{{\tilde\ell}^2} {\phi}_{1, \tau\tau}
+ \iota^4 \frac{4\pi^2}{{\tilde\ell}^2}\phi_{2, \tau\tau}
 \nonumber\\[2mm]
&\qquad
+B_{5}(\iota\phi_1)    + \mathbb{N}_1 (\iota^2 \phi_2)
+ B_4(\iota^2\phi_2)+  B_5 (\iota^2{\phi_2})    \Bigg\}\,{\mathcal Z} \,{\mathrm d}x
\nonumber\\[2mm]
&+\int_\R \eta_\delta^\iota(s) \Bigg\{ \iota^4 \frac{4\pi^2} {{\tilde\ell}^2} a_2\frac{\kappa}{\, \beta^{1+ \mathfrak{p}} \, }  e'' (x+f)\mathcal Z
-\iota^4 \kappa a_2 \frac{ 4\pi^2} {{\tilde\ell}^2} \frac{1}{\, \beta^{1+\mathfrak{p}} \, }  (x+f) f'' e   \mathcal Z'
 \nonumber\\[2mm]
&\qquad\qquad\qquad
 + \frac{\iota^4}{\beta^{1+\mathfrak{p}}} \hat{\mathcal F}_1\left[\tau, x, \beta^{\mathfrak p} f, \frac{\beta^{\mathfrak p} f''}{\tilde\ell^2} \right]
+B_{6}(\iota\phi_1)
+  \frac{\iota^4}{\beta^{1+\mathfrak{p}}} \mathcal F_1\left[\tau, x, \beta^{\mathfrak p} f, \frac{\beta^{\mathfrak p} f'}{\tilde\ell}\right]
 \nonumber\\[2mm]
&\qquad\qquad\qquad
+\frac{\iota^4}{\beta^{1+\mathfrak{p}}}  \mathcal G_1\left[\tau, x, \beta^{\mathfrak p} f, \frac{\beta^{\mathfrak p} f'}{\tilde\ell},   e, \frac{e'}{\tilde\ell}\right]
+  B_6 (\iota^2{\phi_2}) \Bigg\} \,{\mathcal Z} \,{\mathrm d}x.
\end{align}

\smallskip\noindent$\clubsuit$
 From the expression of $\iota^3 S_{42}$ as in \eqref{barS41}, we have
\begin{align} \label{S42 on Z}
 \int_{\R}  \iota^2 S_{42}\mathcal Z  \,{\mathrm d}x
\,=\,&\iota^2  \int_{\R} \Bigg\{
-\frac{1}{2}\, \frac{ \iota^2W_{tt}(0, \iota z)}{\beta^4}|f|^2 U
+
\Big(\frac{4\pi^2}{{\tilde\ell}^2}|f'|^2
-\frac{4\pi}{{\tilde\ell}}\frac{\beta'}{\beta^2}ff'
+\frac{|\beta'|^2}{\beta^4} |f|^2
\Big)U''
\nonumber\\[2mm]
 & \qquad\qquad
+ \frac1{\Lambda}  \frac{\kappa}{\, \beta^{1+ \mathfrak{p}} \, } fe{\mathcal Z}
+\frac{(p_\epsilon-1) (p_\epsilon- 2)}{2\beta^{2\mathfrak{p}}}e^2U^{p_\epsilon- 3} {\mathcal Z}^2
\Bigg\}\mathcal Z  \,{\mathrm d}x
\nonumber\\[2mm]
\,=\,& -\frac{1}{2}\, \iota^2 \frac{ \iota^2W_{tt}(0, \iota z)}{\beta^4}|f|^2  b_2
+ \iota^2\Big(\frac{4\pi^2}{{\tilde\ell}^2}|f'|^2
-\frac{4\pi}{{\tilde\ell}}\frac{\beta'}{\beta^2}ff'
+\frac{|\beta'|^2}{\beta^4} |f|^2
\Big) b_3
\nonumber\\[2mm]
& +  \iota^2  \frac1{\Lambda}\frac{\kappa}{\, \beta^{1+\mathfrak{p}}\, } fe b_0
 +  \iota^2\frac{(p_\epsilon-1) (p_\epsilon- 2)}{2\beta^{2\mathfrak{p}}}e^2b_1 ,
\end{align}
where
\[b_0 = \int_{\R}{\mathcal Z}^2  \,{\mathrm d}x=1,
\qquad  b_1= \int_{\R}U^{p_\epsilon- 3}{\mathcal Z}^3  \,{\mathrm d}x,
\quad b_2 = \int_{\R}  U\mathcal Z  \,{\mathrm d}x,
\quad b_3 = \int_{\R} U''{\mathcal Z}  \,{\mathrm d}x.   \]

  \smallskip\noindent $\clubsuit$
 From the expression of $\iota^3 S_{6} $ as in \eqref{S6}, we have
\begin{align} \label{S6 on Z}
 \int_{\R}  \iota^3 S_{6} \mathcal Z  \,{\mathrm d}x
 =&\frac{\iota^3}{\beta^{1+\mathfrak{p}}}\tilde{\mathcal K}_2\left[\tau, {\beta^\mathfrak{p} f}, \frac{\beta^\mathfrak{p}{ f}'}{\tilde\ell} \right]
+\iota^3\frac{8\pi^2}{{\tilde\ell}^2}\frac{\kappa}{\beta}   f''b_4,
\end{align}
where
$b_4 =  \int_\R xU' \mathcal Z \,{\mathrm d}x,
$
and
\begin{align*}
& \frac{\iota^3}{\beta^{1+\mathfrak{p}}}\tilde{\mathcal K}_2\left[\tau, {\beta^\mathfrak{p} f}, \frac{\beta^\mathfrak{p}{ f}'}{\tilde\ell} \right]
\nonumber\\[2mm]
 &  =  \iota^3\Bigg\{
 - \frac{1}{2} \frac{\,\iota^2 W_{ttt}(0, \iota z)}{\beta^5}f \int_\R x^2 U \mathcal Z\,{\mathrm d}x
\ +\
 \Big[-\frac{\, \, \iota^2 W_{ttt}(0, \iota z)}{6\beta^5}f^3 -2\kappa\mathfrak{p}\big[\frac{\beta''}{\beta^4} + (\mathfrak{p}-1)  \frac{|\beta'|^2}{\beta^5}  \big]   f  + \mathfrak{p} \frac{\kappa'\beta'} {\beta^4} f\Big]\int_\R U\mathcal Z\,{\mathrm d}x
\nonumber\\[2mm]
 & \quad +  \Big[
\Big( \frac{2\kappa^3}{\beta^3}
 -4\mathfrak{p} \frac{\kappa|\beta'|^2}{\beta^5}
 -\frac{2\kappa\beta''}{\beta^4}f +\frac{2\kappa'\beta'}{\beta^4} \Big)f
 -\frac{8\pi^2}{{\tilde\ell}^2}\frac{\kappa}{\beta}   f''
+\Big(\big[8\mathfrak{p} +4 \big]\frac{\pi }{ {\tilde\ell}}\frac{\kappa\beta'}{\beta^3} - \frac{2\pi}{{\tilde\ell}} \frac{\kappa'}{\beta^2} \Big) f'
 \Big]\int_\R xU'\mathcal Z\,{\mathrm d}x
\nonumber\\[2mm]
 & \quad     +  \Big[ -\frac{6\kappa|\beta'|^2}{\beta^5}f
 +\frac{8\pi}{ {\tilde\ell}}\frac{\kappa\beta'}{\beta^3}   f'\Big]  \int_\R x^2 U''\mathcal Z\,{\mathrm d}x
  + \Big[ -\frac{2\kappa|\beta'|^2}{\beta^5}f^3
 +\frac{8\pi}{ {\tilde\ell}} \frac{\kappa\beta'}{\beta^3} f' f^2
 -\frac{8\pi^2} {{\tilde\ell}^2}\frac{\kappa}{\beta} |f'|^2 f \Big]\int_\R U''\mathcal Z\,{\mathrm d}x
 \Bigg\}.
\end{align*}

\smallskip\noindent $\clubsuit$
 From the expression of $\iota^3 S_8$ as in \eqref{S8}, we have
\begin{align} \label{S8 on Z}
\int_{\R} \iota^3 S_8{\mathcal Z}  \,{\mathrm d}x
\,=\,& \frac{\iota^3}{\beta^{\mathfrak{p}}}\int_{\R} \Bigg\{
 \frac{ 2\pi}{\tilde\ell} \frac{\beta'}{\beta^2}  e'{\mathcal Z}
+\frac{4\pi}{{\tilde\ell}}\frac{\beta'}{\beta^2} e' x {\mathcal Z}'
-\frac{1}{2}\,  \frac{\iota^2W_{tt}(0, \iota z)}{\beta^4}f^2  e  {\mathcal Z}
\nonumber\\[2mm]
&\qquad \qquad
+\Big[\frac{4\pi^2}{{\tilde\ell}^2}|f'|^2
-\frac{4\pi}{ {\tilde\ell}}\frac{\beta'}{\beta^2}ff'
+\frac{|\beta'|^2}{\beta^4} f^2
\Big]e{\mathcal Z}'' -\frac{1}{2}\, \frac{\iota^2W_{tt}(0, \iota z)}{\beta^4}  e x^2{\mathcal Z}
\nonumber\\[2mm]
&\qquad \qquad
+\Big[\frac{\beta''} {\beta^3}
\,-\,\frac{\kappa^2}{\beta^2}
\Big]e x{\mathcal Z}'
+\frac{|\beta'|^2}{\beta^4} e x^2{\mathcal Z}''
\Bigg\} {\mathcal Z} \,{\mathrm d}x
\nonumber\\[2mm]
\,=\,&\iota^3
\frac{1}{\beta^{1+\mathfrak{p}} } \tilde{\mathcal K}_3\left[\tau,  \beta^{\mathfrak p} f, \frac{\beta^{\mathfrak p} { f}'}{\tilde\ell}\right] e,
\end{align}
where we have used the fact that
\begin{align*}
\int_\R x{\mathcal Z}'{\mathcal Z}  \,{\mathrm d}x=  - \frac12\int_\R |\mathcal Z|^2 \,{\mathrm d}x,
\end{align*}

 \smallskip\noindent $\clubsuit$
 By the expression of $N_{0\e} \left(\iota\frac{e} {\beta^{\mathfrak{p}}} \mathcal Z\right)$ given in \eqref{N0epsilonZ}, we have
\begin{align}  \label{N0epsilon on Z}
  \int_\R
  N_{0\e} \left(\iota\frac{e} {\beta^{\mathfrak{p}}} \mathcal Z\right)  {\mathcal Z} \,{\mathrm d}x  =  \frac{\iota^3}{\beta^{3\mathfrak{p}} } \tilde{\mathcal K}_4\left[\tau,   e\right] e^3,
 \end{align}
  where  $\tilde{\mathcal K}_4\left[\tau,   e\right] $ is uniformly bounded.

\smallskip\noindent $\clubsuit$
 From the expression of $ \iota^2 {\bf B}_{42}$ as in \eqref{B42}, we have
\begin{align}\label{B4emathcal Z2}
 \int_{\R}\iota^2{\bf B}_{42}\mathcal Z  \,{\mathrm d}x
 =\iota^2\int_{\R}   {\mathbf b}_1^2 \Big[
\frac{1}{\Lambda} f\varpi_{2,x}
+ \frac{1}{\Lambda^2} f x\varpi_2
+\frac{1}{\Lambda} f\varpi_1
+ \frac{1}{\Lambda^2} f^2\varpi_2
\Big]\mathcal Z  \,{\mathrm d}x = \iota^2 \frac{1}{\Lambda^2} f^2  {\mathbf b}_1^2 b_5,
\end{align}
where
$$
b_5=  \int_\R \varpi_2\mathcal Z  \,{\mathrm d}x.
$$

\smallskip\noindent $\clubsuit$
 From the expressions of $ {\bf N}_{02}, {\bf N}_{03},  {\bf N}_{04}$ as in \eqref{N02}-\eqref{N03}, we have
\begin{align}\label{nonlinear term N02emathcal Z}
  \int_{\R}  \iota^2 {\bf N}_{02} {\mathcal Z}  \,{\mathrm d}x
  \,=\, &\iota^2  (p_\epsilon-1)(p_\epsilon- 2)  \Bigg[  \frac 1{2\Lambda^2} {\mathbf b}_1^2f^2  \int_{\R} U^{p_\epsilon- 3} \varpi_2^2{\mathcal Z}  \,{\mathrm d}x
+  \frac{e} {\beta^{\mathfrak{p}}}  \frac 1{\Lambda} {\mathbf b}_1f   \int_{\R}  U^{p_\epsilon- 3} \varpi_2 |\mathcal Z|^2  \,{\mathrm d}x         \Bigg]
\nonumber\\[2mm]
  \,=\, &  \iota^2
 (p_\epsilon-1)(p_\epsilon- 2)  \Big[  \frac 1{2\Lambda^2} {\mathbf b}_1^2f^2b_6
+  \frac{e} {\beta^{\mathfrak{p}}}  \frac 1{\Lambda} {\mathbf b}_1f  b_7\Big],
\end{align}
and
\begin{align}\label{nonlinear term N03 emathcal Z}
  & \int_{\R}  \left[ \iota^3 {\bf N}_{03} +  \iota^4 {\bf N}_{04} \right]  {\mathcal Z}  \,{\mathrm d}x  =  \frac{\iota^3} {\beta^{1+ \mathfrak{p}}} \tilde{\mathcal K}_5\left[\tau,  {\beta^{\mathfrak{p}} f},  e\right],
\end{align}
where
\begin{align*}
b_6 =\int_{\R} U^{p_\epsilon- 3} \varpi_2^2{\mathcal Z}  \,{\mathrm d}x, \quad b_7 =  \int_{\R}  U^{p_\epsilon- 3} \varpi_2 |\mathcal Z|^2  \,{\mathrm d}x.
\end{align*}
 Moreover  the term $\tilde{\mathcal K}_5\left[\tau,  {\beta^{\mathfrak p} f},  e\right] $ is bounded and satisfies
\begin{align*}
 \left| \tilde{\mathcal K}_5\left[\tau,  {\beta^{\mathfrak p} f_1},  e_1\right]
 -  \tilde{\mathcal K}_5\left[\tau,  {\beta^{\mathfrak p} f_2},  e_2\right] \right|
 \le C\Big( \beta^{\mathfrak{p}} |f_1- f_2| + |e_1- e_2|\Big).
 \end{align*}

\smallskip\noindent $\clubsuit$
Next, we have
\begin{align} \label{main term mathccal Z}
 & \int_{\R} \Bigg[\frac 1{\beta^{\mathfrak{p}}} \iota^3\frac{4\pi^2}{{\tilde\ell}^2}  e''\mathcal Z
 +\frac 1{\beta^{\mathfrak{p}} } \iota\lambda_0 e{\mathcal Z}   \Bigg] {\mathcal Z}  \,{\mathrm d}x
  =  \frac 1{\beta^{\mathfrak{p}}} \Bigg[\iota^3\frac{4\pi^2}{{\tilde\ell}^2}  e''  +  \iota\lambda_0 e \Bigg].
\end{align}

\smallskip\noindent $\clubsuit$
   Let's review the expression for $\phi_{1, \tau\tau} $ in \eqref{expression of phi1tautau}, then we have
\begin{align} \label{phitautau mathcal Z}
\iota^3 \frac{4\pi^2}{{\tilde\ell}^2}\int_\R\phi_{1, \tau\tau}  {\mathcal Z} \,{\mathrm d}x
\,=\,&\iota^3  \int_\R \Big[  {\mathbf b}_3 \varpi_{1}
+
\frac{1}{\Lambda} {\mathbf b}_3 f \varpi_{2}
+
\frac{2}{\Lambda} \frac{2\pi}{{\tilde\ell}}{\mathbf b}_2f' \varpi_{2}
+
\frac{1}{\Lambda} \frac{4\pi^2}{{\tilde\ell}^2} {\mathbf b}_1 f''\varpi_{2}  \Big]
\mathcal Z \,{\mathrm d}x
\nonumber\\[2mm]
 =\, &
\frac{\iota^3}{\beta^{1+\mathfrak p}}\tilde{\mathcal K}_5\left[\tau,  \beta^{\mathfrak p} f, \frac{\beta^{\mathfrak p} {f}'}{\tilde\ell}\right]
+  \frac{\iota^3}{\beta^{\mathfrak p}}\frac{4\pi^2} {\Lambda}{\mathbf b}_1 \frac{\beta^{\mathfrak p} f''}{\tilde \ell^2}b_3,
\end{align}
where
\begin{align*}
\frac{\iota^3}{\beta^{1+\mathfrak p}}\tilde{\mathcal K}_5\left[\tau,  \beta^{\mathfrak p} f, \frac{\beta^{\mathfrak p} {f}'}{\tilde\ell}\right]  = \iota^3  \Big[
 \frac{1}{\Lambda} {\mathbf b}_3 f
+
\frac{2}{\Lambda} \frac{2\pi}{{\tilde\ell}}{\mathbf b}_2f'
 \Big] \int_\R
\varpi_{2}  \mathcal Z \,{\mathrm d}x.
\end{align*}

\smallskip\noindent $\clubsuit$
By the definition of $\iota^2\phi_2$ given in \eqref{definition of e2phi2}, we  directly get
\begin{align} \label{phi2tautau mathcal Z}
\iota^4 \frac{4\pi^2}{{\tilde\ell}^2}\int_\R\phi_{2, \tau\tau}  {\mathcal Z} \,{\mathrm d}x =\frac{\iota^4}{\beta^{1+\mathfrak{p}}} \tilde{\mathcal D}_5(\theta).
\end{align}

\smallskip\noindent
$\clubsuit$ From the expression of $B_{5}(\iota \phi_1)$ as in \eqref{expression of B5phi1}, we have
\begin{align}\label{B5phi1 intergral Z}
 \int_{\R} B_{5}(\iota \phi_1)\,{\mathcal Z} \,{\mathrm d}x
\,=\,&
\iota^3
\int_\R
\Bigg\{\big(4\mathfrak{p}+2 \big)\frac{\pi}{{\tilde\ell}}\frac{\beta'} {\beta^2}\Big[- f' \phi_{1,x} +\phi_{1,\tau}\Big]
 +\frac{4\pi^2}{{\tilde\ell}^2}\Big[|f'|^2 \phi_{1,xx}  -f''\phi_{1,x}  -2f' \phi_{1,x\tau}\Big]
\nonumber\\[2mm]
&\qquad\quad
+ \frac{4\pi}{{\tilde\ell}} \frac{\beta'}{\beta^2} (x+f)\big[- f' \phi_{1,xx}  +\phi_{1,x\tau}\big]
\nonumber\\[2mm]
&\qquad\quad+\Big[\frac{\beta''}{\beta^3} +2\mathfrak{p} \frac{|\beta'|^2}{\beta^4}  \,-\,\frac{\kappa^2}{\beta^2} \Big]
(x+f) \phi_{1,x}
 +\Big[\mathfrak{p} \frac{\beta''}{\beta^3}  +  \mathfrak{p}  (\mathfrak{p} -1) \frac{|\beta'|^2}{\beta^4} \Big] \phi_1
 \nonumber\\[2mm]
&\qquad\quad+ \frac{|\beta'|^2}{\beta^4}  (x+f)^2 \phi_{1,xx}
-\frac{1}{2}\, \frac{\iota^2W_{tt}(0, \iota z)}{\beta^4} (x+f)^2\phi_{1}
\Bigg\}  {\mathcal Z} \,{\mathrm d}x
\nonumber\\[2mm]
\,=\,&  -  \iota^3 \frac{4\pi^2}{{\tilde\ell}^2} f'' \int_{\R}  \left[{\mathbf b}_1\varpi_{1,x} +\frac{1}{\Lambda}{\mathbf b}_1 f \varpi_{2,x} \right]\,  \mathcal Z \,{\mathrm d}x +  \frac{\iota^3}{\beta}\tilde{\mathcal K}_7\left[\tau,  \beta^{\mathfrak p} f, \frac{\beta^{\mathfrak p} {f}'}{\tilde\ell}\right]
\nonumber\\[2mm]
\,=\,& -  \iota^3 {\mathbf b}_1\frac{4\pi^2}{{\tilde\ell}^2} f''b_8
+  \frac{\iota^3}{\beta^{1+\mathfrak{p}}}\tilde{\mathcal K}_7\left[\tau,  \beta^{\mathfrak p} f, \frac{\beta^{\mathfrak p} {f}'}{\tilde\ell}\right],
\end{align}
where
\begin{align*}
b_8 =   \int_{\R} \varpi_{1,x} \,  \mathcal Z \,{\mathrm d}x,
\end{align*}
and
\begin{align*}
\frac{\iota^3}{\beta^{1+\mathfrak{p}}}\tilde{\mathcal K}_7\left[\tau,  \beta^{\mathfrak p} f, \frac{\beta^{\mathfrak p} {f}'}{\tilde\ell}\right]
\,=\,&
 \iota^3
\int_\R
\Bigg\{\big(4\mathfrak{p}+2 \big)\frac{\pi}{{\tilde\ell}}\frac{\beta'} {\beta^2}\Big[- f' \phi_{1,x} +\phi_{1,\tau}\Big]
 +\frac{4\pi^2}{{\tilde\ell}^2}\Big[|f'|^2 \phi_{1,xx}    -2f' \phi_{1,x\tau}\Big]
\nonumber\\[2mm]
&\qquad\quad
+ \frac{4\pi}{{\tilde\ell}} \frac{\beta'}{\beta^2} (x+f)\big[- f' \phi_{1,xx}  +\phi_{1,x\tau}\big]
+\Big[\frac{\beta''}{\beta^3} +2\mathfrak{p} \frac{|\beta'|^2}{\beta^4}  \,-\,\frac{\kappa^2}{\beta^2} \Big]
(x+f) \phi_{1,x}
\nonumber\\[2mm]
&\qquad\quad
 +\Big[\mathfrak{p} \frac{\beta''}{\beta^3}  +  \mathfrak{p}  (\mathfrak{p} -1) \frac{|\beta'|^2}{\beta^4} \Big] \phi_{1}+ \frac{|\beta'|^2}{\beta^4}  (x+f)^2 \phi_{1,xx}
\nonumber\\[2mm]
&\qquad\quad-\frac{1}{2}\, \frac{\iota^2W_{tt}(0, \iota z)}{\beta^4} (x+f)^2\phi_{1}
\Bigg\}  {\mathcal Z} \,{\mathrm d}x.
\end{align*}

\smallskip\noindent
$\clubsuit$
Since $\phi_2$ is even of $x$, from the expression of $ B_4$ as in \eqref{new opreator B4},  then we have
\begin{align}\label{B4phi2 intergral Z}
 \int_{\R} B_4(\iota^2\phi_2)\,{\mathcal Z} \,{\mathrm d}x
 \,=\,& \iota^3 \frac{\kappa}{\, \beta\, }\int_{\R} \Big[ \phi_{2,x} +\frac1{\Lambda} (x+f)\phi_2\Big]\mathcal Z \,{\mathrm d}x
\nonumber\\[2mm]
 \,=\,&\iota^3 \frac{\kappa}{\, \beta\, }\int_{\R} \frac1{\Lambda} f\phi_2\mathcal Z \,{\mathrm d}x
\,=\,
  \iota^3 \frac1{\Lambda}\frac{\kappa}{\, \beta^3\, } f\,\tilde {D}_6 (\theta) ,
\end{align}
 where
\begin{align*}
 \frac1{\beta^2}\tilde {D}_6 (\theta) = \int_{\R} \phi_2\mathcal Z \,{\mathrm d}x.
 \end{align*}

\smallskip\noindent
$\clubsuit$
And  by  the expression of the operator $B_5$ given in \eqref{new opreator B5}, we have

\begin{align}   \label{B5phi2 intergral Z}
 \int_{\R}\, B_5 (\iota^4{\phi_2}) \,{\mathcal Z} \,{\mathrm d}x
\,=\,&\iota^4 \int_\R\Bigg\{\frac{4\pi^2}{{\tilde\ell}^2} |f'|^2 \phi_{2, xx}
+
\big(4\mathfrak{p}+2 \big) \frac{\pi}{{\tilde\ell}}\frac{\beta'} {\beta^2}\phi_{2, \tau}  +  4 \frac{\pi}{{\tilde\ell}} \frac{\beta'}{\beta^2} (x+f) \big(-f' \phi_{2,xx} \big)
\nonumber\\[2mm]
&  \qquad\quad
 + \Big[
\frac{\beta''}{\beta^3}\,+\, 2\mathfrak{p}\frac{|\beta'|^2}{\beta^4} \,-\,\frac{\kappa^2}{\beta^2} \Big]  x\phi_{2, x}
  +  \Big[\mathfrak{p} \frac{\beta''}{\beta^3}  +  \mathfrak{p}  (\mathfrak{p} -1) \frac{|\beta'|^2}{\beta^4} \Big] \phi_2
\nonumber\\[2mm]
&   \qquad\quad  + \frac{|\beta'|^2}{\beta^4}(x^2+f^2) \phi_{2, xx}
-\frac{1}{2}\, \iota^2\frac{W_{tt}(0, \iota z)}{\beta^4} (x^2+f^2)\phi_2 \Bigg\} \mathcal Z \,{\mathrm d}x
 \nonumber\\[2mm]
\,=\,&    \frac{\iota^4}{\beta^{1+\mathfrak{p}}}\tilde{\mathcal K}_8\left[\tau,  \beta^{\mathfrak p} f, \frac{\beta^{\mathfrak p} {f}'}{\tilde\ell}\right].
\end{align}
Here we note we used the identities
\begin{align*}
\int_\R  \phi_{2, x}  \mathcal Z \,{\mathrm d}x=0, \quad \int_\R  \phi_{2, x\tau}  \mathcal Z \,{\mathrm d}x=0.
\end{align*}

\smallskip\noindent
$\clubsuit$
By \eqref{definition of mathbbN1}, we have
\begin{align}   \label{N1phi2 intergral Z}
 \int_\R\, \mathbb{N}_1 (\iota^2 \phi_2)      \mathcal Z \,{\mathrm d}x
= & \int_\R\,  \left[ \iota^3 {\bf N}_{11}+\iota^4 {\bf N}_{12}\right] \mathcal Z \,{\mathrm d}x
\nonumber\\[2mm]
= & \int_\R \Bigg\{
  \iota^3 (p_\epsilon-1) (p_\epsilon- 2) U^{p_\epsilon- 3}\Big(\phi_1+ \frac{e} {\beta^{\mathfrak{p}}} \mathcal Z \Big)  \phi_2  \Bigg\} \mathcal Z \,{\mathrm d}x
\nonumber\\[2mm]
&  \quad+   \int_\R \Bigg\{   \left( w_1+ \iota^2\phi_2\right)^{p_\epsilon-1}  -w_1^{p_\epsilon-1}   - \iota^2 (p_\epsilon-1) w_1^{p_\epsilon-2} \phi_2
 \nonumber\\[2mm]
 &\quad +     \iota^2 (p_\epsilon-1)\phi_2 \Big[  (w_1^{p_\epsilon-2} -U^{p_\epsilon-2})  -  (p_\epsilon- 2) \iota U^{p_\epsilon- 3} \Big(\phi_1+ \frac{e} {\beta^{\mathfrak{p}}} \mathcal Z\Big) \Big]
  \Bigg\} \mathcal Z \,{\mathrm d}x
\nonumber\\[2mm]
=  & \frac{\iota^3}{\beta^{1+ \mathfrak{p} }}\tilde{\mathcal K}_9\left[\tau, e,  \beta^{\mathfrak p} f\right].
\end{align}

\smallskip\noindent
$\clubsuit$
By the expression of $a_2(\iota s, \iota z) $ given  in \eqref{a1a2a3} in the region
\begin{align}   \label{second term of e intergral Z}
I_{12}
&  =  - \iota^4 \frac{8\pi^2} {{\tilde\ell}^2} \frac{\kappa}{\, \beta^{1+ \mathfrak{p} }\, }  e''  f \int_\R |\mathcal Z|^2 \,{\mathrm d}x
\nonumber\\[2mm]
&  \quad
+\underbrace{ \int_\R \eta_\delta^\iota(s) \Bigg[ \iota^4 \frac{4\pi^2} {{\tilde\ell}^2} \Big(\iota \frac{3\kappa}{\beta}(x+f) + \iota^2  {\tilde a}_2 \Big) \frac{\kappa}{\, \beta^{1+ \mathfrak{p} }\, }  e'' (x+f)\mathcal Z \Bigg] \mathcal Z \,{\mathrm d}x +O(e^{-c\frac{|\beta \delta|} \iota} )e''}_{:= \frac{\iota^3}{\beta^{1+\mathfrak{p}}}\tilde{\mathcal K}_{10}  \left[\tau,  \beta^{\mathfrak p} f, \frac{\iota^2 e''}{\tilde \ell^2}\right]}
\nonumber\\[2mm]
 &\quad +  \iota^4 \kappa  \frac{ 8\pi^2} {{\tilde\ell}^2} \frac{1}{\, \beta^{1+ \mathfrak{p}}\, }    f'' e   \int_\R     x \mathcal Z' \mathcal Z \,{\mathrm d}x
\nonumber\\[2mm]
&  \quad
+  \underbrace{\int_\R \eta_\delta^\iota(s) \Bigg[  -\iota^4 \kappa  \Big(\iota \frac{3\kappa}{\beta}(x+f) + \iota^2  {\tilde a}_2 \Big)  \frac{ 4\pi^2} {{\tilde\ell}^2} \frac{1}{\, \beta^{1+ \mathfrak{p}}\, }  (x+f) f'' e   \mathcal Z' \Bigg] \mathcal Z \,{\mathrm d}x  + O(e^{-c\frac{|\beta \delta|} \iota} )f''}_{: = \frac{\iota^5}{\beta^{1+\mathfrak{p}}}\tilde{\mathcal K}_{11} \left[\tau,  \beta^{\mathfrak p} f, \frac{\beta^{\mathfrak p}f''}{\tilde \ell^2}\right]}
\nonumber\\[2mm]
&  =
- \iota^4 \frac{8\pi^2} {{\tilde\ell}^2} \frac{\kappa}{\, \beta^{1+ \mathfrak{p}}\, }  e''  f
-\iota^4  \frac{ 4\pi^2} {{\tilde\ell}^2} \frac{\kappa }{\, \beta^{1+ \mathfrak{p}}\, }    f'' e
+  \frac{\iota^3}{\beta^{1+\mathfrak{p}}}\tilde{\mathcal K}_{10} \left[\tau,  \beta^{\mathfrak p} f, \frac{\iota^2 e''}{\tilde \ell^2}\right]
+  \frac{\iota^5}{\beta^{1+\mathfrak{p}}}\tilde{\mathcal K}_{11} \left[\tau,  \beta^{\mathfrak p} f, \frac{\beta^{\mathfrak p}f''}{\tilde \ell^2}\right].
\end{align}

\smallskip\noindent$\clubsuit$
 We note that the expressions of the following terms
$$
 \frac{\iota^4}{\beta^{1+\mathfrak{p}}}  \mathcal F_1\left[\tau, x, \beta^{\mathfrak p} f, \frac{\beta^{\mathfrak p} f'}{\tilde\ell} \right],
 \qquad
\frac{\iota^4}{\beta^{1+\mathfrak{p}}} \hat{\mathcal F}_1\left[\tau, x, \beta^{\mathfrak p} f, \frac{\beta^{\mathfrak p} f''}{\tilde\ell^2} \right],
 \qquad
 \frac{\iota^4}{\beta^{1+\mathfrak{p}}}  \mathcal G_1\left[\tau, x, \beta^{\mathfrak p} f, \frac{\beta^{\mathfrak p} f'}{\tilde\ell},   e, \frac{e'}{\tilde\ell}\right],
$$
have been  given in \eqref{mathcal F1}, \eqref{mathcal G1}, \eqref{hat mathcal F1}.  We denote
\begin{align} \label{mathcal G1 on Z}
& \frac{\iota^4}{\beta^{1+\mathfrak{p}}}  \int_\R \eta_\delta^\iota(s)     \mathcal G_1\left[\tau, x, \beta^{\mathfrak p} f, \frac{\beta^{\mathfrak p} f'}{\tilde\ell},   e, \frac{e'}{\tilde\ell}\right]  \mathcal Z \,{\mathrm d}x
  =  \frac{\iota^3}{\beta^{1+\mathfrak{p}}}\tilde{\mathcal K}_{12}\left[\tau, \beta^{\mathfrak p} f, \frac{\beta^{\mathfrak p} {f}'}{\tilde\ell},  e, \frac{\iota e'}{\tilde\ell}\right],
 \end{align}
\begin{align} \label{mathcal F 1 on Z}
 \frac{\iota^4}{\beta^{1+\mathfrak{p}}}  \int_\R \eta_\delta^\iota(s)   \mathcal F_1\left[\tau, x, \beta^{\mathfrak p} f, \frac{\beta^{\mathfrak p} f'}{\tilde\ell} \right] \mathcal Z \,{\mathrm d}x
 = \frac{\iota^4}{\beta^{1+\mathfrak{p}}}\tilde{\mathcal K}_{14}\left[\tau,  \beta^{\mathfrak p} f, \frac{\beta^{\mathfrak p} f'}{\tilde\ell}\right],
 \end{align}
 and
\begin{align} \label{hat mathcal F 1 on Z}
  & \frac{\iota^4}{\beta^{1+\mathfrak{p}}} \int_\R \eta_\delta^\iota(s)   \hat{\mathcal F}_1\left[\tau, x, \beta^{\mathfrak p} f, \frac{\beta^{\mathfrak p} f''}{\tilde\ell^2} \right]  \mathcal Z \,{\mathrm d}x
  =  \frac{\iota^4}{\beta^{1+\mathfrak{p}}}\tilde{\mathcal K}_{13}\left[\tau,  \beta^{\mathfrak p} f, \frac{\beta^{\mathfrak p} f''}{\tilde\ell^2}\right].
 \end{align}

\smallskip\noindent$\clubsuit$
At last,  by the decompositions \eqref{expression of B6phi1}, \eqref{decomposition of B6phi2} and the estimate \eqref{estimate of mathcal F3-x}, we have
\begin{align} \label{B6phi1 + B6phi2-intergral Z}
  & \int_\R \eta_\delta^\iota(s)   \Big[B_6(\iota \phi_1) + B_6(\iota^2 \phi_2) \Big]  \mathcal Z \,{\mathrm d}x
   \nonumber\\[2mm]
&  =     \int_\R  \eta_\delta^\iota(s)   \Bigg\{ \frac{\iota^4}{\beta^{1+\mathfrak{p}}}  \mathcal F_2\left[\theta, x, f, \frac{\beta^{\mathfrak p} f'}{\tilde\ell}\right]+ \frac{\iota^4}{\beta^{1+\mathfrak{p}}} \hat{\mathcal F}_2\left[\tau, x, \beta^{\mathfrak p} f, \frac{\beta^{\mathfrak p} f''}{\tilde\ell^2} \right]  \Bigg\} \mathcal Z \,{\mathrm d}x
 \nonumber\\[2mm]
&  \quad+     \int_\R \eta_\delta^\iota(s)    \Bigg\{ -\iota^5 \kappa(\iota z) a_2(\iota s, \iota z) \frac{ 4\pi^2} {{\tilde\ell}^2}\frac{1}{\beta} (x+f)f'' \phi_{2, x}+   \frac{\iota^5}{\beta^3} \mathcal F_3\left[\tau, x, \beta^{\mathfrak p} f, \frac{\beta^{\mathfrak p} f'}{\tilde\ell}\right] \Bigg\} \mathcal Z \,{\mathrm d}x
 \nonumber\\[2mm]
&  = \frac{\iota^4}{\beta^{1+\mathfrak p}}\tilde{\mathcal K}_{15}\left[\tau,  \beta^{\mathfrak p} f, \frac{\beta^{\mathfrak p} f''}{\tilde\ell^2}\right]
+ \frac{\iota^4}{\beta^{1+\mathfrak p}}\tilde{\mathcal K}_{16}\left[\tau,  \beta^{\mathfrak p} f, \frac{\beta^{\mathfrak p} f'}{\tilde\ell}\right],
 \end{align}
 where
\begin{align*}
  \frac{\iota^4}{\beta^{1+\mathfrak p}}\tilde{\mathcal K}_{15}\left[\tau,  \beta^{\mathfrak p} f, \frac{\beta^{\mathfrak p} f''}{\tilde\ell^2}\right]
\,=\,& \int_\R  \eta_\delta^\iota(s)   \Bigg\{  -\iota^5 \kappa(\iota z) a_2(\iota s, \iota z) \frac{ 4\pi^2} {{\tilde\ell}^2}\frac{1}{\beta} (x+f)f'' \phi_{2, x}
\\\nonumber
&\qquad\qquad\quad
+ \frac{\iota^4}{\beta^{1+\mathfrak{p}}} \hat{\mathcal F}_2\left[\tau, x, \beta^{\mathfrak p} f, \frac{\beta^{\mathfrak p} f''}{\tilde\ell^2} \right]  \Bigg\} \mathcal Z \,{\mathrm d}x,
 \end{align*}
 and
\begin{align*}
 \frac{\iota^4}{\beta^{1+\mathfrak p}}\tilde{\mathcal K}_{16}\left[\tau,  \beta^{\mathfrak p} f, \frac{\beta^{\mathfrak p} f'}{\tilde\ell}\right]
 =   \int_\R  \eta_\delta^\iota(s)   \Bigg\{   \frac{\iota^4}{\beta^{1+\mathfrak{p}}}  \mathcal F_2\left[\theta, x, f, \frac{\beta^{\mathfrak p} f'}{\tilde\ell}\right] + \frac{\iota^5}{\beta^3} \mathcal F_3\left[\tau, x, \beta^{\mathfrak p} f, \frac{\beta^{\mathfrak p} f'}{\tilde\ell}\right]  \Bigg\} \mathcal Z \,{\mathrm d}x.
 \end{align*}

\medskip
Finally, we combine \eqref{projected eqs for enew},
 \eqref{S42 on Z},  \eqref{S6 on Z},  \eqref{S8 on Z},  \eqref{N0epsilon on Z},   \eqref{B4emathcal Z2},  \eqref{nonlinear term N02emathcal Z},  \eqref{nonlinear term N03 emathcal Z},  \eqref{main term mathccal Z}, \eqref{phitautau mathcal Z}, \eqref{phi2tautau mathcal Z},   \eqref{B5phi1 intergral Z},   \eqref{B4phi2 intergral Z}, \eqref{B5phi2 intergral Z}, \eqref{N1phi2 intergral Z},  \eqref{second term of e intergral Z},  \eqref{mathcal G1 on Z}, \eqref{mathcal F 1 on Z}, \eqref{hat mathcal F 1 on Z},   \eqref{B6phi1 + B6phi2-intergral Z},   and then get
\begin{align}   \label{Sw2 integral 2}
 \frac{1}{\beta^{\mathfrak{p}+2}}\int_\R\eta_\delta^\iota(s){\mathcal E}\,{\mathcal Z} \,{\mathrm d}x
 & = \frac 1{\beta^{\mathfrak{p}}} \Bigg[\iota^3\frac{4\pi^2}{{\tilde\ell}^2}  e''
  +\iota\lambda_0 e \Bigg]  -\frac{1}{2}\, \iota^2 \frac{ \iota^2W_{tt}(0, \iota z)}{\beta^4}|f|^2  b_2
  \nonumber\\[2mm]
 & \quad
+ \iota^2\Big(\frac{4\pi^2}{{\tilde\ell}^2}|f'|^2
-\frac{4\pi}{{\tilde\ell}}\frac{\beta'}{\beta^2}ff'
+\frac{|\beta'|^2}{\beta^4} |f|^2
\Big) b_3
+ \iota^2  \frac1{\Lambda}\frac{\kappa}{ \beta^{1+\mathfrak p} } fe b_0
+ \iota^2 \frac{(p_\epsilon-1) (p_\epsilon- 2)}{2\beta^{2\mathfrak{p}}}e^2b_1
\nonumber\\[2mm]
 & \quad
 + \frac{\iota^3}{\beta^{1+\mathfrak p} }\tilde{\mathcal K}_2\left[\tau, {\beta f}, \frac{\beta{ f}'}{\tilde\ell} \right]
 +\iota^3\frac{8\pi^2}{{\tilde\ell}^2}\frac{\kappa}{\beta}   f''b_4
 + \iota^3  \frac{1}{\beta^{1+\mathfrak{p}} } \tilde{\mathcal K}_3\left[\tau,  \beta^{\mathfrak p} f, \frac{\beta^{\mathfrak p} { f}'}{\tilde\ell}\right] e
+ \iota^2 \frac{1}{\Lambda^2} f^2  {\mathbf b}_1^2 b_5
\nonumber\\[2mm]
 & \quad
 + \iota^2  (p_\epsilon-1)(p_\epsilon- 2)  \Big[  \frac 1{2\Lambda^2} {\mathbf b}_1^2f^2b_6
+  \frac{e} {\beta^{\mathfrak{p}}}  \frac 1{\Lambda} {\mathbf b}_1f  b_7\Big]
+ \iota^3\frac{1} {\beta^{1+ \mathfrak{p}}} \tilde{\mathcal K}_5\left[\tau,  {\beta f},  e\right]
 \nonumber\\[2mm]
 & \quad
  -  \iota^3 {\mathbf b}_1\frac{4\pi^2}{{\tilde\ell}^2} f''b_8
  +  \iota^3\frac{1}{\beta^{1+\mathfrak p} }\tilde{\mathcal K}_7\left[\tau,  \beta^{\mathfrak p} f, \frac{\beta^{\mathfrak p} {f}'}{\tilde\ell}\right]
  +\iota^3 \frac1{\Lambda}  \frac{\kappa}{\, \beta^3\, }  \tilde {D}_4 (\theta) f
  +  \iota^4\frac{1}{\beta^{1+\mathfrak p} }\tilde{\mathcal K}_8\left[\tau,  \beta^{\mathfrak p} f, \frac{\beta^{\mathfrak p} {f}'}{\tilde\ell}\right]
 \nonumber\\[2mm]
 & \quad
 +  \iota^3 \frac{1}{\beta^{1+ \mathfrak{p} }}\tilde{\mathcal K}_9\left[\tau, e, { \beta^{\mathfrak p} f}\right]
  - \iota^4 \frac{8\pi^2} {{\tilde\ell}^2} \frac{\kappa}{\, \beta^{1+ \mathfrak{p}}\, }  e''  f
   +  \iota^3\frac{1}{\beta^{1+\mathfrak p} }\tilde{\mathcal K}_{10} \left[\tau,  \beta^{\mathfrak p} f, \frac{\iota^2 e''}{\tilde \ell^2}\right]
   - \iota^4 \frac{ 4\pi^2} {{\tilde\ell}^2} \kappa  \frac{1}{\, \beta^{1+ \mathfrak{p}}\, }    f'' e
 \nonumber\\[2mm]
 & \quad
 +  \frac{\iota^5}{\beta^{1+\mathfrak p} }\tilde{\mathcal K}_{11} \left[\tau,  \beta^{\mathfrak p} f, \frac{f''}{\tilde \ell^2}\right]
 +  \frac{\iota^3}{\beta^{1+\mathfrak p} }\tilde{\mathcal K}_{12}\left[\tau, \beta^{\mathfrak p} f, \frac{\beta^{\mathfrak p} {f}'}{\tilde\ell},  e, \frac{e'}{\tilde\ell}\right]
 +\frac{\iota^4}{\beta^{1+\mathfrak p} }\tilde{\mathcal K}_{13}\left[\tau,  \beta^{\mathfrak p} f, \frac{\beta^{\mathfrak p} f''}{\tilde\ell^2}\right]
 \nonumber\\[2mm]
 & \quad
 +\frac{\iota^4}{\beta^{1+\mathfrak p} }\tilde{\mathcal K}_{14}\left[\tau, \beta^{\mathfrak p} f, \frac{\beta^{\mathfrak p} {f}'}{\tilde\ell}\right]
 +\frac{\iota^4}{\beta^{1+\mathfrak p} }\tilde{\mathcal K}_{15}\left[\tau,  \beta^{\mathfrak p} f, \frac{\beta^{\mathfrak p} f''}{\tilde\ell^2}\right]
 + \frac{\iota^4}{\beta^{1+\mathfrak p} }\tilde{\mathcal K}_{16}\left[\tau,  \beta^{\mathfrak p} f, \frac{\beta^{\mathfrak p} f'}{\tilde\ell}\right],
  \nonumber\\[2mm]
 &  =  \frac 1{\beta^{\mathfrak{p}}} \Bigg[\iota^3\frac{4\pi^2}{{\tilde\ell}^2}  e''  +  \iota\lambda_0 e \Bigg]
  -\frac{1}{2}\, \iota^2 \frac{ \iota^2W_{tt}(0, \iota z)}{\beta^4}|f|^2  b_2
 \nonumber\\[2mm]
 & \quad
 + \iota^2\Big[\frac{4\pi^2}{{\tilde\ell}^2}|f'|^2
-\frac{4\pi}{{\tilde\ell}}\frac{\beta'}{\beta^2}ff'
+\frac{|\beta'|^2}{\beta^4} |f|^2
\Big] b_3
+ \iota^2\frac1{\Lambda} \frac{\kappa}{\, \beta^2\, } fe b_0
+  \iota^2\frac{(p_\epsilon-1) (p_\epsilon- 2)}{2\beta^{2\mathfrak{p}}}e^2b_1
\nonumber\\[2mm]
 & \quad
 + \iota^2  \frac{1}{\Lambda^2}f^2  {\mathbf b}_1^2 b_5
 +  \iota^2 (p_\epsilon-1)(p_\epsilon- 2)  \Big[  \frac 1{2\Lambda^2} {\mathbf b}_1^2f^2b_6
+  \frac{e} {\beta^{\mathfrak{p}}}  \frac 1{\Lambda} {\mathbf b}_1f  b_7\Big]
    - \iota^4 \frac{8\pi^2} {{\tilde\ell}^2} \frac{\kappa}{\, \beta^{1+ \mathfrak{p}}\, }  e''  f
     \nonumber\\[2mm]
 & \quad
 +  \frac{\iota^3}{\beta^{1+\mathfrak p}}\tilde{\mathcal K}^*\left[\tau, \beta^{\mathfrak p} f, \frac{\beta^{\mathfrak p} {f}'}{\tilde\ell},  e, \frac{e'}{\tilde\ell}\right]
  +  \frac{\iota^4}{\beta^{1+\mathfrak p}}\tilde{\mathcal K}^{***}  \left[\tau, \beta^{\mathfrak p} f, \frac{\beta^{\mathfrak p} {f}'}{\tilde\ell}, e, \frac{\iota e'}{\tilde\ell}\right]
  \nonumber\\[2mm]
 & \quad  +   \frac{\iota^3}{\beta^{1+\mathfrak p}}\tilde{\mathcal K}^{**}   \left[\tau,  \beta^{\mathfrak p} f, \frac{\beta^{\mathfrak p} f''}{\tilde\ell^2}, \frac{\iota^2 e''}{\tilde \ell^2}\right],
 \end{align}
where
\begin{align}
 \frac{\iota^3}{\beta^{1+\mathfrak p}} \tilde{\mathcal K}^* \left[\tau, \beta^{\mathfrak p} f, \frac{\beta^{\mathfrak p} {f}'}{\tilde\ell}, e, \frac{\iota e'}{\tilde\ell}\right]
\,=\, & \frac{\iota^3}{\beta}\tilde{\mathcal K}_2\left[\tau, {\beta f}, \frac{\beta{ f}'}{\tilde\ell} \right]   + \iota^3
\frac{1}{\beta^{1+\mathfrak{p}} } \tilde{\mathcal K}_3\left[\tau,  \beta^{\mathfrak p} f, \frac{\beta^{\mathfrak p} { f}'}{\tilde\ell}\right] e  + \frac{\iota^3} {\beta^{1+ \mathfrak{p}}} \tilde{\mathcal K}_5\left[\tau,  \beta^{\mathfrak p} f,  e\right]
 \nonumber\\[2mm]
 &  + \frac1{\Lambda} \iota^3 \frac{\kappa}{\, \beta^3\, }  \tilde {D}_4 (\theta) f
 + \frac{\iota^3}{\beta^{1+\mathfrak{p}}}\tilde{\mathcal K}_7\left[\tau,  \beta^{\mathfrak p} f, \frac{\beta^{\mathfrak p} {f}'}{\tilde\ell}\right]
  +   \frac{\iota^3}{\beta^{1+ \mathfrak{p} }}\tilde{\mathcal K}_9\left[\tau, e, { f}\right] ,
 \end{align}
\begin{align}
\iota^4 \tilde{\mathcal K}^{***}  \left[\tau, \beta^{\mathfrak p} f, \frac{\beta^{\mathfrak p} {f}'}{\tilde\ell}, e, \frac{\iota e'}{\tilde\ell}\right]
 = & \frac{\iota^4}{\beta^{1+\mathfrak{p}}}\tilde{\mathcal K}_8\left[\tau,  \beta^{\mathfrak p} f, \frac{\beta^{\mathfrak p} {f}'}{\tilde\ell}\right]
 +\frac{\iota^4}{\beta^{1+\mathfrak{p}}}\tilde{\mathcal K}_{14}\left[\tau, \beta^{\mathfrak p} f, \frac{\beta^{\mathfrak p} {f}'}{\tilde\ell}\right]
  \nonumber\\[2mm]
 &  +\frac{\iota^4}{\beta^{1+\mathfrak{p}}}\tilde{\mathcal K}_{16}\left[\tau,  \beta^{\mathfrak p} f, \frac{\beta^{\mathfrak p} f'}{\tilde\ell}\right],
\end{align}
and
\begin{align}
\iota^3 \tilde{\mathcal K}^{**}   \left[\tau,  \beta^{\mathfrak p} f, \frac{\beta^{\mathfrak p} f''}{\tilde\ell^2}, \frac{\iota^2 e''}{\tilde \ell^2}\right]
\,=\,&     \frac{\iota^3}{\beta^{1+\mathfrak{p}}}\tilde{\mathcal K}_{10} \left[\tau,  \beta^{\mathfrak p} f, \frac{\iota^2 e''}{\tilde \ell^2}\right]
-\iota^3\frac{8\pi^2}{{\tilde\ell}^2}\frac{\kappa}{\beta}   f''b_4  -  \iota^3 {\mathbf b}_1\frac{4\pi^2}{{\tilde\ell}^2} f''b_8
   \nonumber\\[2mm]
 &  - \iota^4 \frac{ 4\pi^2} {{\tilde\ell}^2} \kappa  \frac{1}{\, \beta^{1+ \mathfrak{p}}\, }    f'' e   +\frac{\iota^5}{\beta^{1+\mathfrak{p}}}\tilde{\mathcal K}_{11} \left[\tau,  \beta^{\mathfrak p} f, \frac{\beta^{\mathfrak p}f''}{\tilde \ell^2}\right]
  \nonumber\\[2mm]
 &
 +\frac{\iota^4}{\beta^{1+\mathfrak{p}}}\tilde{\mathcal K}_{13}\left[\tau,  \beta^{\mathfrak p} f, \frac{\beta^{\mathfrak p} f''}{\tilde\ell^2}\right]
  +\frac{\iota^4}{\beta^{1+\mathfrak{p}}}\tilde{\mathcal K}_{15}\left[\tau,  \beta^{\mathfrak p} f, \frac{\beta^{\mathfrak p} f''}{\tilde\ell^2}\right].
   \end{align}
Moreover, we have
\begin{align}
\left\|\frac{\iota^3}{\beta^{1+\mathfrak{p}}} \tilde{\mathcal K}^* \left[\tau, \beta^{\mathfrak p} f, \frac{\beta^{\mathfrak p} {f}'}{\tilde\ell}, e, \frac{\iota e'}{\tilde\ell}\right]   \right\|_{L^2(0, 2\pi)}  \le C \iota^3  ( \|f\|_* + \|e\|_{**} ),
\end{align}
\begin{align}
\left\|\frac{\iota^4}{\beta^{1+\mathfrak{p}}}\tilde{\mathcal K}^{***}  \left[\tau, \beta^{\mathfrak p} f, \frac{\beta^{\mathfrak p} {f}'}{\tilde\ell}, e, \frac{\iota e'}{\tilde\ell}\right]   \right\|_{L^2(0, 2\pi)}  \le C  (  \iota^4+ \iota^4\|f\|_* + \iota^4 \|e\|_{**} ),
\end{align}
\begin{align}
\left\| \frac{\iota^3}{\beta^{1+\mathfrak{p}}} \tilde{\mathcal K}^{**}  \left[\tau,  \beta^{\mathfrak p} f, \frac{\beta^{\mathfrak p} f''}{\tilde\ell^2}, \frac{\iota^2 e''}{\tilde \ell^2}\right]      \right\|_{L^2(0, 2\pi)}  \le C \iota^3  ( \|f\|_* + \|e\|_{**} ),
\end{align}
\begin{align}
& \left\|\frac{\iota^3}{\beta^{1+\mathfrak{p}}}\left( \tilde{\mathcal K}^* \left[\tau, \beta^{\mathfrak p} f_1, \frac{\beta^{\mathfrak p} {f_1'}}{\tilde\ell}, e_1, \frac{\e e_1'}{\tilde\ell}\right]  -  \tilde{\mathcal K}^* \left[\tau, \beta^{\mathfrak p} f_2, \frac{\beta^{\mathfrak p} {f_2'}}{\tilde\ell}, e_2, \frac{\iota e_2'}{\tilde\ell}\right]   \right) \right\|_{L^2(0, 2\pi)}
 \nonumber\\[2mm]
 &
   \le C \iota^3 (\|f_1- f_2\|_* + \|e_1- e_2\|_{**} ),
\end{align}
and
\begin{align}
&\left\|\frac{\iota^3}{\beta^{1+\mathfrak{p}}} \left(\tilde{\mathcal K}^{**}  \left[\tau,  \beta^{\mathfrak p} f_1, \frac{\beta^{\mathfrak p} f_1''}{\tilde\ell^2}, \frac{\iota^2 e_1''}{\tilde \ell^2}\right]   -  \tilde{\mathcal K}^{**}  \left[\tau,  \beta^{\mathfrak p} f_2, \frac{\beta^{\mathfrak p} f_2''}{\tilde\ell^2}, \frac{\iota^2 e_2''}{\tilde \ell^2}\right]  \right)    \right\|_{L^2(0, 2\pi)}
  \nonumber\\[2mm]
 & \le C \iota^3 (\|f_1- f_2\|_* + \|e_1- e_2\|_{**} ).
\end{align}

\medskip
\medskip

\medskip
\subsection{ Projections of terms involving the perturbation term}\label{section9.2}\

By  \eqref{errorrelationinterior} and \eqref{expression of mathcal L}, we transfer problem \eqref{project 6.13}-\eqref{projectedproblem3} with $c=0, d=0$ to the following
\begin{align}
 L_2(\phi^*) \ +\ \chi(\iota|s|) \big[\, B_4(\phi^*)+ B_5(\phi^*)+ B_6(\phi^*)\,\big] = - \chi(\iota|s|) \left[S(w_2) + {\mathcal N}(\phi^*)\right],
\end{align}

\begin{equation}
{\phi^* }(x, 0) \,=\, {\phi^* }(x, 2\pi), \quad \phi^*_\tau(x, 0) \,=\, \phi^* _\tau(x, 2\pi),
\quad
-\infty<x<+\infty,
\end{equation}

\begin{align}
 \int_{-\infty}^{+\infty} \phi^*(x, \tau){\mathcal Z} \,{\mathrm d}x \,=\, \int_{-\infty}^{+\infty} \phi^*(x, \tau) U'(x)\,{\mathrm d}x \,=\, 0, \quad 0<\tau<2\pi,
\end{align}
 where
\begin{align}\label{Nmathcal}
{\mathcal N}(\phi^*)
\,=\,
\frac{1}{\beta^{2+\mathfrak{p}} }\left[\, \mathbf{N}\big(\eta_{3\delta}^{\iota}{\tilde\phi}+{\tilde\psi}\big)
-(p_\epsilon-1)\, \mathbf{U}^{p_\epsilon-2}  \, {\tilde\psi}(\tilde\phi) \,\right],
\qquad
\tilde\phi(y) = \beta^{\mathfrak{p}}(\iota z)\chi(\iota|s|) \phi^*(x, \tau),
\end{align}
the operator $ L_2(\phi^*) $ is given in \eqref{definition of L2operator},
 and $\chi(r)$ is smooth cut-off function that equals 1 for $r\leq 10\delta$, vanishes identically for $r>20\delta$.
Recall the norm $\|\cdot\|_{H^2_*(\mathfrak S)} $ given in \eqref{norm h2 mathfrak S}.

\medskip
Let
\begin{align} \label{Pi1}
\Pi_1(\tau, f, e) = \int_{{\mathbb R}} L_2(\phi^*) U'(x)\,{\mathrm d}x,
\end{align}
\begin{align} \label{Pi2}
\Pi_2(\tau, f, e) = \int_{{\mathbb R}} \chi(\iota|s|) \Big[B_4(\phi^*)+ B_5(\phi^*)+ B_6(\phi^*)\Big] U'(x)\,{\mathrm d}x,
\end{align}
and
\begin{align} \label{Pi3}
\Pi_3(\tau, f, e) = -\int_{{\mathbb R}}\chi(\iota|s|) {\mathcal N}(\phi^*) U'(x)\,{\mathrm d}x.
\end{align}

 We recall that
\begin{align*}
 \iota\phi_1(x, \tau)
 \,=\, \iota{\phi_{11}}(x, \tau) + \iota{\phi_{12}}(x, \tau)
 = \iota{\mathbf b}_1(\iota z)\varpi_1(x)
\ +\
\iota\frac1{\Lambda} {\mathbf b}_1(\iota z)f(\tau) \varpi_2(x),
\end{align*}
the expression of $\iota^2\phi_2$ given in \eqref{definition of e2phi2}-\eqref{phi21-phi210},
 and there hold
\begin{align*}
 |\iota\phi_1| \leq \iota \frac{C}{\beta} e^{-c|x|}, \quad |\iota^2\phi_2| \leq \iota^2 \frac{C}{\beta^2} e^{-c|x|}.
\end{align*}
 For $\Pi_1(\tau, f, e) $, we have
\begin{align*}
 \Pi_1 (\tau, f, e) & = \int_{{\mathbb R}} L_2(\phi^*) U'(x)\,{\mathrm d}x
\nonumber\\[2mm]
 & =
\int_{{\mathbb R}}\left[ \frac{\iota^2 4\pi^2}{{\tilde\ell}^2} \phi^*_{\tau\tau} + \phi^*_{xx} - \phi^*+ (p_\epsilon-1)w_2^{p_\epsilon-2} \phi^* \right] U'(x)\,{\mathrm d}x
\nonumber\\[2mm]
 & =\int_{{\mathbb R}}\left[\phi^*_{xx} - \phi^*+ (p_\epsilon-1)U^{p_\epsilon-2} \phi^* + (p_\epsilon-1)\left(w_2^{p_\epsilon-2}  - U^{p_\epsilon-2}  \right) \phi^* \right] U'(x)\,{\mathrm d}x
\nonumber\\[2mm]
 & =\int_{{\mathbb R}} (p_\epsilon-1)\left(w_2^{p_\epsilon-2}  - U^{p_\epsilon-2} \right) \phi^* U'(x)\,{\mathrm d}x
\nonumber\\[2mm]
 & =  (p_\epsilon-1) \int_{{\mathbb R}} \left[\, \Big(   U +  \frac1{\beta^{\mathfrak{p}}}\iota e \mathcal{Z}+\iota \phi_1+ \iota^2\phi_2 \Big)^{p_\epsilon-2}
 - U^{p_\epsilon-2}\right] \phi^* U'(x)\,{\mathrm d}x
 \nonumber\\[2mm]
 & =   (p_\epsilon-1)  (p_\epsilon- 2)   \int_{{\mathbb R}}\, \Big[ U +  t \Big(\frac1{\beta^{\mathfrak{p}}}\iota e \mathcal{Z}+\iota \phi_1+ \iota^2\phi_2 \Big) \Big] ^{3-\e}  \Big(\frac1{\beta^{\mathfrak{p}}}\iota e \mathcal{Z}+\iota \phi_1+ \iota^2\phi_2 \Big) \phi^* U'(x)\,{\mathrm d}x
\nonumber\\[2mm]
 & \le C\frac{\iota} {\beta^{2\mathfrak{p}}}  \left(|e| + |f| +1\right) \|\beta^{\mathfrak{p}} \phi^*(\cdot, \tau)\|_{L^2({\mathbb R})},
\end{align*}
  for some $t\in (0, 1)$.

 Next, we consider $\Pi_2(\tau, f, e)$.
 By the expression of operators $B_4, B_5, B_6$ given in \eqref{new opreator B4}, \eqref{new opreator B5}, \eqref{B2w-2222new},
 we have
\begin{align*}
\Pi_2(\tau, f, e)
\,=\, &\int_{{\mathbb R}} \chi(\iota|s|) \big[B_4(\phi^*)+ B_5(\phi^*)+ B_6(\phi^*)\big] U'(x)\,{\mathrm d}x
\nonumber\\[2mm]
\,=\,&\iota \frac{\kappa}{\, \beta\, }\int_{\mathbb R} \Big[ \phi^*_x +\frac1{\Lambda} (x+f)\phi^*\Big] U'(x)\,{\mathrm d}x
\nonumber\\[2mm]
&
+\iota^2\frac{4\pi^2}{{\tilde\ell}^2}\int_{\mathbb R}\Big[|f'|^2 \phi^*_{xx} -f''\phi^*_x -2f' \phi^*_{x\tau}\Big]U'(x)\,{\mathrm d}x
\ +\ \iota^2\big(4\mathfrak{p} +2 \big) \frac{\pi}{{\tilde\ell}}\frac{\beta'} {\beta^2}\int_{\mathbb R}\Big[- f' \phi^*_x +\phi^*_\tau\Big]U'(x)\,{\mathrm d}x
\nonumber\\[2mm]
&+ \iota^2 \frac{4\pi}{{\tilde\ell}} \frac{\beta'}{\beta^2} \int_{\mathbb R} (x+f)\big[- f' \phi^*_{xx} +\phi^*_{x\tau}\big]U'(x)\,{\mathrm d}x
\ +\
\iota^2\Big[\frac{\beta''}{\beta^3} +2 \mathfrak{p}\frac{|\beta'|^2}{\beta^4}\,-\,\frac{\kappa^2}{\beta^2} \Big]\int_{\mathbb R} (x+f) \phi^*_{x} U'(x)\,{\mathrm d}x
\nonumber\\[2mm]
&
+ \iota^2 \frac{|\beta'|^2}{\beta^4}\int_{\mathbb R}(x+f)^2 \phi^*_{xx} U'(x)\,{\mathrm d}x
\ -\
\frac{1}{2}\, \iota^2\frac{\iota^2W_{tt}(0, \iota z)}{\beta^4} \int_{\mathbb R}(x+f)^2 \phi^* U'(x)\,{\mathrm d}x
\nonumber\\[2mm]
&- \frac{\iota^3}{6} \frac{\, \, \iota^2 W_{ttt}(0, \iota z)}{\beta^5} \int_{\mathbb R} (x+f)^3 \phi^*U'(x)\,{\mathrm d}x
\ -\
\frac{\iota^4}{24} \frac{\, \, \iota^2 W_{tttt}\big(\hbar(\iota z)\iota s, \, \iota z\big)}{\beta^6} (x+f)^4\phi^* U'(x)\,{\mathrm d}x
\nonumber\\[2mm]
&
 +\iota^3 \kappa^3(\iota z) a_3(\iota s, \iota z) \frac{1}{\, \beta^3\, }\int_{\mathbb R} (x+f)^2\phi^*_{x} U'(x)\,{\mathrm d}x
\nonumber\\[2mm]
&\,+\,
\iota^3 \kappa(\iota z) a_2(\iota s, \iota z)
\Bigg\{
 \left[ \mathfrak{p}\frac{\beta''} {\beta^4} +( \mathfrak{p}^2- \mathfrak{p})  \frac{|\beta'|^2}{\beta^5} \right] \int_{\mathbb R}(x+f) \phi^*U'(x)\,{\mathrm d}x
\nonumber\\[2mm]
&\qquad \qquad \qquad \qquad \qquad  + \left[ 2\mathfrak{p} \frac{|\beta'|^2}{\beta^5} +  \frac{\beta''}{\beta^4}  \right]\int_{\mathbb R}(x+f)^2 \phi^*_xU'(x)\,{\mathrm d}x
\nonumber\\[2mm]
& \qquad\qquad\qquad\qquad\qquad
+ \big[4\mathfrak{p} +2 \big] \frac{\pi}{ {\tilde\ell}}\frac{\beta'} {\beta^3}
\int_{\mathbb R} (x+f)\big[-f'\phi^*_x +\phi^*_\tau\big]U'(x)\,{\mathrm d}x
\nonumber\\[2mm]
&\qquad\qquad\qquad\qquad\qquad
+\frac{ 4\pi}{ {\tilde\ell}}\frac{\beta'} {\beta^3}
\int_{\mathbb R} (x+f)^2[- f' \phi^*_{xx} + \phi^*_{x\tau}]U'(x)\,{\mathrm d}x
\nonumber\\[2mm]
& \qquad\qquad\qquad\qquad\qquad
+ \frac{|\beta'|^2}{\beta^5} \int_{\mathbb R}(x+f)^3 \phi^*_{xx}U'(x)\,{\mathrm d}x
+ \frac{\beta''}{\beta^4}\int_{\mathbb R} (x+f)^2 \phi^*_{x}U'(x)\,{\mathrm d}x
\nonumber\\[2mm]
&\qquad\qquad\qquad\qquad\qquad
+\frac{ 4\pi^2} {{\tilde\ell}^2}\frac{1}{\beta}
\int_{\mathbb R}(x+f)\big[ |f'|^2 \phi^*_{xx} -f''\phi^*_x +\phi^*_{\tau\tau} -2f' \phi^*_{x\tau}\big]U'(x)\,{\mathrm d}x
\Bigg\}
\nonumber\\[2mm]
&\,+\,
\iota^3 \kappa'(\iota z)a_1(\iota s, \iota z)
\Bigg\{
 \mathfrak{p}\frac{\beta'}{\beta^4} \int_{\mathbb R}(x+f)\phi^*U'(x)\,{\mathrm d}x
+\frac{\beta'}{\beta^4} \int_{\mathbb R} (x+f)^2 \phi^*_{x}U'(x)\,{\mathrm d}x
\nonumber\\[2mm]
&\qquad\qquad\qquad\qquad\qquad+\frac{ 2\pi}{ {\tilde\ell}}\frac{1}{ \beta^2}\int_{\mathbb R} (x+f)\big[-f'\phi^*_x +\phi^*_\tau \big]U'(x)\,{\mathrm d}x
\Bigg\}.
\end{align*}
We select the less regular terms from $\Pi_{2}(\tau, f, e) $
 as $\Pi_{2*}(\tau, f, e) $, where
\begin{align}\label{Pi2star}
\Pi_{2*}(\tau, f, e)
=& -\iota^2f''\frac{4\pi^2}{{\tilde\ell}^2}\int_{\mathbb R}\, \phi^*_x U'(x)\,{\mathrm d}x
\nonumber\\[2mm]
&+ \iota^3 \kappa(\iota z) a_2(\iota s, \iota z) \frac{4\pi^2}{{\tilde\ell}^2}\frac{1}{\beta} \int_{\mathbb R} (x+f)\Big[ -f''\phi^*_x +\phi^*_{\tau\tau} \Big] U'(x)\,{\mathrm d}x.
\end{align}
Moreover, we have
\begin{align*}
& \Big|\Pi_{2}(\tau, f, e) - \Pi_{2*}(\tau, f, e)\Big|
\nonumber\\[2mm]
& \le \iota \frac1{\beta} \left\{ 1+ |\beta^{\mathfrak p} f| + \Big|\frac{\beta  f'} {\tilde \ell}\Big| \right\}
 \times
\Bigg[\|\phi^*(\cdot, \tau)\|_{L^2({\mathbb R})} + \|\phi^*_x(\cdot, \tau)\|_{L^2({\mathbb R})}
+ \|\phi^*_{xx}(\cdot, \tau)\|_{L^2({\mathbb R})}
\nonumber\\[2mm]
& \qquad\qquad \qquad\qquad \qquad\qquad \qquad
+ \frac{\iota}{\tilde\ell}\|\phi^*_\tau(\cdot, \tau)\|_{L^2({\mathbb R})}+ \frac{\iota}{\tilde\ell}\|\phi^*_{x\tau}(\cdot, \tau)\|_{L^2({\mathbb R})} \Bigg],
\end{align*}
and
\begin{align} \label{minus of Pi2star}
&\Pi_{2*} (\theta, f_1, e_1) - \Pi_{2*} (\theta, f_2, e_2)
\nonumber\\[2mm]
& =
 - \iota^2 f_1''\frac{4\pi^2}{{\tilde\ell}^2}\int_{\mathbb R}\,\chi(\iota|s|) \phi^*_x(f_1, e_1) U'(x)\,{\mathrm d}x
 + \iota^2 f_2''\frac{4\pi^2}{{\tilde\ell}^2}\int_{\mathbb R}\, \chi(\iota|s|) \phi^*_x(f_2, e_2) U'(x)\,{\mathrm d}x
\nonumber\\[2mm]
&\quad + \iota^3 \kappa(\iota z) a_2(\iota s, \iota z) \frac{4\pi^2}{{\tilde\ell}^2}\frac{1}{\beta} \int_{\mathbb R}\, \chi(\iota|s|) (x+f_1)\Big[ -f_1''\phi^*_x(f_1, e_1) +\phi^*_{\tau\tau}(f_1, e_1) \Big] U'(x)\,{\mathrm d}x
\nonumber\\[2mm]
&\quad  -\iota^3 \kappa(\iota z) a_2(\iota s, \iota z) \frac{4\pi^2}{{\tilde\ell}^2}\frac{1}{\beta} \int_{\mathbb R}\, \chi(\iota|s|) (x+f_2)\Big[ -f_2''\phi^*_x(f_2, e_2)+\phi^*_{\tau\tau}(f_2, e_2) \Big] U'(x)\,{\mathrm d}x,
\end{align}
which implies
\begin{align}\label{minus of Pi2star1}
&\Big|\Pi_{2*} (\theta, f_1, e_1) - \Pi_{2*} (\theta, f_2, e_2)\Big|
\nonumber\\[2mm]
& \le \iota^2\big|f_1''-f_2''\big|\frac{4\pi^2}{{\tilde\ell}^2}
\int_{\mathbb R}\,  \big|\phi^*_x(f_1, e_1)\big| |U'(x)|\,{\mathrm d}x
+ \iota^2|f_2''| \frac{4\pi^2}{{\tilde\ell}^2}
\int_{\mathbb R}\, \Big|\phi^*_x(f_1, e_1) - \phi^*_x(f_2, e_2)\Big| \big|U'(x)\big|\,{\mathrm d}x
\nonumber\\[2mm]
&\quad + \iota^3 \frac{4\pi^2}{{\tilde\ell}^2}\frac{1}{\beta}
\big|f_1''-f_2''\big|
\int_{\mathbb R}\, x\big|\phi^*_x(f_1, e_1)\big|U'(x)\,{\mathrm d}x
+ \iota^3 \frac{4\pi^2}{{\tilde\ell}^2}\frac{1}{\beta}
\big|f_2''\big| \int_{\mathbb R}\, x\big|\phi^*_x(f_1, e_1)- \phi^*_x(f_2, e_2)\big||U'(x)|\,{\mathrm d}x
\nonumber\\[2mm]
&\quad + \iota^3 \frac{4\pi^2}{{\tilde\ell}^2}\frac{1}{\beta}
 \int_{\mathbb R}\, x\big|\phi^*_{\tau\tau}(f_1, e_1)- \phi^*_{\tau\tau}(f_2, e_2)\big| |U'(x)|\,{\mathrm d}x
\nonumber\\[2mm]
&\quad + \iota^3 \frac{4\pi^2}{{\tilde\ell}^2}\frac{1}{\beta} \Big[\big|f_1- f_2\big|| f_1''| +|f_2| |f_1''- f_2''|\Big]\int_{\mathbb R}\, \big|\phi^*_x(f_1, e_1)\big||U'(x)|\,{\mathrm d}x
\nonumber\\[2mm]
&\quad + \iota^3 \frac{4\pi^2}{{\tilde\ell}^2}\frac{1}{\beta}  |f_2||f_2''| \int_{\mathbb R}\, \Big(\big|\phi^*_x(f_1, e_1)-\phi^*_x(f_2, e_2)\big| \Big)|U'(x)|\,{\mathrm d}x
\nonumber\\[2mm]
&\quad + \iota^3 \frac{4\pi^2}{{\tilde\ell}^2}\frac{1}{\beta}\big|f_1-f_2\big|
 \int_{\mathbb R}\,\big|\phi^*_{\tau\tau}(f_1, e_1)\big| |U'(x)|\,{\mathrm d}x
\nonumber\\[2mm]
&\quad + \iota^3 \frac{4\pi^2}{{\tilde\ell}^2}\frac{1}{\beta}\big|f_2\big|
 \int_{\mathbb R}\, \big|\phi^*_{\tau\tau}(f_1, e_1)- \phi^*_{\tau\tau}(f_2, e_2)\big| |U'(x)|\,{\mathrm d}x.
\end{align}

\medskip
 For  $\Pi_3(\tau, f, e)$, we have
\begin{align} \label{PI3estimate1}
\Big|\Pi_3(\tau, f, e) \Big|\,=\,&\Big|-\int_{{\mathbb R}}\chi(\iota|s|){\mathcal N}(\phi^*) U'(x)\,{\mathrm d}x\Big|
\nonumber\\[2mm]
\,\leq\,& \frac{C}{\beta^{2+  \mathfrak{p}}} \displaystyle\int_{{\mathbb R}}\chi(\iota|s|)
\Bigg[\,\mathbf{U}^{p_\epsilon- 3}  \left(
|{\tilde\phi}|^2+|{\tilde\psi}|^2\right)+\left(
|{\tilde\phi}|^{p_\epsilon-1} +|{\tilde\psi}|^{p_\epsilon-1}\right)
\nonumber\\[2mm]
&\qquad\qquad\qquad\qquad
+(p_\epsilon-1)\, \mathbf{U}^{p_\epsilon-2}  \,| {\tilde\psi}(\tilde \phi)|\Bigg]\,|U'(x)| \,{\mathrm d}x
\nonumber\\[2mm]
\,\leq\,&  C  \int_{{\mathbb R}}\chi(\iota|s|)  \left[U^{p_\epsilon- 3}  \big( |\phi^*|^2+ |\frac{\tilde\psi} {\beta^{\mathfrak{p}}}|^2\big) \right] \,|U'(x)| \,{\mathrm d}x
\nonumber\\[2mm]  &\quad +
C  \int_{{\mathbb R}}\chi(\iota|s|)  \left[ |\phi^*|^{5\pm \e}+ |\frac{\tilde\psi} {\beta^{\mathfrak{p}}}|^2  + (p_\epsilon-1)\,  U^{p_\epsilon-2} |\frac{\tilde\psi} {\beta^{\mathfrak{p}}}|\right] \,|U'(x)| \,{\mathrm d}x
\end{align}
with $\mathbf{U}(y)$ and $\tilde\phi(y)$ as in \eqref{definition of bf U} and \eqref{relation betwen tildephi and phi}.

\medskip
We now obtain the estimates of the above terms.
 We note
\begin{align} \label{inequalities for f}
\|\beta f(\tau)\|_{L^\infty(0, 2\pi)} \le\, C\, \tilde\ell^{\frac12} \|f\|_*,
\qquad
\Big\|\frac{\beta f'(\tau)}{\tilde \ell}\Big\|_{L^\infty(0, 2\pi)} \le \, C\,\tilde\ell^{\frac12} \|f\|_*,
\end{align}
and
\begin{align} \label{7-89}
\sup_{\tau \in (0, 2\pi)} \|\phi^*(\cdot, \tau)\|^2_{L^2({\mathbb R})} \le \frac{\tilde \ell} {\iota} \|\phi^*\|_{H^2_*(\mathfrak S)}^2, \quad
\sup_{\tau \in (0, 2\pi)} \|\phi_x^*(\cdot, \tau)\|^2_{L^2({\mathbb R})} \le \frac{\tilde \ell} {\iota} \|\phi^*\|_{H^2_*(\mathfrak S)}^2.
\end{align}
 For \eqref{7-89},     a similar proof can be found in Appendices  A in \cite{DWY}. We omit the proof here for concise.

\smallskip\noindent$\clubsuit$
We have
\begin{align*}
\int_0^{2\pi} \beta^2\Pi_1^2 \,{\mathrm d}\tau
& \leq C \iota^2 \int_0^{2\pi}  \beta^2  \frac{1}{\beta^{4\mathfrak{p}}} \left( \left|\beta f\right|^2+|e|^2 + 1\right) \|\beta^{\mathfrak{p}} \phi^*(\cdot, \tau)\|^2_{L^2({\mathbb R})}{\mathrm d}\tau
\nonumber\\[2mm]
& \leq
C \iota^2\left(1+\|e\|_{**}^2+ {\tilde\ell} \|f\|_*^2\right)\|\beta^{\mathfrak{p}} \phi^*\|_{H^2_*(\mathfrak S)}^2.
\end{align*}
Note that we used the fact  $4\mathfrak{p} -2= \frac {8}{ p_\epsilon -2} -2>0. $

\smallskip\noindent$\clubsuit$
Similarly, we have
\begin{align*}
 \int_0^{2\pi} \beta^2\Pi_2^2 \,{\mathrm d}\tau
& \leq C \iota^2 \int_0^{2\pi}
 \Bigg[\|\phi^*(\cdot, \tau)\|^2_{L^2({\mathbb R})} + \|\phi_x^*(\cdot, \tau)\|^2_{L^2({\mathbb R})}+\|\phi^*_{xx}(\cdot, \tau)\|^2_{L^2({\mathbb R})}
\nonumber\\[2mm]
& \qquad\qquad\qquad
 + \frac{\iota^2}{\tilde\ell^2}\|\phi^*_{\tau}(\cdot, \tau)\|^2_{L^2({\mathbb R})}+ \frac{\iota^2}{\tilde\ell^2}\|\phi^*_{x\tau}(\cdot, \tau)\|^2_{L^2({\mathbb R})} \Bigg] \left(\left|f\right|^2 + \left|\frac{{\beta^{\mathfrak p} f}'}{\tilde\ell}\right|^2 + 1\right) {\mathrm d}\tau
\nonumber\\[2mm]
& \quad+ \int_0^{2\pi} \beta^2\Pi_{2*}^2 \,{\mathrm d}\tau
\nonumber\\[2mm]
& \leq
 \iota^2\left(1 + {\tilde\ell} \|f\|_*^2\right)\|\phi^*\|_{H^2_*(\mathfrak S)}^2+ \int_0^{2\pi} \beta^2\Pi_{2*}^2 \,{\mathrm d}\tau.
\end{align*}
On the other hand, by \eqref{Pi2star} and \eqref{7-89}, there holds
\begin{align*}
\int_0^{2\pi} \beta^2\Pi_{2*}^2 \,{\mathrm d}\tau
&
\leq \iota^4 \int_0^{2\pi} \left|\frac{{\beta f}''}{\tilde\ell^2}\right|^2 \|\phi^*_x(\cdot, \tau)\|^2_{L^2({\mathbb R})}\,{\mathrm d}\tau
\ +\
\iota^2 \left(1 + {\tilde\ell} \|f\|_*^2\right) \int_0^{2\pi} \frac{\iota^4}{\tilde\ell^4}\|\phi^*_{\tau\tau}(\cdot, \tau)\|^2_{L^2({\mathbb R})} \,{\mathrm d}\tau
\nonumber\\[2mm]
&
\leq \iota^2\left(1 + {\tilde\ell} \|f\|_*^2\right)\|\phi^*\|_{H^2_*(\mathfrak S)}^2.
\end{align*}
Therefore, we have
\begin{align*}
& \int_0^{2\pi} \beta^2\Pi_2^2 \,{\mathrm d}\tau \leq \iota^2\left(1 + {\tilde\ell} \|f\|_*^2\right)\|\phi^*\|_{H^2_*(\mathfrak S)}^2.
\end{align*}

Moreover, by \eqref{minus of Pi2star1} and using the results of Proposition \ref{proposition6.5}, \eqref{inequalities for f} and \eqref{7-89}, we can verify that
\begin{align} \label{minus of Pi2star-in-L2-norm}
&\Big\|\beta\big(\Pi_{2*} (\theta, f_1, e_1) - \Pi_{2*} (\theta, f_2, e_2)\big)\Big\|_{L^2(0, 2\pi)}
\nonumber\\[2mm]
& \le C \iota^2 \sqrt{{\tilde\ell}/\iota\,}\, \big\|\phi^*_{f_1, e_1}\big\|_{H^2_*(\mathfrak S)}\big\|f_1- f_2\big\|_*
+ \iota^2\sqrt{{\tilde\ell}/\iota\,}\, \big\|f_2\big\|_* \,\big\|\phi^*_{f_1, e_1}- \phi^*_{f_2, e_2} \big\|_{H^2_*(\mathfrak S)}
\nonumber\\[2mm]
& \quad
+\iota^3 \sqrt{{\tilde\ell}/\iota\,}\, \big\|\phi^*_{f_1, e_1}\big\|_{H^2_*(\mathfrak S)}\big\|f_1- f_2\big\|_*
+ \iota^3\sqrt{{\tilde\ell}/\iota\,}\, \big\|f_2\big\|_* \,\big\|\phi^*_{f_1, e_1}- \phi^*_{f_2, e_2} \big\|_{H^2_*(\mathfrak S)}
\nonumber\\[2mm]
& \quad
+\iota^3 \sqrt{{\tilde\ell}/\iota\,}\, {\tilde \ell}^{\frac 12}\big\|f_1\big\|_* \big\|f_1- f_2\big\|_* \big\|\phi^*_{f_1, e_1}\big\|_{H^2_*(\mathfrak S)}
+
\iota^3 \sqrt{{\tilde\ell}/\iota\,}\, {\tilde \ell}^{\frac 12}\big\|f_2\big\|_* \big\|f_1- f_2\big\|_* \big\|\phi^*_{f_1, e_1}\big\|_{H^2_*(\mathfrak S)}
\nonumber\\[2mm]
& \quad+ \iota^3\sqrt{{\tilde\ell}/\iota\,}\, {\tilde \ell}^{\frac 12} \big\|f_2\big\|_*^2\,\big\|\phi^*_{f_1, e_1}- \phi^*_{f_2, e_2} \big\|_{H^2_*(\mathfrak S)}
+\iota\big\|\phi^*_{f_1, e_1}- \phi^*_{f_2, e_2} \big\|_{H^2_*(\mathfrak S)}
\nonumber\\[2mm]
&\quad+ \iota \tilde\ell^{\frac{1}{2}} \, \big\|\phi^*_{f_1, e_1}\big\|_{H^2_*(\mathfrak S)}\big\|f_1- f_2\big\|_*
+\iota \tilde\ell^{\frac{1}{2}}  \big\|f_2\big\|_*\big\|\phi^*_{f_1, e_1}- \phi^*_{f_2, e_2} \big\|_{H^2_*(\mathfrak S)}
\nonumber\\[2mm]
& \le C \iota^{3} \left(\big\|f_1- f_2\big\|_*+ \big\|e_1- e_2\big\|_* \right),
\end{align}
where we have used the fact that
\begin{align*}
 {\tilde\ell}^{\frac 12} \|f_i\|_* \ll 1,  \quad i=1,2.
\end{align*}

\smallskip\noindent$\clubsuit$
Next, from the expressions $\mathbf{U}(y)$ and $\tilde\phi(y)$ as in \eqref{definition of bf U}, \eqref{relation betwen tildephi and phi} and \eqref{PI3estimate1}, we have
\begin{align*}
 \int_0^{2\pi} \beta^2\Pi_3^2 \,{\mathrm d}\tau
& \leq
\int_0^{2\pi}  \frac{\beta^2} {\beta^{4\mathfrak{p}}} \Bigg[\int_{{\mathbb R}}  |\beta^{\mathfrak{p}}\phi^*|^2 |U'(x)| \,{\mathrm d}x\Bigg]^2 \,{\mathrm d}\tau
\ +\
\int_0^{2\pi} \frac{\beta^2} {\beta^{(10-2\e)\mathfrak{p}}} \Bigg[\int_{{\mathbb R}} |\beta^{\mathfrak{p}}\phi^*|^{p_\epsilon-1} |U'(x)| \,{\mathrm d}x \Bigg]^2 \,{\mathrm d}\tau
\nonumber\\[2mm]
& \quad + \int_0^{2\pi} \beta^2\left\{\int_{{\mathbb R}}\chi(\iota|s|)\left[  |{\tilde\psi}|^2+   \big| {\tilde\psi}\big|^{p_\epsilon-1} +\,  U^{p_\epsilon-2}  \, \big| {\tilde\psi}\big|\right] |U'(x)| \,{\mathrm d}x \right\}^2\,{\mathrm d}\tau.
\end{align*}
By the estimate
\begin{align*}
\|{\tilde\psi}({\tilde\phi})\|_{L^{\infty}}\, & \leq \, C\, \iota \lf[\|\tilde\phi\|_{L^{\infty}(|s|>\frac{20\delta}{\iota})}\,+\, \|\nabla \tilde\phi\|_{L^{\infty}(|s|>\frac{20\delta}{\iota})}\ri]\,+\,
e^{-\frac{ \upsilon_2 \delta}{\iota}}
\nonumber\\[2mm]
& \leq
C\left( 1 + \|\beta^{\mathfrak{p}} \phi^*\|_{H^2_*(\mathfrak S)}^2 \right) e^{-\frac{ \upsilon_3 \delta}{\iota}},
\end{align*}
we get
\begin{align*}
  \int_0^{2\pi} \beta^2\left\{\int_{{\mathbb R}}\chi(\iota|s|)\left[  |{\tilde\psi}|^2+   \big| {\tilde\psi}\big|^{p_\epsilon-1} +\,  U^{p_\epsilon-2}  \, \big| {\tilde\psi}\big|\right] |U'(x)| \,{\mathrm d}x \right\}^2\,{\mathrm d}\tau \leq
 C\left( 1 + \|\beta^{\mathfrak{p}} \phi^*\|_{H^2_*(\mathfrak S)}^2 \right) e^{-\frac{2\upsilon_3 \delta}{\iota}}.
\end{align*}

\medskip
By \eqref{7-89} and  Sobolev inequality, we can directly get
\begin{align} \label{phisuqareestimate}
\int_0^{2\pi} \left|\int_{{\mathbb R}} U |\beta^{\mathfrak{p}} \phi^*|^2 |U'(x)| \,{\mathrm d}x\right|^2 \,{\mathrm d}\tau
&\leq \int_0^{2\pi} \|\beta^{\mathfrak{p}} \phi^* \|_{L^{\infty}(\mathbb R)}^4 \,{\mathrm d}\tau
\nonumber\\[2mm]
&\leq \int_0^{2\pi} \left(\|\beta^{\mathfrak{p}} \phi^* \|_{L^{2}(\mathbb R)}^4+\|(\beta^{\mathfrak{p}} \phi^*)_{x} \|_{L^{2}(\mathbb R)}^4\right) \,{\mathrm d}\tau
\nonumber\\[2mm]
&\leq C \frac{{\tilde\ell}^2}{\iota^2}\|\beta^{\mathfrak{p}} \phi^*\|^4_{H^2_*(\mathfrak S)}.
\end{align}
 Similarly, we can prove
\begin{align*}
\int_0^{2\pi} \Bigg[\int_{{\mathbb R}} \left( \beta^{\mathfrak{p}} \right)^{p_\epsilon-1}U'(x)\,{\mathrm d}x\Bigg]^2 \,{\mathrm d}\tau
&\leq \int_0^{2\pi} \|\beta^{\mathfrak{p}} \phi^* \|_{L^{\infty}(\mathbb R)}^{10-2\e}  \,{\mathrm d}\tau
\nonumber\\[2mm]
&\leq \int_0^{2\pi} \left(\|\beta^{\mathfrak{p}} \phi^* \|_{L^{2}(\mathbb R)}^{10-2\e}+\|(\beta^{\mathfrak{p}} \phi^*)_{x} \|_{L^{2}(\mathbb R)}^{10-2\e}\right) \,{\mathrm d}\tau
\nonumber\\[2mm]
&\leq C \frac{{\tilde\ell}^{p_\epsilon-1}}{\iota^{p_\epsilon-1}} \|\beta^{\mathfrak{p}} \phi^*\|^{10-2\e}_{H^2_*(\mathfrak S)}.
\end{align*}
 Therefore, we have
\begin{align*}
 \int_0^{2\pi} \beta^2\Pi_3^2 \,{\mathrm d}\tau \leq C \frac{{\tilde\ell}^2}{\iota^2} \|\beta^{\mathfrak{p}} \phi^* \|^4_{H^2_*(\mathfrak S)} + \frac{{\tilde\ell}^{p_\epsilon-1}}{\iota^{p_\epsilon-1}} \|\beta^{\mathfrak{p}} \phi^*\|^{10-2\e}_{H^2_*(\mathfrak S)}.
\end{align*}

\medskip
Finally, by denoting
\begin{align} \label{hatPi1}
\hat\Pi_1(\tau, f, e) = \int_{{\mathbb R}} L_2(\phi^*) {\mathcal Z} \,{\mathrm d}x,
\end{align}
\begin{align} \label{hatPi2}
\hat\Pi_2(\tau, f, e) = \int_{{\mathbb R}} \chi(\iota|s|) \big[B_4(\phi^*)+ B_5(\phi^*)+ B_6(\phi^*)\big]{\mathcal Z} \,{\mathrm d}x,
\end{align}
\begin{align} \label{hatPi3}
\hat\Pi_3(\tau, f, e) = -\int_{{\mathbb R}}\chi(\iota|x|) {\mathcal N}(\phi^*)\mathcal Z \,{\mathrm d}x,
\end{align}
we have
\begin{align*}
 \int_0^{2\pi} \beta^2 \hat\Pi_1^2 \,{\mathrm d}\tau \leq \iota^2\left(1+\|e\|_{**}^2+ {\tilde\ell} \|f\|_*^2\right)\|\beta^{\mathfrak{p}}\phi^*\|_{H^2_*(\mathfrak S)}^2,
\end{align*}
\begin{align*}
 \int_0^{2\pi} \beta^2 \hat\Pi_2^2 \,{\mathrm d}\tau \leq \iota^2\left(1 + {\tilde\ell} \|f\|_*^2\right)\|\phi^*\|_{H^2_*(\mathfrak S)}^2,
\end{align*}
\begin{align*}
 \int_0^{2\pi} \beta^2 \hat\Pi_3^2 \,{\mathrm d}\tau
 \leq
C \frac{{\tilde\ell}^2}{\iota^2} \|\beta^{\mathfrak{p}} \phi^*\|^4_{H^2_*(\mathfrak S)}
+  \frac{{\tilde\ell}^{p_\epsilon-1}}{\iota^{p_\epsilon-1}} \|\beta^{\mathfrak{p}} \phi^*\|^{10-2\e}_{H^2_*(\mathfrak S)}.
\end{align*}
Moreover, we select the less regular terms from $\hat\Pi_1, \hat\Pi_2, \hat\Pi_3$ as $\hat\Pi_{2*}$, where
\begin{align}\label{hat-Pi2star}
\hat \Pi_{2*}(\tau, f, e)
=& -f''\iota^2\frac{4\pi^2}{{\tilde\ell}^2}\int_{\mathbb R}\, \phi^*_x \mathcal Z\,{\mathrm d}x
\nonumber\\[2mm]
&+ \iota^3 \kappa(\iota z) a_2(\iota s, \iota z) \frac{4\pi^2}{{\tilde\ell}^2}\frac{1}{\beta} \int_{\mathbb R} (x+f)\big[ -f''\phi^*_x +\phi^*_{\tau\tau} \big] \mathcal Z\,{\mathrm d}x.
\end{align}

We note that for $f$ and $e$ with constraints in \eqref{constraint-f}-\eqref{constraint-e}
\begin{align*}
 {\tilde\ell}^{\frac 12} \|f\|_* \ll 1, \quad \|e\|_{**} \ll \|f\|_*.
\end{align*}
Now we summarise the above calculations and present the following lemma.
\begin{lemma}
 Recall the definitions of $\Pi_\ell (\tau, f, e)$ an $\hat\Pi_\ell (\tau, f, e)$ with $\ell =1,2,3$ given in \eqref{Pi1}-\eqref{Pi3}, \eqref{hatPi1}-\eqref{hatPi3}. We have
\begin{align*}
 \int_0^{2\pi} \beta^2 \Pi_1^2 \,{\mathrm d}\tau &\leq \iota^2\left(1+\|e\|_{**}^2+ {\tilde\ell} \|f\|_*^2\right)\|\beta^{\mathfrak{p}} \phi^*\|_{H^2_*(\mathfrak S)}^2
 \le C \iota^{8-2\alpha},
\end{align*}
\begin{align*}
 \int_0^{2\pi} \beta^2\Pi_2^2 \,{\mathrm d}\tau & \leq C\iota^2\left(1 + {\tilde\ell} \|f\|_*\right)\|\phi^*\|_{H^2_*(\mathfrak S)}^2 \le C \iota^{8-2\alpha},
\end{align*}
\begin{align*}
 \int_0^{2\pi} \beta^2 \Pi_3^2 \,{\mathrm d}\tau
 & \leq C \frac{{\tilde\ell}^2}{\iota^2} \|\beta^{\mathfrak{p}} \phi^*\|^4_{H^2_*(\mathfrak S)} +  \frac{{\tilde\ell}^{p_\epsilon-1}}{\iota^{p_\epsilon-1}}  \|\beta^{\mathfrak{p}} \phi^*\|^{10-2\e}_{H^2_*(\mathfrak S)}
\nonumber\\[2mm]
&
 \le C \iota^{2\varrho -4} \iota^{12-4\alpha} +C\iota^{3\varrho -6} \iota^{18-6\alpha}
 \le C \iota^{8+2\varrho-4\alpha},
\end{align*}
\begin{align*}
 \int_0^{2\pi} \beta^2 \hat\Pi_1^2 \,{\mathrm d}\tau &\leq \iota^2\left(1+\|e\|_{**}^2+ {\tilde\ell} \|f\|_*^2\right)\| \beta^{\mathfrak{p}} \phi^*\|_{H^2_*(\mathfrak S)}^2
 \le C \iota^{8-2\alpha},
\end{align*}
\begin{align*}
 \int_0^{2\pi} \beta^2 \hat\Pi_2^2 \,{\mathrm d}\tau &\leq C\iota^2\left(1 + {\tilde\ell} \|f\|_*\right)\|\phi^*\|_{H^2_*(\mathfrak S)}^2
 \le \iota^{8-2\alpha},
\end{align*}
\begin{align*}
 \int_0^{2\pi} \beta^2 \hat\Pi_3^2 \,{\mathrm d}\tau
\leq
C \frac{{\tilde\ell}^2}{\iota^2} \|\beta^{\mathfrak{p}} \phi^*\|^4_{H^2_*(\mathfrak S)}
+   \frac{{\tilde\ell}^{p_\epsilon-1}}{\iota^{p_\epsilon-1}} \|\beta^{\mathfrak{p}} \phi^*\|^{10-2\e}_{H^2_*(\mathfrak S)}
 \le \iota^{8+2\varrho-4\alpha}.
\end{align*}
Moreover, the terms $ \Pi(\tau, f, e) = \sum_{\ell =1}^3\Pi_\ell(\tau, f, e)$ and $\hat\Pi(\tau, f, e) = \sum_{\ell =1}^3\hat \Pi_\ell (\tau, f, e)$ can be decomposed in the following way
\begin{align*}
 \Pi(\tau, f, e) = \Pi^*(\tau, f, e) + \Pi_{2*}(\tau, f, e), \qquad \hat\Pi(\tau, f, e) = \hat\Pi^*(\tau, f, e) + \hat\Pi_{2*}(\tau, f, e),
\end{align*}
In the above, the terms $\Pi^*(\tau, f, e)$ and $\hat\Pi^*(\tau, f, e)$ define compact operators of pair $(f,e)$ into $L^2(0, 2\pi)$, while the terms $\Pi_{2*}(\tau, f, e) $ and $ \hat\Pi_{2*}(\tau, f, e)$ satisfy
\begin{align*}
 \|\beta\left(\Pi_{2*} (\theta, f_1, e_1) - \Pi_{2*} (\theta, f_2, e_2) \right)\|_{L^2(0, 2\pi)}
&\leq C \iota^{3} \big[\|f_1- f_2\|_*+ \|e_1- e_2\|_{**}\big],
\end{align*}
\begin{align*}
\|\beta\left(\hat\Pi_{2*} (\theta, f_1, e_1) - \hat\Pi_{2*} (\theta, f_2, e_2) \right)\|_{L^2(0, 2\pi)}
&
\leq C \iota^{3} \big[\|f_1- f_2\|_*+ \|e_1- e_2\|_{**}\big].
\end{align*}
\qed
\end{lemma}

\medskip
\medskip
\subsection{  The reduced system of differential equations}\label{section9.3}

By the calculations in the above, we conclude that \eqref{project 6.13-1} and \eqref{project 6.13-2} will imply the following system
\begin{align*}
0= &
 \iota^2\varrho_1\Bigg\{-  \frac{4\pi^2}{{\tilde\ell}^2}f''  - 4\mathfrak{p}\frac{\pi}{\tilde\ell}  \frac{\beta'}{\beta^2} f'
+ \Big[(2\mathfrak{p}-1) \frac{|\beta'|^2}{\beta^4}
+\frac{\beta''} {\beta^3}
\,-\, \big(\frac 1{\Lambda}+1\big)\frac{\kappa^2}{\beta^2}  +\Lambda\frac{\iota^2 W_{tt}(0, \iota z)}{\beta^4} \Big]   f\Bigg\}
  \nonumber\\[2mm]
&  \quad
+\iota^2  \frac {{\mathbf b}_1(d_1+d_2)e}{\beta^{\mathfrak{p}}}
+ \iota^4  \frac{d_8} {\beta^{\mathfrak{p}}}{\mathbf b}_1 \frac{e''}{\tilde\ell^2}
+  \frac{\iota^3}{\beta^{\mathfrak{p}} }  \Big[ -\frac{8\pi^2}{{\tilde\ell}^2}f'
+ \frac{4\pi}{ {\tilde\ell}}\frac{\beta'}{\beta^2}  f\Big]e'd_7
   +  \frac{ \iota^3}{\beta^3}\tilde{\mathcal D}(\theta)
   \nonumber\\[2mm]
&  \quad
+\frac{\iota^3}{\beta^{1+\mathfrak{p}}}\tilde{\mathcal F}^{*} \left[\tau, \beta^{\mathfrak p} f, \frac{\beta^{\mathfrak p} {f}'}{\tilde\ell},  e, \frac{\iota e'}{\tilde\ell}\right]
 + \sum_{\ell =1}^3\Pi_\ell(\tau, f, e)  -\Pi_2^*(\tau, f, e) + \tilde M(f, e),
\end{align*}
and
\begin{align*}
 0& =  \frac 1{\beta^{\mathfrak{p}}} \Bigg[\iota^3\frac{4\pi^2}{{\tilde\ell}^2}  e''  +  \iota\lambda_0 e \Bigg]
 -\frac{1}{2}\, \iota^2 \frac{ \iota^2W_{tt}(0, \iota z)}{\beta^4}|f|^2  b_2
+ \iota^2\Big[\frac{4\pi^2}{{\tilde\ell}^2}|f'|^2
-\frac{4\pi}{{\tilde\ell}}\frac{\beta'}{\beta^2}ff'
+\frac{|\beta'|^2}{\beta^4} |f|^2
\Big] b_3
 \nonumber\\[2mm]
 & \quad
+ \iota^2 \frac1{\Lambda}\frac{\kappa}{\, \beta^2\, } fe
+  \iota^2 \frac{(p_\epsilon-1) (p_\epsilon- 2)}{2\beta^{2\mathfrak{p}}}e^2b_1
+ \iota^2  \frac{1}{\Lambda^2}f^2  {\mathbf b}_1^2 b_5
\nonumber\\[2mm]
 & \quad
+ \iota^2 (p_\epsilon-1)(p_\epsilon- 2) \Big[  \frac 1{2\Lambda^2} {\mathbf b}_1^2f^2b_6
+ \frac{e} {\beta^{\mathfrak{p}}}  \frac 1{\Lambda} {\mathbf b}_1f  b_7\Big]
- \iota^4 \frac{8\pi^2} {{\tilde\ell}^2} \frac{\kappa}{\, \beta^{1+ \mathfrak{p}}\, }  e''  f
\nonumber\\[2mm]
 & \quad
 +  \frac{\iota^3}{\beta^{1+\mathfrak{p}}}\tilde{\mathcal K}\left[\tau, \beta^{\mathfrak p} f, \frac{\beta^{\mathfrak p} {f}'}{\tilde\ell},  e, \frac{e'}{\tilde\ell}\right]
 +  \frac{\iota^4}{\beta^{1+\mathfrak{p}}} \tilde{\mathcal K}^{***}  \left[\tau, \beta^{\mathfrak p} f, \frac{\beta^{\mathfrak p} {f}'}{\tilde\ell}, e, \frac{\iota e'}{\tilde\ell}\right] +     \sum_{\ell =1}^3\hat\Pi_\ell(\tau, f, e)  -\hat\Pi_2^*(\tau, f, e) + \tilde N(f, e),
\end{align*}
 where
\begin{align}
 \tilde M(f, e) =
 \frac{\iota^3}{\beta^{1+\mathfrak{p}}}\tilde{\mathcal F}^{**} \left[\tau, \beta^{\mathfrak p} f,  e, \frac{\iota^2 e''}{\tilde\ell^2},  \frac{\beta^{\mathfrak p} f''}{\tilde\ell^2}\right]
 +\Pi_{2*}(\tau, f, e),
 \end{align}
\begin{align}
 \tilde N(f, e) =     \frac{\iota^3}{\beta^{1+\mathfrak{p}}} \tilde{\mathcal K}^{**} \left[\tau,  \beta^{\mathfrak p} f, \frac{\beta^{\mathfrak p} f''}{\tilde\ell^2}, \frac{\iota^2 e''}{\tilde \ell^2}\right]+ \hat\Pi_{2*}(\tau, f, e).
 \end{align}
By reducing, we can obtain
\begin{align}   \label{last reduced equation on U'}
 &  -  \frac{4\pi^2}{{\tilde\ell}^2}f''  - 4\mathfrak{p}\frac{\pi}{\tilde\ell}  \frac{\beta'}{\beta^2} f'
+ \Big[(2\mathfrak{p}-1) \frac{|\beta'|^2}{\beta^4}
+\frac{\beta''} {\beta^3}
\,-\, \Big(\frac 1{\Lambda}+1\Big)\frac{\kappa^2}{\beta^2}  +\Lambda\frac{\iota^2 W_{tt}(0, \iota z)}{\beta^4} \Big]   f
 \nonumber\\[2mm]
 &
=  -\frac1{\varrho_1}\left[    \frac { {\mathbf b}_1(d_1+d_2)e}{\beta^{\mathfrak{p}}}   +
  \iota^2  \frac{d_8} {\beta^{\mathfrak{p}}}{\mathbf b}_1 \frac{e''}{\tilde\ell^2} +   \frac{\iota}{\beta^{\mathfrak{p}} }  \Big(-\frac{8\pi^2}{{\tilde\ell}^2}f'  +\frac{4\pi}{ {\tilde\ell}}\frac{\beta'}{\beta^2}  f\Big) e'd_7  + \frac{ \iota}{\beta^3}\tilde{\mathcal D}(\theta) \right]
   \nonumber\\[2mm]
 &\quad - \mathcal A_{1\e}   (f, e) -  \mathcal K_{1\e}    (f, e),
 \end{align}
 and
 \smallskip
\begin{align}   \label{last reduced equation on mathcal Z}
 \iota^2\frac{4\pi^2}{{\tilde\ell}^2}  e''  + \lambda_0 e
  = &
 - \beta^{\mathfrak{p}}\Bigg\{-\frac{1}{2}\, \iota \frac{ \iota^2W_{tt}(0, \iota z)}{\beta^4}|f|^2  b_2
 + \iota\Big[\frac{4\pi^2}{{\tilde\ell}^2}|f'|^2
-\frac{4\pi}{{\tilde\ell}}\frac{\beta'}{\beta^2}ff'
+\frac{|\beta'|^2}{\beta^4} |f|^2\Big] b_3
\nonumber\\[2mm]
 & \qquad \quad
 +\iota  \frac1{\Lambda}\frac{\kappa}{\, \beta^2\, } fe
 + \iota \frac{(p_\epsilon-1) (p_\epsilon- 2)}{2\beta^{2\mathfrak{p}}}e^2b_1
 +  \iota^2\frac{1}{\Lambda} f^2  {\mathbf b}_1^2 b_5
    \nonumber\\[2mm]
 & \qquad \quad
 +   \iota (p_\epsilon-1)(p_\epsilon- 2) \Big[  \frac 1{2\Lambda^2} {\mathbf b}_1^2f^2b_6
+  \frac{e} {\beta^{\mathfrak{p}}}  \frac 1{\Lambda} {\mathbf b}_1f  b_7\Big]
   - \iota^3 \frac{8\pi^2} {{\tilde\ell}^2} \frac{\kappa}{\, \beta^{1+ \mathfrak{p}}\, }  e''  f        \Bigg\}
      \nonumber\\[2mm]
 &
  - \mathcal A_{2\e}   (f, e) -  \mathcal K_{2\e}    (f, e),
\end{align}
where
\begin{align}
  \mathcal K_{1\e}    (f, e)    & =
\frac{\iota}{\beta\, \varrho_1}\tilde{\mathcal F}^{*} \left[\theta, \beta^{\mathfrak p} f, \frac{\beta^{\mathfrak p}{f}'}{\tilde\ell},  e, \frac{e'}{\tilde\ell}\right]+ \frac{ \sum_{\ell =1}^3\Pi_\ell(\tau, f, e)  -\Pi_2^*(\tau, f, e) }{\iota^2\varrho_1},
\end{align}
\begin{align}
 \mathcal K_{2\e}    (f, e)
 &=    \frac{\iota^2}{\beta}  \tilde{\mathcal K}^* \left[\tau, \beta^{\mathfrak p} f, \frac{\beta^{\mathfrak p} {f}'}{\tilde\ell}, e, \frac{\iota e'}{\tilde\ell}\right]
 +  \frac{\iota^3}{\beta}   \tilde{\mathcal K}^{***} \left[\tau, \beta^{\mathfrak p}  f, \frac{\beta^{\mathfrak p} {f}'}{\tilde\ell}, e, \frac{\iota e'}{\tilde\ell}\right]
\nonumber\\[2mm]
& \quad
+ \frac{\beta^{\mathfrak{p}} \sum_{\ell =1}^3\hat\Pi_\ell(\tau, f, e)  -\beta^{\mathfrak{p}} \hat\Pi_2^*(\tau, f, e)}\iota,
 \end{align}
 and
\begin{align}
  \mathcal A_{1\e}   (f, e)  =  \frac{\tilde M(f, e)} {\iota^2\, \varrho_1},  \quad  \mathcal A_{2\e}   (f, e)  =\frac{\tilde N(f, e)} {\iota}.
  \end{align}

    \medskip
   Based on the above results,  we state the following important lemma.
\begin{lemma}    \label{Lemma7-2}
    To achieve the validity of  problem
    \eqref{project 6.13-1}-\eqref{project 6.13-2},  we need to find $(f, e)\in  \hat{\mathcal F}$
 such that \eqref{last reduced equation on U'}- \eqref{last reduced equation on mathcal Z} hold.
 Here we know that
 $\mathcal K_{\ell \e}  (f, e), \ell =1,2$  define    compact operators of pair $(f,e)$ into  $L^2(0, 2\pi)$ and   $  \mathcal A_{\ell\e}   (f, e), \ell =1,2$
  satisfy  the following      Lipschitz continuity
\begin{align}
  \| \mathcal A_{1\e}   (f, e) -  \mathcal A_{1\e}   (f, e) \|  \le   \iota (\|f_1- f_2\|_* + \|e_1- e_2\|_{**} ),
  \end{align}
\begin{align}
  \| \mathcal A_{2\e}   (f, e) -  \mathcal A_{2\e}   (f, e) \|  \le     \iota^2 (\|f_1- f_2\|_* + \|e_1- e_2\|_{**} ).
  \end{align}
  Moreover,  we have
\begin{align}
  \|\mathcal K_{1 \e}  (f, e) \|_{L^2(0, 2\pi)} +   \|\mathcal A_{1 \e}  (f, e) \|_{L^2(0, 2\pi)}    \le  C \left(\iota + \iota \|f\|_* + \iota\|e\|_{**}  \right),
  \end{align}
\begin{align}
  \|\mathcal K_{2 \e}  (f, e) \|_{L^2(0, 2\pi)}  +  \|\mathcal A_{2 \e}  (f, e) \|_{L^2(0, 2\pi)}   \le    C \left(\iota^3 + \iota^2 \|f\|_* + \iota^2 \|e\|_{**}  \right).
  \end{align}
    \end{lemma}

\bigskip
\section{Solving the reduced equations} \label{section10}

 \smallskip
\subsection{  Auxiliary linear problems} \label{section 8-1}

For the sake of later use, the resolution results for two linear problems will be provided
in Proposition \ref{reduced tilde f boundarycondition-lemma} and Proposition \ref{linear theory for e}.

\subsubsection{}
Consider
\begin{align}
-  \frac{4\pi^2}{{\tilde\ell}^2}{\mathfrak f}_{\tau\tau}  -4 \mathfrak{p}\frac{\pi}{\tilde\ell}  \frac{\beta'}{\beta^2}{\mathfrak f}_{\tau}
+ \Big[\big(2  \mathfrak{p}-1\big)\frac{|\beta'|^2}{\beta^4}
+\frac{\beta''} {\beta^3}
\,-\, \Big( 1+  \frac 1{\Lambda_\e}\Big)\frac{\kappa^2}{\beta^2} \Big]{\mathfrak f}  +  \,  \Lambda\frac{\iota^2 W_{tt}(0, \theta)}{\beta^4} {\mathfrak f}   = {\mathfrak h}(\tau)
\label{problem-auxil}
\end{align}
with
\begin{align}
{\mathfrak f}(0) = {\mathfrak f}(2\pi), \quad  {\mathfrak f}'(0) = {\mathfrak f}'(2\pi),
\label{problem-auxil-boundary}
\end{align}
where the relation between $\tau$ and $\theta=\iota z$ is given in \eqref{ztau}.
Note that  $\beta$ and $\kappa$ are smooth functions on the variable $\theta$, see \eqref{beta} and \eqref{curvature bounded}.

\begin{proposition}  \label{reduced tilde f boundarycondition-lemma}
For any given $\epsilon$ small with the validity of \eqref{resonace-f}, we have  a constant $C>0$ such that
\begin{align} \label{reduced mathfrak f boundarycondition}
\left\|\beta^{\mathfrak p} {\mathfrak f}\right\|_{L^2(0, 2\pi)}^2
+\left\|\frac{\beta^{\mathfrak p}}{\tilde\ell} {\mathfrak f}_\tau\right\|_{L^2(0, 2\pi)}^2
+\left\|\frac{\beta^{\mathfrak p}}{{\tilde\ell}^2} {\mathfrak f}_{\tau\tau}\right\|_{L^2(0, 2\pi)}^2
\leq
\frac{C}{{\mathrm d}_\iota^2}\left\|\beta^{\mathfrak p}  {\mathfrak h}\right\|_{L^2(0, 2\pi)}^2,
\end{align}
where ${\mathrm d}_\iota$ is given in \eqref{resonace-f}.
\qed
\end{proposition}

\medskip
\begin{proof} The proof will be divided into several steps.

{\textbf{Step 1.}}
By denoting
\begin{align}
{\mathfrak f}(\tau) = \beta(\theta) {\tilde f}(\tau), \qquad {\mathfrak h}(\tau) = \beta(\theta){\tilde h}(\tau),
\label{setting1}
\end{align}
we have
\begin{align}
{\mathfrak f}_\tau = \beta {\tilde f}_\tau + \frac{\tilde\ell}{2\pi}\frac{\beta'}{\beta} {\tilde f},
\label{setting2}
\end{align}
\begin{align}
{\mathfrak f}_{\tau\tau} & = \beta {\tilde f}_{\tau\tau} +2 \beta'  {\tilde f}_\tau \frac{\tilde\ell} {2\pi \beta}  + \beta'' {\tilde f}\Big|\frac{\tilde\ell} {2\pi \beta}\Big|^2 -\beta'  {\tilde f} \frac{\tilde\ell} {2\pi \beta^2}\beta'\frac{\tilde\ell} {2\pi \beta}
\nonumber\\[2mm]
& =  \beta {\tilde f}_{\tau\tau}
+ \frac{\tilde\ell}{\pi}\frac{\beta'}{\beta} {\tilde f}_\tau
+ \frac{\tilde\ell^2} {4\pi^2} \frac{\beta\beta''- |\beta'|^2}{\beta^3} {\tilde f}.
\label{setting3}
\end{align}
The above setting will give a new form of \eqref{problem-auxil}
\begin{align*}
& -  \frac{4\pi^2}{{\tilde\ell}^2} \beta\Big[ {\tilde f}_{\tau\tau} + \beta'  {\tilde f}_\tau \frac{\tilde\ell} {\pi \beta^2}  + \beta'' {\tilde f}\frac{\tilde \ell^2} {4\pi^2 \beta^3} -|\beta'|^2  {\tilde f} \frac{\tilde \ell^2} {4\pi^2 \beta^4}\Big]
\nonumber\\[2mm]
& - 4 \mathfrak{p}\frac{\pi}{\tilde\ell}  \frac{\beta'}{\beta^2} \Big[\beta {\tilde f}_\tau + \frac{\tilde\ell}{2\pi}\frac{\beta'}{\beta} {\tilde f} \Big]
+
\Big[(2  \mathfrak{p}-1)\frac{|\beta'|^2}{\beta^4}
+\frac{\beta''} {\beta^3}
\,-\, \Big( 1+  \frac 1{\Lambda}\Big)\frac{\kappa^2}{\beta^2} \Big]\beta{\tilde f}   + \, \Lambda\frac{\iota^2 W_{tt}(0, \theta)}{\beta^4}\beta{\tilde f}
=  \beta{\tilde h},
\end{align*}
i.e,
\begin{align*}
-  \frac{4\pi^2}{{\tilde\ell}^2}\Big[ {\tilde f}_{\tau\tau} + \beta'  {\tilde f}_\tau \frac{\tilde\ell} {\pi \beta^2}  \Big]
- 4 \mathfrak{p}\frac{\pi}{\tilde\ell}  \frac{\beta'}{\beta^2} \Big[{\tilde f}_\tau + \beta'  {\tilde f} \frac{\tilde\ell} {2\pi \beta^2}\Big]
+\Big[2\mathfrak{p}\frac{|\beta'|^2}{\beta^4}\,-\, \Big( 1+  \frac 1{\Lambda}\Big)\frac{\kappa^2}{\beta^2} \Big]{\tilde f}   + \, \Lambda\frac{\iota^2 W_{tt}(0, \theta)}{\beta^4}{\tilde f}
=  {\tilde h},
\end{align*}
i.e,
\begin{align*}
-  \frac{4\pi^2}{{\tilde\ell}^2}\Big[ {\tilde f}_{\tau\tau} +(1+ \mathfrak{p}) \beta'  {\tilde f}_\tau \frac{\tilde\ell} {\pi \beta^2}  \Big]
+
\Big[\Lambda\frac{\iota^2 W_{tt}(0, \theta)}{\beta^4}\,-\,\Big( 1+  \frac 1{\Lambda}\Big)\frac{\kappa^2}{\beta^2} \Big]{\tilde f}
=  {\tilde h},
\end{align*}
i.e,
\begin{align*}
   -  \frac{4\pi^2}{{\tilde\ell}^2}  \frac 1{ \beta^{2+2 \mathfrak{p}} } \p_\tau \left( \beta^{2+2\mathfrak{p}} {\tilde f}_\tau \right)
+
\Big[\Lambda\frac{\iota^2 W_{tt}(0, \theta)}{\beta^4}\,-\,\Big( 1+  \frac 1{\Lambda}\Big)\frac{\kappa^2}{\beta^2} \Big]{\tilde f}
=  {\tilde h}.
\end{align*}

\medskip
{\textbf{Step 2.}}
We consider the  weighted eigenvalue problem
\begin{align}  \label{weighted-reduced tilde f}
\tilde{\mathbb{F}}_\epsilon ({\tilde f}) :=  -  \frac{4\pi^2}{{\tilde\ell}^2} \p_\tau \left(  \beta^{2+ 2\mathfrak p}  {\tilde f}_\tau \right)
+
\Big[\Lambda\frac{\iota^2 W_{tt}(0, \theta)}{\beta^4}\,-\,\Big( 1+  \frac 1{\Lambda}\Big)\frac{\kappa^2}{\beta^2} \Big] \beta^{2+ 2\mathfrak p} {\tilde f} =  \Xi \beta^{2+ 2\mathfrak p} {\tilde h},
\end{align}
subjecting to periodic boundary conditions
\begin{align} \label{newreduced tilde f boundarycondition}
\tilde f(0) = \tilde f(2\pi), \quad  \tilde f'(0) = \tilde f'(2\pi).
\end{align}
Simple calculations will show that problem \eqref{weighted-reduced tilde f}-\eqref{newreduced tilde f boundarycondition} is
equivalent to the eigenvalue problem \eqref{mathbbF}-\eqref{mathbbF-boundary}.
By introducing the weighted inner product $\langle \cdot,  \cdot\rangle_{W}$ as
\begin{align*}
 \langle {\tilde f},  {\tilde h}\rangle_{W} = \int_0^{2\pi} \beta^{2+ 2\mathfrak p} {\tilde f} {\tilde h}{\mathrm d} \tau,
\end{align*}
we first give the following claim, which is well-known.
\\
{\textbf{Claim 1.}} \label{weighted eigenvalue problem lemma}
{\em
There exists  a sequence $\{\Lambda_{\epsilon, k}\}_{k=1}^\infty$ of eigenvalues to the weighted eigenvalue problem  \eqref{weighted-reduced tilde f}-\eqref{newreduced tilde f boundarycondition}.
Moreover, the corresponding eigenfunctions  $\{\omega_{\epsilon, k}\}_{k=1}^\infty$ form  an orthonormal basis of $L^2(0, 2\pi)$ with respect to the weighted inner product $\langle \cdot,  \cdot\rangle_{W}$.
}

\medskip
{\textbf{Step 3.}}
We then concern the problem
\begin{align}
 -  \frac{4\pi^2}{{\tilde\ell}^2} \p_\tau \left(  \beta^{2+ 2\mathfrak p}  {\tilde f}_\tau \right)
+
\Big[\Lambda\frac{\iota^2 W_{tt}(0, \theta)}{\beta^4}\,-\,\Big( 1+  \frac 1{\Lambda}\Big)\frac{\kappa^2}{\beta^2} \Big] \beta^{2+ 2\mathfrak p} {\tilde f}   =   \beta^{2+ 2\mathfrak p}{\tilde h},
\label{weighted-reduced tilde f222222}
\end{align}
with
\begin{align}
\tilde f(0) = \tilde f(2\pi), \quad  \tilde f'(0) = \tilde f'(2\pi),
\label{newreduced tilde f boundarycondition222222}
\end{align}
and obtain the following claim.

\smallskip
\noindent
{\textbf{Claim 2.}}   \label{reduced tilde f boundarycondition-lemma22222}
{\em
Consider \eqref{weighted-reduced tilde f222222}-\eqref{newreduced tilde f boundarycondition222222}.
For any given $\epsilon$ small, we have  a constant $C>0$ such that
\begin{align} \label{reduced tilde f boundarycondition}
\left\| \beta^{1+\mathfrak p}{\tilde f} \right\|_{L^2(0, 2\pi)}^2
+\left\| \frac{\beta^{1+\mathfrak p}}{\tilde\ell}{\tilde f}_\tau \right\|_{L^2(0, 2\pi)}^2
+\left\| \frac{\beta^{1+\mathfrak p}}{{\tilde\ell}^2}{\tilde f}_{\tau\tau} \right\|_{L^2(0, 2\pi)}^2
\leq
\frac{C}{{\mathrm d}_\iota^2}\left\| \beta^{1+\mathfrak p}{\tilde h} \right\|_{L^2(0, 2\pi)}^2.
\end{align}
}
  We note that the proof of  {\textbf{Claim 2}}    is based on  {\textbf{Claim 1}} and   can be found in \cite{DWY}.

\medskip
{\textbf{Step 4.}}
A conclusive remark is that the proof of Proposition \ref{reduced tilde f boundarycondition-lemma} will be easily derived from Claim 2
and the identities \eqref{setting1}-\eqref{setting3}.
In fact, some details will be provided here.
It is obvious that
\begin{align*}
\left\| \beta^{\mathfrak p}{\mathfrak f} \right\|_{L^2(0, 2\pi)}^2
= \left\| \beta^{1+\mathfrak p}{\tilde f} \right\|_{L^2(0, 2\pi)}^2,
\end{align*}

\begin{align*}
\left\| \frac{\beta^{\mathfrak p}}{\tilde\ell} {\mathfrak f}_\tau \right\|_{L^2(0, 2\pi)}^2
\leq
C\left\| \frac{\beta^{1+\mathfrak p}}{\tilde\ell}{\tilde f}_\tau \right\|_{L^2(0, 2\pi)}^2
+
C\left\| \beta^{\mathfrak p}{\tilde f} \right\|_{L^2(0, 2\pi)}^2,
\end{align*}

\begin{align*}
\left\| \frac{\beta^{\mathfrak p}}{{\tilde\ell}^2} {\mathfrak f}_{\tau\tau} \right\|_{L^2(0, 2\pi)}^2
\leq
C\left\| \frac{\beta^{1+\mathfrak p}}{{\tilde\ell}^2}{\tilde f}_{\tau\tau} \right\|_{L^2(0, 2\pi)}^2
+
C\left\| \frac{\beta^{\mathfrak p} }{{\tilde\ell}}{\tilde f}_\tau \right\|_{L^2(0, 2\pi)}^2
+
C\left\|\beta^{\mathfrak p} {\tilde f} \right\|_{L^2(0, 2\pi)}^2,
\end{align*}
due to \eqref{beta-bound}.
The combination of the above estimates and \eqref{reduced tilde f boundarycondition}
will finish the proof of  \eqref{reduced mathfrak f boundarycondition}.
\end{proof}

\medskip
\subsubsection{}
For the linear problem
\begin{align} \label{reduced for e2 new}
\iota^2\frac{4\pi^2}{{\tilde\ell}^2} {\mathfrak e}'' + \lambda_0 {\mathfrak e} \,=\, d\quad\mbox{in }\ (0,2\pi),
\quad {\mathfrak e}(0) \,=\, {\mathfrak e}(2\pi),
\quad {\mathfrak e}'(0) \,=\, {\mathfrak e}' (2\pi),
\end{align}
where the constants ${\tilde\ell}$, $\lambda_0$ are given in \eqref{ztau} and \eqref{lambda0}, we now present the following linear theory.

\begin{proposition} \label{linear theory for e}
Suppose that \eqref{gap condition1} holds.
If $d\in L^2(0,  2\pi)$, there exits a unique solution ${\mathfrak e}\in H^2 (0,2\pi)$ to the linear problem \eqref{reduced for e2 new}, which satisfies the estimate
\begin{align} \label{estimate norm of e1}
 \|{\mathfrak e}\|^2_{**} :\,=\, \|{\mathfrak e}\|^2_{L^\infty(0,  2\pi)} + \frac{\iota^2}{{\tilde\ell}^2}  \|{\mathfrak e}'\|^2_{L^2 (0, 2\pi)}
 + \frac{\iota^4} {{\tilde\ell}^4}\|{\mathfrak e}''\|^2_{L^2(0, 2\pi)}\leq C \frac{{\tilde\ell}^2} {\iota^2} \| d\|^2_{L^2(0, 2\pi)}.
\end{align}

\end{proposition}
\begin{proof}
 The proof of this Proposition can be found in \cite{DWY}.
 \end{proof}

 \medskip
\subsection{  Solving the reduced nonlinear system}

 Note that for any pair $(f, e)$ in $\hat{\mathcal F}$, there hold
$$
 \|f \|_*\leq \iota^{1-\alpha}, \quad \|e \|_{**}\leq  \iota^{\frac 32 \varrho + \frac 12- 2\alpha}.
$$
 We denote $\hat{\mathcal L}  (f, e)= \left( \hat {\mathcal L}_1,   \hat {\mathcal L}_2 \right) (f, e), $ with
\begin{align*}
 \hat {\mathcal L}_1 (f, e):
\,=\, &  -  \frac{4\pi^2}{{\tilde\ell}^2}{\mathfrak f}_{\tau\tau}
-4 \mathfrak{p}\frac{\pi}{\tilde\ell}  \frac{\beta'}{\beta^2}{\mathfrak f}_{\tau}
+ \Big[\big(2  \mathfrak{p}-1\big)\frac{|\beta'|^2}{\beta^4}
+\frac{\beta''} {\beta^3}
\,-\, \Big( 1+  \frac 1{\Lambda_\e}\Big)\frac{\kappa^2}{\beta^2} \Big]{\mathfrak f}  +  \,  \Lambda\frac{\iota^2 W_{tt}(0, \theta)}{\beta^4} {\mathfrak f},
\nonumber\\[2mm]
\hat {\mathcal L}_2 (f, e):
\,=\, &   \iota^2\frac{4\pi^2}{{\tilde\ell}^2}  e''  + \lambda_0 e.
\end{align*}
By  Proposition  \ref{reduced tilde f boundarycondition-lemma} and  Proposition  \ref{linear theory for e} provided in section \ref{section 8-1}, we know that the problem
\begin{equation}\label{liner hd}
\hat {\mathcal L} (f, e) =  \left({\mathfrak h},    d\right)
\end{equation}
with the boundary conditions
\begin{align}
f(0) \,=\, f(2\pi), \quad f'(0) \,=\, f' (2\pi), \quad e(0) \,=\, e(2\pi), \quad e'(0) \,=\, e' (2\pi),
\label{boundaryhatf-e}
\end{align}
is invertible and $(f, e)$  satisfy the estimates
\begin{align}
 \|f\|_* \le \frac{C} {\iota^\alpha}\|\beta^{\mathfrak p} {\mathfrak h}\|_{L^2(0, 2\pi)}, \quad   \|e\|_{**}  \le \frac{C}{\iota^{2-\varrho}}  \|d\|_{L^2(0, 2\pi)}.
\end{align}
Denote
\begin{align*}
 \mathcal Q_{1\e} (f, e) =  -\frac1{\varrho_1}\Bigg\{  \frac { {\mathbf b}_1(d_1+d_2)e}{\beta^{\mathfrak{p}}}   +
  \iota^2  \frac{d_8} {\beta^{\mathfrak{p}}}{\mathbf b}_1 \frac{e''}{\tilde\ell^2} +   \frac{\iota}{\beta^{\mathfrak{p}} }  \Big[ -\frac{8\pi^2}{{\tilde\ell}^2}f'  +\frac{4\pi}{ {\tilde\ell}}\frac{\beta'}{\beta^2}  f\Big]e'd_7  + \frac{ \iota}{\beta^3}\tilde{\mathcal D}(\theta)  \Bigg\},
 \end{align*}
\begin{align*}
 \mathcal Q_{2\e} (f, e)
\,=\,&   - \beta^{\mathfrak{p}}\Bigg\{ -\frac{1}{2}\, \iota \frac{ \iota^2W_{tt}(0, \iota z)}{\beta^4}|f|^2  b_2  + \iota\Big[\frac{4\pi^2}{{\tilde\ell}^2}|f'|^2
-\frac{4\pi}{{\tilde\ell}}\frac{\beta'}{\beta^2}ff'
+\frac{|\beta'|^2}{\beta^4} |f|^2
\Big] b_3
\nonumber\\[2mm]
&\qquad\qquad
+\iota\frac1{\Lambda_\e} \frac{\kappa}{\, \beta^2\, } fe b_0
+ \iota \frac{(p_\epsilon-1) (p_\epsilon- 2)}{2\beta^{2\mathfrak{p}}}e^2b_1
+ \iota \frac{1}{\Lambda_\e}f^2  {\mathbf b}_1^2 b_5
\nonumber\\[2mm]
&\qquad\qquad +   \iota(p_\epsilon-1)(p_\epsilon- 2)  \Big[  \frac 1{2\Lambda_\e^2} {\mathbf b}_1^2f^2b_6
+  \frac{e} {\beta^{\mathfrak{p}}}  \frac 1{\Lambda_\e} {\mathbf b}_1f  b_7\Big]
  - \iota^3 \frac{8\pi^2} {{\tilde\ell}^2} \frac{\kappa}{\, \beta^{1+ \mathfrak{p}}\, }  e''  f        \Bigg\}.  \end{align*}
For nonlinear coupling problems,
\begin{align}    \label{7-138}
 \hat {\mathcal L} (f, e) =  \Big[\left(\mathcal Q_{1\e},   \mathcal Q_{2\e}\right) - \left(\mathcal A_{1\e},   \mathcal A_{2\e}  \right) \Big](f, e)+ \left({\mathfrak h},    d\right),
 \end{align}
we can deduce from the contraction mapping principle that the  problem  is uniquely solvable  if
$$\|\beta^{\mathfrak p} {\mathfrak h}\|_{L^2(0, 2\pi)}  \le {C} \iota, \quad    \|d\|_{L^2(0, 2\pi)} \le  C \iota^{\frac12\varrho+\frac5 2-2 \alpha}.
$$
Define  the set
\[\tilde{\mathcal B} = \Big\{ (f, e)\,\Big|\, (f, e)~\text{satisfies} ~\eqref{boundaryhatf-e}~\text{and}~\|f\|_* \le \bar{C}_1 \iota^{1-\alpha},   \|e\|_{**} \le \bar{C}_2 \iota^{\frac 32 \varrho+ \frac12-2\alpha} \Big\},  \]
with $\bar C_1, \bar C_2$ large enough.
 In fact,   the  problems \eqref{7-138}  can be transformed into the following fixed point problem
\begin{align}   \label{fixed point problem}
  (f,e)  & = \hat {\mathcal L}^{-1}\Bigg( \Big[  \left(\mathcal Q_{1\e},   \mathcal Q_{2\e}\right) - \left(\mathcal A_{1\e},   \mathcal A_{2\e}\right)   \Big]  (f, e)  + \left({\mathfrak h},    d\right)\Bigg)
: =  \mathbb{T} (f, e),
  \end{align}
and
\begin{align*}
\hat {\mathcal L}^{-1}  = \left( \hat {\mathcal L}_1^{-1},   \hat {\mathcal L}_2^{-1}\right).
\end{align*}
 By the results provided in Lemma \ref{Lemma7-2}, and combining  \eqref{constraint-alpha-varrho}, \eqref{Linfty of f},    we can show that
$\mathbb{T} $ is a contradiction map from $\tilde{\mathcal B}$ to itself.    From the contraction mapping principle, problem \eqref{fixed point problem} has   a fixed point in $\tilde{\mathcal B}$.
   Next, we consider  the problem
\begin{align}    \label{7-}
 \hat {\mathcal L} (f, e) = \Big[\left(\mathcal Q_{1\e},   \mathcal Q_{2\e}\right) - \left(\mathcal A_{1\e},   \mathcal A_{2\e}   \right)  - \left(\mathcal K_{1\e},   \mathcal K_{2\e}\right)\Big](f, e),
 \end{align}
 under the boundary condition \eqref{boundaryhatf-e}.
  Denote by
\begin{align*}
  \tilde {\mathcal L}  (e,f)=   \hat {\mathcal L} (f, e) -   \Big[\left(\mathcal Q_{1\e},   \mathcal Q_{2\e}\right) - \left(\mathcal A_{1\e},   \mathcal A_{2\e}\right)\Big] (f,e).
  \end{align*}
  Then we transfer problem \eqref{7-} the following fixed point problem
\begin{align}
  (f,e) =  \underbrace{- \tilde {\mathcal L}^{-1} \left(\mathcal K_{1\e},   \mathcal K_{2\e}\right)}_{\text{compact operator}}   (f,e).
  \end{align}
  By  Schauder's fixed point theorem and the estimates provided in Lemma \ref{Lemma7-2}, we can get the existence of $(f, e)$ for \eqref{7-}, which satisfies \eqref{constraint-f}-\eqref{constraint-e}.
 \qed

\medskip
 \section{Choosing the parameter $\lambda$ from mass constraint (\ref{constraint2})}\label{Section9}

After finding a solution $v$ of \eqref{eq3} in Sections \ref{section-approximatesolution}-\ref{section10}, we now consider the constraint \eqref{constraint2},
which can be formulated as the following
\begin{align} \label{constraint6}
 a\iota^{-2+ \frac4{p_\epsilon-2}}& \,=\, \int_{ {\mathfrak S}_\delta} |v|^2\,{\mathrm d}y
\ +\
\int_{{\mathbb R}^2\setminus{\mathfrak S}_\delta} |v|^2\,{\mathrm d}y
\nonumber\\[2mm]
 & \,=\, \int_{{\mathfrak S}_\delta } |v|^2 \Big[1+ \iota\kappa(\iota z) s \Big] {\mathrm d} s\,{\mathrm d} z
 \ +\
\int_{{\mathbb R}^2\setminus {\mathfrak S}_\delta } |v|^2\,{\mathrm d}y
\nonumber\\[2mm]
 & \,=\,
\frac1\iota\int_{{\mathfrak S}_\delta} |v|^2 \Bigg[1+ \kappa \frac{\iota (x+f)} {\beta} \Bigg] \frac {\tilde \ell} {2\pi \beta}{\mathrm d} \left(\frac{ x+f} {\beta}\right)\,{\mathrm d}\tau
\ +\
\int_{{\mathbb R}^2\setminus{\mathfrak S}_\delta} |v|^2\,{\mathrm d}y,
\end{align}
where
\begin{align*}
{\mathfrak S}_\delta=\Bigg\{ (x, \tau): \quad |x| \le \frac {\beta \delta} \iota, \quad 0<\tau< 2\pi\Bigg\}.
\end{align*}
 Since $|x| \le  \frac {\beta \delta} \iota$, then we have
\begin{align} \label{out of delta}
 -\frac\delta \iota + \frac{\iota f} {\beta} < s= \frac{(x+f)} {\beta} \le \frac \delta \iota+ \frac{\iota f} {\beta}.
\end{align}
 Following the analysis in Section \ref{section8}, we know that
\begin{align*}
 v(y) & = \mathbf{U}(y) \,+\, \eta_{3\delta}^{\iota}(\iota |s|){\tilde\phi}(y)
 \,+\, {\tilde\psi}(y)
\nonumber\\[2mm]
 &
 = \eta_{3\delta}^\iota(s) \beta^{\mathfrak p}(\iota z)w_2(x, \tau)
 \,+\, \eta_{3\delta}^{\iota}(\iota |s|){\tilde\phi}(y)
 \,+\, {\tilde\psi}(y).
\end{align*}
In the above, $ \eta_{3\delta}^\iota(s) $ is the cut-off function given in Section \ref{section 6-1}, and the function $w_2(x, \tau) $ is given in \eqref{newthirdapproximatesolution} with the form
\begin{align*}
w_2(x, \tau) = U(x) \,+\, \iota \frac{e(\tau)}{\beta^{\mathfrak p}(\iota z)} \mathcal Z(x) \,+\, \iota \phi_1(x, \tau) \,+\, \iota^2\phi_2(x, \tau).
\end{align*}
where $\phi_1, \phi_2$ are given in \eqref{definition of phi1}, \eqref{definition of e2phi2}.

In the region ${\mathbb R}^2\setminus{\mathfrak S}_\delta$, which means far from the curve $\Gamma_\iota$, the functions
\[ w_2(x, \tau), \quad \eta_{3\delta}^{\iota}(\iota |s|){\tilde\phi}(y), \quad {\tilde\psi}(y), \]
exhibit exponential decay, and by \eqref{out of delta}, thus the integral
\begin{align}
 \int_{{\mathbb R}^2\setminus{\mathfrak S}_\delta} |v|^2\,{\mathrm d}y \in O\left(e^{-\frac{\tilde \delta}\iota} \right),
\label{estimate1}
\end{align}
where $\tilde \delta$ is small positive constant.
On the other hand, we have
\begin{align}
 & \frac1\iota\int_{{\mathfrak S}_\delta} |v|^2 \Bigg[1+ \kappa \frac{\iota (x+f)} {\beta} \Bigg] \frac {\tilde \ell} {2\pi \beta}{\mathrm d} \left(\frac{ x+f} {\beta}\right)\,{\mathrm d}\tau
\nonumber\\[2mm]
 &
 = \frac {\tilde \ell} {2\pi \iota} \int_{{\mathfrak S}_\delta} |v|^2 \frac {1} {\beta^2}\, {\mathrm d} x \,{\mathrm d}\tau
\nonumber\\[2mm]
 &
 = \frac {\tilde \ell} {2\pi \iota} \int_{{\mathfrak S}_\delta} \frac {1} {\beta^2} \left[ |\beta|^{2\mathfrak p} w_2^2+ |{\tilde\phi}|^2+ |{\tilde\psi}|^2 + 2\beta^{\mathfrak p} w_2 \big( {\tilde\phi} + {\tilde\psi} \big) + 2 \tilde \phi\tilde \psi \right] \, {\mathrm d} x \,{\mathrm d}\tau
\nonumber\\[2mm]
 &
 = \frac {\tilde \ell} {2\pi \iota} \int_{{\mathfrak S}_\delta} \beta^{2\mathfrak p -2}w_2^2 \, {\mathrm d} x \,{\mathrm d}\tau
\ +\
\frac {\tilde \ell} {2\pi \iota} \int_{{\mathfrak S}_\delta} \frac {1} {\beta^2} \left[ |{\tilde\phi}|^2+ |{\tilde\psi}|^2 + 2\beta^{\mathfrak p} w_2 \big( {\tilde\phi} + {\tilde\psi} \big) + 2 \tilde \phi\tilde \psi \right] \, {\mathrm d} x \,{\mathrm d}\tau.
\label{estimate2}
\end{align}
The terms in the above formula will be estimated in the following way.

\smallskip
\noindent $\spadesuit$\
For the first term, we have
\begin{align*}
  \frac {\tilde \ell} {2\pi \iota} \int_{{\mathfrak S}_\delta} \beta^{2\mathfrak p -2}w_2^2 \, {\mathrm d} x \,{\mathrm d}\tau
\,=\,&\frac {\tilde \ell} {2\pi \iota} \int_{{\mathfrak S}_\delta} \beta^{2\mathfrak p -2}\big|U(x) + \iota \frac e{\beta^{\mathfrak p}}  \mathcal Z + \iota \phi_1 + \iota^2\phi_2 \big|^2\, {\mathrm d} x \,{\mathrm d}\tau
 \nonumber\\[2mm]
\,=\,&\frac {\tilde \ell} {2\pi \iota} \int_{{\mathfrak S}} \beta^{2\mathfrak p -2}\big|U(x) + \iota \frac e{\beta^{\mathfrak p}}  \mathcal Z + \iota \phi_1 + \iota^2\phi_2 \big|^2\, {\mathrm d} x \,{\mathrm d}\tau
  +O\left(e^{-\frac{\tilde \delta}\iota} \right)
\nonumber\\[2mm]
\,=\,& \frac {\tilde \ell} {2\pi \iota} \int_{{\mathfrak S}}
\beta^{2\mathfrak p -2}\Big[ U^2 + \iota^2\frac{e^2}{\beta^{2\mathfrak p}}|\mathcal Z|^2 + \iota^2|\phi_1|^2+ \iota^4|\phi_2|^2\Big]\, {\mathrm d} x \,{\mathrm d}\tau
\nonumber\\[2mm]
& + \frac {\tilde \ell} {\pi \iota} \int_{{\mathfrak S}}  \beta^{2\mathfrak p -2}U \left(\iota \frac e{\beta^{\mathfrak p}}  \mathcal Z + \iota \phi_1 + \iota^2\phi_2 \right)\, {\mathrm d} x \,{\mathrm d}\tau
\nonumber\\[2mm]
 &
+\iota\frac {\tilde \ell} {\pi} \int_{{\mathfrak S}} \beta^{2\mathfrak p -2}\frac e{\beta^{\mathfrak p}} \mathcal Z \phi_1\, {\mathrm d} x \,{\mathrm d}\tau
 + \iota^2\frac {\tilde \ell} {\pi} \int_{{\mathfrak S}} \beta^{2\mathfrak p -2}\left(\frac e{\beta^{\mathfrak p}} \mathcal Z \phi_2 + \phi_1 \phi_2 \right) \, {\mathrm d} x \,{\mathrm d}\tau
 +O\left(e^{-\frac{\tilde \delta}\iota} \right).
\end{align*}
 We have
\begin{align*}
 & \frac {\tilde \ell} {2\pi \iota} \int_{{\mathfrak S}}
 \beta^{2\mathfrak p -2}\Big[U^2 + \iota^2 \frac{e^2}{\beta^{2 \mathfrak p}}|\mathcal Z|^2 + \iota^2|\phi_1|^2+ \iota^4|\phi_2|^2\Big]\, {\mathrm d} x \,{\mathrm d}\tau
\nonumber\\[2mm]
 & = \frac {\tilde \ell} {2\pi \iota} \Bigg[ \int_0^{2\pi} \beta^{2\mathfrak p -2} \int_{\mathbb R}  U^2 \, {\mathrm d} x \,{\mathrm d}\tau
 + \iota^2 \int_0^{2\pi}  \beta^{2\mathfrak p -2}\int_{\mathbb R} \frac{e^2}{\beta^{2\mathfrak p}} |\mathcal Z|^2 \, {\mathrm d} x \,{\mathrm d}\tau
 \nonumber\\[2mm]
 & \qquad\qquad + \iota^2 \int_0^{2\pi}  \beta^{2\mathfrak p -2} \int_{\mathbb R} |\phi_1|^2 \, {\mathrm d} x \,{\mathrm d}\tau + \iota^4 \int_0^{2\pi} \beta^{2\mathfrak p -2}\int_{\mathbb R} |\phi_2|^2 \, {\mathrm d} x \,{\mathrm d}\tau \Bigg]
\nonumber\\[2mm]
& =
 \frac {\tilde \ell} {2\pi \iota} \Bigg[   \int_0^{2\pi} \beta^{2\mathfrak p -2} \,{\mathrm d}\tau \varrho_0 + \iota^2 \int_0^{2\pi} \beta^{2\mathfrak p -2} \,{\mathrm d}\tau O\left(\| e\|^2_{L^\infty(0, 2\pi)}\right)
  \nonumber\\[2mm]
 & \qquad\qquad + \iota^2 \int_0^{2\pi} \beta^{2\mathfrak p -2}\int_{\mathbb R} \left( {\mathbf b}_1^2\varpi_1^2+ \frac1{\Lambda^2}f^2 {\mathbf b}_1^2\varpi_2^2 \right) \, {\mathrm d} x \,{\mathrm d}\tau + \iota^4 \bar{D}_{1\iota} \Bigg]
\nonumber\\[2mm]
& =
  \frac {1} { \iota}  \int_0^{\ell} \beta^{2\mathfrak p -1} \,{\mathrm d}\theta \left[  \varrho_0  + \iota^2O\left( 1\right)\right].
\end{align*}
 Moreover, there hold
\begin{align*}
  \frac {\tilde \ell} {\pi \iota} \int_{{\mathfrak S}} \beta^{2\mathfrak p -2}U \left(\iota \frac e{\beta^{\mathfrak p}}  \mathcal Z + \iota \phi_1 + \iota^2\phi_2 \right)\, {\mathrm d} x \,{\mathrm d}\tau
  =
   \int_0^{\ell} \beta^{2\mathfrak p -1} \,{\mathrm d}\theta\, O\Big(\|f\|_* + \|e\|_{**}  +1\Big),
\end{align*}
and
\begin{align*}
& \iota\frac {\tilde \ell} {\pi} \int_{{\mathfrak S}} \beta^{2\mathfrak p -2} \frac e{\beta^{\mathfrak p}} \mathcal Z \phi_1\, {\mathrm d} x \,{\mathrm d}\tau
 + \iota^2\frac {\tilde \ell} {\pi} \int_{{\mathfrak S}} \beta^{2\mathfrak p -2}\left( \frac e{\beta^{\mathfrak p}} \mathcal Z \phi_2 + \phi_1 \phi_2 \right) \, {\mathrm d} x \,{\mathrm d}\tau
=
   \int_0^{\ell} \beta^{2\mathfrak p -1} \,{\mathrm d}\theta\, O\Big(\iota^2\|f\|_* +\iota \|e\|_{**} \Big),
\end{align*}
 Then we conclude that
\begin{align} \label{w2integral}
\frac {\tilde \ell} {2\pi \iota} \int_{{\mathfrak S}_\delta} w_2^2 \, {\mathrm d} x \,{\mathrm d}\tau
  = \frac {1} { \iota}  \int_0^{\ell} \beta^{2\mathfrak p -1} \,{\mathrm d}\theta \Big[  \varrho_0  + \iota O\left( 1\right) + \iota \|f\|_* + \iota \|e\|_{**} \Big].
\end{align}


\smallskip
\noindent $\spadesuit$\
 By the relation $\tilde \phi = \beta^{\mathfrak p} \phi^* $, and \eqref{contraction-varphi},  we have
\begin{align} \label{tilde phi+tilde psi}
& \frac {\tilde \ell} {2\pi \iota} \int_{{\mathfrak S}_\delta} \frac {1} {\beta^2}
\left[
|{\tilde\phi}|^2+ |{\tilde\psi}|^2
+ 2\beta^{\mathfrak p} w_2 \big( {\tilde\phi} + {\tilde\psi} \big)
+ 2 \tilde \phi\tilde \psi
\right] \, {\mathrm d} x \,{\mathrm d}\tau
\nonumber\\[2mm]
&  \le    \frac {\tilde \ell} {2\pi \iota}     \int_0^{2\pi} \beta^{2\mathfrak p -2} \,{\mathrm d}\tau \sup_{\tau \in (0, 2\pi)} \|\phi^*(\cdot, \tau)\|^2_{L^2({\mathbb R})}
\nonumber\\[2mm]
&  \quad
+
O\left(
\iota \lf[\|{\tilde\phi}\|_{L^{\infty}(|s|>\frac{\delta}{\iota})}
\,+\, \|\nabla {\tilde\phi}\|_{L^{\infty}(|s|>\frac{\delta}{\iota})}\ri]
\,+\,
e^{-\frac{ \upsilon_2 \delta}{\iota}}
\right)
\frac {\tilde \ell} {2\pi \iota} \int_{{\mathfrak S}_\delta} \frac {1} {\beta^2} \, {\mathrm d} x \,{\mathrm d}\tau
\nonumber\\[2mm]
&  =    O\left(\frac{\tilde \ell} {\iota^2}\right)  \int_0^{\ell} \beta^{2\mathfrak p -1} \,{\mathrm d}\theta  \|\phi^*\|_{H^2_*(\mathfrak S)}^2
+
\frac{\tilde \ell} \iota O\left(e^{-\frac{\tilde \delta}\iota} \right)
\nonumber\\[2mm]
&  =    \frac 1{\iota}\int_0^{\ell} \beta^{2\mathfrak p -1} \,{\mathrm d}\theta    O\left(\frac{\tilde \ell} {\iota}\right)    \left( \iota^6+ \iota^4\|f\|_*^2 + \iota^4\|e\|_{**}^2\right)
\end{align}
 Here we note that we can show
$
\iota \lf[\|{\tilde\phi}\|_{L^{\infty}(|s|>\frac{\delta}{\iota})}
 \,+\, \|\nabla {\tilde\phi}\|_{L^{\infty}(|s|>\frac{\delta}{\iota})}\ri]
\in O\left(e^{-\frac{\tilde \delta}\iota} \right),
$
 by a barrier argument.

\medskip
 Combining \eqref{constraint6}, \eqref{estimate1} and \eqref{estimate2} together with \eqref{w2integral}-\eqref{tilde phi+tilde psi}, we get
\begin{align} \label{a-equality}
  a\iota^{-2+ \frac4{p_\epsilon-2}} &=     \frac {1} { \iota}  \int_0^{\ell} \beta^{2\mathfrak p -1} \,{\mathrm d}\theta \left[  \varrho_0  + \iota O\left( 1\right) + \iota \|f\|_* + \iota \|e\|_{**} \right]
+ \frac{\tilde \ell} \iota O\left(e^{-\frac{\tilde \delta}\iota} \right).
\end{align}
due to \eqref{constraint-f}, \eqref{constraint-e} and \eqref{constraint-alpha-varrho}.
By \eqref{a-equality}, we conclude that
\begin{align*}
1&=     \iota^{1- \frac4{p_\epsilon-2}} \frac{\hat \ell   \varrho_0}a \left(1+ O\left(\iota\right) \right),
  \end{align*}
  where $\hat \ell =\int_0^{\ell} \beta^{2\mathfrak p -1} \,{\mathrm d}\theta. $

In the following, we discuss from two cases.
\\
{\em Case 1}:     If  $ \hat \ell   \varrho_0 <a$,
we  consider
\begin{align} \label{identity from mass constraint2}
1&=     \iota^{1- \frac4{p_\epsilon -2}} \frac{\hat \ell   \varrho_0}a \left(1+ O\left(\iota\right) \right).
  \end{align}
  We  recall that,  in this case, we have
  \[  \frac{\ell_0  \varrho_0}a <\frac{\hat \ell   \varrho_0}a < c_1(\e)<1.   \]
  By Lemma \ref{lemma 7-1}, we know that  there exists a solution  $\iota (\epsilon)$ to \eqref{identity from mass constraint2} with
  \[\iota (\epsilon) \in \left(e^{\frac{6} \e \ln  {\frac{\ell_0   \varrho_0}a}},   e^{\frac{2} \e \ln  c_1(\e)}   \right). \]
 \\
{\em Case 2}:  If  $ \hat \ell   \varrho_0 >a$, we  consider
\begin{align}  \label{identity from mass constraint}
1&=     \iota^{1- \frac4{4+\e}} \frac{\hat \ell   \varrho_0}a \left(1+ O\left(\iota\right) \right).
  \end{align}
 We  recall that,  in this case, we have
  \[   1<c_2 (\epsilon) <\frac{\hat \ell   \varrho_0}a <c_3(\e).   \]
 Then by Lemma \ref{lemma 7-1}, we can show that  there exists
\[\iota (\epsilon) \in  \left( e^{-\frac{2} \e \ln c_2(\e)},   e^{-\frac{6} \e \ln  c_3(\e)}  \right)  \]
 solves \eqref{identity from mass constraint}.

\begin{lemma}    \label{lemma 7-1}
 Consider  the
\begin{align}  \label{identity for iota}
1&=     \iota^{1- \frac4{p_\epsilon-2}} \frac{\hat \ell   \varrho_0}a \left(1+ O\left(\iota\right) \right),
  \end{align}
 where
\begin{align}
  \hat \ell =\int_0^{\ell} \beta^{2\mathfrak p -1} \,{\mathrm d}\theta  =  \int_0^{\ell} \frac1{\beta^{\frac{\pm \e} {4\pm \e}}}  \,{\mathrm d}\theta.
  \end{align}
  We have the following results to the existence of \eqref{identity for iota}.
\begin{itemize}

\item[(i).]
{\em
If  $\frac{\hat \ell   \varrho_0}a <1, $  we take $p_\epsilon = 6- \e$.
Suppose there exist $\tilde c_1(\e)$ and $\tilde c_2(\e)$  independent of $\iota$,
\begin{align}
0< \tilde c_1(\e)     <\frac{\hat \ell \varrho_0} {a}  < \tilde c_2(\e)<1.
\end{align}
Then there exists a solution $\iota(\epsilon)$ to \eqref{identity for iota}  which satisfies $\iota(\epsilon) \in ( e^{\frac{6} \e \ln \tilde c_1(\e)},   e^{\frac{2} \e \ln \tilde c_2(\e)}  )$
  }
   \\[0.5mm]
   \item[(ii).]
{\em
  If  $\frac{\hat \ell   \varrho_0}a >1, $  we take $p_\epsilon = 6+ \e$.
  Suppose there exist $\tilde c_3(\e)$ and $\tilde c_4(\e)$  independent of $\iota$,
\begin{align}
1< \tilde c_3(\e)     <\frac{\hat \ell \varrho_0} {a}  < \tilde c_4(\e).
\end{align}
Then there exists a solution $\iota(\epsilon)$ to \eqref{identity for iota}  which satisfies $\iota(\epsilon) \in ( e^{\frac{2} \e \ln \tilde c_3(\e)},   e^{\frac{6} \e \ln \tilde c_4(\e)}  )$

  }

     \end{itemize}

\end{lemma}
\begin{proof}

  {\em Case 1:}
 If $\frac{\hat \ell   \varrho_0}a <1$, we take $p_\epsilon = 6- \e$.  In this case, we consider
\begin{align*}
1&=     \iota^{1- \frac4{4-\e}} \frac{\hat \ell   \varrho_0}a \left(1+ O\left(\iota\right) \right),  \quad \text{with}~ \hat \ell =  \int_0^{\ell}  \beta^{\frac{ \e} {4- \e}} \,{\mathrm d}\theta.
  \end{align*}
  Denote
  \[  \mathfrak w_1 (\iota) : = \iota^{1- \frac4{4-\e}} \frac{\hat \ell   \varrho_0}a \left(1+ O\left(\iota\right) \right). \]
  Suppose there exist $\tilde c_1(\e)$ and $\tilde c_2(\e)$  independent of $\iota$,
\begin{align}
0< \tilde c_1(\e)     <\frac{\hat \ell \varrho_0} {a}  < \tilde c_2(\e)<1.
\end{align}
We  give the following claim: \\
{\bf Claim 1:}  {\em
 We take $\iota_1 (\e) = e^{\frac{6} \e \ln \tilde c_1(\e)}  $ and   $\iota_2 (\e)= e^{\frac{2} \e \ln \tilde c_2(\e)}  $.  Obviously,  $ \iota_2 (\e) >\iota_1 (\e).  $ We can verify that
\begin{align}
 \mathfrak w_1 (\iota_1) >1, \quad   \mathfrak w_1 (\iota_2) <1.
 \end{align}
}

\noindent{\em Proof of {\bf Claim 1}}:  We have
\begin{align}
 \iota^{ \frac\e{4-\e}} =  \frac{\hat \ell   \varrho_0}a \left(1+ O\left(\iota\right) \right),
\end{align}
i.e,
\begin{align}
 \frac\e{4-\e} \ln \iota =  \ln{\frac{\hat \ell   \varrho_0}a \left(1+ O\left(\iota\right) \right)},
\end{align}
\begin{align}
  \ln \iota =  \frac{4-\e}\e \ln{\frac{\hat \ell   \varrho_0}a \left(1+ O\left(\iota\right) \right)},
\end{align}
\begin{align}
  \iota =  e^{ \frac{4-\e}\e \ln{\frac{\hat \ell   \varrho_0}a \left(1+ O\left(\iota\right) \right)}}.
\end{align}
We have
\begin{align}
 \frac\e{4-\e} \ln \iota_1 =   \frac\e{4-\e} \ln {e^{\frac{6} \e \ln \tilde c_1(\e)} } =   \frac6{4-\e}  \ln \tilde c_1(\e) < \frac6{4-\e}  \ln{\frac{\hat \ell   \varrho_0}a \left(1+ O\left(\iota_1\right) \right)},
 \end{align}
 i.e,
\begin{align}
 \iota_1^{ \frac\e{4-\e}} <  \left(\frac{\hat \ell   \varrho_0}a \left(1+ O\left(\iota_1\right) \right) \right)^{\frac6{4-\e}} <\frac{\hat \ell   \varrho_0}a \left(1+ O\left(\iota_1\right) \right).
\end{align}
   The last inequality hold due to the fact $\frac{\hat \ell \varrho_0} {a}<1$.  Therefore we have $ \mathfrak w_1 (\iota_1) >1$.

 On the other hand, we have
\begin{align}
 \frac\e{4-\e} \ln \iota_2 =   \frac\e{4-\e} \ln {e^{\frac{2} \e \ln \tilde c_2(\e)} } =   \frac2{4-\e}  \ln \tilde c_2(\e) > \frac2{4-\e}  \ln{\frac{\hat \ell   \varrho_0}a \left(1+ O\left(\iota_2\right) \right)},
 \end{align}
  i.e,
\begin{align}
 \iota_2^{ \frac\e{4-\e}}
> \left(\frac{\hat \ell   \varrho_0}a \left(1+ O\left(\iota_2\right) \right) \right)^{\frac2{4-\e}} >\frac{\hat \ell   \varrho_0}a \left(1+ O\left(\iota_2\right) \right).
\end{align}
Therefore we have $ \mathfrak w_1 (\iota_2) <1$.

 \medskip

  {\em Case 2:}
 If $\frac{\hat \ell   \varrho_0}a >1$, we take $p_\epsilon = 6+ \e$.  In this case, we consider
\begin{align*}
1&=     \iota^{1- \frac4{4+\e}} \frac{\hat \ell   \varrho_0}a \left(1+ O\left(\iota\right) \right),  \quad \text{with}~ \hat \ell =  \int_0^{\ell}  \beta^{\frac{ \e} {4+\e}} \,{\mathrm d}\theta.
  \end{align*}
   Denote
  \[  \mathfrak w_2 (\iota) : = \iota^{1- \frac4{4+\e}} \frac{\hat \ell   \varrho_0}a \left(1+ O\left(\iota\right) \right). \]
   Suppose there exist $\tilde c_3(\e)$ and $\tilde c_4(\e)$  independent of $\iota$,
\begin{align}
1< \tilde c_3(\e)     <\frac{\hat \ell \varrho_0} {a}  < \tilde c_4(\e).
\end{align}
We  give the following claim: \\
{\bf Claim 2:}  {\em
 We take $\iota_3(\e) = e^{-\frac{2} \e \ln \tilde c_3(\e)}  $ and   $\iota_4 (\e)= e^{-\frac{6} \e \ln \tilde c_4(\e)}  $.  Obviously,  $ 0<\iota_3 (\e) <\iota_4(\e).  $ We can verify that
\begin{align}
 \mathfrak w_2 (\iota_3) >1, \quad   \mathfrak w_2 (\iota_4) <1.
 \end{align}
}

\noindent{\em Proof of {\bf Claim 2}}:  We have
\begin{align}
 \iota^{ \frac\e{4+\e}} =  \frac{\hat \ell   \varrho_0}a \left(1+ O\left(\iota\right) \right),
\end{align}
i.e,
\begin{align}
 \frac\e{4+\e} \ln \iota =  \ln{\frac{\hat \ell   \varrho_0}a \left(1+ O\left(\iota\right) \right)},
\end{align}
\begin{align}
  \ln \iota =  \frac{4+\e}\e \ln{\frac{\hat \ell   \varrho_0}a \left(1+ O\left(\iota\right) \right)},
\end{align}
\begin{align}
  \iota =  e^{ \frac{4+\e}\e \ln{\frac{\hat \ell   \varrho_0}a \left(1+ O\left(\iota\right) \right)}}.
\end{align}
We have
\begin{align}
- \frac\e{4+\e} \ln \iota_3 =  - \frac\e{4+\e} \ln {e^{\frac{2} \e \ln \tilde c_3(\e)} } = -   \frac2{4+\e}  \ln \tilde c_3(\e) >- \frac2{4+\e}  \ln{\frac{\hat \ell   \varrho_0}a \left(1+ O\left(\iota_1\right) \right)},
 \end{align}
 i.e,
\begin{align}
 \iota_3^{ \frac\e{4+\e}} <  \left(\frac{\hat \ell   \varrho_0}a \left(1+ O\left(\iota_3\right) \right) \right)^{\frac2{4+\e}} <\frac{\hat \ell   \varrho_0}a \left(1+ O\left(\iota_1\right) \right).
\end{align}
   The last inequality hold due to the fact $\frac{\hat \ell \varrho_0} {a}>1$.  Therefore we have $ \mathfrak w_2(\iota_3) >1$.

 On the other hand, we have
\begin{align}
 -\frac\e{4+\e} \ln \iota_4 =  - \frac\e{4+\e} \ln {e^{\frac{6} \e \ln \tilde c_4(\e)} } =   -\frac6{4+\e}  \ln \tilde c_4(\e) <- \frac6{4+\e}  \ln{\frac{\hat \ell   \varrho_0}a \left(1+ O\left(\iota_2\right) \right)},
 \end{align}
  i.e,
\begin{align}
 \iota_4^{ \frac\e{4+\e}}
> \left(\frac{\hat \ell   \varrho_0}a \left(1+ O\left(\iota_2\right) \right) \right)^{\frac6{4+\e}} >\frac{\hat \ell   \varrho_0}a \left(1+ O\left(\iota_2\right) \right).
\end{align}
Therefore we have $ \mathfrak w_2 (\iota_4) <1$.

\end{proof}

\medskip
\medskip
{\bf Acknowledgements:}
L. Duan and S. Wei was supported by NSFC grant (No. 12201140, No.12001203), Guangdong Basic and Applied Basic Research Foundation (No. 2022A1515111131, No. 2026A1515010521). J. Yang was supported by NSFC (No. 12171109).

\end{document}